\documentclass[english]{article} 
\usepackage[absolute,overlay]{textpos} 
\usepackage{amsthm}
\usepackage{amsmath}
\usepackage{amssymb}
\usepackage[utf8]{inputenc}
\usepackage[english]{babel}
\usepackage{graphicx}
\usepackage{geometry}
\usepackage{caption}
\usepackage{subcaption}
\usepackage{pdfpages}
\usepackage{comment}
\usepackage{bm}
\usepackage{etoolbox}
\usepackage{mathtools}
\usepackage{tikz}
\usetikzlibrary{calc}

\usepackage{algpseudocode}
\usepackage{algorithm}
\usepackage{mathrsfs} 
\usepackage{multirow}

\usepackage{pgfplots}
\usepackage{pgfplotstable}
\usepackage{stmaryrd} 

\usepackage{hyperref}

\usepackage{booktabs}

\usepackage[numbers]{natbib}  

\usepackage[noabbrev,capitalise]{cleveref}
\crefformat{appendix}{#2#1#3}

\usetikzlibrary{fillbetween}
\usepgfplotslibrary{fillbetween}

\usetikzlibrary{shapes.misc}
\usetikzlibrary{math} 

\usetikzlibrary{arrows}

\definecolor{darkspringgreen}{rgb}{0., 0.55, 0.3}
\definecolor{dartmouthgreen}{rgb}{0.05, 0.5, 0.06}
\definecolor{etonblue}{rgb}{0.59, 0.78, 0.64}
\definecolor{airforceblue}{rgb}{0., 0.4, 0.66}
\definecolor{arylideyellow}{rgb}{0.91, 0.84, 0.42}
\definecolor{emerald}{rgb}{0.31, 0.78, 0.47}
\definecolor{uclagold}{rgb}{1.0, 0.7, 0.0}
\definecolor{cadmiumorange}{rgb}{0.93, 0.53, 0.18}

\theoremstyle{thmstyleone}

\theoremstyle{thmstyletwo}

\newtheorem{remark}{Remark}

\theoremstyle{thmstylethree}

\newcommand{\uvec}[2][3]{\boldsymbol{#2\mkern-#1mu}\mkern#1mu}

\newcommand\abs[1]{\left\lvert#1\right\rvert}

\newcommand{\R}{\mathbb{R}}

\newcommand{\bu}{\uvec{u}}

\newcommand{\bG}{{\uvec{g}}}

\newcommand{\Temp}{T}

\newcommand{\E}{\boldsymbol{E}}

\newcommand{\RIcolor}[1]{{\leavevmode\color{black} #1}}
\newcommand{\RIIcolor}[1]{{\leavevmode\color{black} #1}}

\newcommand{\bF}{\uvec{f}}
\newcommand{\diff}[1]{{\mathrm{d}{#1}}}

\newcommand{\ubar}{\overline{\uvec{u}}}

\newcommand{\xip}{x_{i+\frac{1}{2}}}
\newcommand{\xin}{x_{i-\frac{1}{2}}}
\newcommand{\yjp}{y_{j+\frac{1}{2}}}
\newcommand{\yjn}{y_{j-\frac{1}{2}}}
\newcommand{\iip}{i+\frac{1}{2}}
\newcommand{\iin}{i-\frac{1}{2}}

\newcommand{\jjp}{j+\frac{1}{2}}
\newcommand{\jjn}{j-\frac{1}{2}}

\newcommand{\jph}{{j+\frac{1}{2}}}
\newcommand{\jmh}{{j-\frac{1}{2}}}

\newcommand{\iph}{{i+\frac{1}{2}}}
\newcommand{\imh}{{i-\frac{1}{2}}}

\newcommand{\qbar}{\overline{q}}

\begin{document}
\author{Lorenzo Micalizzi\footnote{Affiliation: Department of Mathematics, North Carolina State University, SAS Hall, 2108, 2311 Stinson Dr, Raleigh, NC 27607, United States. Email: lmicali@ncsu.edu}, and Eleuterio F. Toro\footnote{Affiliation: Laboratory of Applied Mathematics, DICAM, University of Trento, Via Mesiano 77, 38123 Trento, Italy. Email: eleuterio.toro@unitn.it}}
\title{FORCE-$\alpha$ Numerical Fluxes within the Arbitrary High Order Semidiscrete WENO-DeC Framework: A Competitive Alternative to Upwind Fluxes} 


\maketitle

\abstract{
	This work systematically investigates the performance of FORCE--$\alpha$ numerical fluxes within an arbitrary high order semidiscrete finite volume (FV) framework for hyperbolic partial differential equations (PDEs).
	Such numerical fluxes have been recently introduced by Toro, Saggiorato, Tokareva, and Hidalgo (Journal of Computational Physics, 416, 2020), and constitute a family of centred fluxes obtained from a suitable modification of First--Order Centred (FORCE) numerical fluxes.
	In contrast with upwind fluxes, such as Rusanov, Harten--Lax--van Leer (HLL) or the exact Riemann solver (RS) numerical flux, centred ones do not consider in any way the structure of the Riemann problem at cell interfaces.
	Adopting centred numerical fluxes leads to a high level of flexibility of the resulting numerical schemes, for example in the context of complicated hyperbolic systems, for which RSs may be impossible to construct or computationally expensive.
	
	The baseline framework adopted in this investigation is a FV semidiscrete approach with Weighted Essentially Non--Oscillatory (WENO) spatial reconstruction and Deferred Correction (DeC) time discretization, and results are reported up to order 7.
	Previous investigations involving the same framework have established that increasing the order of accuracy tends to decrease the differences in the results obtained through different numerical fluxes.
	The goal of this paper is to show that the employment of FORCE--$\alpha$ numerical fluxes within such a framework is a competitive alternative to the adoption of more classical upwind fluxes.
	The hyperbolic system considered for this investigation is the ideal Euler equations in one and two space dimensions.
	\RIIcolor{
	The results show that moderate values of $\alpha$ are generally the most efficient and provide the best overall compromise between resolution, robustness, and computational cost: $\alpha\approx 1$--$2$ in one space dimension and $\alpha\approx 2$--$3$ in two space dimensions.
	Overall, they confirm that, in very high order settings, FORCE--$\alpha$ numerical fluxes are able to achieve performance comparable to upwind fluxes, while retaining the simplicity, flexibility, and low computational cost of centred fluxes.	
	}
}

\section{Introduction}
Numerical fluxes constitute a key ingredient required by most of the existing approaches for the numerical solution of hyperbolic partial differential equations (PDEs), such as finite volume  (FV)~\cite{hirsch2007numerical,ToroBook,leveque2002finite,godlewski2021numerical,toro2024computational}, finite difference~\cite{leveque2007finite} and Discontinuous Galerkin finite element methods~\cite{reed1973triangular,cockburn2000development,cockburn2001runge}.
As a matter of fact, they represent an approximation of the information exchanged by confining regions with discontinuous data across the interface.
Giving a complete overview of all existing numerical fluxes is difficult, and we refer the reader to classical literature on the topic~\cite{hirsch2007numerical,ToroBook,leveque2002finite,godlewski2021numerical,toro2024computational}.
It is worth remarking that also interesting works devoted to the comparison of the performance of different numerical fluxes are available, see~\cite{leidi2024performance,qiu2007numerical,qiu2008development,qiu2006numerical,hongxia2020numerical,san2015evaluation,micalizzitoro2024}.
In fact, the choice of the numerical flux influences the properties of the resulting numerical scheme.
For example, the naive adoption of a simple average of the fluxes computed in the data from both sides of the interface (usually referred to as ``central flux'') is well known to lead to unstable schemes.
Moreover, it is also well known that the diffusion of the schemes is highly influenced by the specific numerical flux adopted, with  the Godunov flux from the exact Riemann solver~(RS)~\cite{Godunov} being the most precise numerical flux and the least diffusive monotone one (for scalar problems).
Further, the ability of schemes to sharply capture specific flow features is strictly related to the properties of the employed numerical flux.
For example, the HLLC numerical flux~\cite{toro1992restoration,toro1994restoration} has been specifically designed, as a modification of the Harten--Lax--van Leer (HLL) numerical flux~\cite{harten1983upstream}, to cure its inability in capturing contact discontinuities, shear waves and material interfaces.
Such a fundamental difference between HLLC and HLL lies in the construction itself of these numerical fluxes, with HLLC being a complete upwind numerical flux and HLL being an incomplete upwind one according to the characterization described in the following.

In fact, all numerical fluxes are built upon the Riemann problem~\cite{riemann1860fortpflanzung,Godunov,ToroBook}, involving the governing PDE in the normal direction to the interface and the (reconstructed) states from both sides of
the interface.
Based on this fact, one can define two main families of numerical fluxes, namely, upwind fluxes and centred (or non--upwind) ones.
The first ones make explicit use of the structure of the underlying Riemann problem in terms of wave propagation information; the latter ones do not, and only require simple flux evaluations in their computation.
For example, exact RS, HLL, HLLC and Rusanov~\cite{Rusanov1961} are amongst the most popular upwind fluxes, while, Lax--Friedrichs~\cite{lax1954weak}, First--Order Centred (FORCE)~\cite{Toro1996,toro2000centred,chen2003centred} and FORCE--$\alpha$ are centred ones.
Upwind fluxes are constructed via the (exact or approximate) solution of the underlying Riemann problem through a wave model containing $A$ waves.
Based on this, one can further divide them into complete and incomplete (upwind) fluxes:
complete ones are characterized by $A = E$, with $E$ being the number of characteristic fields in the actual Riemann problem, while, incomplete ones are characterized by $A < E$.
For example, for the Euler equations (in any number of space dimensions), HLLC and exact RS are complete and thus able to capture characteristic fields associated with contact discontinuities, shear waves and material interfaces~\cite{ToroBook}; instead, Rusanov and HLL are incomplete, and this determines their inability to capture such flow features.
The same inability is shared by the Central--Upwind~\cite{kurganov2001semidiscrete} numerical flux, which is as a matter of fact the HLL numerical flux with Davis' speed estimates~\cite{davis1988simplified}, see~\cite{micalizzitoro2024}, and related low--dissipation versions~\cite{KLin} with the exception of the one proposed in~\cite{CKX_Ustar}.

In general, the deeper is the reliance employed in its construction on the correct physical description of the solution of the Riemann problem, the more accurate is the numerical flux in terms of ability of the resulting scheme to capture complex features. 
In fact, due to their particular attention on the underlying physics, upwind fluxes commonly have a higher resolution with respect to centred ones~\cite{micalizzitoro2024}. 
However, the simplicity of the latter ones makes them more suitable for complex applications, for which Riemann structures may not be directly accessible, e.g., for two (or more)--layer shallow-water systems the eigenstructure is not known in closed--form~\cite{abgrall2009two,kim2009two}.
Moreover, one of the main conclusions from~\cite{micalizzitoro2024,micalizzi2025algorithms} suggests that the impact of the adopted numerical flux tends to become less important as the order of accuracy increases to very high values. 
This is the main motivation behind this investigation.
Here, we are concerned with a particular family of centred fluxes, namely, the FORCE--$\alpha$ numerical fluxes~\cite{toro2020low}, as an alternative to classical upwind ones, within a high order framework which may make--up for their higher level of diffusion. 
For this purpose, we consider an arbitrary high order FV semidiscrete setting where a Weighted Essentially Non--Oscillatory (WENO)~\cite{liu1994weighted,shu1998essentially,shu1989efficient} spatial reconstruction is adopted, along with Deferred Correction (DeC)~\cite{micalizzi2023new,Decremi,minion2003semi,Decoriginal,lore_phd_thesis} time integration, and we report results up to order 7 in space and time.
The same setting has been adopted in~\cite{ciallella2022arbitrary,ciallella2023arbitrary,ciallella2025high,micalizzitoro2024,micalizzi2025algorithms} for investigations in various directions.
We remark that one of the distinctive characters of the considered framework is the adoption of the same order of accuracy for both the discretizations in space and time. 
In fact, very often, high order space discretizations are coupled with lower order time discretizations, typically strong stability preserving Runge--Kutta schemes, see for example~\cite{evstigneev2016construction,Evstigneev2016OnTC,gerolymos2009very,balsara2000monotonicity,shi2003resolution,hermes2012linear,gao2020seventh}.
The inefficiency of such a choice has been clearly demonstrated in~\cite{micalizzi2025algorithms}, even (and especially) when the time step is chosen to be proportional to a power of the characteristic mesh size to guarantee a matching between the orders of accuracy in space and time.
We also remark that other frameworks allow for the construction of schemes of arbitrarily high order in space and time, e.g., the fully--discrete ADER method~\cite{toro2001towards,dumbser2008unified} or the semidiscrete ADER--Spectral Difference method~\cite{velasco2023spectral,veiga2024improving}.

The results obtained with FORCE--$\alpha$ numerical fluxes are compared with the ones obtained with Rusanov, HLL and exact RS in order to provide valuable references for the reader.

\RIIcolor{Summarizing, the main goal of this work is to investigate the performance of centred FORCE--$\alpha$ numerical fluxes within a very high order setting. Our numerical results show that:
\begin{itemize}
	\item moderate values of $\alpha$, i.e., $\alpha\approx 1$--$2$ and $\alpha\approx 2$--$3$ in one and two space dimensions respectively, provide the best compromise between resolution, robustness, and computational cost;
	\item FORCE--$\alpha$ numerical fluxes constitute an efficient, computationally inexpensive, and robust alternative to classical upwind fluxes at high orders of accuracy.
\end{itemize}}

The paper is structured as follows.
In Section~\ref{sec:WENO_DeC}, we describe the FV WENO--DeC framework for hyperbolic PDEs adopted in this investigation.
In Section~\ref{sec:FORCE_alpha}, we recall the FORCE--$\alpha$ numerical fluxes and their properties.
Numerical results are reported in Section~\ref{sec:numerical_results}.
Conclusions are finally drawn, along with further perspectives, in Section~\ref{sec:conclusions}.

\section{Semidiscrete Finite Volume WENO--DeC schemes}\label{sec:WENO_DeC}
In this section, we will present the governing equations and the FV WENO--DeC framework adopted for the investigation of the FORCE--$\alpha$ numerical fluxes.
In particular, the governing equations will be described in Section~\ref{sec:governing_equations}, while, the numerical discretization will be detailed in Sections~\ref{sec:space} and~\ref{sec:time}.
Let us remark that the considered framework is a semidiscrete one, meaning that the discretizations in space and in time are essentially decoupled.
In fact, we will start by describing the space discretization in Section~\ref{sec:space}, and we will continue with the time discretization in Section~\ref{sec:time}.

\subsection{Governing equations}\label{sec:governing_equations}
Let us start by defining the analytical problem under investigation. We focus on the two--dimensional case, with the restriction to one space dimension being straightforward. 
We are interested in the numerical solution of hyperbolic systems of PDEs in the form
\begin{equation}\label{eq:sys}
	\frac{\partial}{\partial t} \uvec{u}(x,y,t) + \frac{\partial}{\partial x}\uvec{f}(\uvec{u}(x,y,t))+\frac{\partial}{\partial y}\uvec{g}(\uvec{u}(x,y,t)) = \uvec{0},
\end{equation}
where $(x,y)\in\Omega$ and $t\in[0,T_f]$, with $\Omega\subseteq \mathbb{R}^2$ being a two--dimensional space domain and $T_f\in \mathbb{R}^+$ being the final time.
More in detail, $\uvec{u}:\Omega\times (0,T_f)\longrightarrow \mathbb{R}^{N_c}$ represents the unknown solution in conserved variables, where $N_c\in \mathbb{N}^+$ is the number of equations in the system, while, $\uvec{f},\uvec{g}:\mathbb{R}^{N_c}\longrightarrow\mathbb{R}^{N_c}$ are the fluxes in the $x$- and $y$-directions.
Hyperbolicity is guaranteed by the following assumption: 
the matrix $\omega_1 \frac{\partial\uvec{f}}{\partial \uvec{u}}(\uvec{u})+\omega_2 \frac{\partial\uvec{g}}{\partial \uvec{u}}(\uvec{u})$ admits $N_c$ real eigenvalues and a corresponding set of linearly independent eigenvectors for any $\omega_1,\omega_2\in \mathbb{R}$.
It is useful to introduce the matrix
\begin{equation}
	J_{\uvec{\nu}}(\uvec{u}):=\nu_1 \frac{\partial\uvec{f}}{\partial \uvec{u}}(\uvec{u})+\nu_2 \frac{\partial\uvec{g}}{\partial \uvec{u}}(\uvec{u}),
\end{equation}
for any two--dimensional unit vector $\uvec{\nu}:=(\nu_1,\nu_2)^T\in \mathbb{R}^2$, which is the normal Jacobian matrix in direction $\uvec{\nu}$.
Its eigenvalues are the wave speeds of the system along the direction $\uvec{\nu}$.
To be more specific, here, we will focus on the Euler equations described in the following.

\subsubsection{Euler equations}
The Euler equations in two space dimensions are a hyperbolic system of PDEs in the form~\eqref{eq:sys} with
\begin{equation}
	\uvec{u}:=\begin{pmatrix}\rho\\ \rho u\\ \rho v\\E\end{pmatrix},\quad\uvec{f}(\uvec{u}):=\begin{pmatrix}\rho u\\\rho u^2+p\\\rho uv \\(E+p)u\end{pmatrix},\quad\uvec{g}(\uvec{u}):=\begin{pmatrix}\rho v\\\rho uv\\\rho v^2+p \\(E+p)v\end{pmatrix}.
	\label{eq:conservative_variables_2D}
\end{equation}
In particular, $\rho$ represents the fluid density, $u$ and $v$ the speed components along the $x$- and $y$-directions, $p$ represents the fluid pressure, while, $E$ is the total energy.
The system is closed by a further equation
\begin{align}
	E&=\rho \left[e + \frac{1}{2}(u^2+v^2) \right],
\end{align}
with $e$ being the specific internal energy given by an equation of state $e:=e(\rho,p)$ as a function of density and pressure. As we consider the case of ideal fluid in this work, we have $e(\rho,p):=\frac{p}{\rho(\gamma-1)}$ where $\gamma$ is the adiabatic coefficient, in the numerical experiments set to be $\gamma:=1.4$.
The wave speeds along the $\uvec{\nu}:=(\nu_1,\nu_2)^T$ direction are $v_{\uvec{\nu}}+c$, $v_{\uvec{\nu}}-c$ and $v_{\uvec{\nu}}$, with $v_{\uvec{\nu}}:=u \nu_1+v \nu_2$. The sound speed $c$ is given by
$c:=\sqrt{\frac{\frac{p}{\rho^2}-\frac{\partial}{\partial \rho}e}{\frac{\partial}{\partial p}e}}$, reducing to $c:=\sqrt{\gamma \frac{p}{\rho}}$ in the ideal case.

\subsection{Space discretization}\label{sec:space}
The FV method~\cite{Godunov,hirsch2007numerical,godlewski2021numerical,toro2024computational,shu1988efficient,shu1989efficient,shu1998essentially,leveque2002finite,ToroBook,abgrall1994essentially} discretizes the solution through cell averages over suitable control volumes covering the space domain.
Therefore, we consider a Cartesian setting, and we cover $\Omega$, which is assumed to be rectangular, through uniform non--overlapping rectangles $C_{i,j}:=[x_{\imh},x_{\iph}]\times [y_{\jmh},y_{\jph}]$, with $\xip-\xin=\Delta x$ and $\yjp-\yjn=\Delta y$ for any $i$ and $j$. We further introduce $\ubar_{i,j}(t)$, representing the vector of the cell averages of the conserved variables over the generic cell $C_{i,j}$.
Following a semidiscrete approach, we start by performing a discretization in space only, namely, we integrate the governing PDE~\eqref{eq:sys} over $C_{i,j}$ and apply Gauss's theorem, leading to 
\begin{equation}\label{eq:FV_semidiscretization}
	\frac{d}{dt}\ubar_{i,j}(t) + \frac{1}{\Delta x}(\uvec{f}_{\iip,j}-\uvec{f}_{\iin,j}) + \frac{1}{\Delta y}(\uvec{g}_{i,\jjp}-\uvec{g}_{i,\jjn}) = \uvec{0},
\end{equation}
with 
\begin{align}
	\uvec{f}_{\iip,j} &:= \frac{1}{\Delta y} \sum_{q=1}^{N_w} w_q \widehat\bF(\bu^-(x_{\iip},y_q),\bu^+(x_{\iip},y_q))\approx \frac{1}{\Delta y}\int_{\yjn}^{\yjp}\bF(\bu(\xip,y,t))\;dy, \label{eq:flux_F}   \\ 
	\uvec{g}_{i,\jjp} &:= \frac{1}{\Delta x} \sum_{q=1}^{N_w} w_q \widehat\bG(\bu^-(x_q,y_{\jjp}),\bu^+(x_q,y_{\jjp}))\approx \frac{1}{\Delta x}\int_{\xin}^{\xip}\uvec{g}(\bu(x,\yjp,t))\;dx, \label{eq:flux_G}
\end{align}
representing the semidiscretization in space of the analytical problem.
Notice that the flux surface integrals at cell boundaries have been discretized via the employment of sufficiently accurate quadrature formulas with points and weights $(y_q,w_q)$ and $(x_q,w_q)$.
Further, the physical fluxes, $\bF(\cdot)$ and $\bG(\cdot)$, in the quadrature points have been approximated through numerical fluxes, $\widehat\bF(\cdot,\cdot)$ and $\widehat\bG(\cdot,\cdot)$, taking in input reconstructions of the conserved variables, $\bu^-$ and $\bu^+$, in the surface quadrature points from both the cells sharing the edges.
In fact, two final ingredients must be specified to complete the discretization in space: the employed numerical flux and space reconstruction. 
%
%
Being numerical fluxes the main focus of this work, we leave a more complete description of this point for Section~\ref{sec:FORCE_alpha}.
Concerning the space reconstruction, we aim at achieving two main requirements: (arbitrary) high order of accuracy and ability to handle discontinuities which naturally arise in the context of hyperbolic PDEs even from smooth initial conditions. 
Let us remark that the order of accuracy of the space discretization only depends on the order of the space reconstruction (provided that quadrature formulas of the same order are assumed for the surface integrals and for the initialization of the cell averages). 
In principle, constructing (arbitrary) high order accurate linear space discretizations based on the available cell averages is not difficult. However, recalling Godunov's theorem~\cite{Godunov}, we have that, for scalar problems, there are no monotone linear schemes of order of accuracy greater than one.
Therefore, nonlinearity is a key ingredient, necessary but not sufficient, to be incorporated in a robust space reconstruction.
For the above reasons, we adopt a WENO~\cite{liu1994weighted} space discretization, able to conjugate (arbitrary) high order of accuracy with good shock--capturing performance.
Here, we only sketch the rough idea, referring the interested reader to~\cite{ciallella2022arbitrary,micalizzitoro2024,micalizzi2025algorithms,lore_phd_thesis}, where the approach is thoroughly described, and to classical works on the topic~\cite{shu1998essentially,jiang1996efficient,shu1988efficient,shu1989efficient,abgrall1994essentially}.
The WENO reconstruction of order $2r-1$, for a scalar quantity $q$, is performed locally in each cell $C_{i,j}$ considering a stencil of $(2r-1)\times(2r-1)$ cells 
\begin{align}
C_{\ell,m},\quad\ell=i-(r-1),\dots,i+(r-1), \quad m=j-(r-1),\dots,j+(r-1).
\end{align}
Given a surface quadrature point $(x_{\iip},y_q)$ (or $(x_{i-\frac{1}{2}},y_q)$) in the cell $C_{i,j}$, the reconstruction is performed via two orthogonal one--dimensional sweeps, one in the $x$-direction and one in the $y$-direction (the other way around for points $(x_q,y_{\jjp})$ and $(x_q,y_{j-\frac{1}{2}})$).
In the first one--dimensional sweep, one reconstructs averages of $q$ at cell interfaces with  $x=x_{i+\frac{1}{2}}$ with respect to the $y$-direction
\begin{align}
	\qbar^-_{m} \approx \frac{1}{\Delta y}\int_{y_{m-\frac{1}{2}}}^{y_{m+\frac{1}{2}}} q(\xip,y)dy,
\end{align}
for each $m=j-(r-1),\dots,j+(r-1)$ out of the cell averages $\qbar_{\ell, m}$ $\ell=i-(r-1),\dots,i+(r-1)$.
After this, another sweep in the $y$-direction is performed reconstructing, out of $\qbar^-_{m}$ $m=j-(r-1),\dots,j+(r-1)$, the required values at quadrature points, $q^-(x_{\iip},y_q)$.

Using $\xi$ as a handle for a generic space variable, we will now sketch the idea behind each one--dimensional sweep aiming at reconstructing an approximation of $q(\xi)$, out of cell averages $\qbar_i$ over uniform cells $C_i:=[\xi_{\imh},\xi_{\iph}]$ of size $\Delta \xi$.
Given $\xi^*\in C_i$, in order to reconstruct $q_h^{WENO}(\xi^*)\approx q(\xi^*)$, one considers 
\begin{itemize}
	\item a high order approximation, $q_h^{HO}(\xi^*)$, of order $2r-1$ (associated to a polynomial degree $2r-2$) constructed out of the cell averages of a big stencil, $\mathcal{S}^{HO}:=\left\lbrace C_{i-(r-1)},\dots,C_{i+(r-1)} \right\rbrace$, of $2r-1$ cells;
	\item $r$ low order approximations, $q_h^{\ell,LO}(\xi^*)$ $\ell=0,\dots,r-1$, of order $r$ (associated to a polynomial degree $r-1$) constructed out of the cell averages of smaller stencils, $\mathcal{S}^{LO}_\ell:=\left\lbrace C_{i-(r-1)+\ell},\dots,C_{i+\ell} \right\rbrace$, of $r$ cells.
\end{itemize}
One then looks for the so--called ``linear weights'', $d_\ell^{\xi^*}$, allowing to express the high order approximation as a linear combination of the low order ones
\begin{equation}
	q_h^{HO}(\xi^*)=\sum_{\ell=0}^{r-1}d_\ell^{\xi^*} q_h^{\ell,LO}(\xi^*).
	\label{eq:linear_weights}
\end{equation}
The final WENO approximation is obtained upon replacement of such weights through ``nonlinear'' ones, $\omega_\ell^{\xi^*}$, able to recover $q_h^{HO}(\xi^*)$ in smooth cases and to select the approximations associated to the smoothest stencils when discontinuities are present,
\begin{equation}
	q_h^{WENO}(\xi^*)=\sum_{\ell=0}^{r-1}\omega_\ell^{\xi^*} q_h^{\ell,LO}(\xi^*).
	\label{eq:WENO_approximation}
\end{equation}
%

\RIIcolor{More specifically, in this work, we employ the classical nonlinear weights from~\cite{jiang1996efficient,shu1998essentially}. 
For each candidate low order stencil $\mathcal{S}^{LO}_\ell$, we have
\begin{equation}
	\omega_\ell^{\xi^*} := \frac{\alpha_\ell^{\xi^*}}{\sum^{r-1}_{k=0}\alpha_k^{\xi^*}},\quad \alpha_\ell^{\xi^*} := \frac{d_\ell^{\xi^*}}{(\beta_\ell+\epsilon_{{\small \text{WENO}}})^2},
\end{equation}
with $\epsilon_{{\small \text{WENO}}}$ being a small constant to prevent divisions by zero, herein set to be $10^{-6}$, and $\beta_\ell$ being the smoothness indicator associated with the stencil
\begin{equation}
	\beta_\ell := \sum_{k=1}^{r-1} \Delta\xi^{2k-1}\int_{\xi_{i-1/2}}^{\xi_{i+1/2}} \left(\frac{d^k}{d\xi^k} q_h^{\ell,LO}(\xi)\right)^2 \diff{\xi},\quad \ell = 0,\ldots,r-1.
\end{equation}
More details can be found in~\cite{ciallella2022arbitrary,micalizzitoro2024,micalizzi2025algorithms,lore_phd_thesis}.}

Notice that the WENO algorithm has been described for a scalar quantity. For systems, the same reconstruction is meant to be applied componentwise.
Actually, many works~\cite{qiu2002construction,miyoshi2020short,peng2019adaptive,ghosh2012compact,micalizzitoro2024,micalizzi2025algorithms} put in evidence how, in the case of systems, the simple reconstruction of conserved variables causes spurious oscillations in high order schemes.
Hence, as there suggested, we apply the reconstruction to characteristic variables.
More in detail, in each one--dimensional sweep along the generic $\xi$-direction, we apply the reconstruction to the components of the vectors $L_{\xi}\ubar_{i,j}$, with $L_{\xi}$ being the matrix of the left eigenvectors of the flux Jacobian $J_{\xi}$ in the $\xi$-direction, thus multiplying the reconstructed variables by $R_{\xi}:=L_{\xi}^{-1}$.
Following~\cite{micalizzitoro2024,micalizzi2025algorithms}, for the local reconstruction in the cell $C_{i,j}$, the matrices are ``frozen'' and fixed for the whole stencil: $L_{\xi}:=L_{\xi}(\ubar_{i,j})$ and $R_{\xi}:=R_{\xi}(\ubar_{i,j})$.
\RIIcolor{A sketch of the procedure to obtain the reconstructed value in the cell $C_{i,j}$ at a single vertical edge quadrature point $(x_{i+\frac12},y_q)$ is reported in Algorithm~\ref{alg:char_WENO_single_point}.
The other reconstructed values, corresponding to right states and horizontal edges, are obtained through suitable adaptations of the same algorithm.}

\begin{algorithm}
	\footnotesize
	\caption{\RIIcolor{Local characteristic WENO reconstruction of order $2r-1$ in the cell $C_{i,j}$ at a single vertical edge quadrature point $(x_{i+\frac12},y_q)$}}
	\label{alg:char_WENO_single_point}
	\begin{algorithmic}[1]
		\Require Cell averages $\overline{\uvec{u}}_{\ell,m}$ on the stencil of $C_{i,j}$, scalar WENO reconstruction of order $2r-1$, quadrature point $(x_{i+\frac12},y_q)$
		\State Compute the matrices $L_{x}:=L_{x}(\overline{\uvec{u}}_{i,j})$ and $R_{x}:=L_{x}^{-1}$ associated with the Jacobian $J_{x}(\overline{\uvec{u}}_{i,j}):=\frac{\partial \uvec{f}}{\partial \uvec{u}}(\overline{\uvec{u}}_{i,j})$
		\For{$m=j-(r-1),\dots,j+(r-1)$}
		\State Project the cell averages along the row $m$ onto local characteristic variables:
		\[
		\uvec{w}_{\ell,m}^{x}:=L_x\overline{\uvec{u}}_{\ell,m},
		\qquad
		\ell=i-(r-1),\dots,i+(r-1).
		\]
		\For{$k=1,\dots,N_c$}
		\State Apply the scalar WENO reconstruction in the $x$-direction to the scalar component $k$, $\{w_{\ell,m}^{x,k}\}_{\ell=i-(r-1)}^{i+(r-1)}$, and obtain
		\[
		w_{m}^{x,k,-}\approx
		\frac{1}{\Delta y}\int_{y_{m-\frac12}}^{y_{m+\frac12}}
		w^{x,k}(x_{i+\frac12}^-,y)\,dy .
		\]
		\EndFor
		\State Collect $\uvec{w}_m^{x,-}:=(w_m^{x,1,-},\dots,w_m^{x,N_c,-})^T$ and project back onto conserved variables:
		\[
		\uvec{u}_m^-:=R_x\uvec{w}_m^{x,-}.
		\]
		\EndFor
		\State Compute the matrices $L_y:=L_y(\overline{\uvec{u}}_{i,j})$ and $R_y:=L_y^{-1}$ associated with the Jacobian $J_y(\overline{\uvec{u}}_{i,j}):=\frac{\partial \uvec{g}}{\partial \uvec{u}}(\overline{\uvec{u}}_{i,j})$
		\State Project the edge averages onto local characteristic variables:
		\[
		\uvec{w}_m^{y,-}:=L_y\uvec{u}_m^-,
		\qquad
		m=j-(r-1),\dots,j+(r-1).
		\]
		\For{$k=1,\dots,N_c$}
		\State Apply the scalar WENO reconstruction in the $y$-direction to the scalar component $k$, $\{w_m^{y,k,-}\}_{m=j-(r-1)}^{j+(r-1)}$, and obtain
		\[
		w^{y,k,-}(x_{i+\frac12},y_q).
		\]
		\EndFor
		\State Collect $\uvec{w}^{y,-}:=(w^{y,1,-},\dots,w^{y,N_c,-})^T$ and project back onto conserved variables:
		\[
		\uvec{u}^{\text{WENO}}(x_{i+\frac12}^-,y_q):=R_y\uvec{w}^{y,-}.
		\]
		\State Return $\uvec{u}^{\text{WENO}}(x_{i+\frac12}^-,y_q)$
	\end{algorithmic}
\end{algorithm}

\subsection{Time discretization}\label{sec:time}

Once fixed the space discretization (in terms of quadrature formulas, space reconstruction and numerical flux), Equations~\eqref{eq:FV_semidiscretization}--\eqref{eq:flux_F}--\eqref{eq:flux_G} constitute a system of ordinary differential equations (ODEs) which must be numerically solved in time.
A time integration scheme is hence necessary to complete the discretization.
Aiming at an arbitrary high order framework, we adopt here a DeC scheme, also considered in~\cite{han2021dec,micalizzi2023new,micalizzitoro2024,micalizzi2025algorithms,lore_phd_thesis}, and we present it for a generic ODE initial value problem of the type
\begin{equation}
	\label{eq:ODE}
	\begin{cases}
		\frac{d}{dt}\uvec{y}(t) = \uvec{G}(t,\uvec{y}(t)),\quad t\in[0,T_f], \\
		\uvec{y}(0)=\uvec{z},
	\end{cases}
\end{equation}
where $\uvec{y}:[0,T_f] \rightarrow \mathbb{R}^{N_c}$ is the unknown solution, with $N_c\in \mathbb{N}^+$ being the number of its components, $\uvec{G}: [0,T_f] \times \R^{N_c} \to \R^{N_c}$ is the right-hand side function, smooth enough to guarantee well--posedness, $T_f\in \mathbb{R}^+$ is the final time, and $\uvec{z} \in \R^{N_c}$ is the initial condition. 
The DeC scheme considered herein is a one--step method, therefore, as customary in such a context, we focus on the generic time interval $[t_n,t_{n+1}]$, where $\Delta t:=t_{n+1}-t_n$, with the main goal of computing $\uvec{y}_{n+1}\approx \uvec{y}(t_{n+1})$ from a known approximation $\uvec{y}_{n}\approx \uvec{y}(t_{n})$.
The scheme is based on introducing $M+1$ subtimenodes $t^m\in[t_n,t_{n+1}]$ $m=0,\dots,M$, with $M$ being proportional to the desired order of accuracy as described in the following, such that $t_n=:t^0<t^1<\dots<t^M:=t_{n+1}$, and in considering an approximation $\uvec{y}^{m}\approx \uvec{y}(t^{m})$ of the exact solution in each of those.
In the first subtimenode, thanks to the assumption of one--step method setting, we can set $\uvec{y}^{0}:=\uvec{y}_{n}$, while in the other ones the approximated solution has to be determined.
In order to do this, we consider the high order implicit discretization of the integral version of the ODE over each time interval $[t^0,t^m]$
\begin{equation}
	\label{eq:apprint}
	\uvec{y}^m=\uvec{y}_n+\Delta t \sum_{\ell=0}^{M} \theta^m_\ell \uvec{G}(t^\ell,\uvec{y}^\ell), \quad m=1,\dots,M,
\end{equation}
where $\theta^m_\ell$ $\ell=0,\dots,M$ are the (normalized) quadrature weights of the quadrature formula associated to all subtimenodes over the time interval $[t^0,t^m]$.
Such an implicit discretization consists in a nonlinear system of equations in the unknown approximated values $\uvec{y}^m$ $m=1,\dots,M$, in fact, corresponding to a high order fully--implicit Runge--Kutta scheme.
Rather than considering a computationally costly direct solution, requiring a nonlinear solver, the DeC scheme employs an explicit fixed--point iterative procedure over the unknown approximated values
\begin{align}\label{eq:DeCODE_Remi}
	\uvec{y}^{m,(p)} := \uvec{y}_n+\Delta t \sum_{\ell=0}^{M} \theta^m_\ell \uvec{G}(t^\ell,\uvec{y}^{\ell,(p-1)}), \quad m=1,\dots,M,\quad p>0,
\end{align}
where the index $m$ refers to the subtimenode, while, the index $(p)$ is the iteration index.
Let us remark that we set $\uvec{y}^{m,(p)}:=\uvec{y}_n$ whenever $m=0$ or $p=0$, meaning that the solution in the first subtimenode is never updated and always kept equal to $\uvec{y}_n$.

The following facts hold~\cite{han2021dec,micalizzitoro2024,lore_phd_thesis}:
\begin{itemize}
	\item The order of accuracy of the implicit discretization~\eqref{eq:apprint}, i.e., the order of accuracy of $\uvec{y}^{M}$ with respect to $\uvec{y}(t_{n+1})$, depends on the number and distribution of the subtimenodes. More in detail, for equispaced ones, we have order $M+1$; for Gauss--Lobatto ones, we have order $2M$.
	\item Each iteration~\eqref{eq:DeCODE_Remi} corresponds to an increase in the order of accuracy with respect to the solution of the implicit discretization~\eqref{eq:apprint}.
\end{itemize}

In view of these facts, assuming Gauss--Lobatto subtimenodes, for desired order $P$, we set $M:=\left \lceil \frac{P}{2}\right \rceil$ and we perform $P$ iterations, finally setting $\uvec{y}_{n+1}:=\uvec{y}^{M,(P)}$.
\RIIcolor{A sketch of the generic DeC iteration is reported in Algorithm~\ref{alg:DeC_time_integration}. Note that the computation of the coefficients $\theta_\ell^m$ can be performed once at the beginning of the simulation.}

\begin{algorithm}
	\footnotesize
	\caption{\RIIcolor{Explicit DeC time integration for \eqref{eq:ODE}}}
	\label{alg:DeC_time_integration}
	\begin{algorithmic}[1]
		\Require Initial approximation $\uvec{y}_n$ at time $t_n$, time step $\Delta t$, desired order $P$
		\State Set $M:=\left\lceil P/2\right\rceil$
		\State Fix $M+1$ Gauss--Lobatto subtimenodes $t^m$, $m=0,\dots,M$ such that
		\[
		t_n=:t^0<t^1<\dots<t^M:=t_{n+1}.
		\]
		\State Construct the Lagrange basis functions $\psi^\ell$, $\ell=0,\dots,M$, associated with the subtimenodes $t^0,\dots,t^M$
		\State Compute the normalized quadrature weights
		\[
		\theta_\ell^m
		:=
		\frac{1}{\Delta t}
		\int_{t^0}^{t^m}
		\psi^\ell(t)\,dt,
		\qquad
		m=1,\dots,M,\quad \ell=0,\dots,M .
		\]
		\State Initialize
		\[
		\uvec{y}^{m,(0)}:=\uvec{y}_n,
		\qquad m=0,\dots,M .
		\]
		\For{$p=1,\dots,P$}
		\State Set $\uvec{y}^{0,(p)}:=\uvec{y}_n$
		\For{$m=1,\dots,M$}
		\State Compute
		\[
		\uvec{y}^{m,(p)}
		:=
		\uvec{y}_n
		+
		\Delta t
		\sum_{\ell=0}^{M}
		\theta_\ell^m
		\uvec{G}\!\left(t^\ell,\uvec{y}^{\ell,(p-1)}\right).
		\]
		\EndFor
		\EndFor
		\State Set $\uvec{y}_{n+1}:=\uvec{y}^{M,(P)}$
	\end{algorithmic}
\end{algorithm}

The described scheme is based on Abgrall's DeC formalism introduced in~\cite{Decremi} to construct high order numerical schemes for Continuous Galerkin finite element discretizations overcoming the burden related to the mass matrix.
Such a formalism provided a precious tool to design arbitrary high order schemes with various applications, e.g., to positivity--preservation~\cite{offner2020arbitrary,ciallella2022arbitrary,ciallella2025high}, to well--balancing~\cite{micalizzi2024novel,ciallella2023arbitrary}, and to conservative schemes for primitive formulations~\cite{abgrall2024staggered}.
Adaptive schemes based on the same formalism have been proposed in~\cite{micalizzi2023efficient,micalizzi2023new,veiga2024improving}.
Connections between DeC, Runge--Kutta and ADER, within Abgrall's formalism, have been established in~\cite{han2021dec,micalizzi2023new,lore_phd_thesis}.
Let us remark that the DeC has been originally introduced in~\cite{fox1949some} and, since then, several DeC formulations and schemes have been proposed across the years; we refer the interested reader to~\cite{boscarino2016error,boscarino2018implicit,hamon2019multi,franco2018multigrid,minion2015interweaving,benedusi2021experimental,liu2008strong,ong2020deferred,ketcheson2014comparison,christlieb2009comments,christlieb2010integral} and references therein.

\section{FORCE-$\alpha$ numerical fluxes}\label{sec:FORCE_alpha}

FORCE-$\alpha$ numerical fluxes, which constitute the main focus of this work, have been introduced in~\cite{toro2020low} and belong to the family of centred numerical fluxes.
Such numerical fluxes do not make any use of the structure of the underlying Riemann problem, only requiring simple flux evaluations for their computation.
Such a simplicity, makes them suitable for complex applications in which (approximate or exact) RSs (and related structures) are unavailable or expensive. 
In particular, focusing on the $x$-direction, FORCE-$\alpha$ numerical fluxes are defined as
\begin{equation}
	\label{eq:FORCE_alpha}
	\widehat{\bF}^{\text{FORCE}-\alpha}({\bu}^L,{\bu}^R) := \frac{1}{2}\left[\widehat{\bF}^{LxF-\alpha}({\bu}^L,{\bu}^R) + \widehat{\bF}^{\text{Richtm}-\alpha}({\bu}^L,{\bu}^R)\right], 
\end{equation}
where 
\begin{equation}
	\label{eq:LxF_alpha}
	\widehat{\bF}^{\text{LxF}-\alpha}({\bu}^L,{\bu}^R) := \frac{1}{2}\left(\bF({\bu}^R) + \bF({\bu}^L)\right) - \frac{1}{2}\frac{\Delta x}{\alpha\Delta t}\left({\bu}^R - {\bu}^L\right),
\end{equation}
\begin{align}
	\label{eq:Richtm_alpha}
	\widehat{\bF}^{\text{Richtm}-\alpha}({\bu}^L,{\bu}^R)&:=\bF(\uvec{u}^*({\bu}^L,{\bu}^R)), \\ \uvec{u}^*({\bu}^L,{\bu}^R)&:=\frac{1}{2}( {\bu}^L+{\bu}^R) - \frac{1}{2}\frac{\alpha\Delta t}{\Delta x}\left(\bF({\bu}^R) - \bF({\bu}^L) \right), \label{eq:ustar_alpha}
\end{align}
with $\alpha\geq 1$ being a real parameter to be fixed.
The definition is analogous for the $y$-direction, with suitable replacement of $\uvec{f}$ and $\Delta x$ by $\uvec{g}$ and $\Delta y$ respectively.
Let us observe that, for $\alpha=1$, $\widehat{\bF}^{\text{FORCE}-\alpha}(\cdot,\cdot)$, $\widehat{\bF}^{\text{LxF}-\alpha}(\cdot,\cdot)$ and $\widehat{\bF}^{\text{Richtm}-\alpha}(\cdot,\cdot)$ yield the FORCE numerical flux~\cite{Toro1996,toro2000centred,chen2003centred}, the Lax--Friedrichs numerical flux~\cite{lax1954weak} and the Richtmyer (or two--step Lax--Wendroff) numerical flux~\cite{richtmyer1967difference} respectively.
For benchmarking purposes, the performance of FORCE-$\alpha$ numerical fluxes will be compared to the one of three upwind numerical fluxes: Rusanov~\cite{Rusanov1961}, HLL~\cite{harten1983upstream} and exact RS~\cite{Godunov}.
All of them are thoroughly described in classical references, e.g., in~\cite{ToroBook}.
For the implementation details, we refer the reader to~\cite{micalizzitoro2024} and we remark that 
\begin{itemize}
	\item in Rusanov, the maximum local wave speed in absolute value is computed through simple Davis' speed estimates~\cite{davis1988simplified} involving the states ${\bu}^L$ and ${\bu}^R$;
	\item in HLL, for the slowest and fastest local wave speeds estimates, we consider the rigorous analytical bounds $\text{TMS}_\text{a}$ from~\cite{toro2020bounds};
	\item for the exact RS, we consider the RS for the Euler equations from~\cite{toro1989fast}, also described in \cite[Chapter 4]{ToroBook}.
\end{itemize}

\subsection{Stability considerations}\label{sec:stability}
Let us remark that stability is strictly related to the adoption of a small enough time step, whose size is proportional to the mesh size over the maximum wave speed in absolute value, according to the Courant--Friedrichs--Lewy (CFL) stability condition.
Namely, we choose the time step size, in one- and two--dimensional experiments respectively, as
\begin{equation}
	\Delta t:= C_{CFL} \frac{\Delta x}{\max{(s^x)}}, \quad
	\Delta t:= C_{CFL} \min{\left(\frac{\Delta x}{\max{(s^x)}},\frac{\Delta y}{\max{(s^y)}}\right)}, 
	\label{eq:CFL}
\end{equation}
with $C_{CFL}$ being a constant, and $s^x$ and $s^y$ being estimates of the maximum local wave speeds in absolute value along the $x$- and $y$-directions, \RIIcolor{where $s^y$ is defined only in the two--dimensional setting}.
In particular, see~\cite{ToroBook}, for first order FV schemes with Rusanov, HLL, exact RS, FORCE, i.e., FORCE--1, and Lax--Friedrichs numerical fluxes, $C_{CFL}\leq 1$ ensures stability in one space dimension, provided that rigorous speed estimates for $s^x$ and $s^y$ are adopted, bounding the actual wave speeds, such as the ones from~\cite{toro2020bounds}.
In multiple space dimensions, the situation is much more delicate.
For example, the linear stability analysis conducted in~\cite{toro2000centred} reveals that FORCE, i.e., FORCE--1, and Lax--Friedrichs numerical fluxes lead to unconditionally unstable first order schemes, even though we remark that increasing the order seems to provide a stabilizing effect, see~\cite{micalizzitoro2024}.
On the other hand, under the assumption of rigorous speed estimates bounding the actual wave speeds, $C_{CFL}\leq 0.5$ is sufficient to guarantee stability in two space dimensions for first order FV schemes with Rusanov, HLL and exact RS numerical fluxes.
The stability of FORCE--$\alpha$ numerical fluxes in one, two and three space dimensions has been investigated in~\cite{toro2020low}, and precise time step constraints have been given.
They are conveniently reported, for different values of $\alpha$, in Table~\ref{tab:CFL_constraints} for one and two space dimensions, which are the cases of interest in this work.
More in detail, in the one--dimensional case, a closed expression is available for the $C_{CFL}$ stability limit, reading
\begin{equation}
	C_{CFL}^{\text{max}}:=\frac{\sqrt{2\alpha-1}}{\alpha}.
\end{equation}

It is worth remarking that, according to the systematic investigation conducted in~\cite{toro2020low}, the stability of FORCE--$\alpha$ numerical fluxes in first order FV schemes tends to decrease, and so also $C_{CFL}^{\text{max}}$, as $\alpha$ increases, on the other hand, also the viscosity decreases.
The first effect implies longer computational times for high values of $\alpha$ in order to guarantee stability.
However, the latter one makes high values of $\alpha$ preferable.
It is thus important to find the optimal balance between these two effects, which is one of the main goals of this investigation.

\begin{table}[h!]
	\centering
	\begin{tabular}{c c c }
		\hline
		$\alpha$ & $1-D$  &  $2-D$  \\
		\hline
		1.0  & 1.000 &  $-$ \\
		2.0  & 0.866 &  0.498 \\
		3.0  & 0.745 &  0.470 \\
		4.0  & 0.662 &  0.433 \\
		5.0  & 0.600 &  0.399 \\
		6.0  & 0.554 &  0.371 \\
		7.0  & 0.516 &  0.348 \\
		8.0  & 0.484 &  0.328 \\
		9.0  & 0.457 &  0.314 \\
		10.0 & 0.435 &  0.299 \\
		\hline
	\end{tabular}
	\caption{$C_{CFL}$ stability limits in one and two space dimensions for first order FV schemes with FORCE-$\alpha$ numerical fluxes}
	\label{tab:CFL_constraints}
\end{table}


\section{Numerical results}\label{sec:numerical_results}
In this section, we will present the results of the numerical tests conducted to assess the properties of FORCE--$\alpha$ numerical fluxes within the WENO--DeC framework.
More in detail, in Section~\ref{sec:Euler_1d}, we will present the results obtained for one--dimensional tests, while, in Section~\ref{sec:Euler_2d}, we will focus on the two--dimensional ones.
Before starting, we report in the following some implementation guidelines adopted in the context of this work.

We will consider FORCE--$\alpha$ numerical fluxes, given by Equations~\eqref{eq:FORCE_alpha}--\eqref{eq:LxF_alpha}--\eqref{eq:Richtm_alpha}--\eqref{eq:ustar_alpha}, for $\alpha=$1, 2, 3, 5 and 10.
As already remarked, valuable comparisons will be made with Rusanov, HLL, and exact RS within the same numerical framework.
The time step will be chosen according to~\eqref{eq:CFL}, where we set $C_{CFL}:=\sigma_{CFL} C_{CFL}^{\text{max}}$, with $\sigma_{CFL}$ being a suitable empirical safety factor which multiplies $C_{CFL}^{\text{max}}$ representing the maximum value of $C_{CFL}$ guaranteeing linear stability.
For FORCE--$\alpha$ numerical fluxes, the values of $C_{CFL}^{\text{max}}$ in one and two space dimensions can be found in the Table~\ref{tab:CFL_constraints}.
Since in two dimensions $\alpha=1$ leads to an unstable scheme, in such a multidimensional context we will consider only FORCE--$\alpha$ numerical fluxes with $\alpha\geq 2$.
For Rusanov, HLL and exact RS, instead, we have $C_{CFL}^{\text{max}}=1$ and 0.5 in one and two space dimensions respectively.
Furthermore, concerning the speed estimates $s^x$ and $s^y$, we consider 
direct Davis' estimates~\cite{davis1988simplified}, i.e., $s^x:=\abs{u}+\sqrt{\gamma \frac{p}{\rho}}$ and $s^y:=\abs{v}+\sqrt{\gamma \frac{p}{\rho}}$, where $\rho,u,v,p$ are obtained from the cell averages at time $t_n$. 
%
We remark that, despite their simplicity, they do not constitute rigorous bounds, see~\cite{toro2020bounds}, and they may require smaller values of $C_{CFL}$ to avoid violating the CFL stability condition.

As stated at the end of Section~\ref{sec:space}, we apply the WENO reconstruction to characteristic variables.
In all simulations, the accuracy of the time discretization is selected in such a way to match the spatial accuracy.
We employ Gauss--Legendre quadrature formulas for the computation of the required integrals, considering the minimal number of points to achieve the desired order of accuracy, with only one exception:
as in~\cite{ciallella2022arbitrary,micalizzitoro2024}, we adopt a four--point Gauss--Legendre quadrature rule for WENO5--DeC5 in the two--dimensional case to avoid negative linear weights occurring for the three--point one.
Negative linear weights can cause stability issues and special care is needed in their presence, see~\cite{shi2002technique}.


\RIIcolor{For the sake of readability, the following compact notations are used in the plots: F-$\alpha$, Rus and Ex.RS denote the FORCE-$\alpha$, Rusanov and exact RS numerical fluxes, respectively, while WD$P$ denotes the WENO--DeC scheme of order $P$.}

\subsection{One--dimensional tests}\label{sec:Euler_1d}
We start with a smooth test from~\cite{toro2005tvd}, involving the advection of a smooth density profile in Section~\ref{sec:Euler_1d_smooth_advection}.
\RIIcolor{We continue, in Section~\ref{sec:Euler_1d_shock_turbulence}, with the shock--turbulence interaction problem from~\cite{titarev2004finite}, which is a modification of the original test introduced in~\cite{shu1989efficient}.
Finally, we consider some Riemann problems from~\cite{ToroBook}, specifically designed to test the robustness of numerical methods.}

\subsubsection{Advection of smooth density profile}\label{sec:Euler_1d_smooth_advection}
This test, taken from~\cite{toro2005tvd}, is meant to verify the order of accuracy of our one--dimensional implementation, and to establish a first performance evaluation over a smooth problem.
The initial condition is prescribed as
\begin{align}
	\begin{cases}
		\rho(x,0)&:=2+\sin^4{\left(\pi x\right)},\\
		u(x,0)&:=u_\infty, \\
		p(x,0)&:=p_\infty,
	\end{cases}
\end{align}
where $u_\infty:=1$ and $p_\infty:=1,$ on the computational domain $\Omega:=[-1,1]$ with periodic boundary conditions.
The test corresponds to the advection of the density profile with speed $u_\infty$, and the exact solution is $\uvec{u}(x,t):=\uvec{u}(x-u_{\infty}t,0)$.
We run our simulations until the final time $T_f:=2$, considering a safety CFL coefficient $\sigma_{CFL}:=0.9$. 

\begin{figure}[htbp]
	\centering
	
	\begin{subfigure}[t]{0.49\textwidth}
		\centering
		\includegraphics[width=\textwidth, trim={0 0 140 0}, clip]{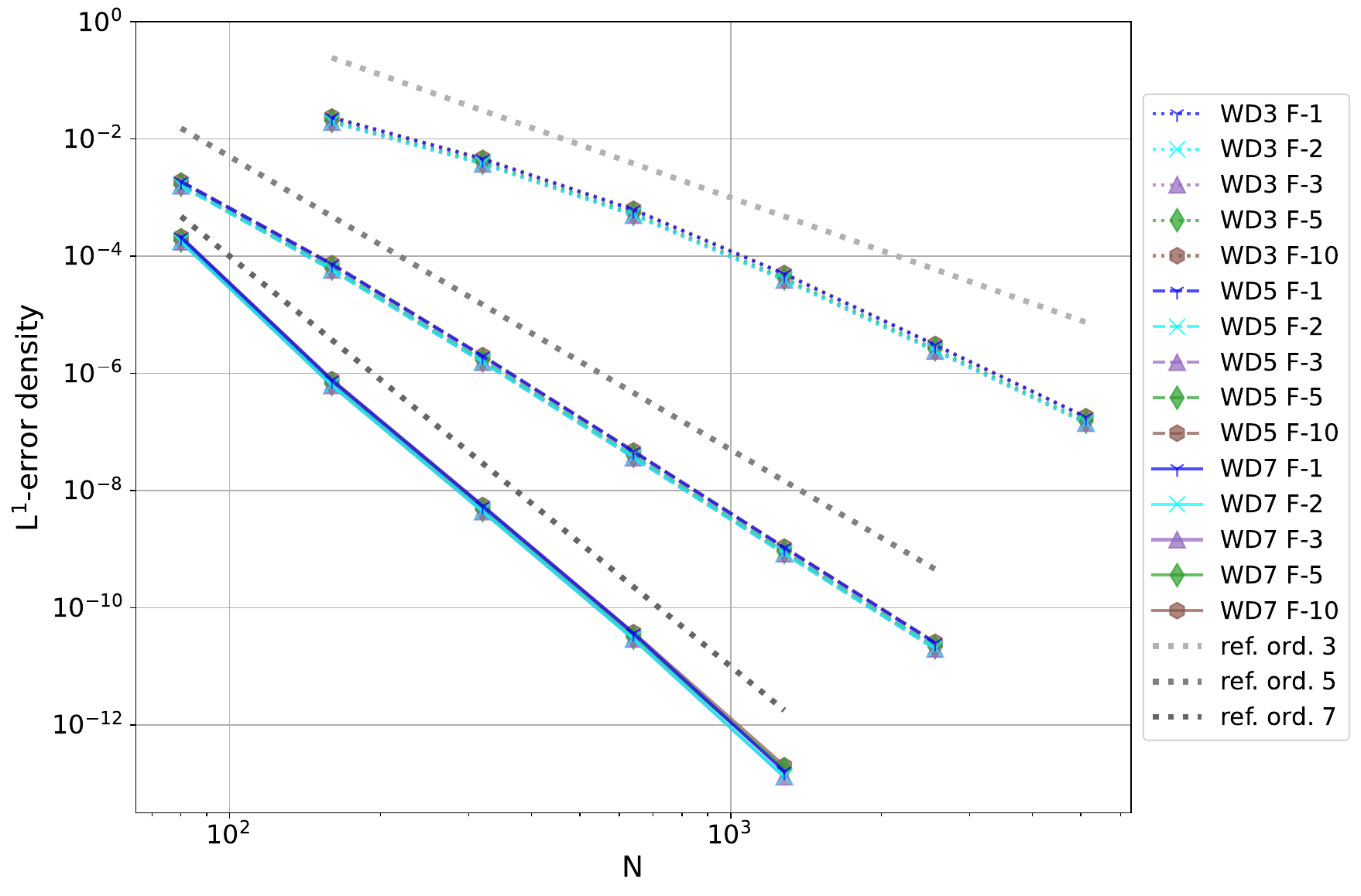}
		\caption{Convergence analysis}
	\end{subfigure}
	\hfill
	\begin{subfigure}[t]{0.49\textwidth}
		\centering
		\includegraphics[width=\textwidth, trim={0 0 150 0}, clip]{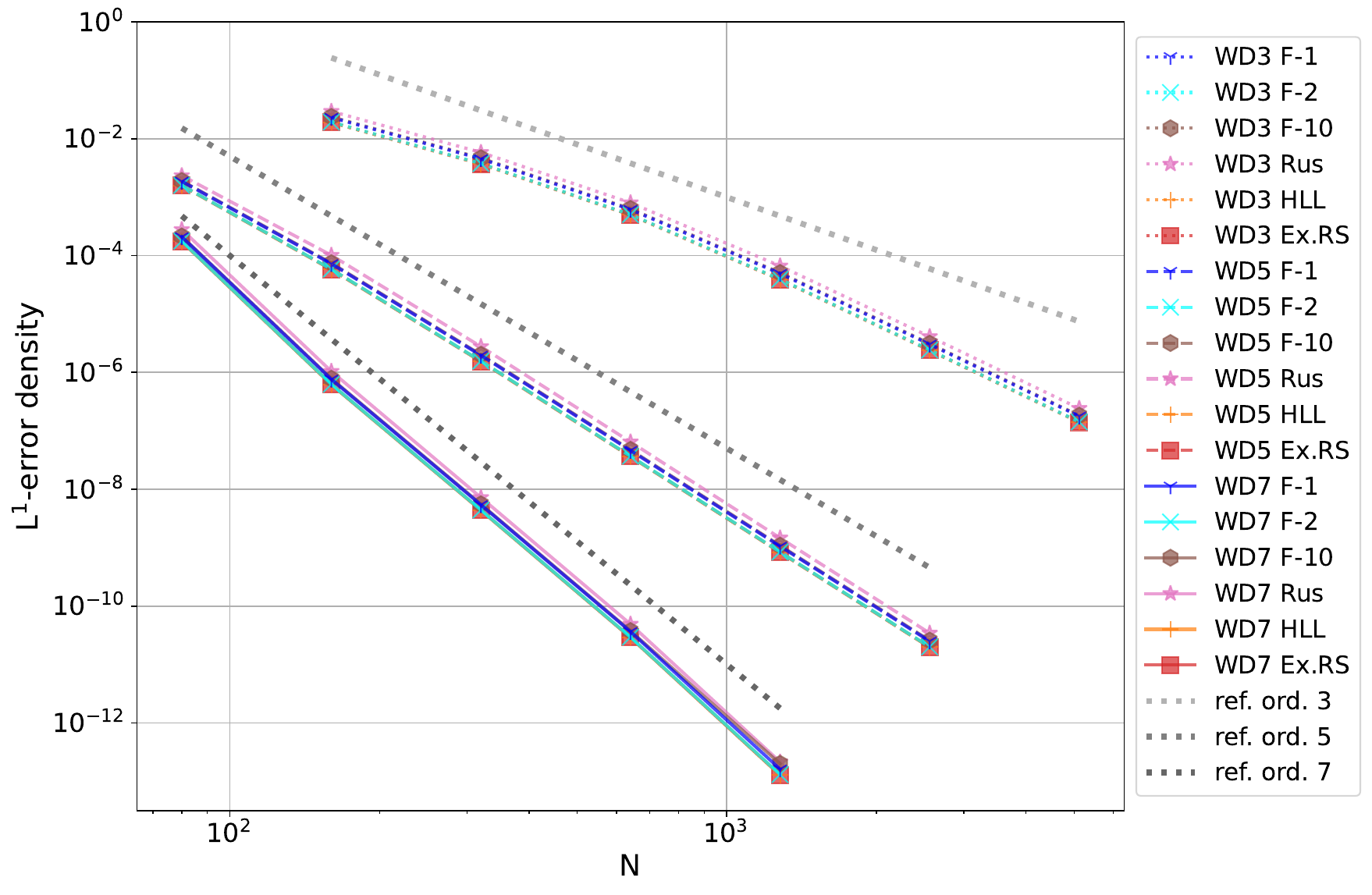}
		\caption{Convergence comparison with upwind fluxes}
	\end{subfigure}
	
	\vspace{0.3cm}
	
	\begin{subfigure}[t]{0.49\textwidth}
		\centering
		\includegraphics[width=\textwidth, trim={0 0 140 0}, clip]{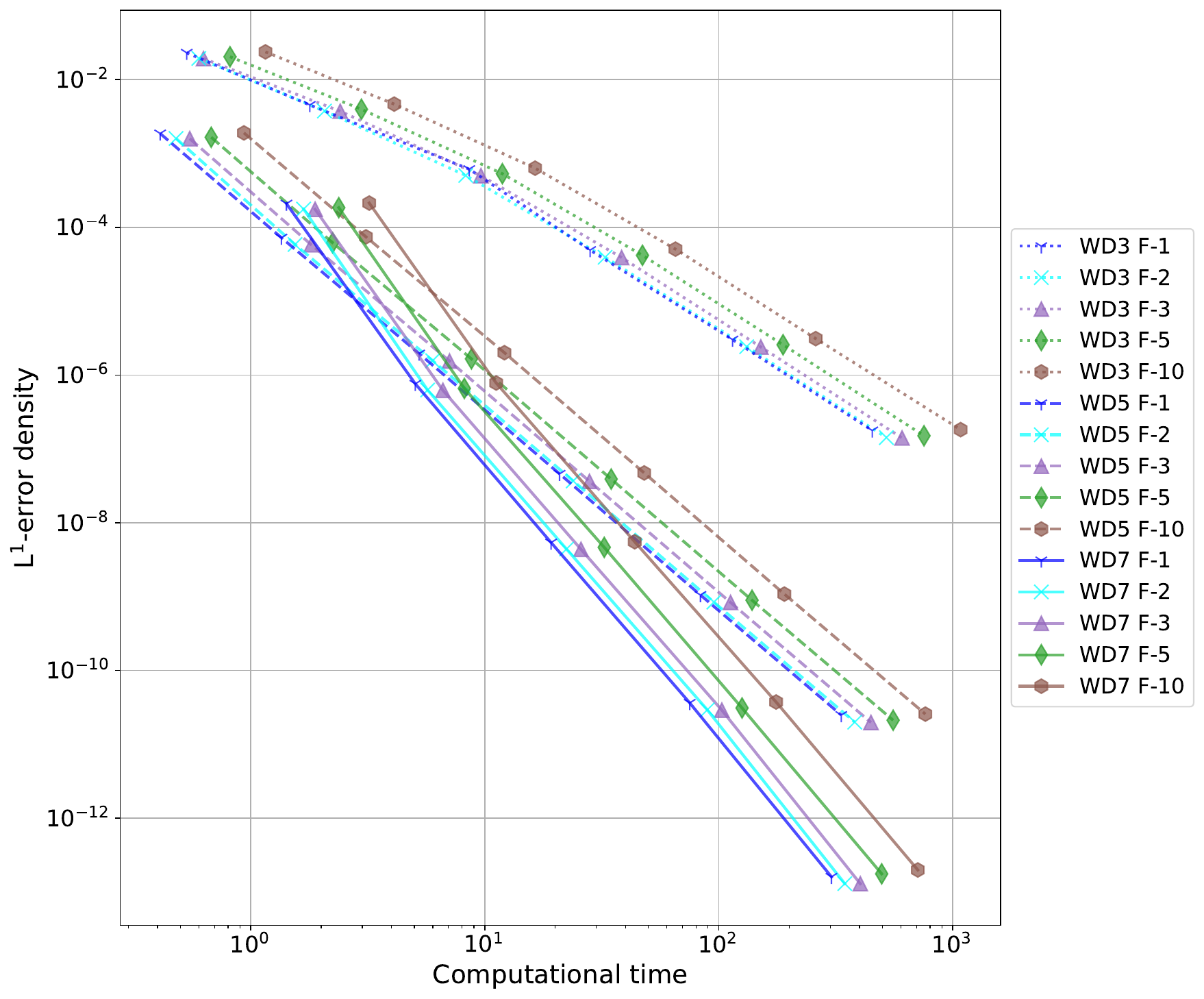}
		\caption{Efficiency analysis}
	\end{subfigure}
	\hfill
	\begin{subfigure}[t]{0.49\textwidth}
		\centering
		\includegraphics[width=\textwidth, trim={0 0 150 0}, clip]{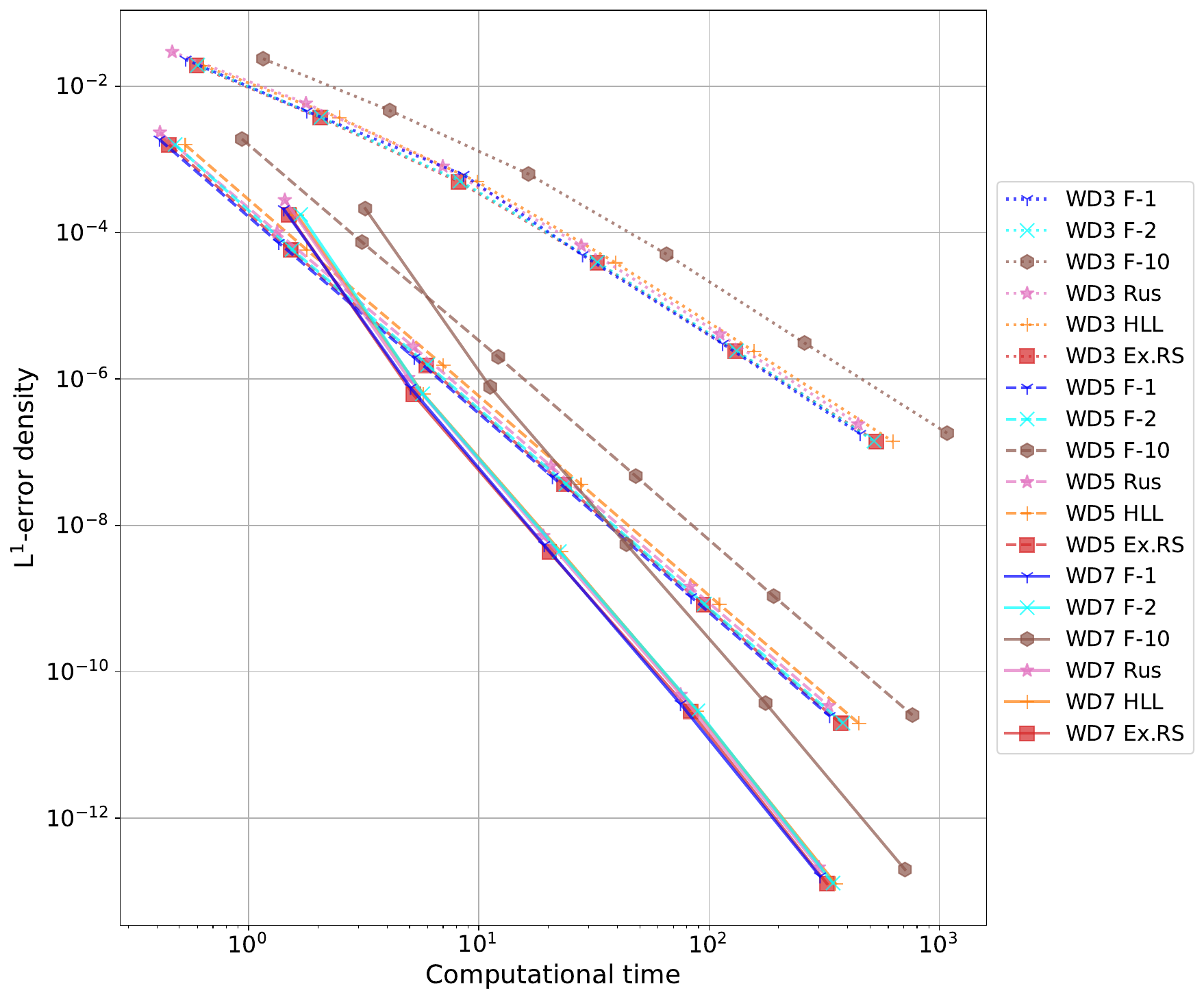}
		\caption{Efficiency comparison with upwind fluxes}
	\end{subfigure}

	\vspace{0.3cm}

	\begin{subfigure}[t]{1.0\textwidth}
	\centering
	\includegraphics[width=0.9\textwidth]{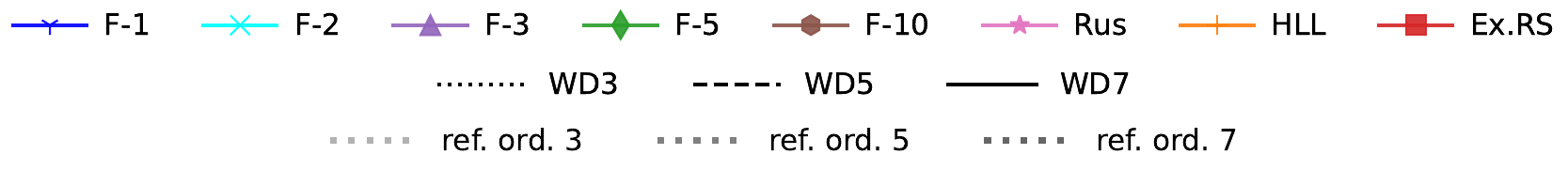}
	\end{subfigure}

	\caption{\RIcolor{Advection of smooth density profile: Convergence and efficiency analyses. The left column illustrates the influence of the parameter $\alpha$ within the FORCE-$\alpha$ family, while the right column compares representative FORCE-$\alpha$ numerical fluxes ($\alpha=1$, $2$ and $10$) with standard upwind fluxes (Rusanov, HLL and exact RS).  On the bottom a common legend for numerical fluxes (distinguished by color and marker) and order (distinguished by linestyle) is reported.}}
	\label{fig:Euler_1d_sin4_WENODeC}
\end{figure}

\begin{figure}[htbp]
	\centering
	\includegraphics[width=1.0\linewidth]{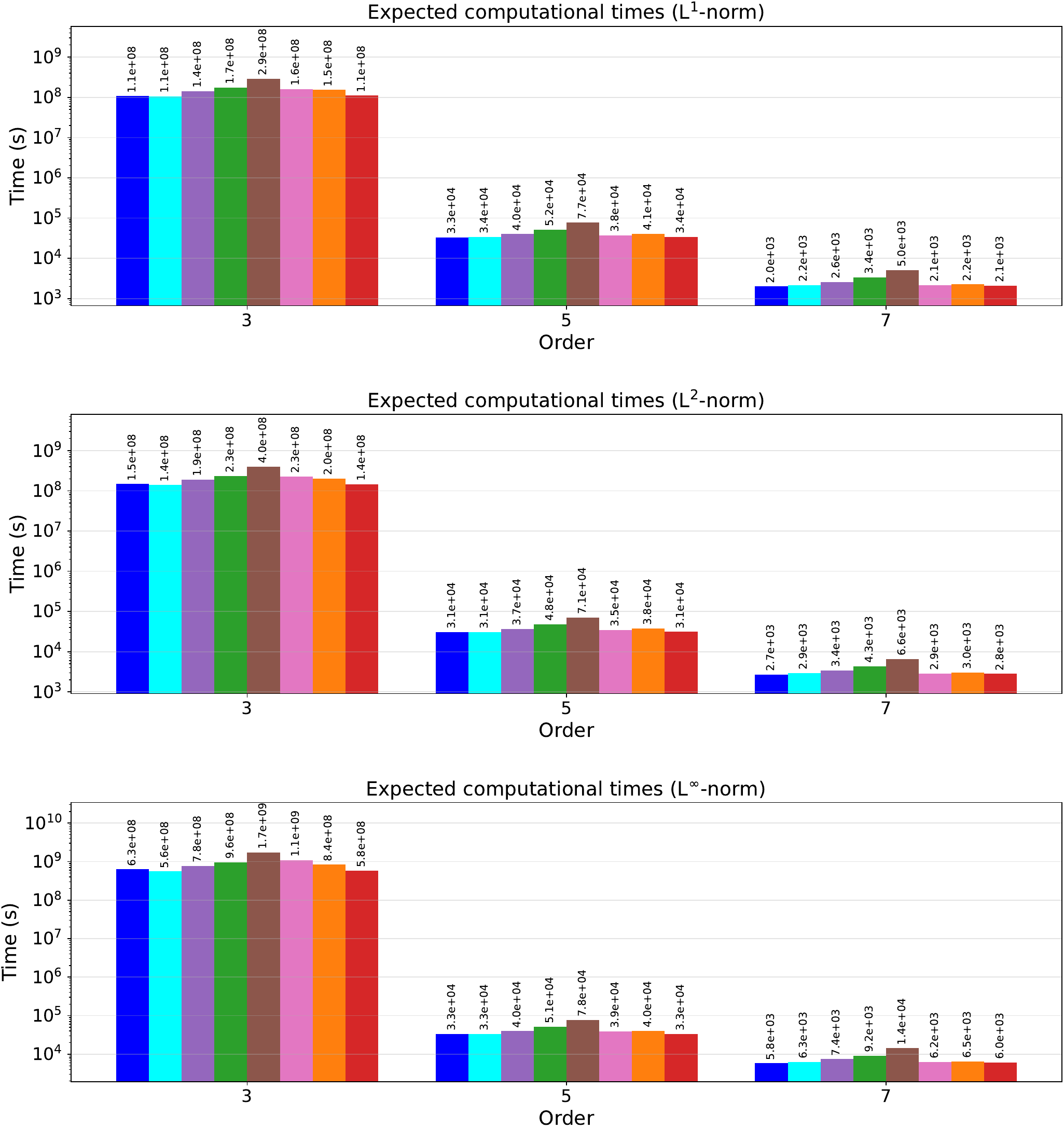}\\
	\includegraphics[width=1.0\linewidth]{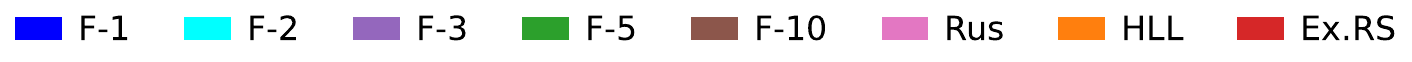}	
	\caption{\RIcolor{Advection of smooth density profile: Estimated computational times (in seconds) required to achieve a density error of $10^{-16}$ in the $L^1$-, $L^2$- and $L^\infty$-norms for the different numerical fluxes and orders considered. The estimates are obtained from the efficiency curves by linear regression in logarithmic scale. A common legend for the numerical fluxes (distinguished by color) is reported at the bottom.}}
	\label{fig:expected_time_Euler_1d}
\end{figure}

The results of the convergence analysis for the density are reported in Tables~\ref{tab:Euler_1d_convergence_tables_WENO_DeC_char_order3},~\ref{tab:Euler_1d_convergence_tables_WENO_DeC_char_order5}, and~\ref{tab:Euler_1d_convergence_tables_WENO_DeC_char_order7}.
As one can see, the expected order of accuracy is obtained for all orders and all numerical fluxes in the three norms $L^1$, $L^2$ and $L^{\infty}$.
Same results have been obtained with respect to the other variables and are hence omitted.
For convenience, the results of the convergence analysis in the $L^1$--norm are also graphically displayed in Figure~\ref{fig:Euler_1d_sin4_WENODeC} along with the efficiency analysis in the same norm.
\RIIcolor{In particular, the plots on the left column consider only the FORCE-$\alpha$ numerical fluxes, while, in the plots on the right, we report comparisons between FORCE-$\alpha$ fluxes with references values $\alpha=1$, 2 and 10 with Rusanov, HLL and exact RS.}
No big differences among the numerical fluxes can be appreciated from the convergence analysis in this test: namely, for a given refinement, all numerical fluxes yield similar results, and the performance of FORCE--$\alpha$ numerical fluxes is comparable to the one of Rusanov, HLL, and exact RS. Interestingly, Rusanov has a higher level of diffusion with respect to FORCE--$\alpha$ for any considered value of $\alpha$, despite being an upwind numerical flux. This is in line with what obtained in~\cite{micalizzitoro2024}, where only FORCE, i.e., FORCE--1, has been considered.
\RIIcolor{Focusing on the efficiency analysis plots at the bottom}, we can see that FORCE--1 is the most efficient numerical flux.
It is followed by FORCE-2, exact RS, HLL and Rusanov, whose performance is difficult to be distinguished.
The least efficient numerical fluxes are indeed FORCE-5 and FORCE-10.
%
%
The performance of FORCE-3 is intermediate between the two aforementioned groups of most and least efficient numerical fluxes.
Let us remark that, coherently with what obtained in~\cite{micalizzitoro2024}, differences among numerical fluxes tend to decrease for increasing order, and for order 7 the results of FORCE-1, FORCE-2, exact RS, HLL and Rusanov in the efficiency plot are practically indistinguishable.
To provide a quantitative comparison, we report in Figure~\ref{fig:expected_time_Euler_1d} the estimated computational times needed to obtain an accuracy tolerance of $10^{-16}$ in the three considered norms, computed through the linear regression in logarithmic scale of the curves for each order and numerical flux from the efficiency plots.
Results confirm what previously described, leading to the following conclusion: in the considered test, the employment of FORCE-1 and FORCE-2 is a competitive alternative to upwind numerical fluxes, especially for high order.

\begin{table}[htbp]
	\centering
	\caption{Advection of smooth density profile: convergence tables for order 3}
	\label{tab:Euler_1d_convergence_tables_WENO_DeC_char_order3}
	\scalebox{0.65}{ 
		\begin{tabular}{c c c c c c c c}
			\toprule
			\multirow{2}{*}{$N$} & \multicolumn{2}{c}{$L^1$-error $\rho$} & \multicolumn{2}{c}{$L^2$-error $\rho$} & \multicolumn{2}{c}{$L^{\infty}$-error $\rho$} & \multirow{2}{*}{CPU Time} \\
			\cmidrule(lr){2-3} \cmidrule(lr){4-5} \cmidrule(lr){6-7}
			& Error & Order & Error & Order & Error & Order & \\
			\midrule
			
			\multicolumn{8}{c}{\textbf{FORCE-1}} \\ 
			\midrule
			160  &   2.313e-02  &  $-$  &   2.633e-02  &  $-$  &   5.639e-02  &  $-$  &   5.328e-01 \\ 
			320  &   4.541e-03  &  2.349  &   6.829e-03  &  1.947  &   1.879e-02  &  1.585  &   1.789e+00 \\ 
			640  &   6.148e-04  &  2.885  &   1.172e-03  &  2.542  &   4.346e-03  &  2.112  &   8.586e+00 \\ 
			1280  &   4.924e-05  &  3.642  &   9.708e-05  &  3.594  &   4.781e-04  &  3.184  &   2.821e+01 \\ 
			2560  &   3.028e-06  &  4.024  &   4.708e-06  &  4.366  &   2.033e-05  &  4.555  &   1.145e+02 \\ 
			5120  &   1.778e-07  &  4.090  &   2.155e-07  &  4.449  &   6.806e-07  &  4.901  &   4.532e+02 \\ 
			\midrule

			\multicolumn{8}{c}{\textbf{FORCE-2}} \\ 
			\midrule
			160  &   1.941e-02  &  $-$  &   2.266e-02  &  $-$  &   4.991e-02  &  $-$  &   6.001e-01 \\ 
			320  &   3.768e-03  &  2.365  &   5.834e-03  &  1.957  &   1.658e-02  &  1.590  &   2.067e+00 \\ 
			640  &   5.048e-04  &  2.900  &   9.850e-04  &  2.566  &   3.771e-03  &  2.136  &   8.297e+00 \\ 
			1280  &   3.940e-05  &  3.680  &   7.860e-05  &  3.648  &   3.961e-04  &  3.251  &   3.270e+01 \\ 
			2560  &   2.425e-06  &  4.022  &   3.772e-06  &  4.381  &   1.632e-05  &  4.601  &   1.317e+02 \\ 
			5120  &   1.424e-07  &  4.090  &   1.725e-07  &  4.450  &   5.448e-07  &  4.905  &   5.206e+02 \\ 
			\midrule

			\multicolumn{8}{c}{\textbf{FORCE-3}} \\ 
			\midrule
			160  &   1.929e-02  &  $-$  &   2.255e-02  &  $-$  &   4.970e-02  &  $-$  &   6.284e-01 \\ 
			320  &   3.744e-03  &  2.365  &   5.804e-03  &  1.958  &   1.650e-02  &  1.591  &   2.413e+00 \\ 
			640  &   5.015e-04  &  2.900  &   9.792e-04  &  2.567  &   3.752e-03  &  2.137  &   9.635e+00 \\ 
			1280  &   3.906e-05  &  3.683  &   7.794e-05  &  3.651  &   3.932e-04  &  3.254  &   3.847e+01 \\ 
			2560  &   2.403e-06  &  4.022  &   3.737e-06  &  4.383  &   1.617e-05  &  4.604  &   1.513e+02 \\ 
			5120  &   1.410e-07  &  4.091  &   1.708e-07  &  4.451  &   5.393e-07  &  4.906  &   6.074e+02 \\ 
			\midrule

			\multicolumn{8}{c}{\textbf{FORCE-5}} \\ 
			\midrule
			160  &   2.043e-02  &  $-$  &   2.366e-02  &  $-$  &   5.167e-02  &  $-$  &   8.158e-01 \\ 
			320  &   3.975e-03  &  2.362  &   6.101e-03  &  1.955  &   1.717e-02  &  1.589  &   2.972e+00 \\ 
			640  &   5.332e-04  &  2.898  &   1.035e-03  &  2.560  &   3.923e-03  &  2.130  &   1.190e+01 \\ 
			1280  &   4.185e-05  &  3.671  &   8.322e-05  &  3.636  &   4.170e-04  &  3.234  &   4.723e+01 \\ 
			2560  &   2.574e-06  &  4.023  &   4.000e-06  &  4.379  &   1.730e-05  &  4.592  &   1.884e+02 \\ 
			5120  &   1.508e-07  &  4.093  &   1.827e-07  &  4.452  &   5.768e-07  &  4.906  &   7.538e+02 \\ 
			\midrule

			\multicolumn{8}{c}{\textbf{FORCE-10}} \\ 
			\midrule
			160  &   2.380e-02  &  $-$  &   2.698e-02  &  $-$  &   5.747e-02  &  $-$  &   1.156e+00 \\ 
			320  &   4.671e-03  &  2.349  &   6.998e-03  &  1.947  &   1.915e-02  &  1.586  &   4.104e+00 \\ 
			640  &   6.338e-04  &  2.882  &   1.204e-03  &  2.539  &   4.439e-03  &  2.109  &   1.642e+01 \\ 
			1280  &   5.077e-05  &  3.642  &   9.990e-05  &  3.591  &   4.908e-04  &  3.177  &   6.535e+01 \\ 
			2560  &   3.118e-06  &  4.025  &   4.842e-06  &  4.367  &   2.090e-05  &  4.553  &   2.597e+02 \\ 
			5120  &   1.826e-07  &  4.094  &   2.212e-07  &  4.452  &   6.981e-07  &  4.904  &   1.081e+03 \\ 
			\midrule

			\multicolumn{8}{c}{\textbf{Rusanov}} \\ 
			\midrule
			160  &   2.950e-02  &  $-$  &   3.208e-02  &  $-$  &   6.616e-02  &  $-$  &   4.658e-01 \\ 
			320  &   5.851e-03  &  2.334  &   8.378e-03  &  1.937  &   2.211e-02  &  1.581  &   1.775e+00 \\ 
			640  &   8.029e-04  &  2.865  &   1.470e-03  &  2.511  &   5.223e-03  &  2.082  &   6.975e+00 \\ 
			1280  &   6.648e-05  &  3.594  &   1.274e-04  &  3.528  &   6.080e-04  &  3.103  &   2.788e+01 \\ 
			2560  &   4.105e-06  &  4.018  &   6.272e-06  &  4.344  &   2.695e-05  &  4.495  &   1.112e+02 \\ 
			5120  &   2.422e-07  &  4.083  &   2.882e-07  &  4.444  &   9.045e-07  &  4.897  &   4.440e+02 \\ 
			\midrule

			\multicolumn{8}{c}{\textbf{HLL}} \\ 
			\midrule
			160  &   1.911e-02  &  $-$  &   2.237e-02  &  $-$  &   4.942e-02  &  $-$  &   6.374e-01 \\ 
			320  &   3.711e-03  &  2.364  &   5.762e-03  &  1.957  &   1.642e-02  &  1.590  &   2.483e+00 \\ 
			640  &   4.975e-04  &  2.899  &   9.722e-04  &  2.567  &   3.732e-03  &  2.137  &   9.839e+00 \\ 
			1280  &   3.883e-05  &  3.679  &   7.754e-05  &  3.648  &   3.911e-04  &  3.254  &   3.927e+01 \\ 
			2560  &   2.391e-06  &  4.021  &   3.722e-06  &  4.381  &   1.611e-05  &  4.602  &   1.569e+02 \\ 
			5120  &   1.406e-07  &  4.088  &   1.704e-07  &  4.449  &   5.383e-07  &  4.903  &   6.283e+02 \\ 
			\midrule

			\multicolumn{8}{c}{\textbf{exact RS}} \\ 
			\midrule
			160  &   1.911e-02  &  $-$  &   2.237e-02  &  $-$  &   4.942e-02  &  $-$  &   5.941e-01 \\ 
			320  &   3.711e-03  &  2.364  &   5.762e-03  &  1.957  &   1.642e-02  &  1.590  &   2.050e+00 \\ 
			640  &   4.975e-04  &  2.899  &   9.722e-04  &  2.567  &   3.732e-03  &  2.137  &   8.164e+00 \\ 
			1280  &   3.883e-05  &  3.679  &   7.754e-05  &  3.648  &   3.911e-04  &  3.254  &   3.266e+01 \\ 
			2560  &   2.391e-06  &  4.021  &   3.722e-06  &  4.381  &   1.611e-05  &  4.602  &   1.303e+02 \\ 
			5120  &   1.406e-07  &  4.088  &   1.704e-07  &  4.449  &   5.383e-07  &  4.903  &   5.325e+02 \\ 
			\midrule

			\bottomrule
	\end{tabular}}
\end{table}

\begin{table}[htbp]
	\centering
	\caption{Advection of smooth density profile: convergence tables for order 5}
	\label{tab:Euler_1d_convergence_tables_WENO_DeC_char_order5}
	\scalebox{0.65}{ 
		\begin{tabular}{c c c c c c c c}
			\toprule
			\multirow{2}{*}{$N$} & \multicolumn{2}{c}{$L^1$-error $\rho$} & \multicolumn{2}{c}{$L^2$-error $\rho$} & \multicolumn{2}{c}{$L^{\infty}$-error $\rho$} & \multirow{2}{*}{CPU Time} \\
			\cmidrule(lr){2-3} \cmidrule(lr){4-5} \cmidrule(lr){6-7}
			& Error & Order & Error & Order & Error & Order & \\
			\midrule
			
			\multicolumn{8}{c}{\textbf{FORCE-1}} \\ 
			\midrule
			80  &   1.869e-03  &  $-$  &   1.576e-03  &  $-$  &   2.533e-03  &  $-$  &   4.102e-01 \\ 
			160  &   7.248e-05  &  4.689  &   7.189e-05  &  4.454  &   1.475e-04  &  4.103  &   1.353e+00 \\ 
			320  &   1.956e-06  &  5.212  &   1.891e-06  &  5.249  &   4.255e-06  &  5.115  &   5.244e+00 \\ 
			640  &   4.602e-08  &  5.409  &   3.978e-08  &  5.571  &   7.233e-08  &  5.878  &   2.089e+01 \\ 
			1280  &   1.053e-09  &  5.449  &   8.634e-10  &  5.526  &   1.169e-09  &  5.952  &   8.357e+01 \\ 
			2560  &   2.496e-11  &  5.399  &   2.186e-11  &  5.304  &   3.264e-11  &  5.162  &   3.341e+02 \\ 
			\midrule

			\multicolumn{8}{c}{\textbf{FORCE-2}} \\ 
			\midrule
			80  &   1.601e-03  &  $-$  &   1.382e-03  &  $-$  &   2.400e-03  &  $-$  &   4.794e-01 \\ 
			160  &   5.873e-05  &  4.768  &   5.859e-05  &  4.559  &   1.228e-04  &  4.288  &   1.547e+00 \\ 
			320  &   1.566e-06  &  5.229  &   1.521e-06  &  5.268  &   3.478e-06  &  5.142  &   5.992e+00 \\ 
			640  &   3.687e-08  &  5.409  &   3.188e-08  &  5.576  &   5.804e-08  &  5.905  &   2.390e+01 \\ 
			1280  &   8.438e-10  &  5.449  &   6.915e-10  &  5.527  &   9.362e-10  &  5.954  &   9.533e+01 \\ 
			2560  &   2.000e-11  &  5.399  &   1.751e-11  &  5.304  &   2.623e-11  &  5.158  &   3.814e+02 \\ 
			\midrule

			\multicolumn{8}{c}{\textbf{FORCE-3}} \\ 
			\midrule
			80  &   1.592e-03  &  $-$  &   1.375e-03  &  $-$  &   2.395e-03  &  $-$  &   5.489e-01 \\ 
			160  &   5.821e-05  &  4.774  &   5.807e-05  &  4.566  &   1.218e-04  &  4.297  &   1.815e+00 \\ 
			320  &   1.552e-06  &  5.229  &   1.507e-06  &  5.268  &   3.448e-06  &  5.143  &   7.071e+00 \\ 
			640  &   3.653e-08  &  5.409  &   3.158e-08  &  5.576  &   5.751e-08  &  5.906  &   2.803e+01 \\ 
			1280  &   8.360e-10  &  5.449  &   6.852e-10  &  5.526  &   9.276e-10  &  5.954  &   1.123e+02 \\ 
			2560  &   1.982e-11  &  5.399  &   1.735e-11  &  5.304  &   2.592e-11  &  5.162  &   4.473e+02 \\ 
			\midrule

			\multicolumn{8}{c}{\textbf{FORCE-5}} \\ 
			\midrule
			80  &   1.661e-03  &  $-$  &   1.428e-03  &  $-$  &   2.438e-03  &  $-$  &   6.779e-01 \\ 
			160  &   6.197e-05  &  4.744  &   6.177e-05  &  4.531  &   1.288e-04  &  4.243  &   2.239e+00 \\ 
			320  &   1.659e-06  &  5.223  &   1.610e-06  &  5.262  &   3.666e-06  &  5.134  &   8.776e+00 \\ 
			640  &   3.910e-08  &  5.407  &   3.380e-08  &  5.574  &   6.149e-08  &  5.898  &   3.471e+01 \\ 
			1280  &   8.951e-10  &  5.449  &   7.336e-10  &  5.526  &   9.931e-10  &  5.952  &   1.388e+02 \\ 
			2560  &   2.122e-11  &  5.398  &   1.858e-11  &  5.303  &   2.795e-11  &  5.151  &   5.566e+02 \\ 
			\midrule

			\multicolumn{8}{c}{\textbf{FORCE-10}} \\ 
			\midrule
			80  &   1.903e-03  &  $-$  &   1.599e-03  &  $-$  &   2.544e-03  &  $-$  &   9.355e-01 \\ 
			160  &   7.442e-05  &  4.676  &   7.376e-05  &  4.439  &   1.508e-04  &  4.076  &   3.108e+00 \\ 
			320  &   2.010e-06  &  5.210  &   1.943e-06  &  5.247  &   4.362e-06  &  5.112  &   1.214e+01 \\ 
			640  &   4.734e-08  &  5.408  &   4.091e-08  &  5.569  &   7.437e-08  &  5.874  &   4.804e+01 \\ 
			1280  &   1.084e-09  &  5.449  &   8.883e-10  &  5.525  &   1.202e-09  &  5.951  &   1.908e+02 \\ 
			2560  &   2.570e-11  &  5.399  &   2.250e-11  &  5.303  &   3.371e-11  &  5.157  &   7.642e+02 \\ 
			\midrule

			\multicolumn{8}{c}{\textbf{Rusanov}} \\ 
			\midrule
			80  &   2.347e-03  &  $-$  &   1.954e-03  &  $-$  &   2.660e-03  &  $-$  &   4.120e-01 \\ 
			160  &   1.013e-04  &  4.534  &   1.015e-04  &  4.267  &   2.013e-04  &  3.724  &   1.330e+00 \\ 
			320  &   2.769e-06  &  5.193  &   2.697e-06  &  5.234  &   6.033e-06  &  5.060  &   5.192e+00 \\ 
			640  &   6.457e-08  &  5.422  &   5.632e-08  &  5.582  &   1.048e-07  &  5.847  &   2.058e+01 \\ 
			1280  &   1.463e-09  &  5.464  &   1.194e-09  &  5.559  &   1.557e-09  &  6.073  &   8.268e+01 \\ 
			2560  &   3.434e-11  &  5.413  &   2.980e-11  &  5.325  &   4.355e-11  &  5.159  &   3.306e+02 \\ 
			\midrule

			\multicolumn{8}{c}{\textbf{HLL}} \\ 
			\midrule
			80  &   1.586e-03  &  $-$  &   1.371e-03  &  $-$  &   2.390e-03  &  $-$  &   5.305e-01 \\ 
			160  &   5.795e-05  &  4.775  &   5.781e-05  &  4.567  &   1.214e-04  &  4.299  &   1.794e+00 \\ 
			320  &   1.545e-06  &  5.229  &   1.500e-06  &  5.268  &   3.434e-06  &  5.144  &   7.002e+00 \\ 
			640  &   3.637e-08  &  5.409  &   3.145e-08  &  5.576  &   5.727e-08  &  5.906  &   2.780e+01 \\ 
			1280  &   8.323e-10  &  5.449  &   6.822e-10  &  5.527  &   9.235e-10  &  5.955  &   1.110e+02 \\ 
			2560  &   1.973e-11  &  5.399  &   1.727e-11  &  5.304  &   2.577e-11  &  5.163  &   4.479e+02 \\ 
			\midrule

			\multicolumn{8}{c}{\textbf{exact RS}} \\ 
			\midrule
			80  &   1.586e-03  &  $-$  &   1.371e-03  &  $-$  &   2.390e-03  &  $-$  &   4.514e-01 \\ 
			160  &   5.795e-05  &  4.775  &   5.781e-05  &  4.567  &   1.214e-04  &  4.299  &   1.520e+00 \\ 
			320  &   1.545e-06  &  5.229  &   1.500e-06  &  5.268  &   3.434e-06  &  5.144  &   5.896e+00 \\ 
			640  &   3.637e-08  &  5.409  &   3.145e-08  &  5.576  &   5.727e-08  &  5.906  &   2.340e+01 \\ 
			1280  &   8.323e-10  &  5.449  &   6.822e-10  &  5.527  &   9.235e-10  &  5.955  &   9.443e+01 \\ 
			2560  &   1.973e-11  &  5.399  &   1.727e-11  &  5.304  &   2.577e-11  &  5.163  &   3.727e+02 \\ 
			\midrule

			\bottomrule
	\end{tabular}}
\end{table}

\begin{table}[htbp]
	\centering
	\caption{Advection of smooth density profile: convergence tables for order 7}
	\label{tab:Euler_1d_convergence_tables_WENO_DeC_char_order7}
	\scalebox{0.65}{ 
		\begin{tabular}{c c c c c c c c}
			\toprule
			\multirow{2}{*}{$N$} & \multicolumn{2}{c}{$L^1$-error $\rho$} & \multicolumn{2}{c}{$L^2$-error $\rho$} & \multicolumn{2}{c}{$L^{\infty}$-error $\rho$} & \multirow{2}{*}{CPU Time} \\
			\cmidrule(lr){2-3} \cmidrule(lr){4-5} \cmidrule(lr){6-7}
			& Error & Order & Error & Order & Error & Order & \\
			\midrule
			
			\multicolumn{8}{c}{\textbf{FORCE-1}} \\ 
			\midrule
			80  &   2.096e-04  &  $-$  &   2.281e-04  &  $-$  &   4.032e-04  &  $-$  &   1.422e+00 \\ 
			160  &   7.620e-07  &  8.103  &   9.008e-07  &  7.984  &   2.055e-06  &  7.617  &   5.053e+00 \\ 
			320  &   5.403e-09  &  7.140  &   9.927e-09  &  6.504  &   5.034e-08  &  5.351  &   1.922e+01 \\ 
			640  &   3.646e-11  &  7.211  &   9.056e-11  &  6.776  &   6.100e-10  &  6.367  &   7.519e+01 \\ 
			1280  &   1.609e-13  &  7.824  &   3.131e-13  &  8.176  &   1.860e-12  &  8.357  &   3.033e+02 \\ 
			\midrule

			\multicolumn{8}{c}{\textbf{FORCE-2}} \\ 
			\midrule
			80  &   1.770e-04  &  $-$  &   1.938e-04  &  $-$  &   3.445e-04  &  $-$  &   1.684e+00 \\ 
			160  &   6.277e-07  &  8.139  &   7.531e-07  &  8.008  &   1.754e-06  &  7.618  &   5.701e+00 \\ 
			320  &   4.413e-09  &  7.152  &   8.249e-09  &  6.512  &   4.263e-08  &  5.362  &   2.236e+01 \\ 
			640  &   2.936e-11  &  7.232  &   7.377e-11  &  6.805  &   5.043e-10  &  6.402  &   8.938e+01 \\ 
			1280  &   1.303e-13  &  7.815  &   2.488e-13  &  8.212  &   1.549e-12  &  8.346  &   3.457e+02 \\ 
			\midrule

			\multicolumn{8}{c}{\textbf{FORCE-3}} \\ 
			\midrule
			80  &   1.759e-04  &  $-$  &   1.927e-04  &  $-$  &   3.426e-04  &  $-$  &   1.885e+00 \\ 
			160  &   6.225e-07  &  8.142  &   7.473e-07  &  8.011  &   1.742e-06  &  7.620  &   6.628e+00 \\ 
			320  &   4.377e-09  &  7.152  &   8.186e-09  &  6.513  &   4.234e-08  &  5.362  &   2.578e+01 \\ 
			640  &   2.909e-11  &  7.233  &   7.314e-11  &  6.806  &   5.003e-10  &  6.403  &   1.031e+02 \\ 
			1280  &   1.299e-13  &  7.807  &   2.455e-13  &  8.219  &   1.527e-12  &  8.356  &   4.027e+02 \\ 
			\midrule

			\multicolumn{8}{c}{\textbf{FORCE-5}} \\ 
			\midrule
			80  &   1.849e-04  &  $-$  &   2.020e-04  &  $-$  &   3.587e-04  &  $-$  &   2.379e+00 \\ 
			160  &   6.596e-07  &  8.131  &   7.884e-07  &  8.002  &   1.827e-06  &  7.617  &   8.186e+00 \\ 
			320  &   4.651e-09  &  7.148  &   8.653e-09  &  6.510  &   4.451e-08  &  5.359  &   3.242e+01 \\ 
			640  &   3.109e-11  &  7.225  &   7.784e-11  &  6.797  &   5.301e-10  &  6.392  &   1.258e+02 \\ 
			1280  &   1.767e-13  &  7.459  &   2.762e-13  &  8.138  &   1.606e-12  &  8.366  &   4.972e+02 \\ 
			\midrule

			\multicolumn{8}{c}{\textbf{FORCE-10}} \\ 
			\midrule
			80  &   2.135e-04  &  $-$  &   2.324e-04  &  $-$  &   4.110e-04  &  $-$  &   3.210e+00 \\ 
			160  &   7.814e-07  &  8.094  &   9.218e-07  &  7.978  &   2.097e-06  &  7.615  &   1.120e+01 \\ 
			320  &   5.541e-09  &  7.140  &   1.016e-08  &  6.503  &   5.141e-08  &  5.350  &   4.375e+01 \\ 
			640  &   3.747e-11  &  7.208  &   9.296e-11  &  6.772  &   6.249e-10  &  6.362  &   1.759e+02 \\ 
			1280  &   1.994e-13  &  7.554  &   3.300e-13  &  8.138  &   1.958e-12  &  8.318  &   7.094e+02 \\ 
			\midrule

			\multicolumn{8}{c}{\textbf{Rusanov}} \\ 
			\midrule
			80  &   2.800e-04  &  $-$  &   3.045e-04  &  $-$  &   5.314e-04  &  $-$  &   1.436e+00 \\ 
			160  &   1.042e-06  &  8.070  &   1.205e-06  &  7.981  &   2.542e-06  &  7.707  &   4.937e+00 \\ 
			320  &   7.211e-09  &  7.175  &   1.275e-08  &  6.563  &   6.275e-08  &  5.340  &   1.906e+01 \\ 
			640  &   4.890e-11  &  7.204  &   1.184e-10  &  6.751  &   7.818e-10  &  6.327  &   7.527e+01 \\ 
			1280  &   2.140e-13  &  7.836  &   4.171e-13  &  8.149  &   2.511e-12  &  8.282  &   3.000e+02 \\ 
			\midrule

			\multicolumn{8}{c}{\textbf{HLL}} \\ 
			\midrule
			80  &   1.751e-04  &  $-$  &   1.919e-04  &  $-$  &   3.411e-04  &  $-$  &   1.619e+00 \\ 
			160  &   6.198e-07  &  8.142  &   7.444e-07  &  8.010  &   1.736e-06  &  7.618  &   5.749e+00 \\ 
			320  &   4.359e-09  &  7.152  &   8.156e-09  &  6.512  &   4.220e-08  &  5.362  &   2.274e+01 \\ 
			640  &   2.897e-11  &  7.233  &   7.285e-11  &  6.807  &   4.985e-10  &  6.404  &   8.877e+01 \\ 
			1280  &   1.277e-13  &  7.826  &   2.460e-13  &  8.210  &   1.496e-12  &  8.380  &   3.547e+02 \\ 
			\midrule

			\multicolumn{8}{c}{\textbf{exact RS}} \\ 
			\midrule
			80  &   1.751e-04  &  $-$  &   1.919e-04  &  $-$  &   3.411e-04  &  $-$  &   1.495e+00 \\ 
			160  &   6.198e-07  &  8.142  &   7.444e-07  &  8.010  &   1.736e-06  &  7.618  &   5.185e+00 \\ 
			320  &   4.359e-09  &  7.152  &   8.156e-09  &  6.512  &   4.220e-08  &  5.362  &   2.022e+01 \\ 
			640  &   2.897e-11  &  7.233  &   7.285e-11  &  6.807  &   4.985e-10  &  6.404  &   8.300e+01 \\ 
			1280  &   1.277e-13  &  7.826  &   2.460e-13  &  8.210  &   1.496e-12  &  8.380  &   3.254e+02 \\ 
			\midrule

			\bottomrule
	\end{tabular}}
\end{table}

\subsubsection{Shock--turbulence interaction}\label{sec:Euler_1d_shock_turbulence}
Here, we consider the test problem proposed in~\cite{titarev2004finite}, which is a modification of that in~\cite{shu1989efficient}.
The initial condition is prescribed as
\begin{align}
	\begin{pmatrix}
		\rho\\
		u\\
		p
	\end{pmatrix}(x,0):=\begin{cases}
		(1.515695,0.523346,1.80500)^T, \quad  &\text{if}~x< -4.5,\\
		(1.0+0.1\sin{(20 \pi x)},0,1)^T, \quad &\text{otherwise},
	\end{cases} 
	\label{eq:Euler_1d_titarev_toro_IC}
\end{align}
over $\Omega:=[-5,5]$, in which inflow boundary conditions are assumed for the left boundary, while, transmissive ones are considered at the right boundary. We consider a final time $T_f:=5.$

\begin{figure}[htbp]
	\centering
	\begin{subfigure}[b]{1.\textwidth}
		\centering
		\includegraphics[width=\textwidth]{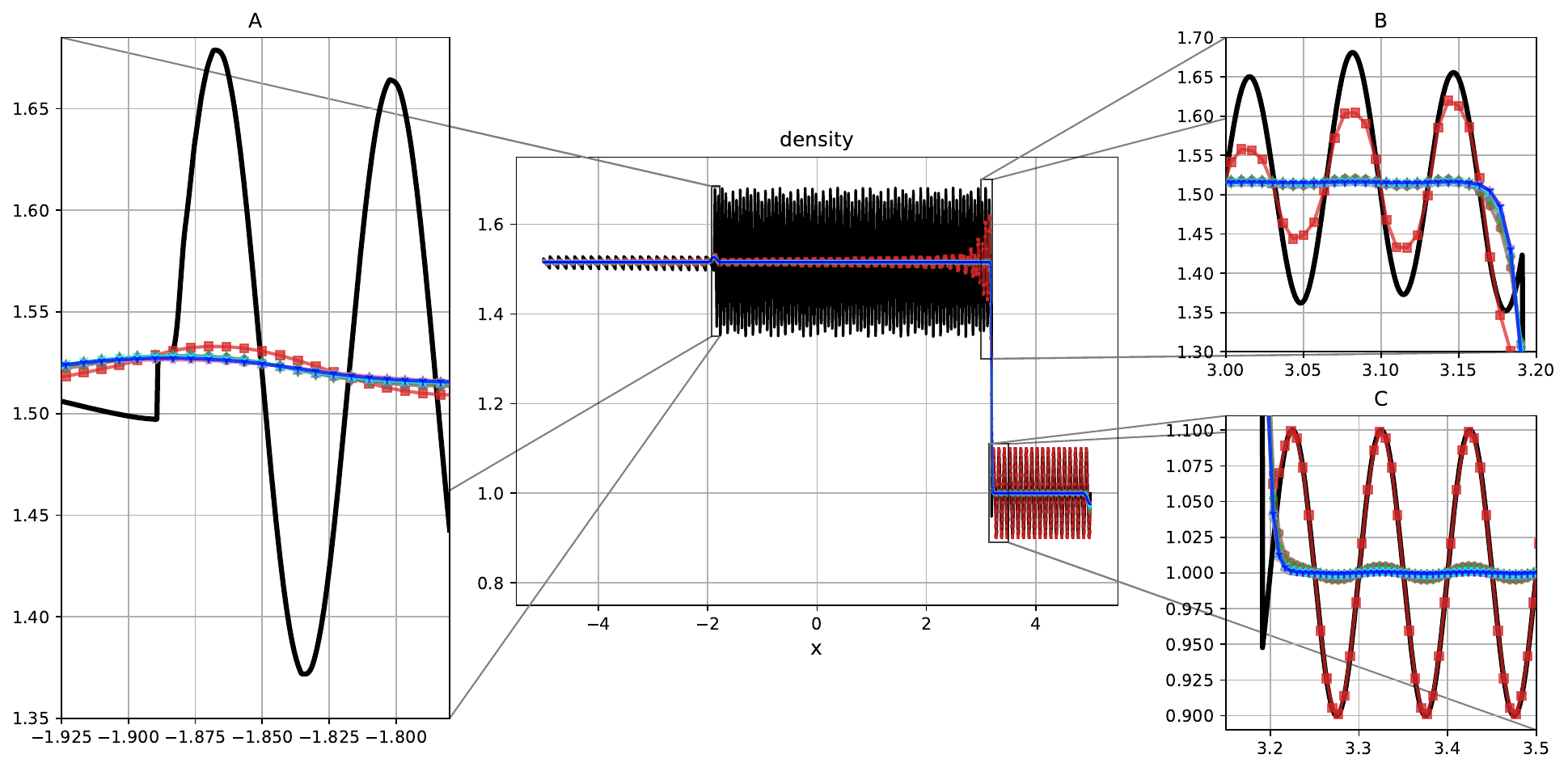}
		\caption{Order 3}
	\end{subfigure}\\
	\begin{subfigure}[b]{1.\textwidth}
		\centering
		\includegraphics[width=\textwidth]{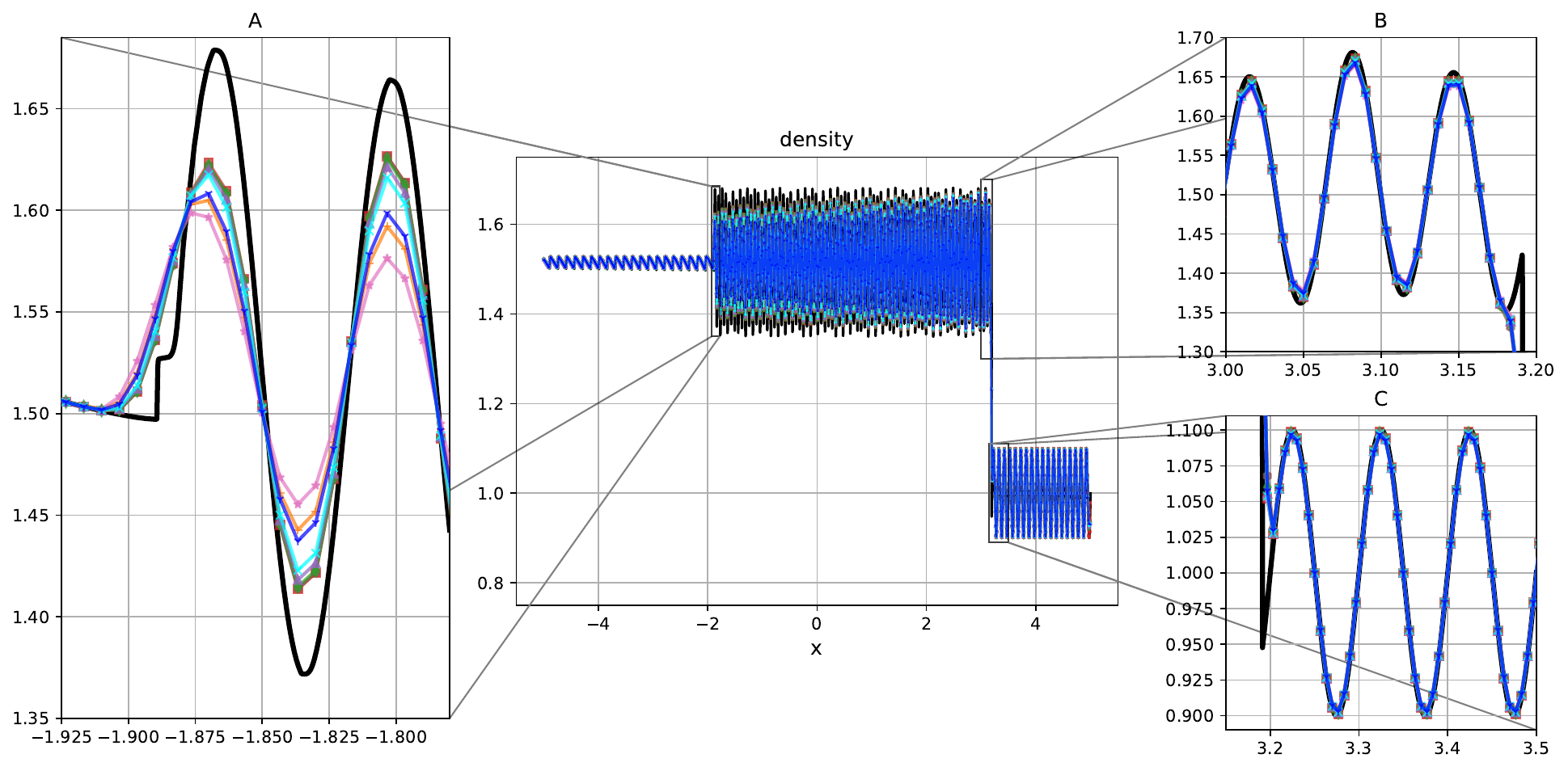}
		\caption{Order 7}
	\end{subfigure}
	\\
	\begin{subfigure}[b]{1.\textwidth}
		\centering
		\includegraphics[width=\textwidth]{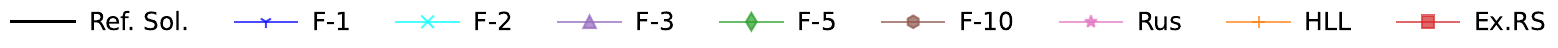}
	\end{subfigure}
	\caption{\RIcolor{Shock--turbulence interaction: Density profile obtained on a mesh with 1500 elements for $\sigma_{CFL}:=0.9$}}
	\label{fig:Euler_1d_shock_turbulence_interaction_titarev_toro_zoom_density}
\end{figure}

Such a test is particularly indicated for the assessment of the performance of (very) high order schemes, as the solution, resulting from the interaction of a shock with a turbulent flow, features several smooth structures to be resolved.
At the same time, the presence of a shock requires a robust shock--capturing character from schemes tackling this problem.
\RIIcolor{We start with a test with fixed mesh resolution and run our simulations with 1500 cells and $\sigma_{CFL}:=0.9$.}
\RIIcolor{The density results are reported in Figure~\ref{fig:Euler_1d_shock_turbulence_interaction_titarev_toro_zoom_density}.
We focus on orders 3 and 7, as they effectively allow us to appreciate how the gap between FORCE--$\alpha$ numerical fluxes and exact RS decreases when moving from a moderate-order to a very-high-order regime. The intermediate order 5 exhibits the same qualitative behavior and yields results lying between those obtained with orders 3 and 7; it is therefore omitted.}
For such a test, no exact solution is available, and the reference solution has been computed, over a mesh of 200,000 elements, through a second order FV scheme with exact RS numerical flux and van Leer's minmod limiter~\cite{AbgrallMishranotes}, SSPRK2 time discretization with
$C_{CFL}:=0.5$, and reconstruction of characteristic variables.

One can see that, for order 3, the results of all numerical fluxes are rather unsatisfactory. The only exception is given by exact RS, which is able to capture some solution details in panel B and to fully resolve the flow structures in panel C. Nonetheless, we must remark how, for the considered mesh refinement, even with exact RS the vast majority of the turbulent part of the flow is essentially unresolved.
The quality of the results increases and the differences amongst the numerical fluxes decrease, as observed in~\cite{micalizzitoro2024}, for higher order of accuracy.
In particular, for order 7, the performance of FORCE--$\alpha$ numerical fluxes is the same as the one of HLL and exact RS in panels B and C.
Slight differences remain in panel A, corresponding to the tail of the turbulent part.
Let us however observe that, in order of increasing quality, we have Rusanov, HLL, the FORCE--$\alpha$ numerical fluxes for increasing $\alpha$ and the exact RS.
It is interesting to notice that the centred FORCE--$\alpha$ fluxes give better results than the incomplete upwind ones, Rusanov and HLL, and that the results of FORCE-3, FORCE-5, FORCE-10 and exact RS are basically indistinguishable.
Again, we remark that FORCE--$\alpha$ numerical fluxes constitute a valid alternative to upwind ones in very-high-order contexts.

\RIIcolor{
Indeed, a fair efficiency comparison should be performed at fixed computational cost rather than at fixed mesh.
Hence, we have run further simulations for all numerical fluxes and orders with different meshes, corresponding to the same computational time as that required by the exact RS with order 7 on a mesh of 1500 elements. Also in this case we have used $\sigma_{CFL}:=0.9$.
The number of elements is reported in Table~\ref{tab:Euler_1d_efficiency_shock_turbulence}.}

\begin{table}[htbp]
	\centering
	\begin{tabular}{|c||c|c|c|}
		\hline
		Numerical flux & Order 3 & Order 5 & Order 7\\
		\hline\hline
		FORCE--1  & 5250 & 2850 & 1750\\
		\hline
		FORCE--2  & 4900 & 2650 & 1650\\
		\hline
		FORCE--3  & 4500 & 2450 & 1500\\
		\hline
		FORCE--5  & 4050 & 2200 & 1350\\
		\hline
		FORCE--10 & 3450 & 1900 & 1150\\
		\hline
		Rusanov   & 5250 & 2850 & 1750\\
		\hline
		HLL       & 3800 & 2150 & 1400\\
		\hline
		exact RS  & 4050 & 2300 & 1500\\
		\hline
	\end{tabular}
	\caption{
		\RIIcolor{Shock--turbulence interaction: Number of elements employed in the fixed-computational-cost analysis of the shock--turbulence interaction problem. The computational cost is chosen equal to the CPU time required by the exact RS with order 7 on a mesh of 1500 elements with $\sigma_{CFL}:=0.9$, corresponding to approximately 19 seconds.}
	}
	\label{tab:Euler_1d_efficiency_shock_turbulence}
\end{table}

\begin{figure}[htbp]
	\centering
	\includegraphics[width=1.0\textwidth]{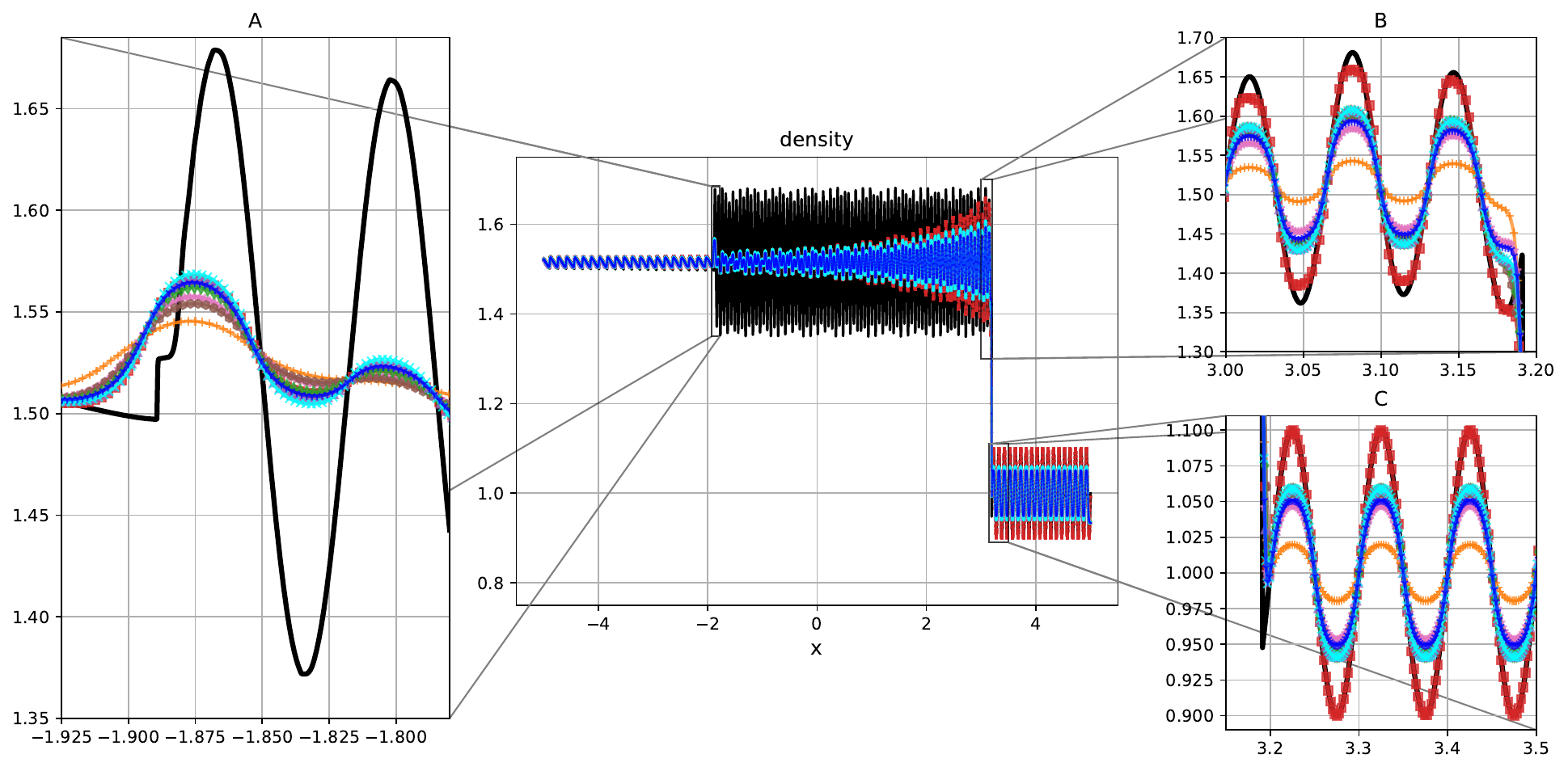}\\
	\includegraphics[width=\textwidth]{figures_new/ref_sol_1d_common_legend.pdf}
	\caption{\RIIcolor{Shock--turbulence interaction: Density profile obtained in the fixed-computational-cost analysis with order 3. The computational budget is chosen equal to the CPU time required by the exact RS with order 7 on a mesh of 1500 elements with $\sigma_{CFL}:=0.9$.}}
	\label{fig:Euler_1d_efficiency_shock_turbulence_order3}
\end{figure}

\begin{figure}[htbp]
	\centering
	\includegraphics[width=1.0\textwidth]{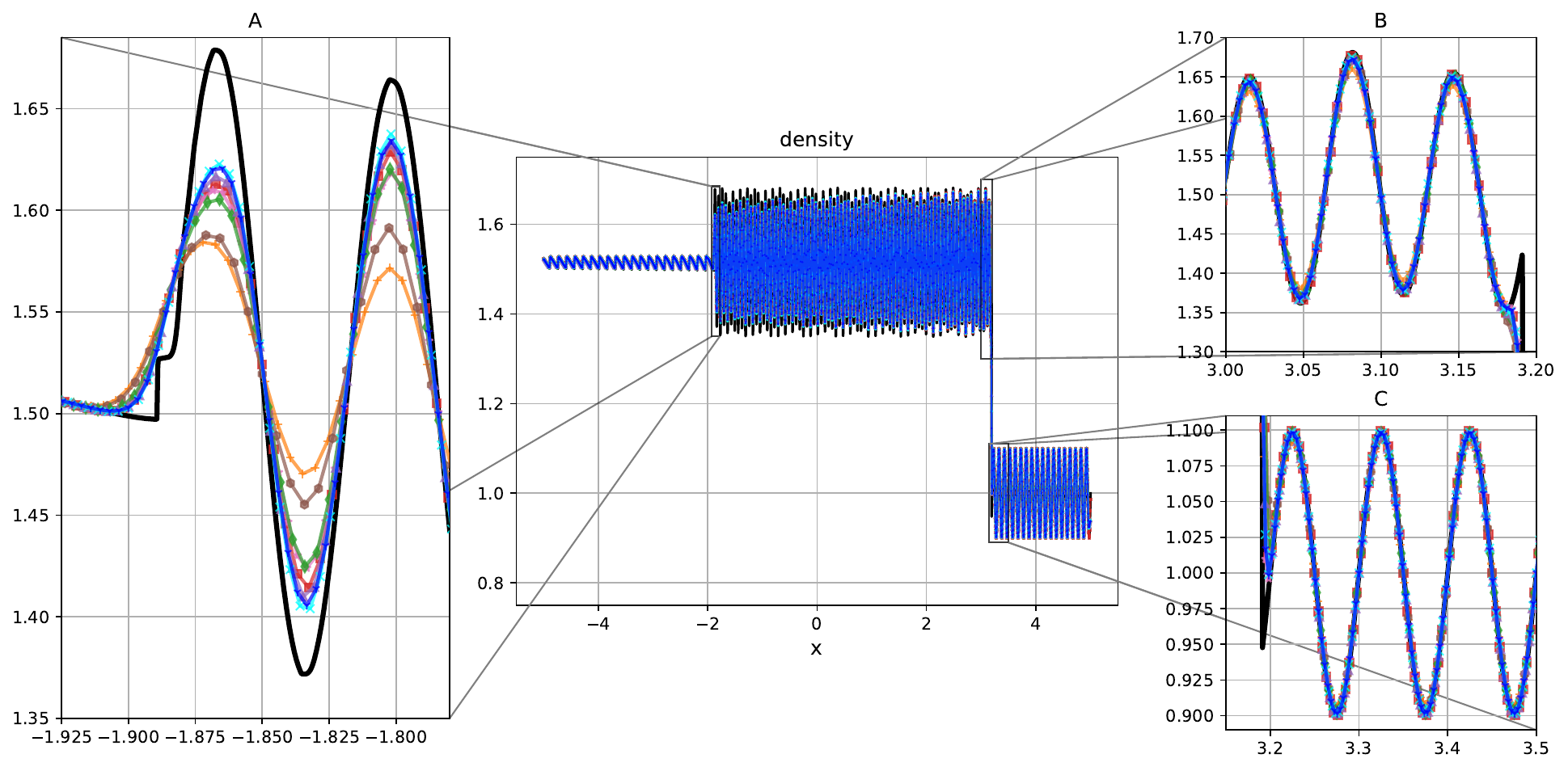}\\
	\includegraphics[width=\textwidth]{figures_new/ref_sol_1d_common_legend.pdf}
	\caption{\RIIcolor{Shock--turbulence interaction: Density profile obtained in the fixed-computational-cost analysis with order 5. The computational budget is chosen equal to the CPU time required by the exact RS with order 7 on a mesh of 1500 elements with $\sigma_{CFL}:=0.9$.}}
	\label{fig:Euler_1d_efficiency_shock_turbulence_order5}
\end{figure}

\begin{figure}[htbp]
	\centering
	\includegraphics[width=1.0\textwidth]{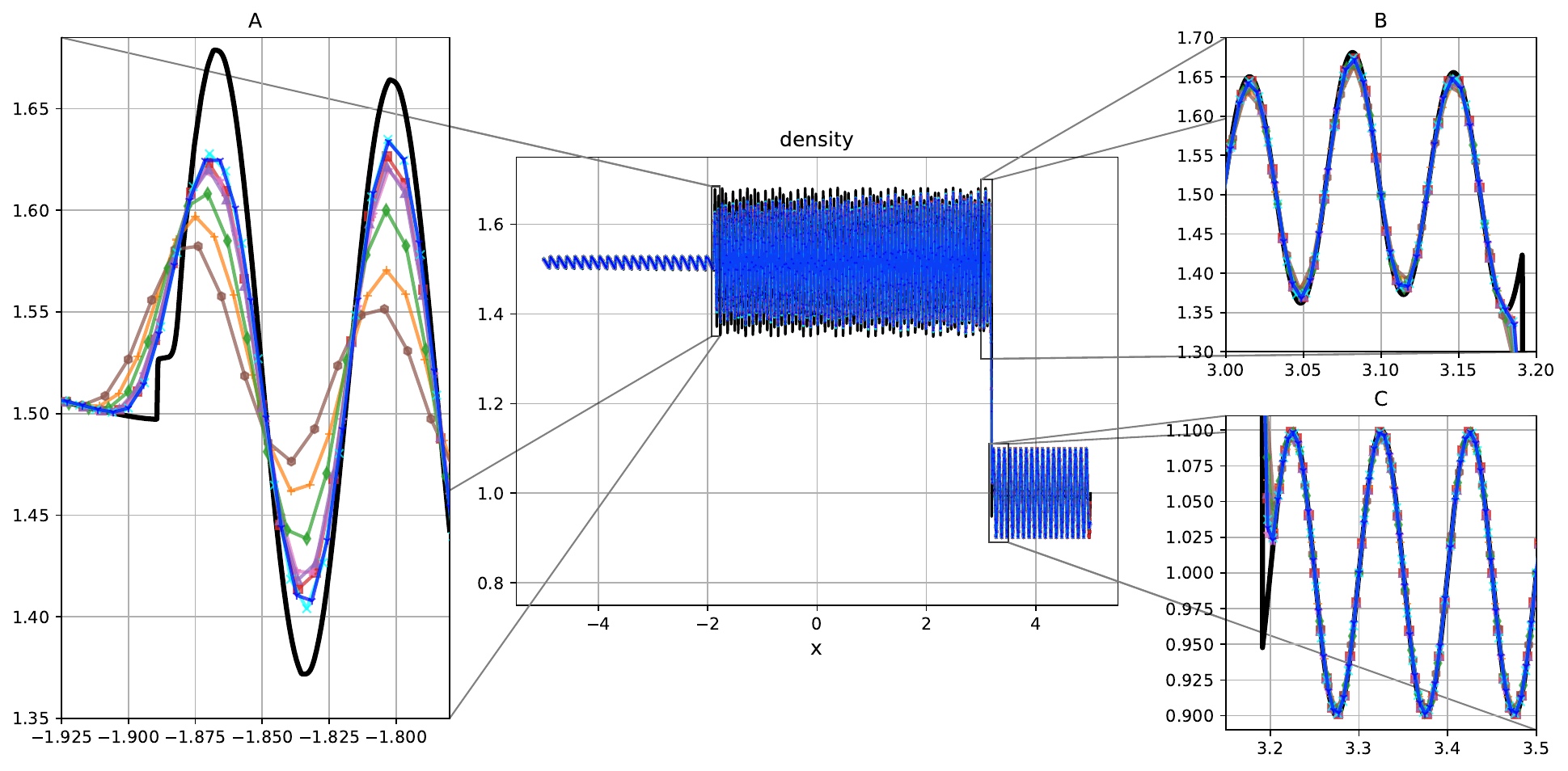}\\
	\includegraphics[width=\textwidth]{figures_new/ref_sol_1d_common_legend.pdf}
	\caption{\RIIcolor{Shock--turbulence interaction: Density profile obtained in the fixed-computational-cost analysis with order 7. The computational budget is chosen equal to the CPU time required by the exact RS with order 7 on a mesh of 1500 elements with $\sigma_{CFL}:=0.9$.}}
	\label{fig:Euler_1d_efficiency_shock_turbulence_order7}
\end{figure}

\RIIcolor{The density results are reported in Figures~\ref{fig:Euler_1d_efficiency_shock_turbulence_order3},~\ref{fig:Euler_1d_efficiency_shock_turbulence_order5} and~\ref{fig:Euler_1d_efficiency_shock_turbulence_order7} for orders 3, 5 and 7 respectively.
We see a huge improvement in the resolution in passing from order 3 to order 5. However, interestingly, the resolution does not improve in passing from order 5 to order 7.
A possible explanation is that, for the considered meshes, the turbulent part is under-resolved by WENO7 with the nonlinear weights degrading the full accuracy.
Indeed, order 7 might gain higher resolution with respect to order 5 by switching from the standard time integration DeC methods considered here to the more efficient versions discussed in~\cite{micalizzi2023new}. There, interpolation processes are introduced to increase the number of subtimenodes across the iteration process, resulting in computational advantages which considerably increase with the order of accuracy; this is however out of the scope of the current investigation.
The numerical results show that, while for order 3  exact RS constitutes by far the best choice (despite FORCE-1, FORCE-2 and FORCE-3 offering similar performances in the tail of the turbulent part in Panel A), for orders 5 and 7 FORCE-$\alpha$ schemes with $\alpha=1$ and $\alpha=2$, characterized by similar resolution, are the best numerical fluxes, overperforming all other competitors.
For such orders, the flow is fully resolved by all numerical fluxes in panels B and C, but differences can be appreciated in Panel A.
FORCE-1 and FORCE-2 are followed by FORCE-3 for similar performance as exact RS. Rusanov is comparable to FORCE-3 for order 7, while it performs slightly worse for order 5.
FORCE-5, FORCE-10 and HLL offer the worst performance.

These results, along with the efficiency analysis on the smooth problem of the previous section, demonstrate that lower values of $\alpha$ are preferable, with $\alpha=1$ and $\alpha=2$ being characterized by very similar performance and constituting the optimal values.
Furthermore, FORCE-1 and FORCE-2 are able to achieve a resolution which is competitive with the one achieved by classical upwind fluxes per fixed computational time.
}

\subsubsection{Riemann problems}\label{sec:Euler_1d_RP}
In this section, we report the results obtained on challenging tests taken from~\cite{ToroBook}.
They are Riemann problems on the interval $\Omega:=[0,1]$, endowed with transmissive boundary conditions, with initial condition
\begin{align}
	(\rho,u,p)(x,0):=\begin{cases}
		(\rho_L,u_L,p_L), \quad x<x_d,\\
		(\rho_R,u_R,p_R), \quad x>x_d.
	\end{cases}
\end{align}
The specific data, i.e., left and right states, discontinuity location, and final time, are reported in Table~\ref{tab:Euler_1d_RP}.
As stated in~\cite{ToroBook}, these tests have been specifically designed to challenge the robustness of numerical methods: oscillations and simulation crashes due to negative values of density and pressure are rather common in the context of the numerical solution of such problems.
\RIIcolor{
	In these tests, we consider a fixed mesh refinement, since their purpose is to assess the ability of the numerical fluxes under investigation to resolve structures such as rarefactions, shocks and contact discontinuities.
	To improve the readability of the manuscript and avoid unnecessary repetitions, in these tests we report detailed results only for orders 3 and 7. These two orders are sufficient to assess the behavior of the considered FORCE--$\alpha$ numerical fluxes in both a moderate-order and a very-high-order regime, while keeping the presentation compact. 
	Like in the shock--turbulence interaction test, for fixed mesh refinement, order 5 is characterized by results of intermediate quality between those obtained with orders 3 and 7.}

\begin{table}
	\centering
	\begin{tabular}{|c||c|c|c||c|c|c||c||c|}
		\hline
		Test& $\rho_L$ & $u_L$ & $p_L$& $\rho_R$& $u_R$ & $p_R$ & $x_d$ & $T_f$\\\hline\hline
		1   &      1.0 &  0.75 &   1.0&    0.125&   0.0 &  0.1  &  0.3  & 0.2 \\
		\hline
		2   &      1.0 &  -2.0 &   0.4&      1.0&   2.0 &  0.4  &  0.5 & 0.15 \\
		\hline
		3   &      1.0 &   0.0 &1000.0&      1.0&   0.0 & 0.01  & 0.5  & 0.012 \\
		\hline
		4   &  5.99924 &19.5975&460.894& 5.99242&-6.19633&46.0950 & 0.4  & 0.035 \\
		\hline
		5   &     1.0  &-19.59745&1000.0& 1.0 & -19.59745 & 0.01 &  0.8 &  0.012\\
		\hline
	\end{tabular}
	\caption{Initial conditions for Riemann problems. The quantities $\rho_L$, $u_L$ and $p_L$ represent the values of the density, of the velocity and of the pressure of the left state, while, the quantities $\rho_R$, $u_R$ and $p_R$ represent the analogous values for the right state. The quantities $x_d$ and $T_f$ indicate the initial position of the discontinuity and the considered final time respectively}\label{tab:Euler_1d_RP}
\end{table}

\RIIcolor{
	
	For the sake of compactness and to focus the discussion on the most representative benchmarks, we omit the detailed results for Riemann problems 2 and 4, and retain only Riemann problems 1, 3 and 5, which better highlight the differences and the similarities among the considered numerical fluxes. We nevertheless summarize below the main conclusions obtained from Riemann problems 2 and 4 on a mesh of 100 elements.

	Riemann problem 2 is characterized by two waves moving in opposite directions, leading to the formation of a low-density and low-pressure region in the middle of the domain. It is essentially a benchmark aimed at assessing positivity-preserving properties rather than the resolution of complex flow features. No significant differences among the considered FORCE--$\alpha$ numerical fluxes can be observed in terms of solution quality in this test, and the already limited discrepancies progressively disappear as the order of accuracy increases. For order 3, FORCE--$\alpha$ fluxes exhibit a performance comparable to Rusanov and HLL, while exact RS leads to simulation failure. Interestingly, FORCE--$\alpha$ numerical fluxes display a higher degree of robustness than the considered upwind fluxes. In fact, for suitable restrictions of $\sigma_{CFL}\in [0.1,0.9]$, all FORCE--$\alpha$ fluxes successfully reach the final simulation time for all investigated orders, whereas Rusanov, HLL and exact RS fail for orders 5 and 7 due to positivity violations. Moreover, lower values of $\alpha$ correspond to a higher level of robustness and more relaxed restrictions on $\sigma_{CFL}$.
		
	%
	%
	%
	%
	%
	%
	%
	%
	%
	%
	
	Riemann problem 4 is characterized by the collision of two high-density and high-pressure states. All numerical fluxes and orders successfully complete the simulation with $\sigma_{CFL}=0.9$. The differences among the considered numerical fluxes are generally small and decrease as the order of accuracy increases. For order 7, the numerical solutions are very similar, confirming that FORCE--$\alpha$ numerical fluxes can provide results comparable to those obtained with classical upwind fluxes when employed within a very-high-order framework.
	
}

\subsubsection*{Riemann problem 1}
This test is a modification of the Sod shock tube problem~\cite{sod1978survey}. 
Such a modified version is indeed more challenging than the original one due to the presence of a sonic point, which constitutes an issue for linearised Riemann solvers, such as Roe--type numerical fluxes~\cite{ToroBook}.

\begin{figure}[htbp]
	\centering
	\begin{subfigure}[b]{1.\textwidth}
		\centering
		\includegraphics[width=\textwidth]{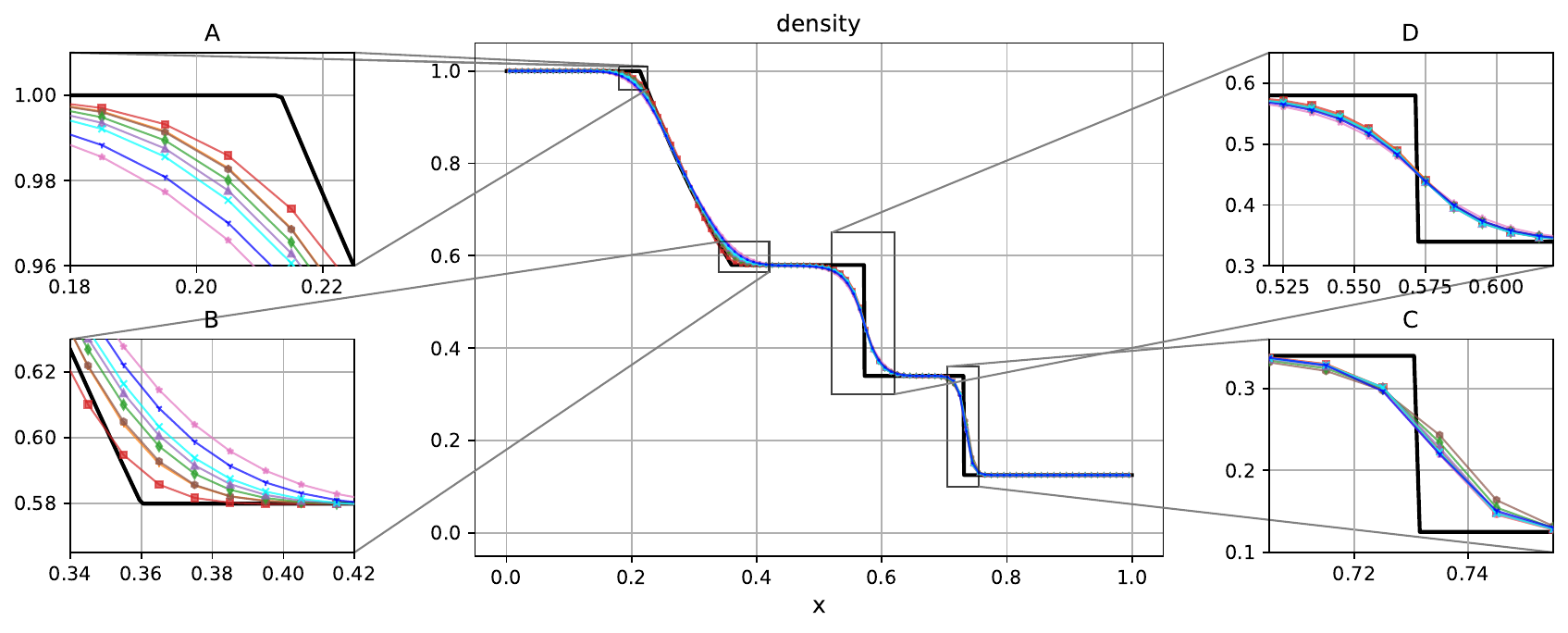}
		\caption{Order 3}
	\end{subfigure}\\
	\begin{subfigure}[b]{1.\textwidth}
		\centering
		\includegraphics[width=\textwidth]{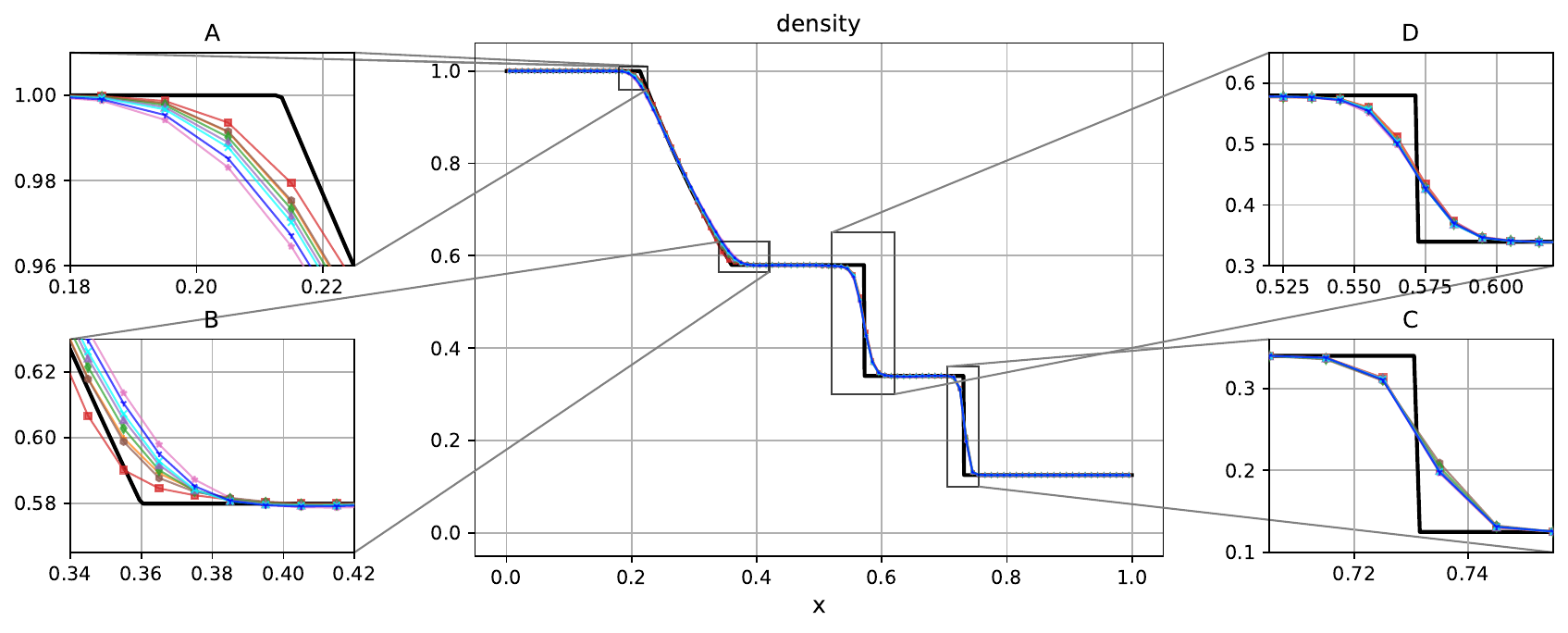}
		\caption{Order 7}
	\end{subfigure}
	\\
	\begin{subfigure}[b]{1.\textwidth}
		\centering
		\includegraphics[width=\textwidth]{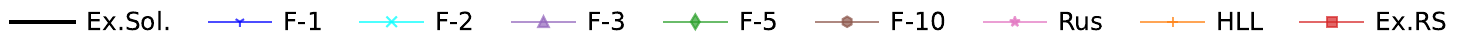}
	\end{subfigure}
	\caption{\RIcolor{Riemann problem 1: Density profile obtained on a mesh with 100 elements for $\sigma_{CFL}:=0.9$}}
	\label{fig:Euler_1d_RP1_zoom_density}
\end{figure}

We have run our simulations with 100 elements and $\sigma_{CFL}:=0.9$, without experiencing any issue for any order and numerical flux. 
The obtained density results are reported in Figure~\ref{fig:Euler_1d_RP1_zoom_density}. 
As one can see, the main differences among the numerical fluxes can be appreciated on the head and on the tail of the rarefaction in panels A and B respectively.
The exact RS gives the best results, while, Rusanov gives the worst ones.
Between these two, we find all FORCE--$\alpha$ numerical fluxes, with increasing quality as $\alpha$ increases, and HLL.
In particular, the performances of FORCE--10 and HLL are comparable, with FORCE--10 being slightly better for order 7, see panel B.
This is interesting as we remark once again that FORCE--$\alpha$ numerical fluxes, being centred fluxes, pay much less attention to the underlying physics with respect to upwind ones, such as Rusanov and HLL. 
Overall, all schemes are able to provide a good description of the numerical solution, showing that FORCE--$\alpha$ numerical fluxes are a competitive option.

\subsubsection*{Riemann problem 3}
This test consists in the left part of the well--known Woodward--Colella problem~\cite{woodward1984numerical}, and it constitutes, due to the huge pressure jump in the initial data, a challenge for any numerical scheme for the Euler equations.

\begin{figure}[htbp]
	\centering
	\begin{subfigure}[b]{1.\textwidth}
		\centering
		\includegraphics[width=\textwidth]{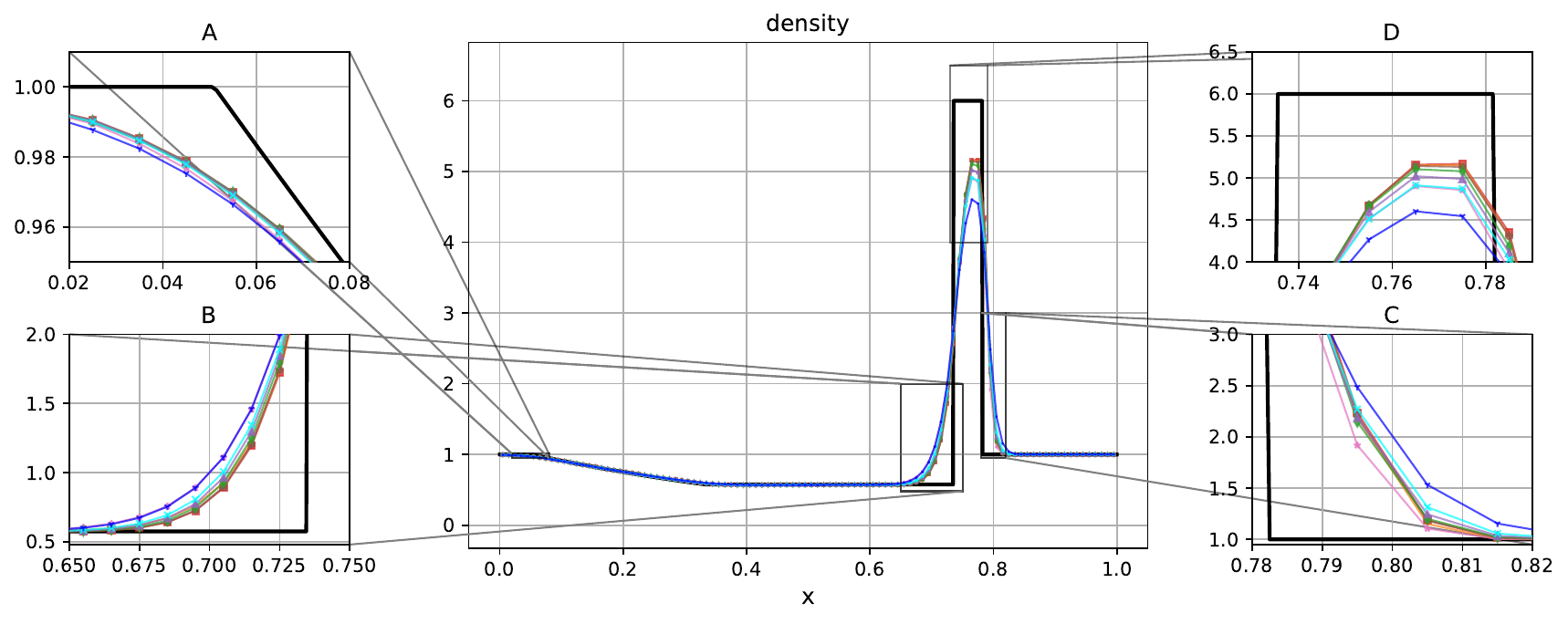}
		\caption{Order 3}
	\end{subfigure}\\
	\begin{subfigure}[b]{1.\textwidth}
		\centering
		\includegraphics[width=\textwidth]{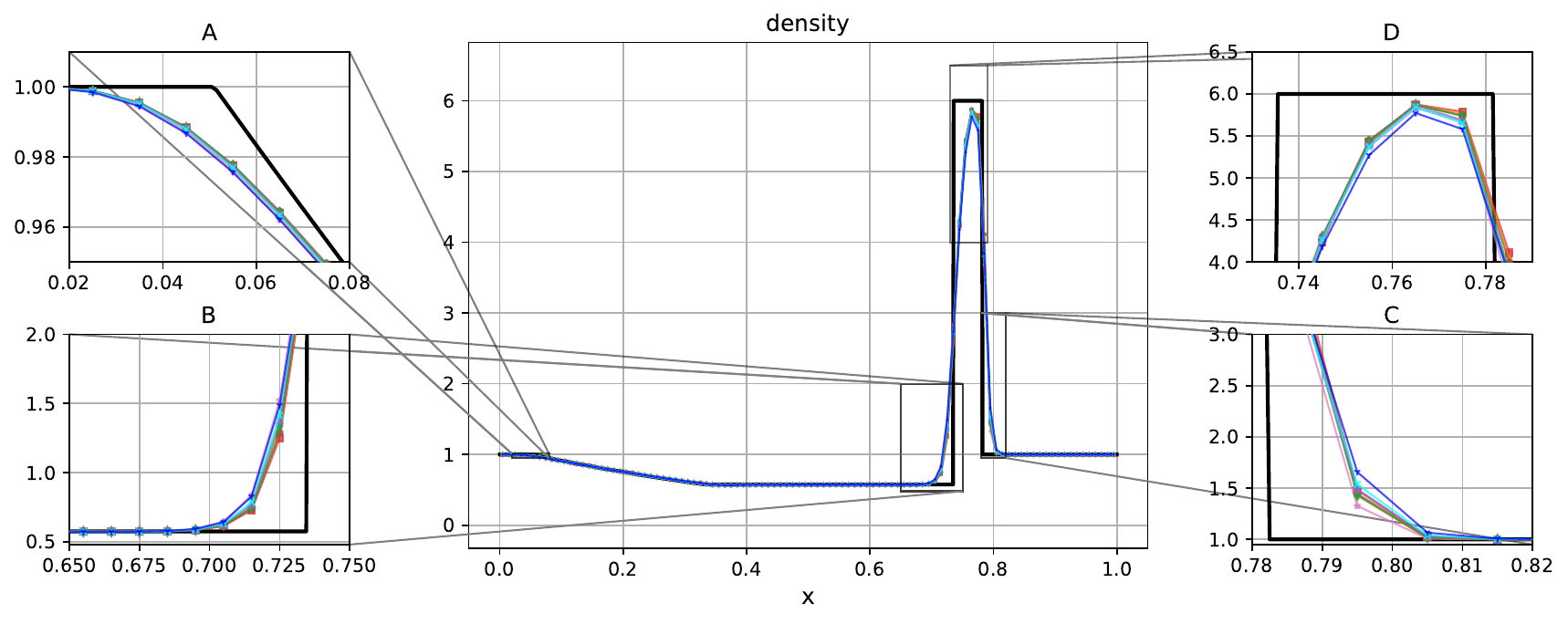}
		\caption{Order 7}
	\end{subfigure}
	\\
	\begin{subfigure}[b]{1.\textwidth}
		\centering
		\includegraphics[width=\textwidth]{figures_new/RP_common_legend.pdf}
	\end{subfigure}
	\caption{\RIcolor{Riemann problem 3: Density profile obtained on a mesh with 100 elements for $\sigma_{CFL}:=0.6$}}
	\label{fig:Euler_1d_RP3_zoom_density}
\end{figure}

For 100 elements, we have found that no simulation crashes occur with FORCE-1, FORCE-2, FORCE-3, FORCE-5 and exact RS, for $\sigma_{CFL}:=0.9$ for any order.
HLL runs with $\sigma_{CFL}\leq 0.9$ for orders 3 and 5, and with $\sigma_{CFL}\leq 0.8$ for order 7.
FORCE-10 and Rusanov require smaller values of $\sigma_{CFL}$.
More in detail, FORCE-10 experiences no issues for $\sigma_{CFL}\leq 0.6$ for all orders, while, Rusanov runs for $\sigma_{CFL}\leq 0.8$ for orders 3 and 5, and for $\sigma_{CFL}\leq 0.7$ for order 7.
The results obtained for the density, with $\sigma_{CFL}:=0.6$, are reported in Figure~\ref{fig:Euler_1d_RP3_zoom_density}. 
The main differences between the numerical fluxes are in the density peak in panel D.
For order 3, one can see that exact RS, HLL and FORCE-10 offer the best performance, and are practically indistinguishable.
Again, we remark that, due to the different nature of the fluxes, this is not a trivial fact:
FORCE-10, which is a centred flux, has the same performance as two upwind fluxes, namely, HLL and exact RS, in particular, with exact RS being also complete.
Still focusing on order 3, in terms of quality of the related results, the aforementioned fluxes are followed by FORCE-5, FORCE-3, FORCE-2, Rusanov and FORCE-1.
In particular, FORCE-2 and Rusanov yield very similar results.
Let us notice that the difference between the best numerical fluxes, exact RS, HLL and FORCE-10, and the worst one, FORCE-1, has an order of magnitude which is around 10\% of the density peak.

The quality of the results increases as the order increases, for all the numerical fluxes.
More in detail, in line with~\cite{micalizzitoro2024}, the difference between the results tends to vanish for increasing order.
For order 7, all numerical fluxes yield negligible differences, making again centred FORCE-$\alpha$ numerical fluxes a valid alternative to upwind ones.

\subsubsection*{Riemann problem 5}


In this test, the two states of density and pressure of Riemann problem 3 move towards left with high speed. The value of such a speed is selected to determine the formation of a (virtually) stationary contact discontinuity in the solution.
Besides the complications of the third Riemann problem, caused by the high pressure jump in the initial conditions, the presence of the stationary contact makes this test even more challenging.
In fact, stationary contact discontinuities are solution features typically captured by complete upwind fluxes, but not from incomplete upwind and centred ones.
Precisely for this reason, this was one of the tests selected in~\cite{micalizzitoro2024} to compare the performance of different numerical fluxes.
For 100 elements, we have determined the following constraints for $\sigma_{CFL}$:

\begin{itemize}
	
	\item FORCE-1, FORCE-2 and FORCE-3 experience no crashes with $\sigma_{CFL} := 0.9$ for any order;

	\item FORCE-5 runs without problems with $\sigma_{CFL} \leq 0.8$ for any order;

	\item FORCE-10 runs without problems with $\sigma_{CFL} \leq 0.6$ for any order;
	
	\item Rusanov runs for $\sigma_{CFL} \leq 0.9$, $\sigma_{CFL} \leq 0.8$ and $\sigma_{CFL} \leq 0.7$, for orders 3, 5 and 7 respectively;
	
	\item HLL runs with $\sigma_{CFL} \leq 0.1$ for order $3$, and it runs with $\sigma_{CFL} \leq 0.2$ for orders 5 and 7; 
	
	\item exact RS requires $\sigma_{CFL}\leq 0.3$ for orders 3 and 7, and $\sigma_{CFL}\leq 0.1$ for order 5.

\end{itemize}

\begin{figure}[htbp]
	\centering
	\begin{subfigure}[b]{1.\textwidth}
		\centering
		\includegraphics[width=\textwidth]{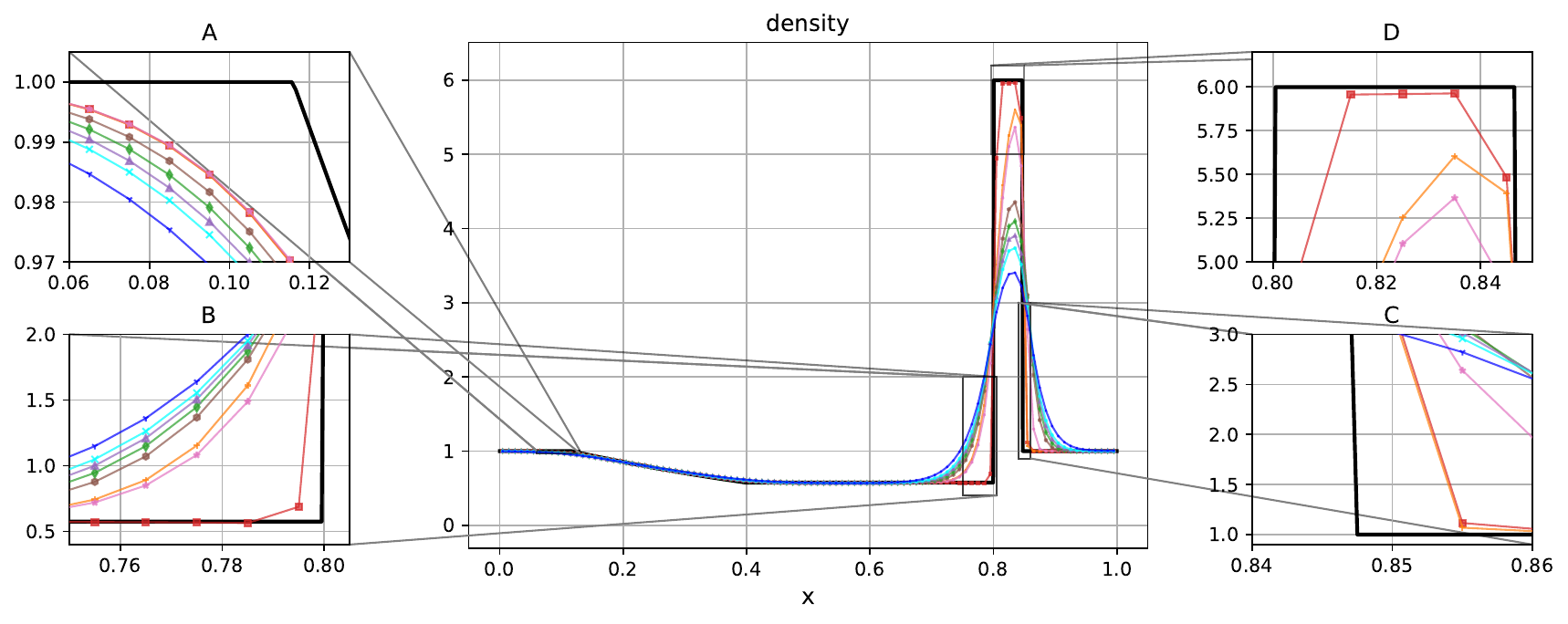}
		\caption{Order 3}
	\end{subfigure}\\
	\begin{subfigure}[b]{1.\textwidth}
		\centering
		\includegraphics[width=\textwidth]{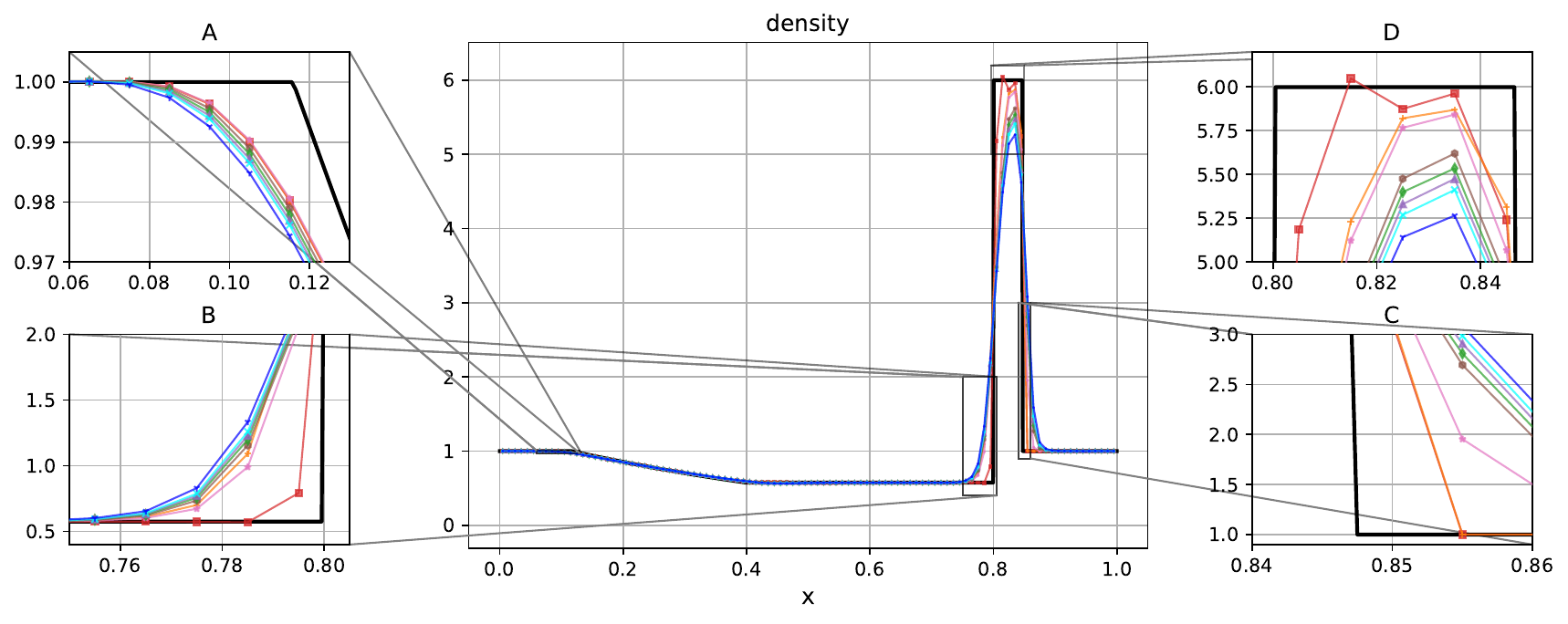}
		\caption{Order 7}
	\end{subfigure}
	\\
	\begin{subfigure}[b]{1.\textwidth}
		\centering
		\includegraphics[width=\textwidth]{figures_new/RP_common_legend.pdf}
	\end{subfigure}
	\caption{\RIcolor{Riemann problem 5: Density profile obtained on a mesh with 100 elements for $\sigma_{CFL}:=0.1$}}
	\label{fig:Euler_1d_RP5_zoom_density}
\end{figure}

The obtained density is displayed, for all investigated orders and  $\sigma_{CFL}:= 0.1$, in Figure~\ref{fig:Euler_1d_RP5_zoom_density}.
Huge differences amongst the numerical fluxes can be appreciated, especially \RIIcolor{in panels B, C and D} for order 3:
the exact RS has by far the best performance, followed by HLL and Rusanov.
FORCE--$\alpha$ numerical fluxes are rather unsatisfactory in the context of this test for such an order of accuracy. While exact RS is able to fully capture the density peak in panel D, FORCE-10 has an error which is more than 25\% with respect to the peak value. As expected, the error increases as $\alpha$ decreases, until becoming around 45\% for FORCE-1.
The situation changes when the order of accuracy increases.
For order 7, the results are much more uniform, even though differences are still evident: the error on the density peak in panel D obtained for FORCE--$\alpha$ numerical fluxes is much smaller, and it is around 6\% for FORCE-10.
Again, the decrease in the difference between the results obtained through different numerical fluxes for increasing order has been pointed out in~\cite{micalizzitoro2024}.

Such a test shows that there are problems in which the adoption of a complete upwind numerical flux results crucial for an accurate approximation of the solution even for high order.
Nonetheless, one must notice that the results obtained through FORCE--$\alpha$ numerical fluxes, unacceptable for order 3, become more satisfactory for order 7.
The main conclusion is that such numerical fluxes become more competitive within a very high order setting, which can make up for their higher level of diffusion.
Also the possibility to run simulations with higher values of $\sigma_{CFL}$ when adopting FORCE--$\alpha$ numerical fluxes in the context of this test must be appreciated.

\subsection{Two--dimensional tests}\label{sec:Euler_2d}
In this section, we test the FORCE--$\alpha$ numerical fluxes on multidimensional benchmarks. 
Initially, in Section~\ref{sec:unsteady_vortex}, we consider a smooth problem to verify the order of accuracy and assess the computational efficiency. 
\RIIcolor{Then, in Section~\ref{sec:shock_vortex}, we report the results obtained for a shock--vortex interaction problem, aimed at verifying the ability of FORCE--$\alpha$ numerical fluxes to capture complex solution features.
In Section~\ref{sec:unsteady_vortex_longer_time}, we assess the diffusion associated with long--time evolution.
Lastly, in Sections~\ref{sec:explosion} and~\ref{sec:2DRP}, we test the robustness and shock-capturing capabilities of FORCE--$\alpha$ numerical fluxes on classical benchmarks involving explosion problems and a two--dimensional Riemann problem.
}

\subsubsection{Smooth isentropic vortex}\label{sec:unsteady_vortex}
The first two--dimensional test consists in a smooth unsteady vortex, and it is meant to verify the order of accuracy of our two--dimensional implementation, and to perform a further assessment of the FORCE--$\alpha$ numerical fluxes on a smooth problem.
Over the computational domain $\Omega := [-10,10]\times [-10,10]$ with periodic boundary conditions, the initial condition is prescribed as
\begin{equation}\label{eq:vortex}
	\begin{cases}
		\rho(x,y,0) := (1+\delta \Temp)^{\frac{1}{\gamma-1}},\\
		\begin{pmatrix}
			u\\v
		\end{pmatrix}(x,y,0):=\begin{pmatrix}
		u_\infty\\v_\infty
	\end{pmatrix}+\frac{\beta}{2\pi}e^{\frac{1-r^2}{2}} \begin{pmatrix}
			-(y-y_c)\\(x-x_c)
		\end{pmatrix},	\\
		p(x,y,0) := (1+\delta \Temp)^{\frac{\gamma}{\gamma-1}},
	\end{cases}
	\quad \delta \Temp := -\frac{(\gamma-1)\beta^2}{8\gamma\pi^2}e^{1-r^2},
\end{equation}
where $\beta:=5$, $u_\infty:=1$, $v_\infty:=1$, and $r(x,y):=\sqrt{(x-x_c)^2+(y-y_c)^2}$ being the distance from the initial center of the vortex $(x_c,y_c)^T:=(0,0)^T$.
Such a set-up corresponds to a vortex translating with background speed $(u_\infty,v_\infty)^T$.
Therefore, at a generic time $t$, the solution is $\uvec{u}(x,y,t)=\uvec{u}(x-u_\infty t,y-v_\infty t,0)$ modulo the periodicity of the domain. 
It is worth pointing out that, although the vortex is $C^{\infty}$, it is not compactly supported, which may lead to spurious boundary effects due to the periodic boundary conditions.
One can verify that our choice of the domain prevents this from happening, with the initial condition near boundaries being equal to the background state $(\rho_{\infty},u_{\infty},v_{\infty},p_{\infty})^T:=(1,1,1,1)^T$ up to machine precision.

Many references make use of this test for the assessment of the accuracy of high order numerical schemes, see inter alia~\cite{shu1998essentially,micalizzi2023efficient,lore_phd_thesis,boscheri2022continuous,boscheri2015direct,micalizzitoro2024,abgrall2025dual}.
This problem is more involved with respect to the simple one--dimensional advection considered in Section~\ref{sec:Euler_1d_smooth_advection}, hence, we expect a higher level of variation from the results obtained through different numerical fluxes.
We have run our convergence analysis over meshes with $N\times N$ elements until final time $T_f:=0.1$ with $\sigma_{CFL}:=0.9$.  
The results are summarized in Tables~\ref{tab:Euler_2d_convergence_tables_WENO_DeC_char_order3},~\ref{tab:Euler_2d_convergence_tables_WENO_DeC_char_order5}, and~\ref{tab:Euler_2d_convergence_tables_WENO_DeC_char_order7}.
Again, we focus only on the density, as analogous results have been obtained for all other variables.

\begin{figure}[htbp]
	\centering
	
	\begin{subfigure}[t]{0.49\textwidth}
		\centering
		\includegraphics[width=\textwidth, trim={0 0 140 0}, clip]{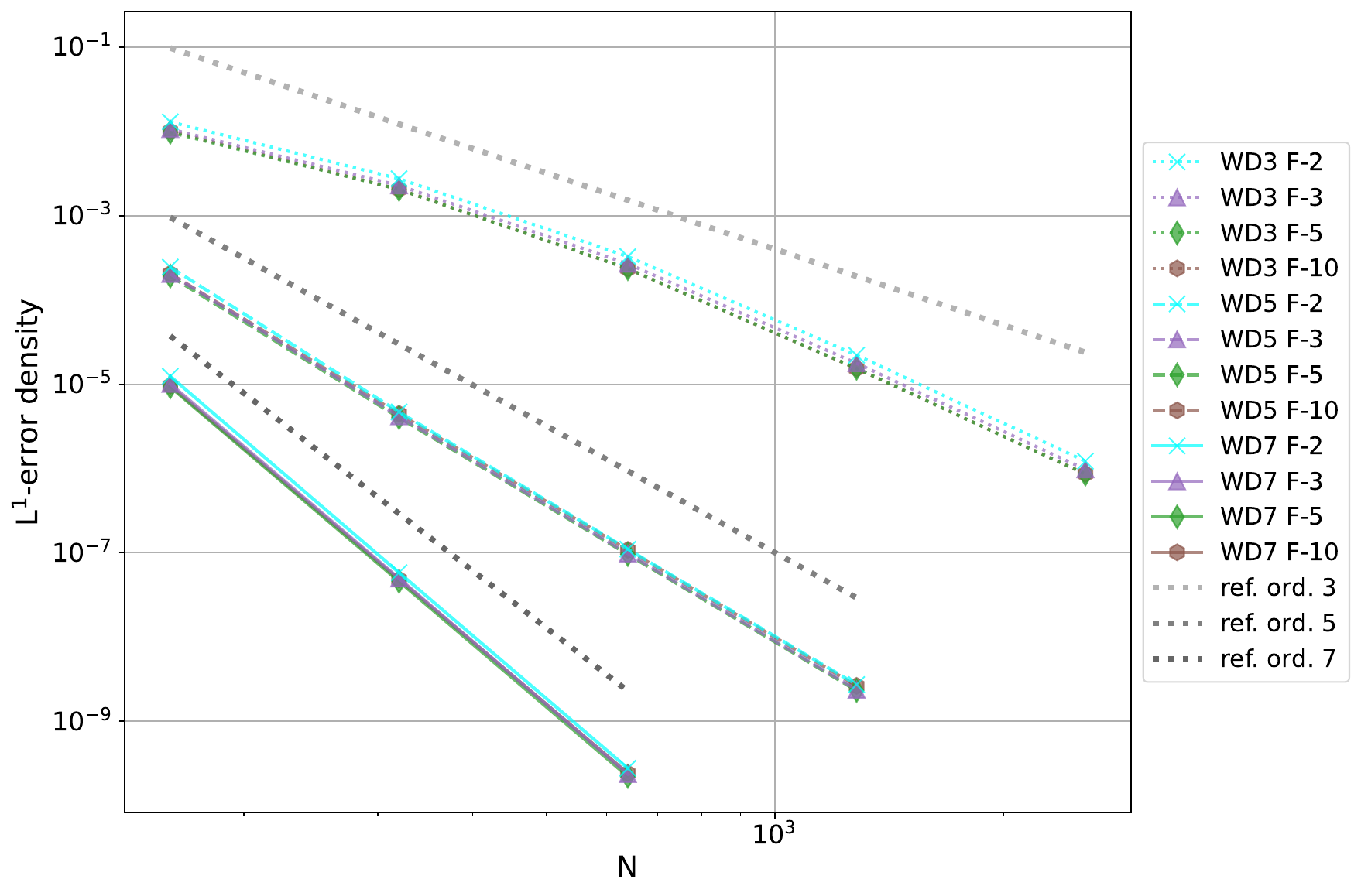}
		\caption{Convergence analysis}
	\end{subfigure}
	\hfill
	\begin{subfigure}[t]{0.49\textwidth}
		\centering
		\includegraphics[width=\textwidth, trim={0 0 150 0}, clip]{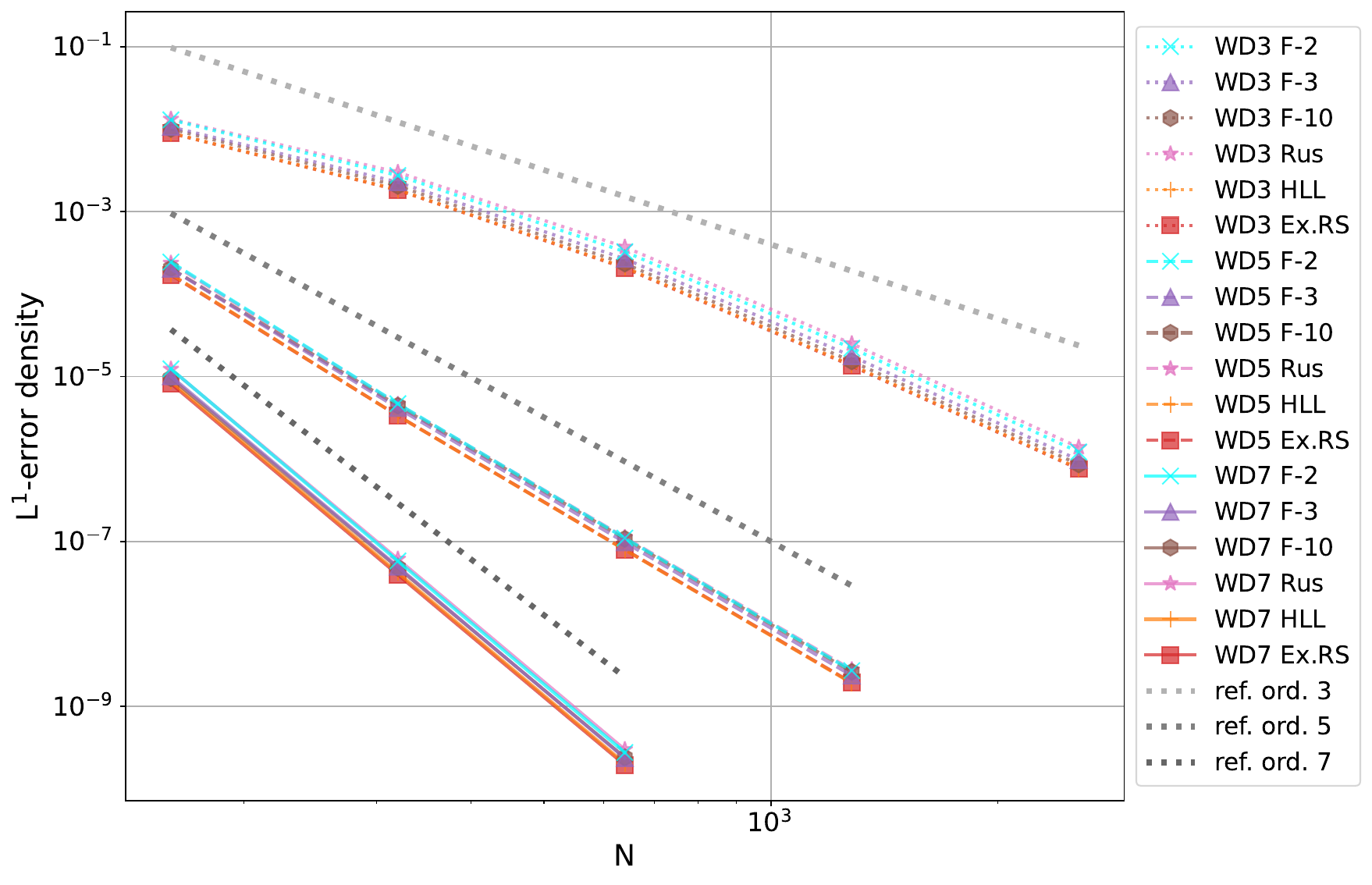}
		\caption{Convergence comparison with upwind fluxes}
	\end{subfigure}
	
	\vspace{0.3cm}
	
	\begin{subfigure}[t]{0.49\textwidth}
		\centering
		\includegraphics[width=\textwidth, trim={0 0 0 0}, clip]{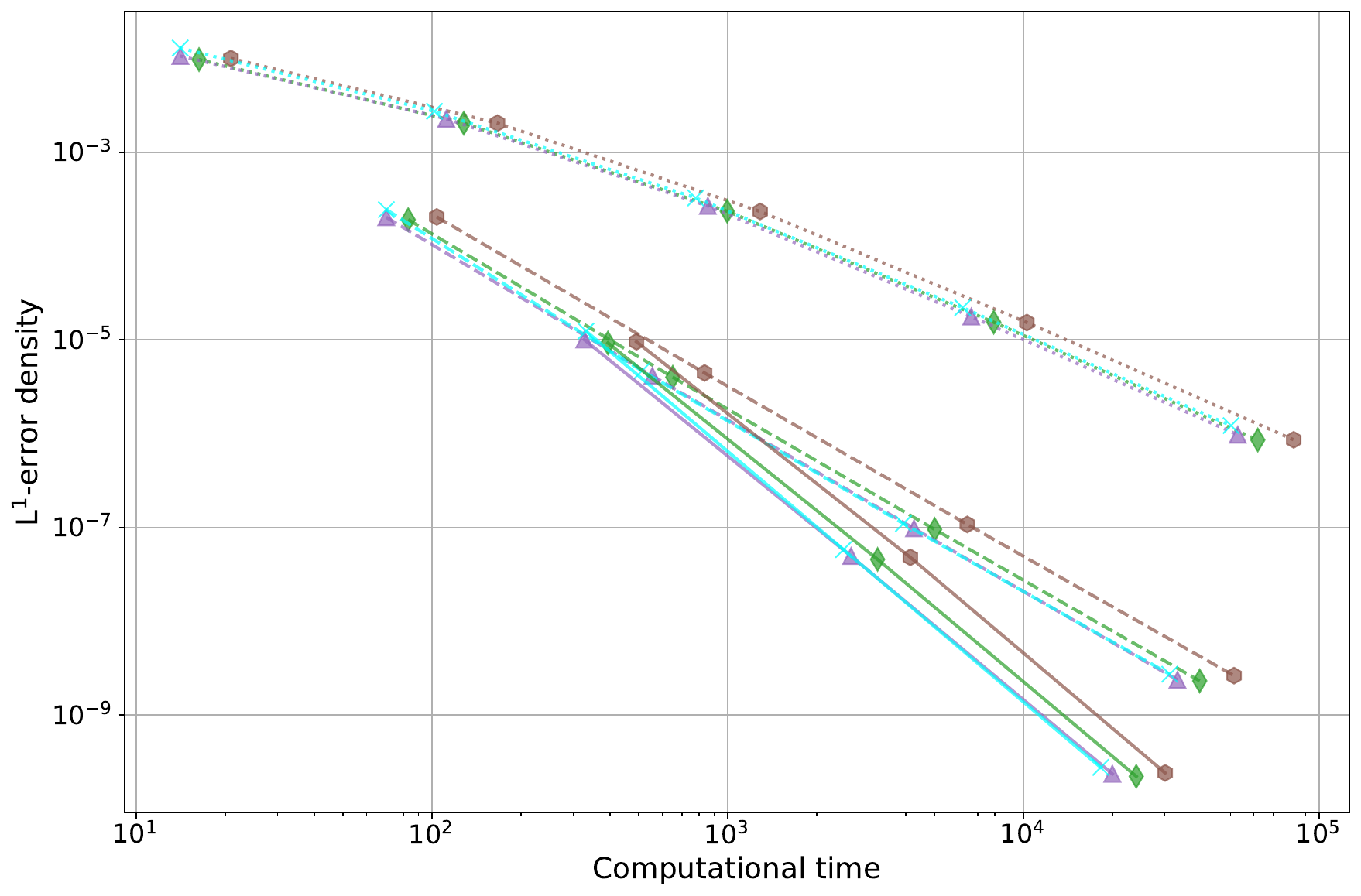}
		\caption{Efficiency analysis}
	\end{subfigure}
	\hfill
	\begin{subfigure}[t]{0.49\textwidth}
		\centering
		\includegraphics[width=\textwidth, trim={0 0 0 0}, clip]{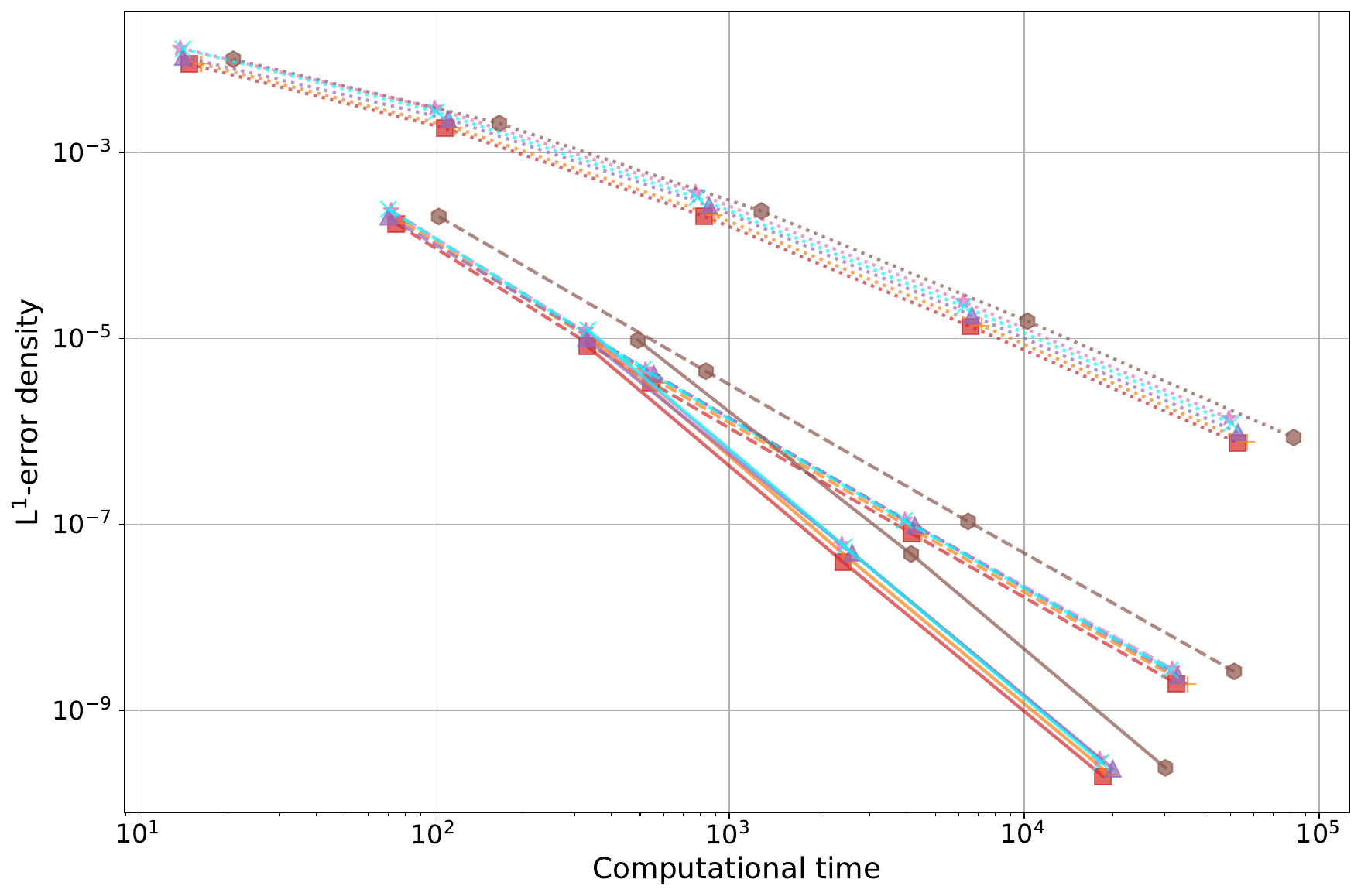}
		\caption{Efficiency comparison with upwind fluxes}
	\end{subfigure}
	
	\vspace{0.3cm}
	
	\begin{subfigure}[t]{1.0\textwidth}
		\centering
		\includegraphics[width=0.9\textwidth]{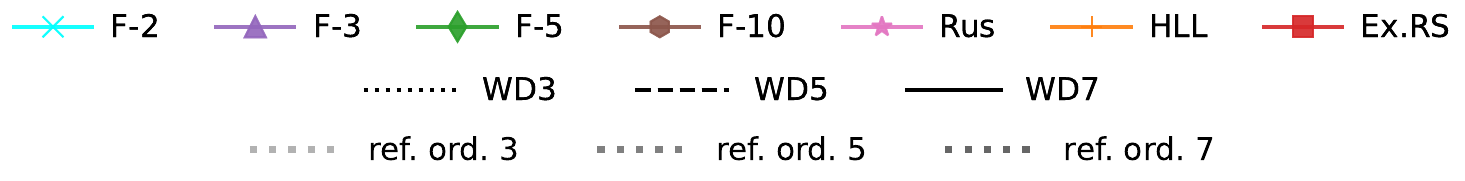}
	\end{subfigure}

	\caption{\RIcolor{Smooth isentropic vortex: Convergence and efficiency analyses. The left column illustrates the influence of the parameter $\alpha$ within the FORCE-$\alpha$ family, while the right column compares representative FORCE-$\alpha$ numerical fluxes ($\alpha=2$, $3$ and $10$) with standard upwind fluxes (Rusanov, HLL and exact RS).  On the bottom a common legend for numerical fluxes (distinguished by color and marker) and order (distinguished by linestyle) is reported.}}
	\label{fig:Euler_2d_unsteady_vortex_WENODeC}
\end{figure}

\begin{figure}[htbp]
	\centering
	\includegraphics[width=0.9\linewidth]{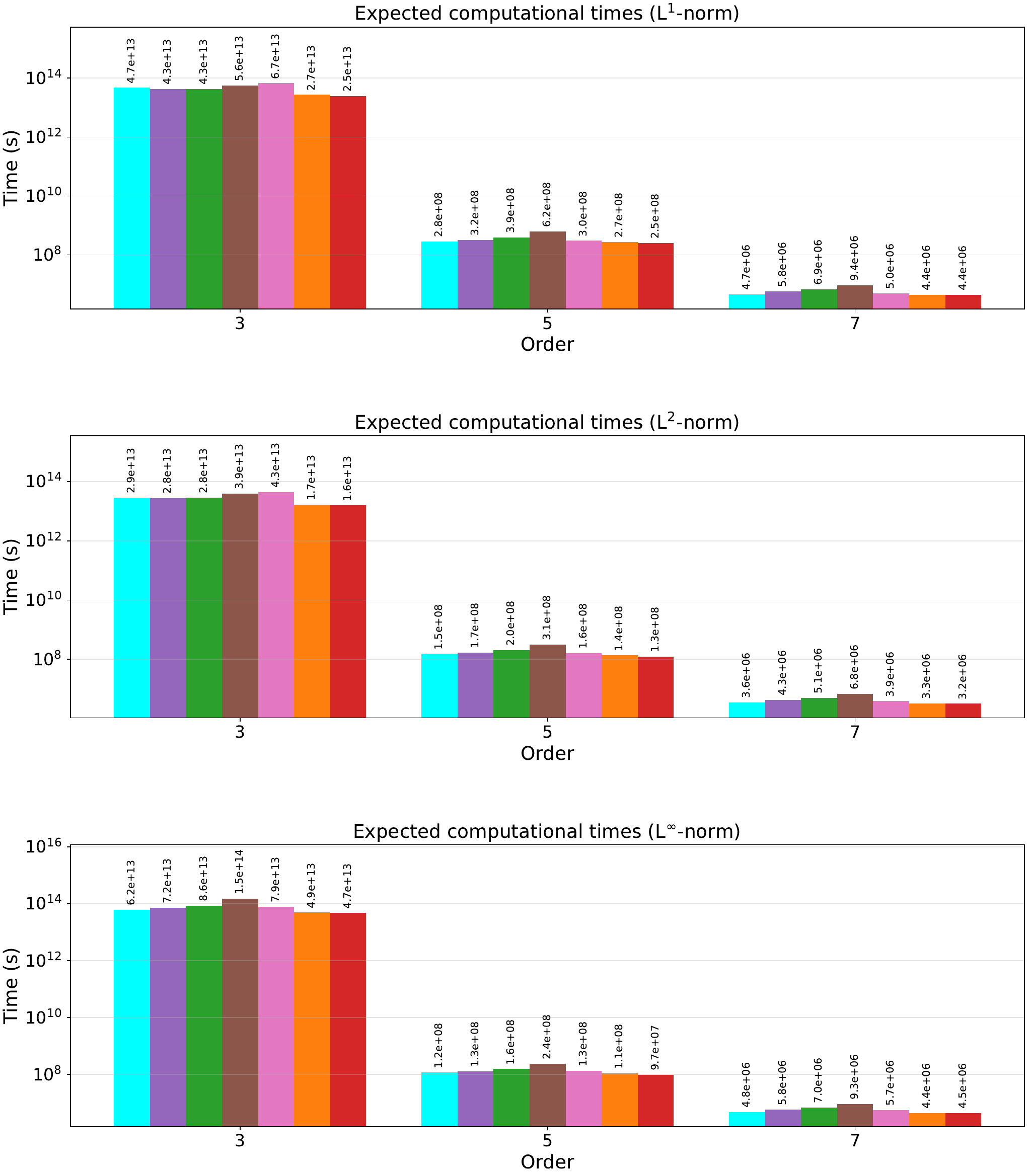}\\
	\includegraphics[width=0.9\linewidth]{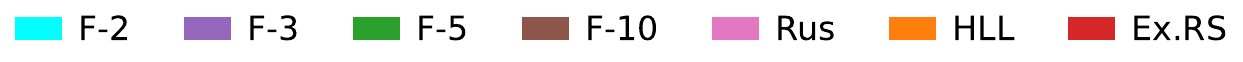}	
	\caption{\RIcolor{Smooth isentropic vortex: Estimated computational times (in seconds) required to achieve a density error of $10^{-16}$ in the $L^1$-, $L^2$- and $L^\infty$-norms for the different numerical fluxes and orders considered. The estimates are obtained from the efficiency curves by linear regression in logarithmic scale. A common legend for the numerical fluxes (distinguished by color) is reported at the bottom.}}
	\label{fig:expected_time_Euler_2d}
\end{figure}

\RIIcolor{Visual representations of the results are reported in Figure~\ref{fig:Euler_2d_unsteady_vortex_WENODeC}.
For readability, just like in the one--dimensional case, we present two sets of plots: the ones on the left column concern only the FORCE-$\alpha$ numerical fluxes, while, the ones on the right concern comparisons between FORCE-$\alpha$ fluxes with references values $\alpha=2,$ 3 and 10 with Rusanov, HLL and exact RS.
Concerning convergence, the expected order of accuracy has been attained for all configurations, and one can once again notice how the difference among the results obtained through different numerical fluxes tends to decrease as the order of accuracy increases, as observed in~\cite{micalizzitoro2024}. For order 7, the lines are rather hard to distinguish, meaning that the error for given mesh refinement is essentially the same.
Concerning efficiency, the best performance is the one obtained through the exact RS, while, the worst one is obtained with FORCE--10.
The exact RS is followed by HLL, however, it is interesting to notice how the related line becomes indistinguishable from the ones of FORCE-2, FORCE-3, and Rusanov for orders 5 and 7.
For such orders, FORCE-5 gives intermediate results between FORCE-10 and the block consisting of HLL, FORCE-2, FORCE-3 and Rusanov.
For order 3, in order of increasing efficiency, we have FORCE-10, Rusanov, FORCE-2, FORCE-5, FORCE-3, HLL and exact RS.
The estimated computational times to reach an error equal to $10^{-16}$, obtained through linear regression of the curves in the efficiency plots, are reported in bar plots in Figure~\ref{fig:expected_time_Euler_2d}.
The results reflect the previous description.
Notice how the relative difference between the estimated computational times of FORCE-2 and of the exact RS goes from 88\% for order 3 to 7\% for order 7.
The difference is supposed to become even smaller for higher order, leading again to the conclusion that FORCE--$\alpha$ numerical fluxes are a very competitive alternative to upwind fluxes within high order frameworks.
In particular, the conclusions are consistent with what obtained in one space dimension: lower values of $\alpha$ are preferable, with $\alpha=2$ and $\alpha=3$ being the best options.
}

\begin{figure}[htbp]
	\centering
	\begin{subfigure}[t]{0.49\textwidth}
		\centering
		\includegraphics[width=\textwidth]{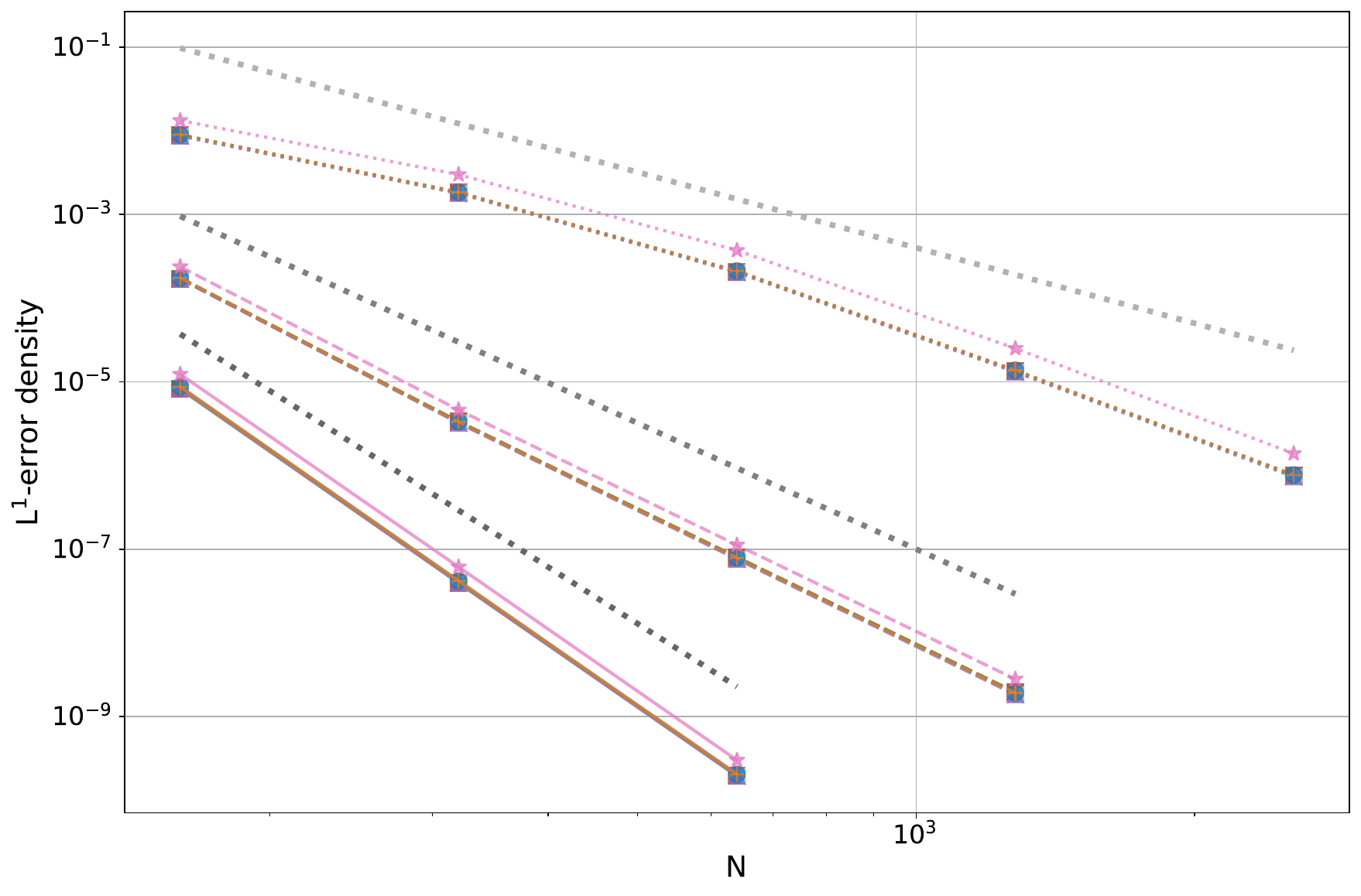}
		\caption{Convergence analysis}
	\end{subfigure}
	\hfill
	\begin{subfigure}[t]{0.49\textwidth}
		\centering
		\includegraphics[width=\textwidth]{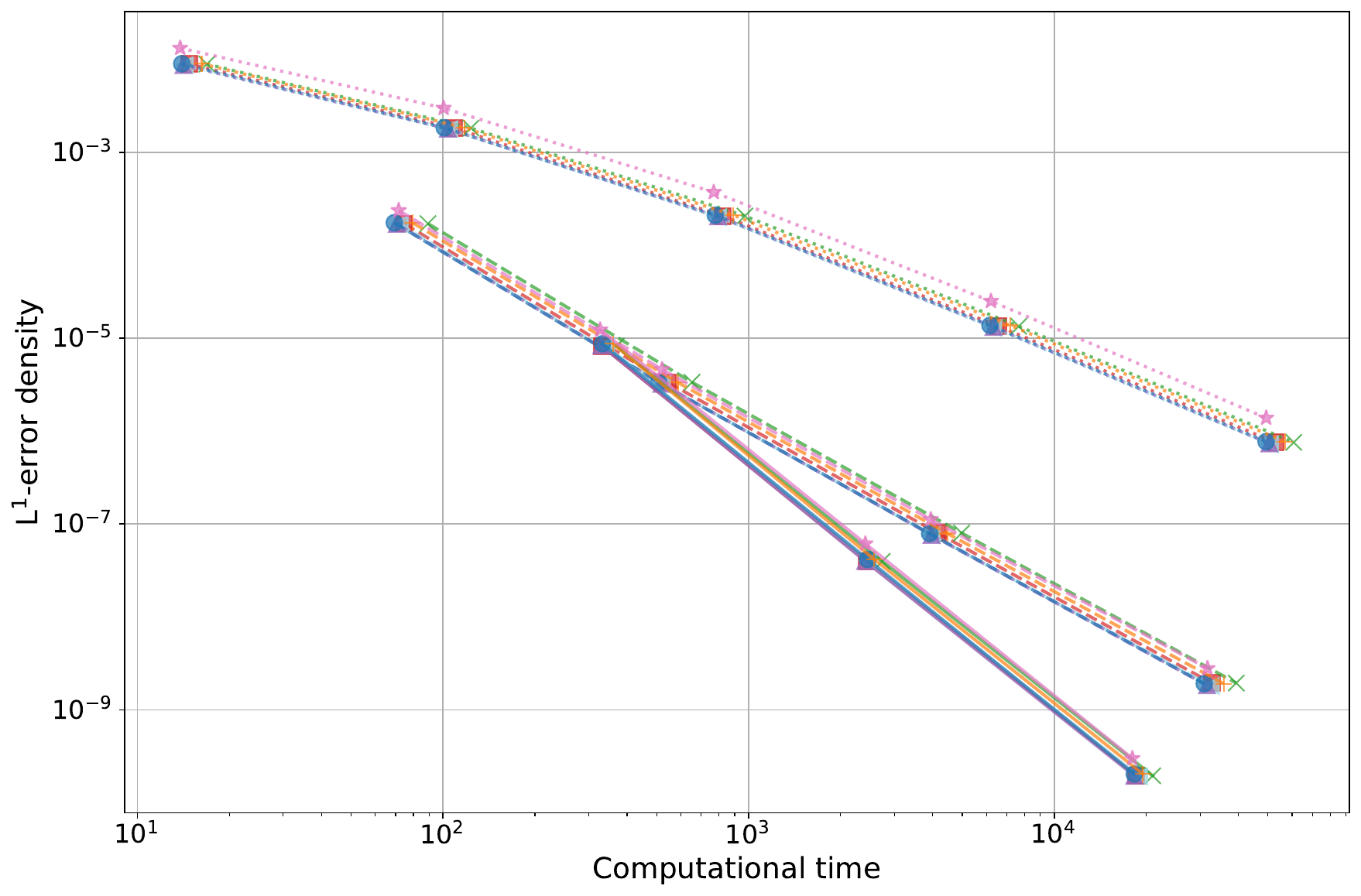}
		\caption{Efficiency analysis}
	\end{subfigure}
	\vspace{0.3cm}
	
	\begin{subfigure}[t]{1.0\textwidth}
		\centering
		\includegraphics[width=0.9\textwidth]{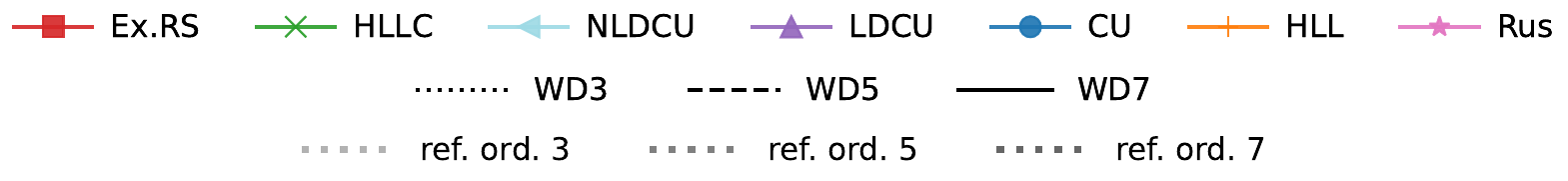}
	\end{subfigure}

	\caption{\RIcolor{Smooth isentropic vortex: Convergence and efficiency analyses for Rusanov~\cite{Rusanov1961}, Harten--Lax--van Leer (HLL)~\cite{harten1983upstream}, Central--Upwind~(CU)~\cite{kurganov2001semidiscrete,kurganov2000new}, Low--Dissipation Central--Upwind~(LDCU)~\cite{KLin}, new Low--Dissipation Central--Upwind~(NLDCU)~\cite{CKX_Ustar}, HLLC~\cite{toro1992restoration,toro1994restoration} and exact Riemann solver~(RS)~\cite{Godunov}.
	On the bottom a common legend for numerical fluxes (distinguished by color and marker) and order (distinguished by linestyle) is reported.}}
	\label{fig:Euler_2d_unsteady_vortex_all_numerical_fluxes_WENODeC}
\end{figure}

\begin{remark}[On the efficiency of the exact RS]
	The results obtained in this work are generally consistent with the ones obtained in~\cite{micalizzitoro2024} as remarked several times throughout the description of the numerical experiments.
	The only registered difference concerns the computational efficiency of the exact RS with respect to HLL in this test.
	As in~\cite{micalizzitoro2024}, we obtain that both exact RS and HLL are more computationally efficient than Rusanov.
	However, here the best efficiency is achieved through the exact RS, while, in~\cite{micalizzitoro2024} it was achieved through HLL.
	We believe that this difference is due to a more careful management of the simulation runs in the present work, which were conducted so as to avoid CPU overload and yield more reliable performance measurements.
	For the sake of completeness, we have carefully re-run the simulations for all numerical fluxes considered in~\cite{micalizzitoro2024}, but FORCE and Lax--Friedrichs, and also for a new low--dissipation Central Upwind (NLDCU) numerical flux recently developed in~\cite{CKX_Ustar}. Lax--Friedrichs and FORCE, i.e., FORCE--1, have been excluded because their first order versions are unstable in multiple space dimensions~\cite{toro2000centred,micalizzitoro2024}.
	Summarizing, we consider the following numerical fluxes: Rusanov~\cite{Rusanov1961}, Harten--Lax--van Leer (HLL)~\cite{harten1983upstream}, Central--Upwind~(CU)~\cite{kurganov2001semidiscrete,kurganov2000new}, Low--Dissipation Central--Upwind~(LDCU)~\cite{KLin}, new Low--Dissipation Central--Upwind~(NLDCU)~\cite{CKX_Ustar}, HLLC~\cite{toro1992restoration,toro1994restoration,toro2019hllc} and exact Riemann solver~(RS)~\cite{Godunov}.
	The results for the $L^1$-norm, obtained with $C_{CFL}:=0.9$, are reported in Figure~\ref{fig:Euler_2d_unsteady_vortex_all_numerical_fluxes_WENODeC}.
	Comparing the results with the ones from~\cite{micalizzitoro2024}, there are no differences in the convergence analysis plot as expected, however, some differences can be noticed in the efficiency analysis plot.
	Much less variation amongst the different numerical fluxes is registered for all orders. In particular, all lines tend to coincide for order 7. This is a further confirmation of the fact that increasing the order of accuracy is associated with a smaller impact of the numerical flux choice, which was the starting point for this work.
	Notice also that exact RS tends to perform as good as HLL and related two-waves incomplete upwind fluxes, such as CU, LDCU and NLDCU, for all orders, despite being classically (and erroneously) considered more computationally expensive.

\end{remark}

\begin{table}[htbp]
	\centering
	\caption{Smooth isentropic vortex: convergence tables for order 3}
	\label{tab:Euler_2d_convergence_tables_WENO_DeC_char_order3}
	\scalebox{0.65}{ 
		\begin{tabular}{c c c c c c c c}
			\toprule
			\multirow{2}{*}{$N$} & \multicolumn{2}{c}{$L^1$-error $\rho$} & \multicolumn{2}{c}{$L^2$-error $\rho$} & \multicolumn{2}{c}{$L^{\infty}$-error $\rho$} & \multirow{2}{*}{CPU Time} \\
			\cmidrule(lr){2-3} \cmidrule(lr){4-5} \cmidrule(lr){6-7}
			& Error & Order & Error & Order & Error & Order & \\
			\midrule
			
			\multicolumn{8}{c}{\textbf{FORCE-2}} \\ 
			\midrule
			160  &   1.298e-02  &  $-$  &   5.668e-03  &  $-$  &   1.108e-02  &  $-$  &   1.408e+01 \\ 
			320  &   2.742e-03  &  2.242  &   1.382e-03  &  2.036  &   3.527e-03  &  1.651  &   1.019e+02 \\ 
			640  &   3.269e-04  &  3.068  &   1.809e-04  &  2.934  &   4.424e-04  &  2.995  &   7.784e+02 \\ 
			1280  &   2.208e-05  &  3.888  &   1.141e-05  &  3.986  &   2.726e-05  &  4.020  &   6.228e+03 \\ 
			2560  &   1.220e-06  &  4.178  &   5.619e-07  &  4.344  &   1.161e-06  &  4.554  &   5.030e+04 \\ 
			\midrule

			\multicolumn{8}{c}{\textbf{FORCE-3}} \\ 
			\midrule
			160  &   1.066e-02  &  $-$  &   4.552e-03  &  $-$  &   8.724e-03  &  $-$  &   1.410e+01 \\ 
			320  &   2.279e-03  &  2.225  &   1.168e-03  &  1.962  &   2.936e-03  &  1.571  &   1.117e+02 \\ 
			640  &   2.680e-04  &  3.088  &   1.507e-04  &  2.954  &   3.916e-04  &  2.907  &   8.561e+02 \\ 
			1280  &   1.770e-05  &  3.921  &   9.214e-06  &  4.032  &   2.311e-05  &  4.082  &   6.654e+03 \\ 
			2560  &   9.730e-07  &  4.185  &   4.506e-07  &  4.354  &   9.553e-07  &  4.597  &   5.301e+04 \\ 
			\midrule

			\multicolumn{8}{c}{\textbf{FORCE-5}} \\ 
			\midrule
			160  &   9.763e-03  &  $-$  &   4.147e-03  &  $-$  &   8.017e-03  &  $-$  &   1.630e+01 \\ 
			320  &   2.048e-03  &  2.253  &   1.053e-03  &  1.977  &   2.695e-03  &  1.573  &   1.280e+02 \\ 
			640  &   2.355e-04  &  3.121  &   1.341e-04  &  2.974  &   3.780e-04  &  2.834  &   9.956e+02 \\ 
			1280  &   1.543e-05  &  3.932  &   8.096e-06  &  4.050  &   2.183e-05  &  4.114  &   7.928e+03 \\ 
			2560  &   8.554e-07  &  4.173  &   3.967e-07  &  4.351  &   8.728e-07  &  4.645  &   6.194e+04 \\ 
			\midrule

			\multicolumn{8}{c}{\textbf{FORCE-10}} \\ 
			\midrule
			160  &   1.005e-02  &  $-$  &   4.164e-03  &  $-$  &   8.060e-03  &  $-$  &   2.090e+01 \\ 
			320  &   2.052e-03  &  2.293  &   1.053e-03  &  1.983  &   2.763e-03  &  1.544  &   1.664e+02 \\ 
			640  &   2.343e-04  &  3.130  &   1.346e-04  &  2.969  &   4.150e-04  &  2.735  &   1.287e+03 \\ 
			1280  &   1.531e-05  &  3.936  &   8.089e-06  &  4.056  &   2.391e-05  &  4.117  &   1.026e+04 \\ 
			2560  &   8.598e-07  &  4.154  &   3.975e-07  &  4.347  &   9.163e-07  &  4.706  &   8.191e+04 \\ 
			\midrule

			\multicolumn{8}{c}{\textbf{Rusanov}} \\ 
			\midrule
			160  &   1.328e-02  &  $-$  &   5.674e-03  &  $-$  &   1.046e-02  &  $-$  &   1.380e+01 \\ 
			320  &   2.999e-03  &  2.147  &   1.494e-03  &  1.926  &   3.512e-03  &  1.575  &   1.005e+02 \\ 
			640  &   3.725e-04  &  3.009  &   2.057e-04  &  2.860  &   4.517e-04  &  2.959  &   7.688e+02 \\ 
			1280  &   2.506e-05  &  3.894  &   1.301e-05  &  3.983  &   2.785e-05  &  4.020  &   6.206e+03 \\ 
			2560  &   1.387e-06  &  4.175  &   6.376e-07  &  4.351  &   1.207e-06  &  4.529  &   4.936e+04 \\ 
			\midrule

			\multicolumn{8}{c}{\textbf{HLL}} \\ 
			\midrule
			160  &   8.995e-03  &  $-$  &   3.640e-03  &  $-$  &   6.567e-03  &  $-$  &   1.623e+01 \\ 
			320  &   1.842e-03  &  2.287  &   9.328e-04  &  1.964  &   2.320e-03  &  1.501  &   1.173e+02 \\ 
			640  &   2.109e-04  &  3.127  &   1.197e-04  &  2.962  &   3.586e-04  &  2.694  &   8.905e+02 \\ 
			1280  &   1.375e-05  &  3.939  &   6.977e-06  &  4.101  &   1.907e-05  &  4.233  &   7.175e+03 \\ 
			2560  &   7.709e-07  &  4.157  &   3.355e-07  &  4.378  &   6.809e-07  &  4.808  &   5.685e+04 \\ 
			\midrule

			\multicolumn{8}{c}{\textbf{exact RS}} \\ 
			\midrule
			160  &   8.925e-03  &  $-$  &   3.635e-03  &  $-$  &   6.555e-03  &  $-$  &   1.480e+01 \\ 
			320  &   1.826e-03  &  2.289  &   9.321e-04  &  1.963  &   2.319e-03  &  1.499  &   1.088e+02 \\ 
			640  &   2.068e-04  &  3.142  &   1.194e-04  &  2.965  &   3.586e-04  &  2.693  &   8.236e+02 \\ 
			1280  &   1.343e-05  &  3.945  &   6.940e-06  &  4.105  &   1.908e-05  &  4.232  &   6.573e+03 \\ 
			2560  &   7.559e-07  &  4.151  &   3.342e-07  &  4.376  &   6.814e-07  &  4.807  &   5.285e+04 \\ 
			\midrule

			\bottomrule
	\end{tabular}}
\end{table}

\begin{table}[htbp]
	\centering
	\caption{Smooth isentropic vortex: convergence tables for order 5}
	\label{tab:Euler_2d_convergence_tables_WENO_DeC_char_order5}
	\scalebox{0.65}{ 
		\begin{tabular}{c c c c c c c c}
			\toprule
			\multirow{2}{*}{$N$} & \multicolumn{2}{c}{$L^1$-error $\rho$} & \multicolumn{2}{c}{$L^2$-error $\rho$} & \multicolumn{2}{c}{$L^{\infty}$-error $\rho$} & \multirow{2}{*}{CPU Time} \\
			\cmidrule(lr){2-3} \cmidrule(lr){4-5} \cmidrule(lr){6-7}
			& Error & Order & Error & Order & Error & Order & \\
			\midrule
			
			\multicolumn{8}{c}{\textbf{FORCE-2}} \\ 
			\midrule
			160  &   2.448e-04  &  $-$  &   1.048e-04  &  $-$  &   1.516e-04  &  $-$  &   7.006e+01 \\ 
			320  &   4.680e-06  &  5.709  &   1.865e-06  &  5.812  &   2.022e-06  &  6.229  &   5.127e+02 \\ 
			640  &   1.099e-07  &  5.412  &   4.140e-08  &  5.493  &   3.952e-08  &  5.677  &   3.924e+03 \\ 
			1280  &   2.717e-09  &  5.338  &   1.040e-09  &  5.315  &   1.073e-09  &  5.202  &   3.119e+04 \\ 
			\midrule

			\multicolumn{8}{c}{\textbf{FORCE-3}} \\ 
			\midrule
			160  &   2.030e-04  &  $-$  &   8.900e-05  &  $-$  &   1.396e-04  &  $-$  &   6.999e+01 \\ 
			320  &   4.138e-06  &  5.617  &   1.664e-06  &  5.741  &   2.032e-06  &  6.103  &   5.558e+02 \\ 
			640  &   9.756e-08  &  5.407  &   3.698e-08  &  5.491  &   3.772e-08  &  5.751  &   4.260e+03 \\ 
			1280  &   2.364e-09  &  5.367  &   9.051e-10  &  5.353  &   9.902e-10  &  5.251  &   3.314e+04 \\ 
			\midrule

			\multicolumn{8}{c}{\textbf{FORCE-5}} \\ 
			\midrule
			160  &   1.915e-04  &  $-$  &   8.630e-05  &  $-$  &   1.453e-04  &  $-$  &   8.308e+01 \\ 
			320  &   4.000e-06  &  5.581  &   1.631e-06  &  5.725  &   2.304e-06  &  5.978  &   6.515e+02 \\ 
			640  &   9.519e-08  &  5.393  &   3.644e-08  &  5.485  &   3.969e-08  &  5.859  &   5.005e+03 \\ 
			1280  &   2.307e-09  &  5.367  &   8.825e-10  &  5.368  &   1.023e-09  &  5.278  &   3.940e+04 \\ 
			\midrule

			\multicolumn{8}{c}{\textbf{FORCE-10}} \\ 
			\midrule
			160  &   2.047e-04  &  $-$  &   9.447e-05  &  $-$  &   1.776e-04  &  $-$  &   1.038e+02 \\ 
			320  &   4.442e-06  &  5.526  &   1.852e-06  &  5.672  &   2.996e-06  &  5.890  &   8.347e+02 \\ 
			640  &   1.076e-07  &  5.367  &   4.176e-08  &  5.471  &   5.095e-08  &  5.878  &   6.454e+03 \\ 
			1280  &   2.624e-09  &  5.358  &   1.002e-09  &  5.381  &   1.218e-09  &  5.386  &   5.148e+04 \\ 
			\midrule

			\multicolumn{8}{c}{\textbf{Rusanov}} \\ 
			\midrule
			160  &   2.374e-04  &  $-$  &   1.008e-04  &  $-$  &   1.429e-04  &  $-$  &   7.161e+01 \\ 
			320  &   4.614e-06  &  5.685  &   1.829e-06  &  5.783  &   2.009e-06  &  6.152  &   5.209e+02 \\ 
			640  &   1.119e-07  &  5.365  &   4.196e-08  &  5.446  &   4.115e-08  &  5.609  &   3.943e+03 \\ 
			1280  &   2.793e-09  &  5.325  &   1.068e-09  &  5.296  &   1.146e-09  &  5.167  &   3.173e+04 \\ 
			\midrule

			\multicolumn{8}{c}{\textbf{HLL}} \\ 
			\midrule
			160  &   1.748e-04  &  $-$  &   7.991e-05  &  $-$  &   1.366e-04  &  $-$  &   7.954e+01 \\ 
			320  &   3.339e-06  &  5.710  &   1.363e-06  &  5.873  &   2.120e-06  &  6.010  &   5.864e+02 \\ 
			640  &   7.905e-08  &  5.401  &   2.982e-08  &  5.515  &   3.594e-08  &  5.882  &   4.473e+03 \\ 
			1280  &   1.912e-09  &  5.369  &   7.124e-10  &  5.387  &   8.342e-10  &  5.429  &   3.583e+04 \\ 
			\midrule

			\multicolumn{8}{c}{\textbf{exact RS}} \\ 
			\midrule
			160  &   1.711e-04  &  $-$  &   7.877e-05  &  $-$  &   1.363e-04  &  $-$  &   7.455e+01 \\ 
			320  &   3.351e-06  &  5.674  &   1.362e-06  &  5.854  &   2.120e-06  &  6.006  &   5.432e+02 \\ 
			640  &   7.991e-08  &  5.390  &   3.014e-08  &  5.498  &   3.594e-08  &  5.882  &   4.151e+03 \\ 
			1280  &   1.942e-09  &  5.363  &   7.257e-10  &  5.376  &   8.342e-10  &  5.429  &   3.276e+04 \\ 
			\midrule

			\bottomrule
	\end{tabular}}
\end{table}

\begin{table}[htbp]
	\centering
	\caption{Smooth isentropic vortex: convergence tables for order 7}
	\label{tab:Euler_2d_convergence_tables_WENO_DeC_char_order7}
	\scalebox{0.65}{ 
		\begin{tabular}{c c c c c c c c}
			\toprule
			\multirow{2}{*}{$N$} & \multicolumn{2}{c}{$L^1$-error $\rho$} & \multicolumn{2}{c}{$L^2$-error $\rho$} & \multicolumn{2}{c}{$L^{\infty}$-error $\rho$} & \multirow{2}{*}{CPU Time} \\
			\cmidrule(lr){2-3} \cmidrule(lr){4-5} \cmidrule(lr){6-7}
			& Error & Order & Error & Order & Error & Order & \\
			\midrule
			
			\multicolumn{8}{c}{\textbf{FORCE-2}} \\ 
			\midrule
			160  &   1.241e-05  &  $-$  &   4.744e-06  &  $-$  &   6.317e-06  &  $-$  &   3.316e+02 \\ 
			320  &   5.796e-08  &  7.742  &   2.421e-08  &  7.614  &   2.918e-08  &  7.758  &   2.461e+03 \\ 
			640  &   2.758e-10  &  7.715  &   1.162e-10  &  7.703  &   1.969e-10  &  7.211  &   1.823e+04 \\ 
			\midrule

			\multicolumn{8}{c}{\textbf{FORCE-3}} \\ 
			\midrule
			160  &   1.007e-05  &  $-$  &   4.020e-06  &  $-$  &   5.738e-06  &  $-$  &   3.272e+02 \\ 
			320  &   4.941e-08  &  7.671  &   2.099e-08  &  7.581  &   2.728e-08  &  7.717  &   2.611e+03 \\ 
			640  &   2.356e-10  &  7.712  &   9.967e-11  &  7.718  &   1.751e-10  &  7.284  &   1.996e+04 \\ 
			\midrule

			\multicolumn{8}{c}{\textbf{FORCE-5}} \\ 
			\midrule
			160  &   9.279e-06  &  $-$  &   3.909e-06  &  $-$  &   5.931e-06  &  $-$  &   3.933e+02 \\ 
			320  &   4.558e-08  &  7.669  &   2.004e-08  &  7.608  &   2.985e-08  &  7.634  &   3.210e+03 \\ 
			640  &   2.212e-10  &  7.687  &   9.489e-11  &  7.723  &   1.746e-10  &  7.417  &   2.405e+04 \\ 
			\midrule

			\multicolumn{8}{c}{\textbf{FORCE-10}} \\ 
			\midrule
			160  &   9.558e-06  &  $-$  &   4.316e-06  &  $-$  &   6.909e-06  &  $-$  &   4.911e+02 \\ 
			320  &   4.798e-08  &  7.638  &   2.230e-08  &  7.597  &   3.783e-08  &  7.513  &   4.141e+03 \\ 
			640  &   2.407e-10  &  7.639  &   1.059e-10  &  7.718  &   1.970e-10  &  7.585  &   3.011e+04 \\ 
			\midrule

			\multicolumn{8}{c}{\textbf{Rusanov}} \\ 
			\midrule
			160  &   1.232e-05  &  $-$  &   4.622e-06  &  $-$  &   6.184e-06  &  $-$  &   3.267e+02 \\ 
			320  &   6.129e-08  &  7.651  &   2.524e-08  &  7.517  &   3.064e-08  &  7.657  &   2.409e+03 \\ 
			640  &   2.999e-10  &  7.675  &   1.272e-10  &  7.633  &   2.375e-10  &  7.011  &   1.801e+04 \\ 
			\midrule

			\multicolumn{8}{c}{\textbf{HLL}} \\ 
			\midrule
			160  &   8.704e-06  &  $-$  &   3.632e-06  &  $-$  &   5.512e-06  &  $-$  &   3.599e+02 \\ 
			320  &   4.165e-08  &  7.707  &   1.818e-08  &  7.642  &   2.653e-08  &  7.699  &   2.577e+03 \\ 
			640  &   2.038e-10  &  7.675  &   8.667e-11  &  7.712  &   1.574e-10  &  7.397  &   1.949e+04 \\ 
			\midrule

			\multicolumn{8}{c}{\textbf{exact RS}} \\ 
			\midrule
			160  &   8.236e-06  &  $-$  &   3.560e-06  &  $-$  &   5.394e-06  &  $-$  &   3.307e+02 \\ 
			320  &   3.956e-08  &  7.702  &   1.765e-08  &  7.656  &   2.653e-08  &  7.668  &   2.434e+03 \\ 
			640  &   1.949e-10  &  7.665  &   8.401e-11  &  7.715  &   1.574e-10  &  7.397  &   1.851e+04 \\ 
			\midrule

			\bottomrule
	\end{tabular}}
\end{table}

\subsubsection{Shock--vortex interaction}\label{sec:shock_vortex}
This test, considered in~\cite{DZLD,geisenhofer2019discontinuous,RCD,abgrall2025dual}, consists of a moving vortex interacting with a stationary shock. Such an interaction causes the formation of several complex structures, making this a perfect benchmark for the validation of (very) high order schemes.
The initialization is described in~\cite{geisenhofer2019discontinuous,RCD}, however, we recall it also here for convenience.
The set-up of the initial condition is outlined in Figure~\ref{fig:shock_vortex_setup}. At the initial time, the spatial domain $\Omega:=[0,2]\times[0,1]$ is divided into four regions, referred to as I, II, III, and IV.
In particular, the constant data in regions III and IV, separated by a vertical line at $x=0.5$, constitute the stationary shock, with the vortex initially occupying regions II and I within region III.
Given the data in region~III, 
\begin{align}
	\rho_{\rm III}:=1,\quad u_{\rm III}:=\sqrt{\gamma}M_s,\quad v_{\rm III}:=0,\quad p_{\rm III}:=1,
\end{align}
and given the Mach number $M_s:=1.5$ of the stationary shock, thanks to the Rankine-Hugoniot conditions~\cite[Section 3.1.3]{ToroBook}, one can explicitly compute the state in region IV
\begin{align}
	\begin{aligned}
		\rho_{\rm IV}&:=\frac{(\gamma+1)M_s^2}{(\gamma-1)M_s^2+2}\rho_{\rm III},\\ 
		u_{\rm IV}&:=\frac{(\gamma-1)M_s^2+2}{(\gamma+1)M_s^2}u_{\rm III},
		\\ 
		v_{\rm IV}&:=0,\\ 
		p_{\rm IV}&:=\frac{2\gamma M_s^2-(\gamma-1)}{\gamma+1}p_{\rm III}.
	\end{aligned}
\end{align}
The vortex, with Mach number $M_v:=0.9$, initially occupies regions I and II, corresponding to a circle centered in $(x_c,y_c):=(0.25,0.5)$.
More in detail, region I is a circle with radius $a:=0.075$, while, region II is a concentrical annulus with inner and outer radii equal to $a$ and $b:=0.175$ respectively.
The velocity field in such regions is prescribed, in terms of the radial coordinates $(r,\vartheta)$ with respect to the center $(x_c,y_c)$, as
\begin{align}
	\begin{aligned}
		u(r,\vartheta)&:=u_{\rm III}-v_\vartheta(r)\sin\vartheta,\\ v(r,\vartheta)&:=v_{\rm III}+v_\vartheta(r)\cos{\vartheta},
	\end{aligned}
\end{align}
with
\begin{align}
	v_\vartheta(r):=\begin{cases}
		v_m\,\frac{r}{a},&r\le a,\\
		v_m\,\frac{a}{a^2-b^2}\left(r-\frac{b^2}{r}\right),&a<r<b,\\
		0,&r\ge b,
	\end{cases}
\end{align}
where $v_m:=M_v\sqrt{\gamma}$.

At this point, the last values to be prescribed in order to complete the initial conditions are density and pressure within the vortex. The balance between pressure gradients and centripetal force~\cite{RCD} gives

\begin{align}
	p&:=p_{\rm III}\left(\frac{T}{T_{\rm III}}\right)^{\frac{\gamma}{\gamma-1}},\\
	\rho&:=\rho_{\rm III}\left(\frac{T}{T_{\rm III}}\right)^{\frac{1}{\gamma-1}},
\end{align}
with $T_{\rm III}=\frac{p_{\rm III}}{\rho_{\rm III}R}$ being the value of the temperature in region III. 
The constant $R:=287$ J/(kg $\cdot$ K) represents the specific gas constant, while, $T(r)$ is the temperature inside the vortex obtained by solving the following ODE
\begin{equation}
	\frac{d}{dr}T(r)=\frac{\gamma-1}{R\gamma}\frac{v_\vartheta^2(r)}{r}.
\end{equation}
Simple computations give
\begin{equation}
	T(r):=\left\{\begin{aligned}
		&A+\frac{\gamma-1}{R\gamma}\frac{v_m^2}{a^2}\frac{r^2}{2},&&r\le a,\\
		&B+\frac{\gamma-1}{R\gamma}v_m^2\frac{a^2}{(a^2-b^2)^2}\left(\frac{r^2}{2}-2b^2\ln r-b^4\frac{r^{-2}}{2}\right),&&a<r<b,\\
		&T_{\rm III},&&r\ge b,
	\end{aligned}\right.
\end{equation}
with
\begin{equation}
	\begin{aligned}
		A&:=B+\frac{\gamma-1}{R\gamma}v_m^2\frac{a^2}{(a^2-b^2)^2}\left(\frac{a^2}{2}-2b^2\ln a-b^4\frac{a^{-2}}{2}\right)-
		\frac{\gamma-1}{R\gamma}\frac{v_m^2}{2},	\\
		B&:=T_{\rm III}-\frac{\gamma-1}{R\gamma}v_m^2\frac{a^2}{(a^2-b^2)^2}\left(\frac{b^2}{2}-2b^2\ln b-b^4\frac{b^{-2}}{2}\right).
	\end{aligned}
\end{equation}
Finally, we set inflow boundary conditions on the left side of the boundary and outflow boundary conditions on the other ones.
For information concerning the phenomenology of the problem and related wave patterns, the interested reader is referred to~\cite{RCD} and references therein.

\begin{figure}[ht!]
	\centering
	\begin{tikzpicture}[scale=8.5]
		\def\angleA{145}   
		\def\angleB{-60}    
		\def\angleR{60}   
		\def\rA{0.09}
		\def\rB{0.18}
		\def\rR{0.16}
		\def\bottom{0.2}
		\def\top{0.8}
		\coordinate (C) at (0.25,0.5);
		\draw[thick] (0,\bottom) rectangle (1.35,\top);
		\draw[very thick] (0.5,\bottom) -- (0.5,\top);
		\draw[line width=0.5mm] (C) circle (\rA); 
		\draw[line width=0.5mm] (C) circle (\rB); 
		\draw[dashed] (C) circle (\rR);  
		\draw[->] (C) -- ++({\angleA}:\rA)
		node[pos=0.4, anchor=south] {\(a\)};
		\draw[->] (C) -- ++({\angleB}:\rB)
		node[pos=0.8, anchor=south] {\(b\)};
		\draw[->] (C) -- ++({\angleR}:\rR)
		node[pos=0.6, anchor=south] {\(r\)};
		\draw[dashed] (C) -- (0.25, \bottom); 
		\draw[dashed] (C) -- (0, 0.5);  
		\node[anchor=east] at (0,\top) {$y=1$};
		\node[anchor=east] at (0,\bottom) {$y=0$};
		\node[anchor=north] at (0,\bottom) {$x=0$};
		\node[anchor=north] at (0.5,\bottom) {$x=0.5$};
		\node[anchor=north] at (1.35,\bottom) {$x=2$};
		\node[anchor=south] at (0.5,\top+0.03) {stationary};
		\node[anchor=south] at (0.5,\top) {shock};
		\node[anchor=north] at (0.25,\bottom) {$x=0.25$};
		\node[anchor=east] at (0,0.5) {$y=0.5$};
		\node[scale=0.8] at (0.219, 0.465) {\textbf{I}};
		\node[scale=0.8] at (0.175, 0.395) {\textbf{II}};
		\node[scale=0.8] at (0.11, 0.29) {\textbf{III}};
		\node[scale=0.8] at (0.9, 0.5) {\textbf{IV}};
		\draw[dashed] (0.25,0.5) -- (0.25+\rB,0.5); 
		\def\rTheta{0.11} 
		\draw[->] ($(C)+(\rTheta,0)$) arc[start angle=0, end angle=\angleR, radius=\rTheta];
		\node at ($(C)+({\angleR/2}:1.2*\rTheta)$) {\(\vartheta\)};		
	\end{tikzpicture}
	\caption{\sf Shock--vortex interaction: Sketch of the initial condition set-up}
	\label{fig:shock_vortex_setup}
\end{figure}

\begin{figure}[htbp]
	\centering
	\begin{subfigure}{0.49\textwidth}
		\includegraphics[width=\linewidth]{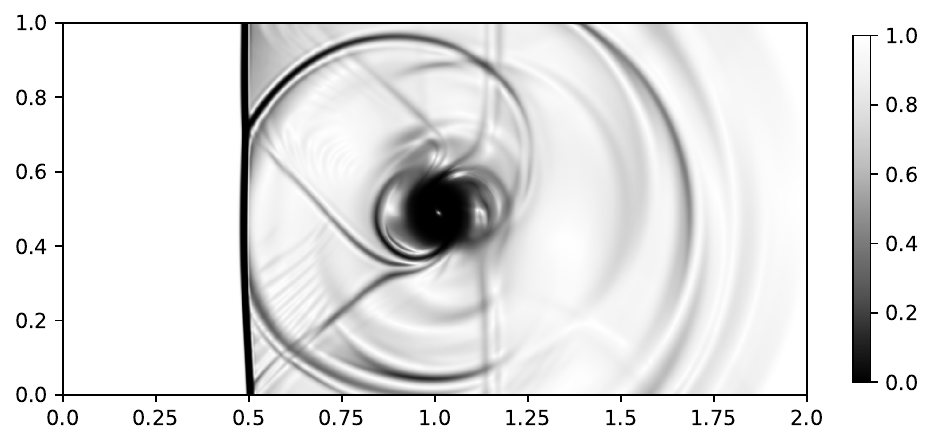}
		\caption{FORCE--$2$}
	\end{subfigure}
	\begin{subfigure}{0.49\textwidth}
		\includegraphics[width=\linewidth]{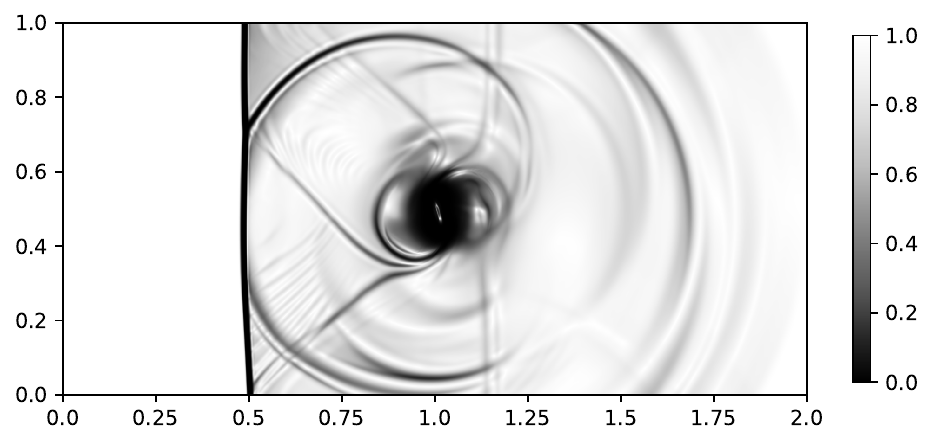}
		\caption{FORCE--$3$}
	\end{subfigure}
	\begin{subfigure}{0.49\textwidth}
		\includegraphics[width=\linewidth]{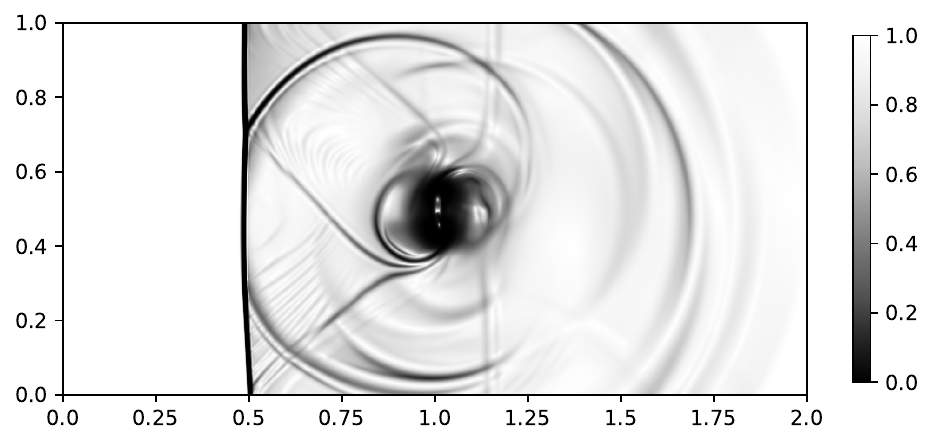}
		\caption{FORCE--$5$}
	\end{subfigure}
	\begin{subfigure}{0.49\textwidth}
		\includegraphics[width=\linewidth]{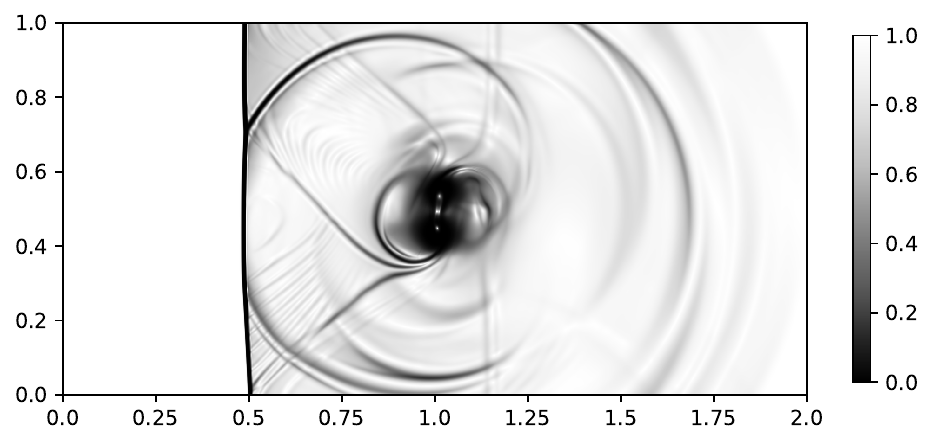}
		\caption{FORCE--$10$}
	\end{subfigure}
	
	\begin{subfigure}{0.49\textwidth}
		\includegraphics[width=\linewidth]{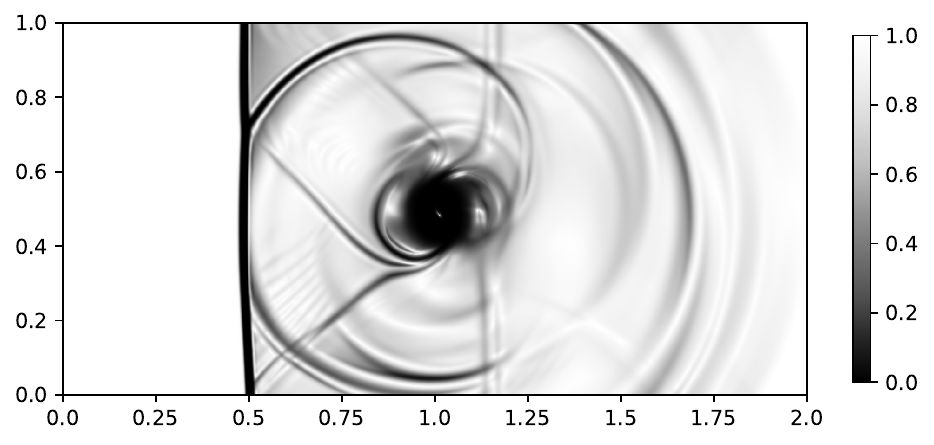}
		\caption{Rusanov}
	\end{subfigure}
	\begin{subfigure}{0.49\textwidth}
		\includegraphics[width=\linewidth]{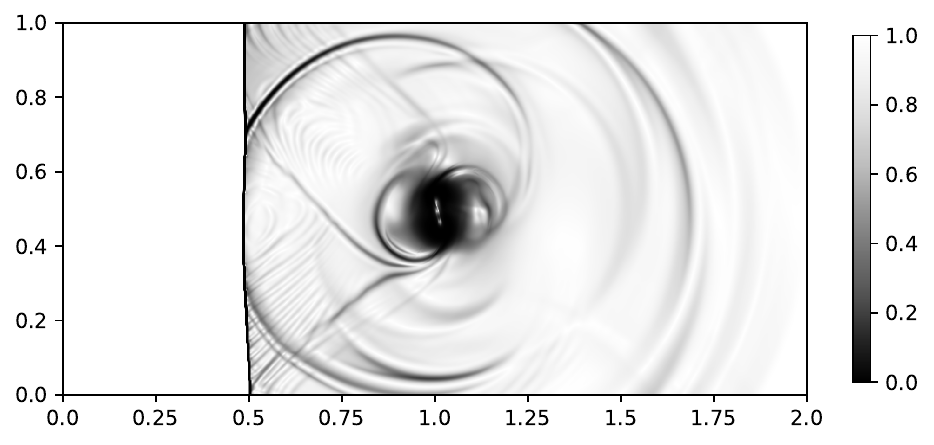}
		\caption{HLL}
	\end{subfigure}
	\begin{subfigure}{0.49\textwidth}
		\includegraphics[width=\linewidth]{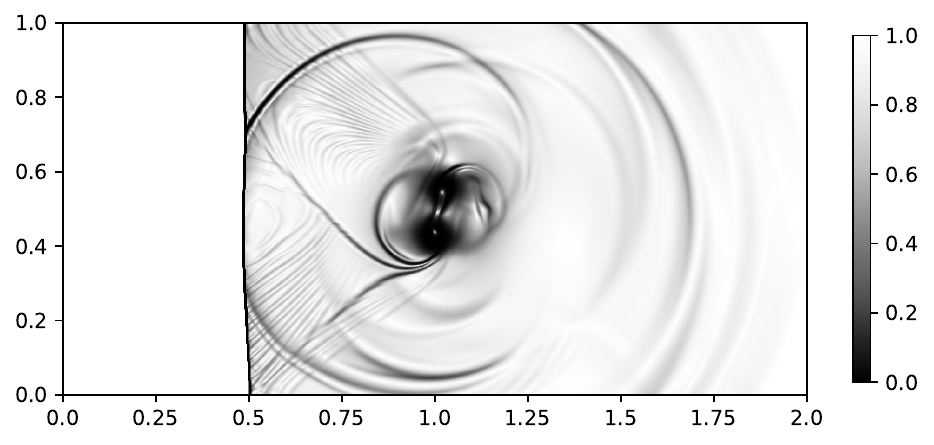}
		\caption{Exact RS}
	\end{subfigure}
	
	\caption{Shock--vortex interaction: Schlieren plots of the density profile obtained over a mesh with $800\times 401$ elements for order 3 with $\sigma_{CFL}:=0.9$}
	\label{fig:shockvortex_order3}
\end{figure}

\begin{figure}[htbp]
	\centering
	\begin{subfigure}{0.49\textwidth}
		\includegraphics[width=\linewidth]{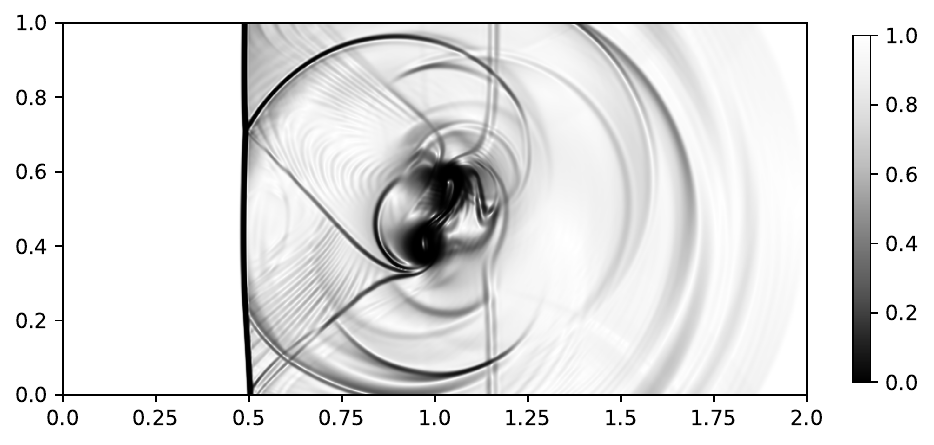}
		\caption{FORCE--$2$}
	\end{subfigure}
	\begin{subfigure}{0.49\textwidth}
		\includegraphics[width=\linewidth]{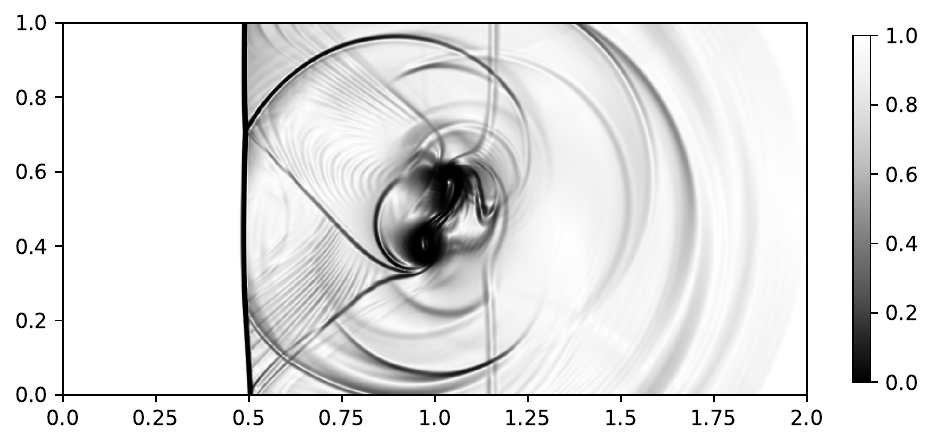}
		\caption{FORCE--$3$}
	\end{subfigure}
	\begin{subfigure}{0.49\textwidth}
		\includegraphics[width=\linewidth]{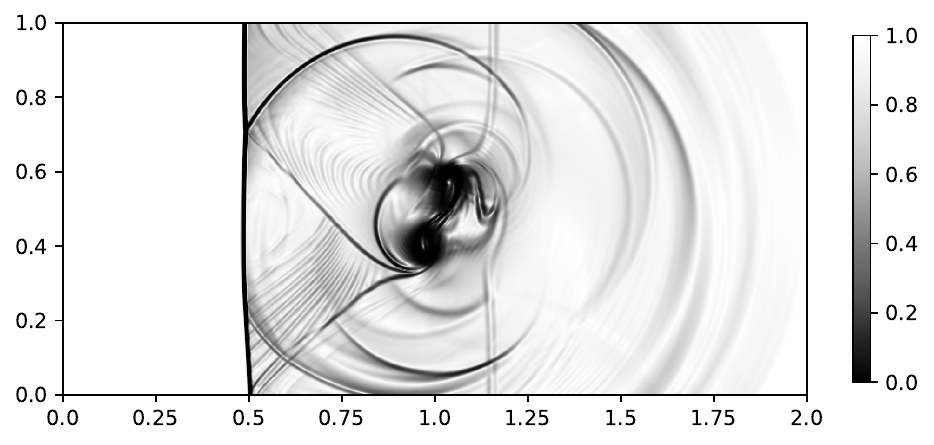}
		\caption{FORCE--$5$}
	\end{subfigure}
	\begin{subfigure}{0.49\textwidth}
		\includegraphics[width=\linewidth]{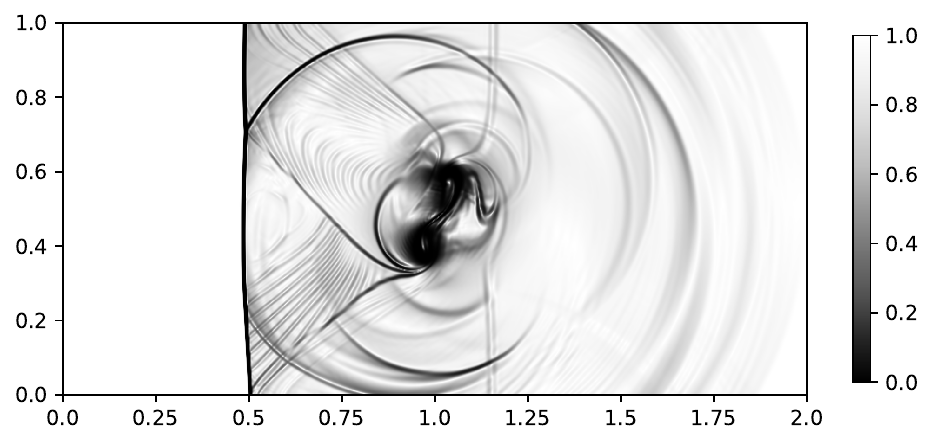}
		\caption{FORCE--$10$}
	\end{subfigure}
	
	\begin{subfigure}{0.49\textwidth}
		\includegraphics[width=\linewidth]{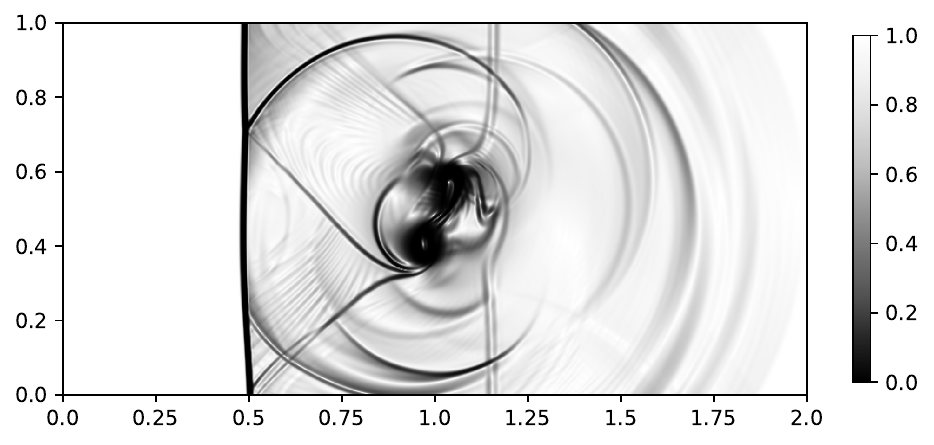}
		\caption{Rusanov}
	\end{subfigure}
	\begin{subfigure}{0.49\textwidth}
		\includegraphics[width=\linewidth]{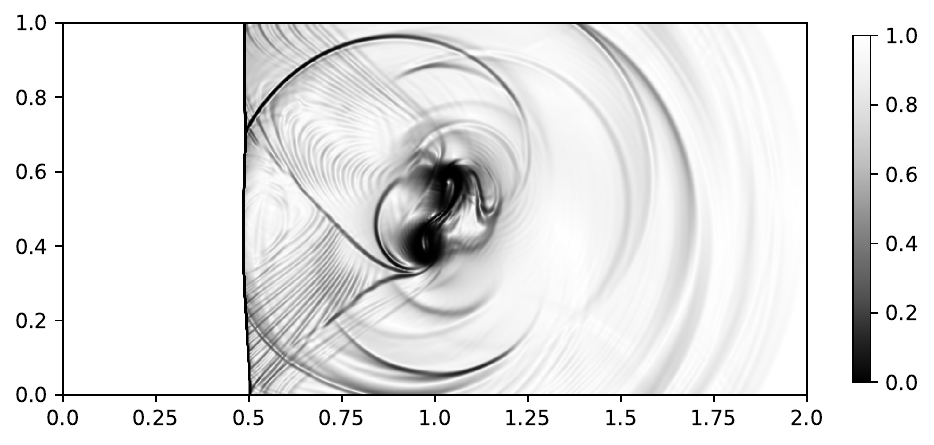}
		\caption{HLL}
	\end{subfigure}
	\begin{subfigure}{0.49\textwidth}
		\includegraphics[width=\linewidth]{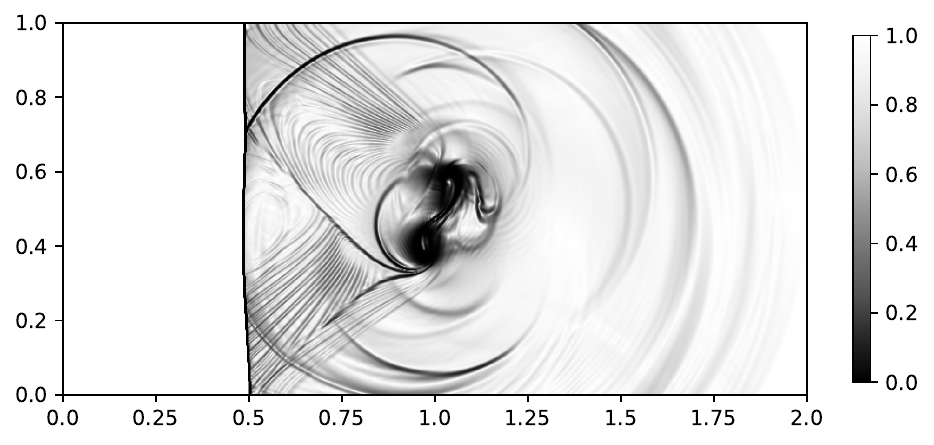}
		\caption{Exact RS}
	\end{subfigure}
	
	\caption{Shock--vortex interaction: Schlieren plots of the density profile obtained over a mesh with $800\times 401$ elements for order 5 with $\sigma_{CFL}:=0.9$}
	\label{fig:shockvortex_order5}
\end{figure}

\begin{figure}[htbp]
	\centering
	\begin{subfigure}{0.49\textwidth}
		\includegraphics[width=\linewidth]{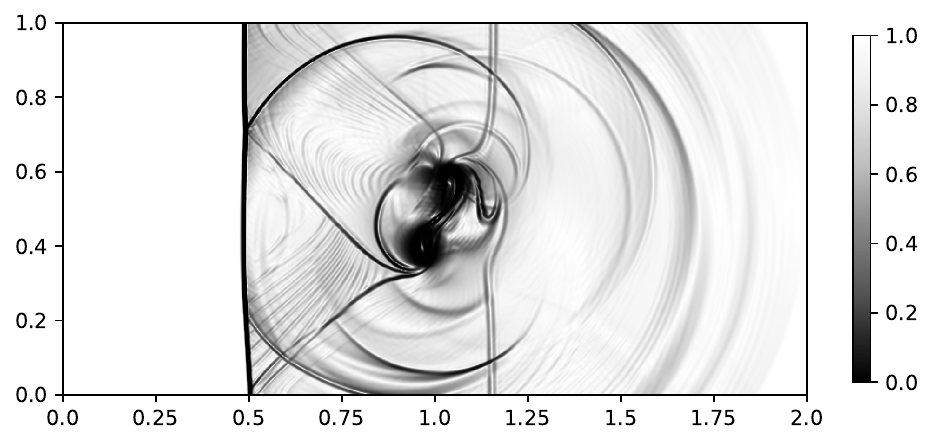}
		\caption{FORCE--$2$}
	\end{subfigure}
	\begin{subfigure}{0.49\textwidth}
		\includegraphics[width=\linewidth]{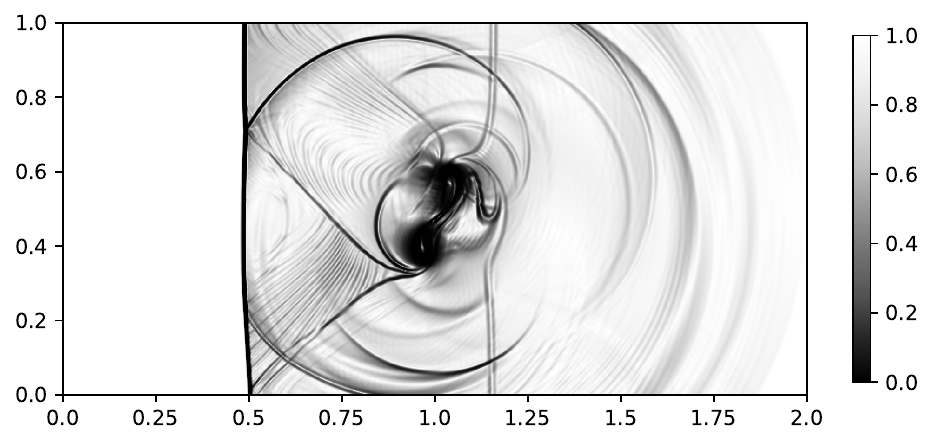}
		\caption{FORCE--$3$}
	\end{subfigure}
	\begin{subfigure}{0.49\textwidth}
		\includegraphics[width=\linewidth]{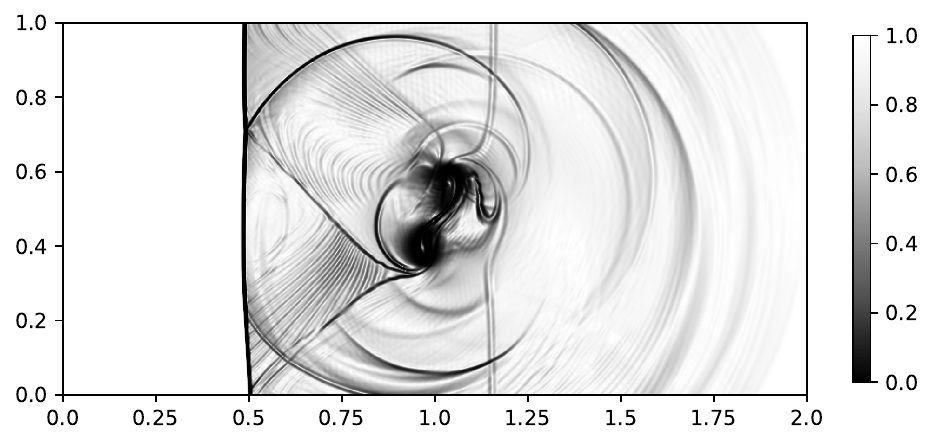}
		\caption{FORCE--$5$}
	\end{subfigure}
	\begin{subfigure}{0.49\textwidth}
		\includegraphics[width=\linewidth]{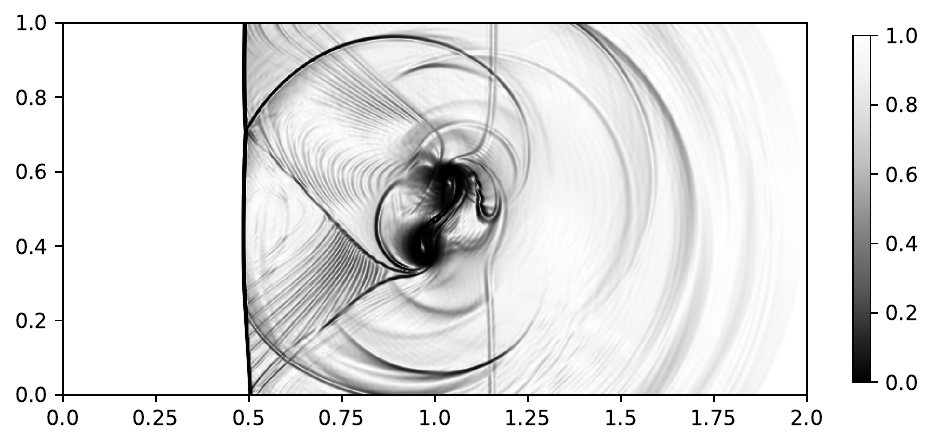}
		\caption{FORCE--$10$}
	\end{subfigure}
	
	\begin{subfigure}{0.49\textwidth}
		\includegraphics[width=\linewidth]{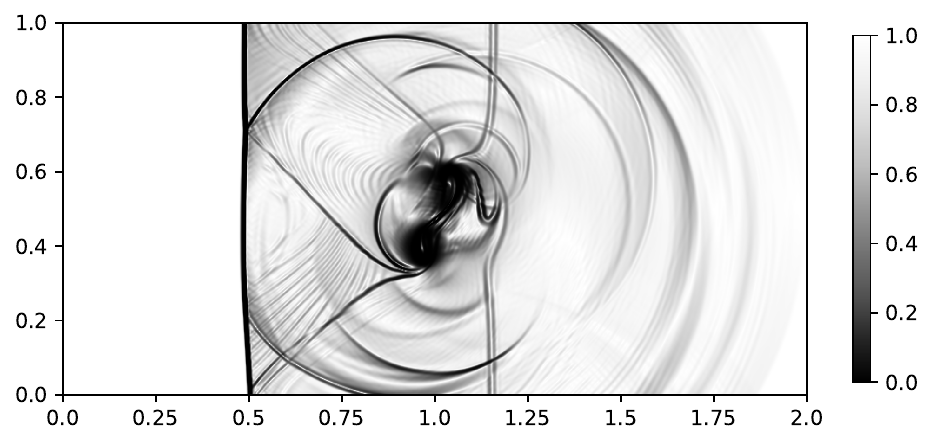}
		\caption{Rusanov}
	\end{subfigure}
	\begin{subfigure}{0.49\textwidth}
		\includegraphics[width=\linewidth]{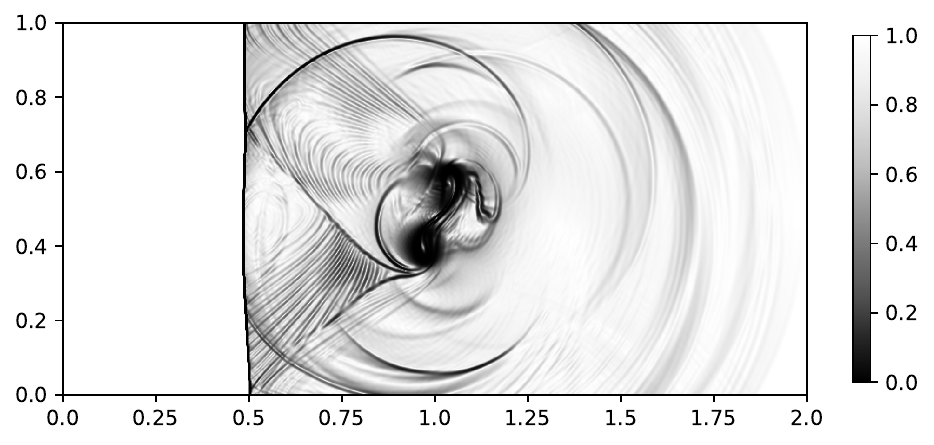}
		\caption{HLL}
	\end{subfigure}
	\begin{subfigure}{0.49\textwidth}
		\includegraphics[width=\linewidth]{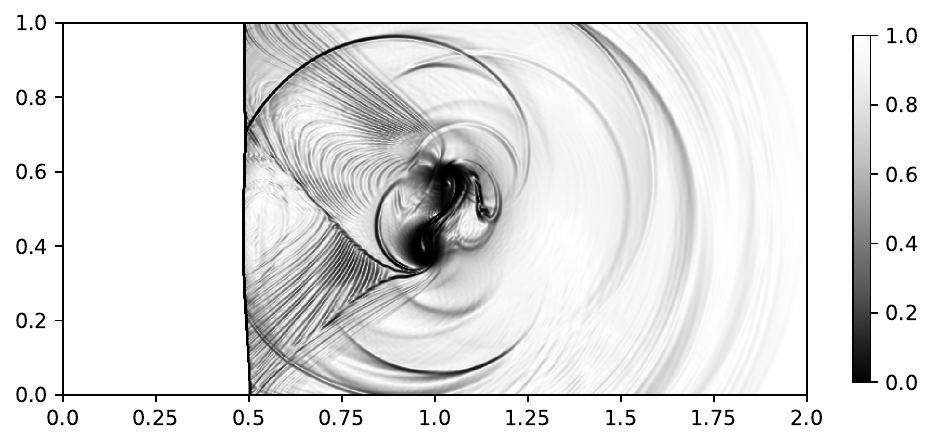}
		\caption{Exact RS}
	\end{subfigure}
	
	\caption{Shock--vortex interaction: Schlieren plots of the density profile obtained over a mesh with $800\times 401$ elements for order 7 with $\sigma_{CFL}:=0.9$}
	\label{fig:shockvortex_order7}
\end{figure}

The two--dimensional density results obtained with $\sigma_{CFL}:=0.9$ at final time $T_f:=0.69$ over meshes of $800\times 401$ elements with FORCE-$\alpha$, Rusanov, HLL and exact RS numerical fluxes are displayed in Figures~\ref{fig:shockvortex_order3},~\ref{fig:shockvortex_order5} and~\ref{fig:shockvortex_order7} for orders 3, 5 and 7 respectively.
In particular, we report Schlieren plots of the density considering the following quantity
\begin{equation}
	\exp{\left(-\frac{K\|\nabla_{\uvec{x}}\rho\|_2}{\max\|\nabla_{\uvec{x}}\rho\|_2}\right)},
\end{equation}
with $K:=80$.

\begin{figure}[htbp]
	\centering
	\begin{subfigure}[b]{1.0\textwidth}
		\centering
		\includegraphics[width=0.7\textwidth]{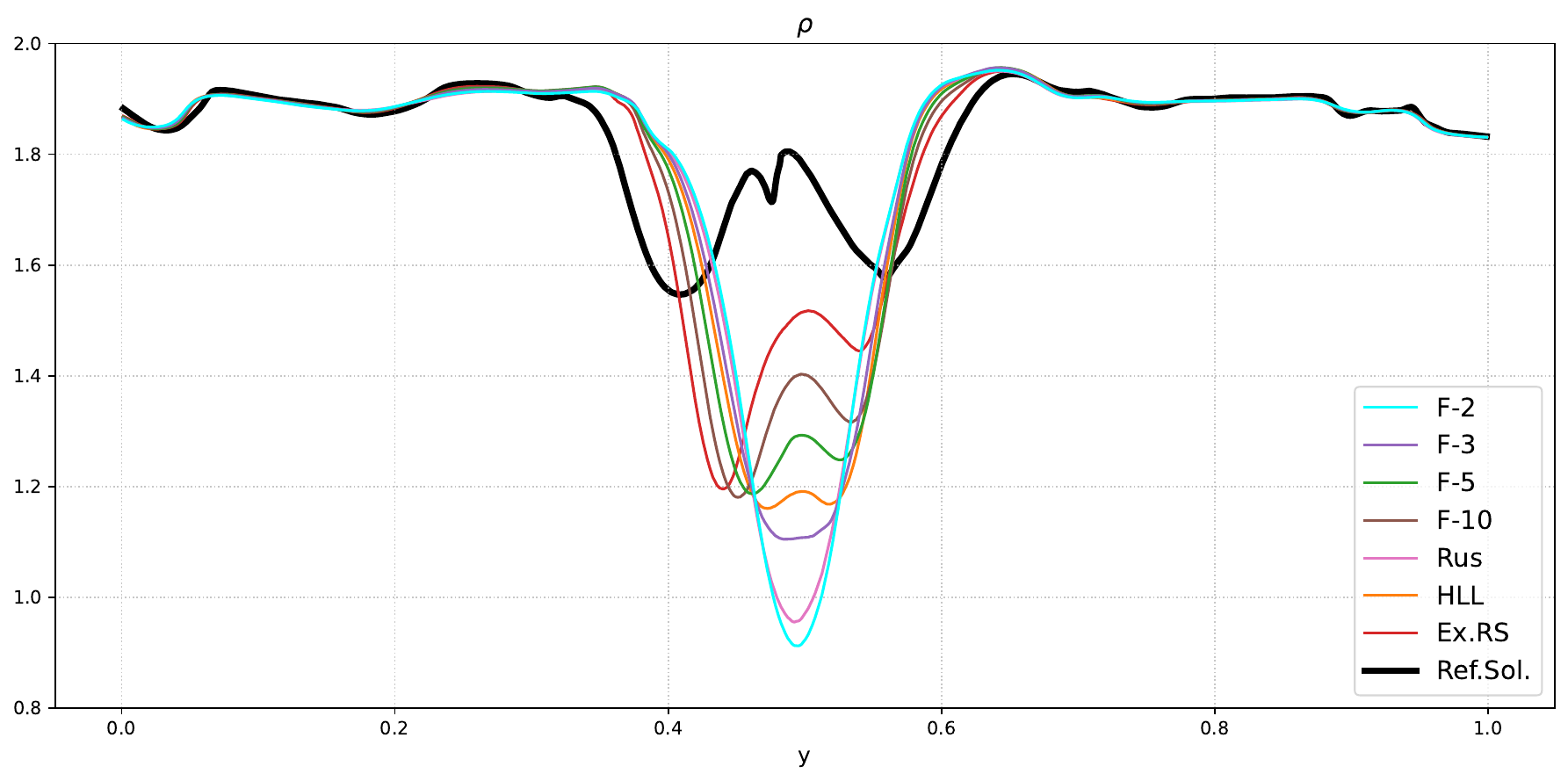}
		\caption{Order 3}
	\end{subfigure}\\
	\begin{subfigure}[b]{1.0\textwidth}
		\centering
		\includegraphics[width=0.7\textwidth]{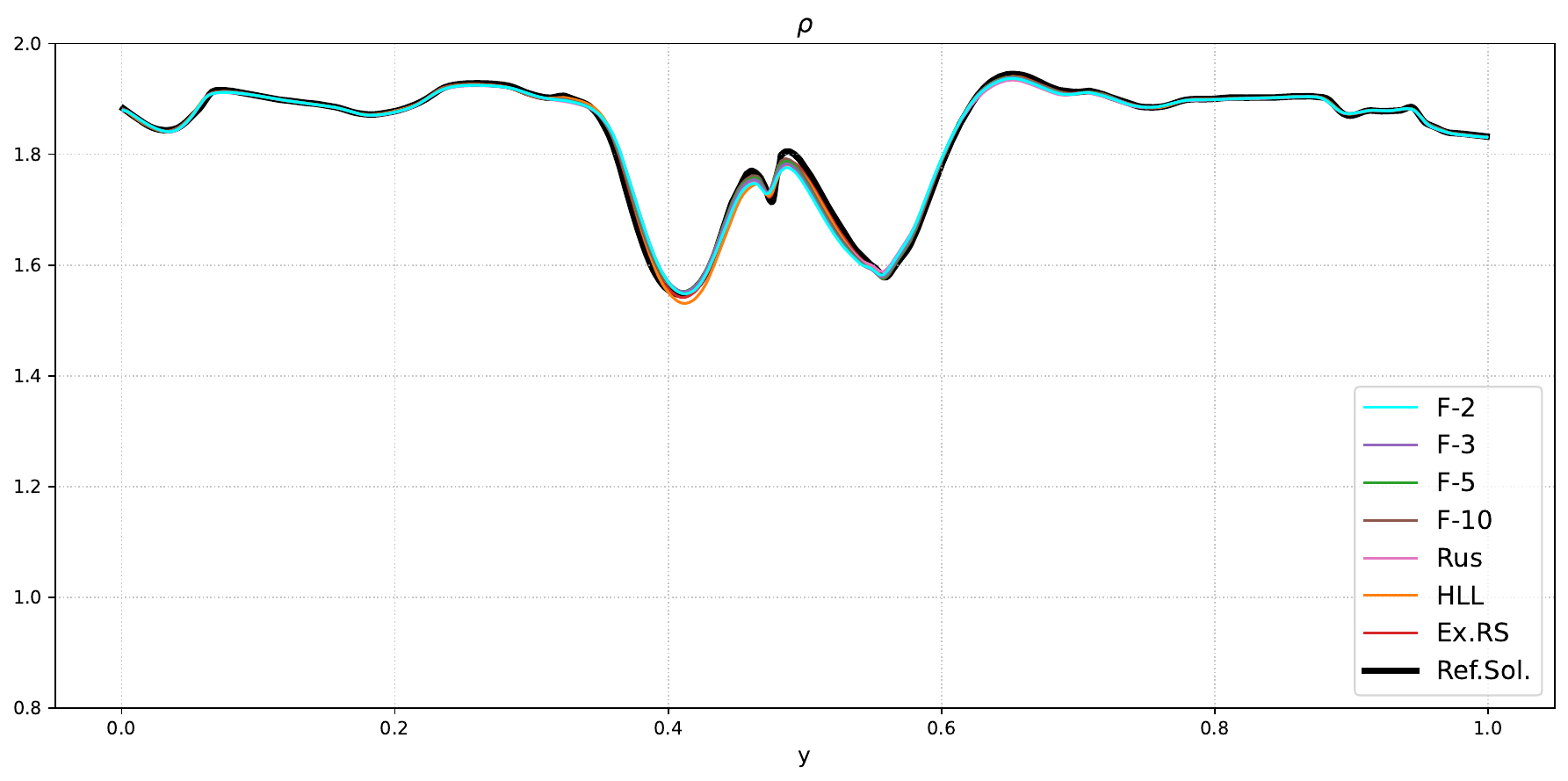}
		\caption{Order 5}
	\end{subfigure}
	\\
	\begin{subfigure}[b]{1.0\textwidth}
		\centering
		\includegraphics[width=0.7\textwidth]{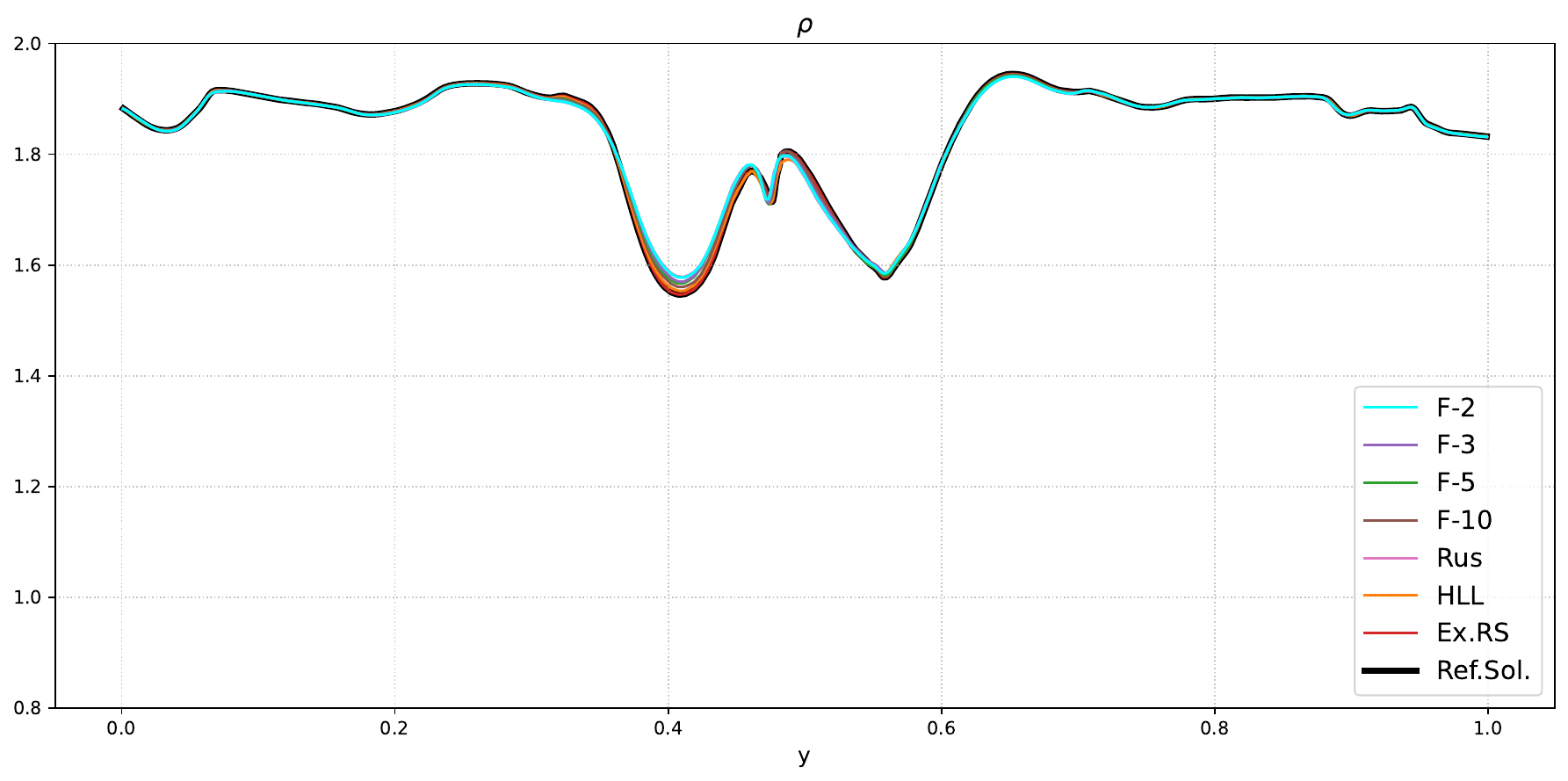}
		\caption{Order 7}
	\end{subfigure}
	\caption{Shock--vortex interaction: Density profile along $x=1$ obtained over a mesh with $800\times 401$ elements with $\sigma_{CFL}:=0.9$}
	\label{fig:density_slices_x1}
\end{figure}

\begin{figure}[htbp]
	\centering
	\begin{subfigure}[b]{1.0\textwidth}
		\centering
		\includegraphics[width=0.7\textwidth]{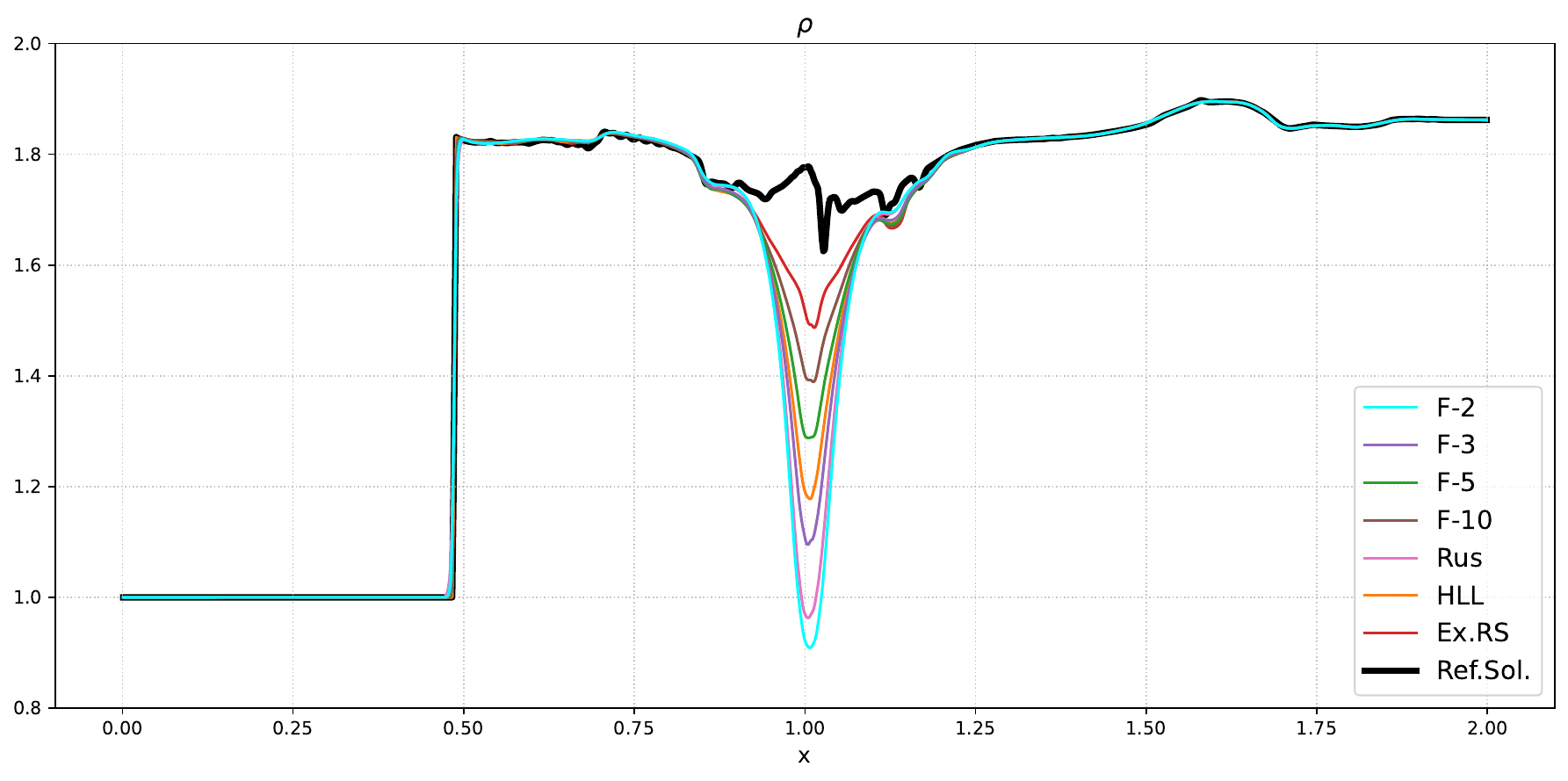}
		\caption{Order 3}
	\end{subfigure}\\
	\begin{subfigure}[b]{1.0\textwidth}
		\centering
		\includegraphics[width=0.7\textwidth]{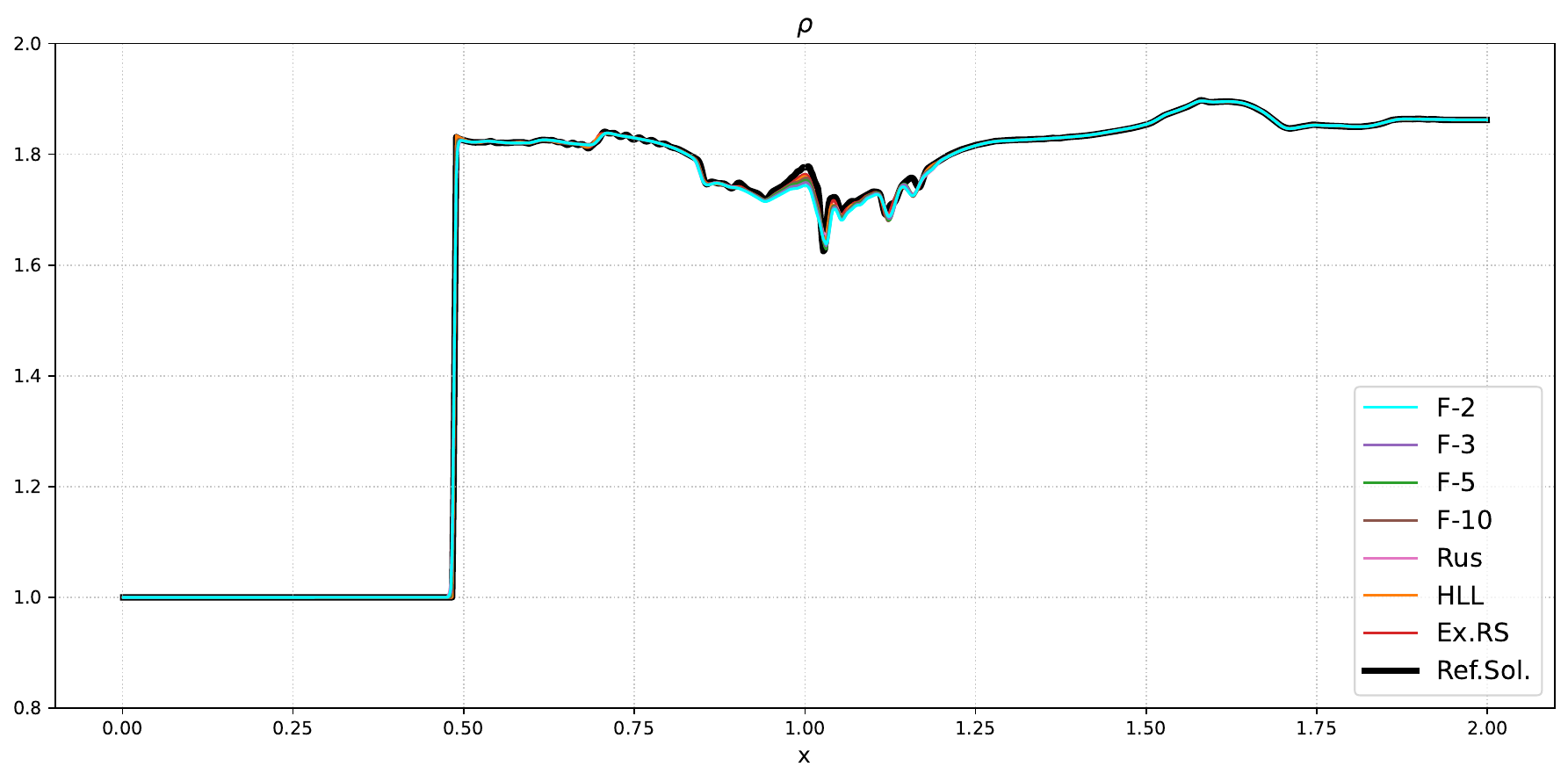}
		\caption{Order 5}
	\end{subfigure}
	\\
	\begin{subfigure}[b]{1.0\textwidth}
		\centering
		\includegraphics[width=0.7\textwidth]{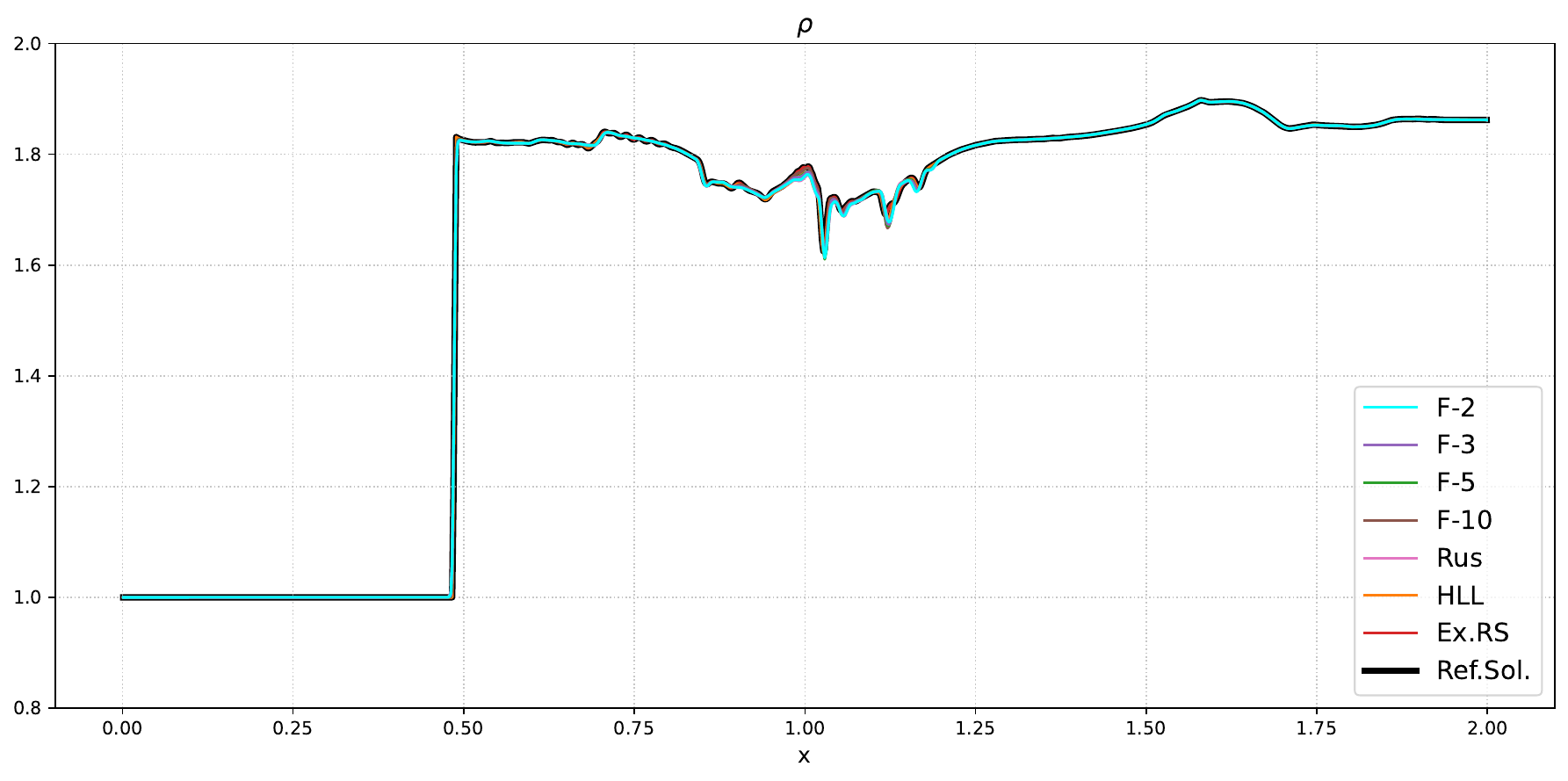}
		\caption{Order 7}
	\end{subfigure}
	\caption{Shock--vortex interaction: Density profile along $y=0.5$ obtained over a mesh with $800\times 401$ elements with $\sigma_{CFL}:=0.9$}
	\label{fig:density_slices_y0.5}
\end{figure}

\begin{figure}[htbp]
	\centering
	\begin{subfigure}[b]{1.0\textwidth}
		\centering
		\includegraphics[width=0.7\textwidth]{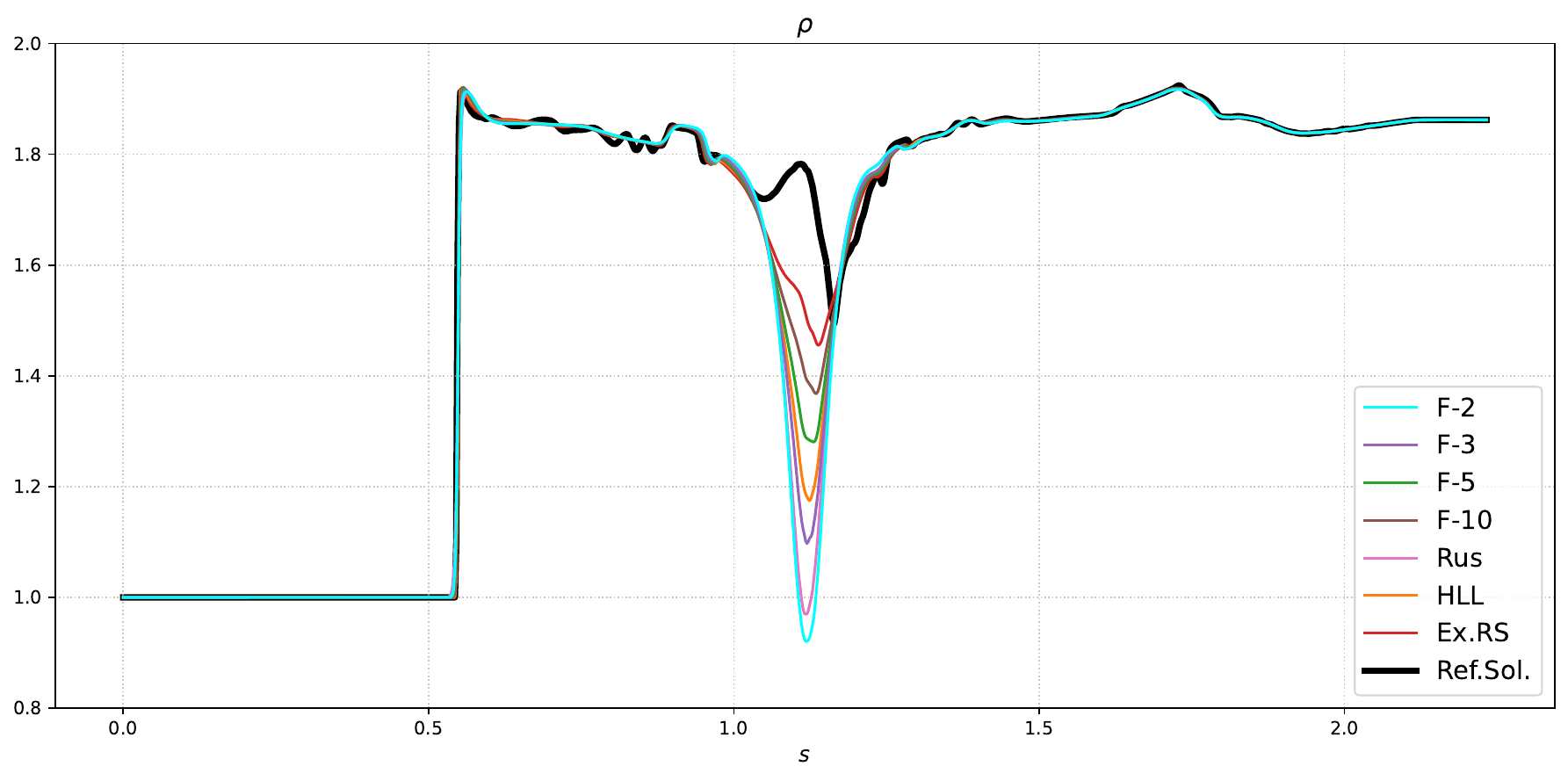}
		\caption{Order 3}
	\end{subfigure}\\
	\begin{subfigure}[b]{1.0\textwidth}
		\centering
		\includegraphics[width=0.7\textwidth]{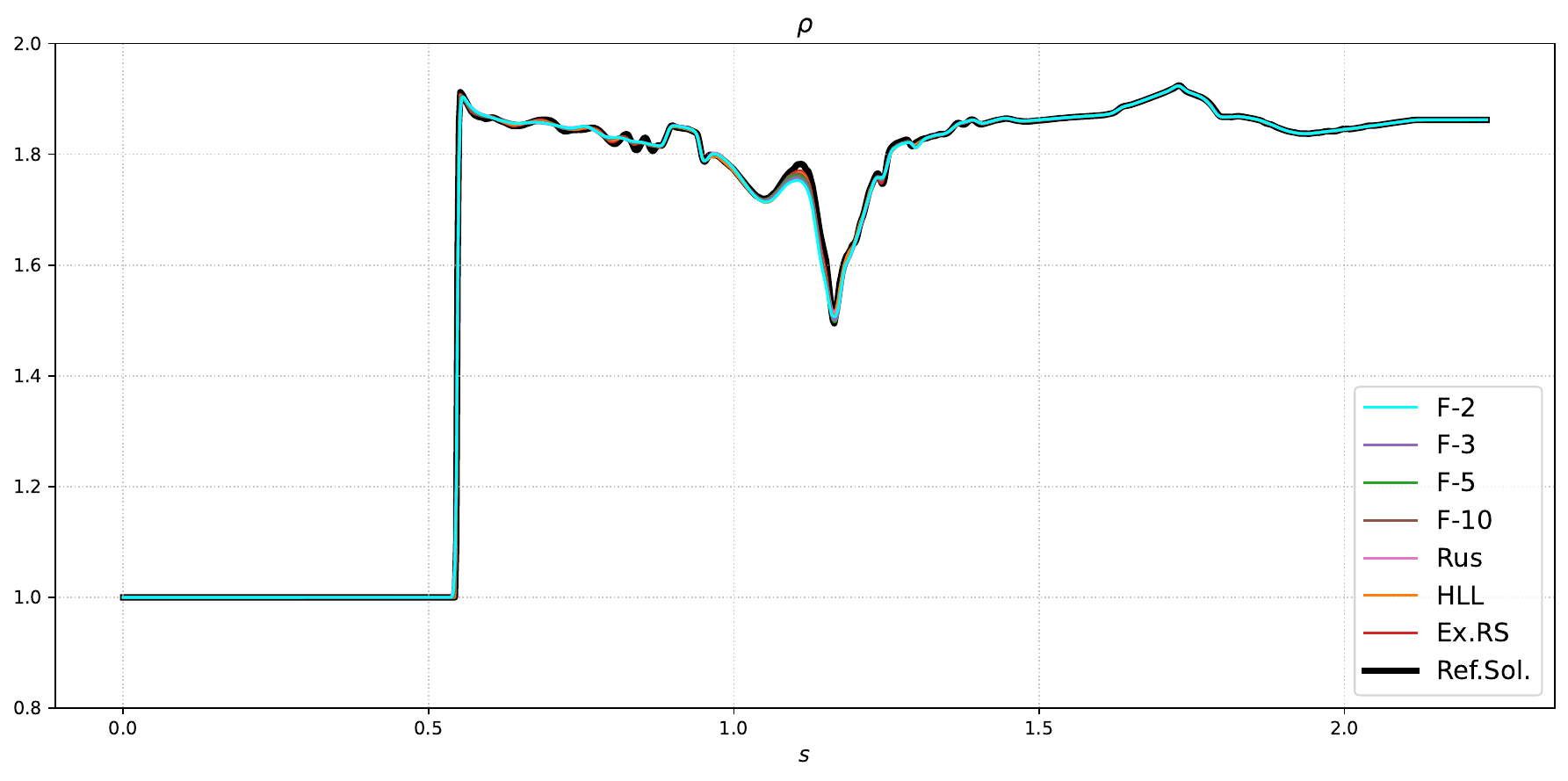}
		\caption{Order 5}
	\end{subfigure}
	\\
	\begin{subfigure}[b]{1.0\textwidth}
		\centering
		\includegraphics[width=0.7\textwidth]{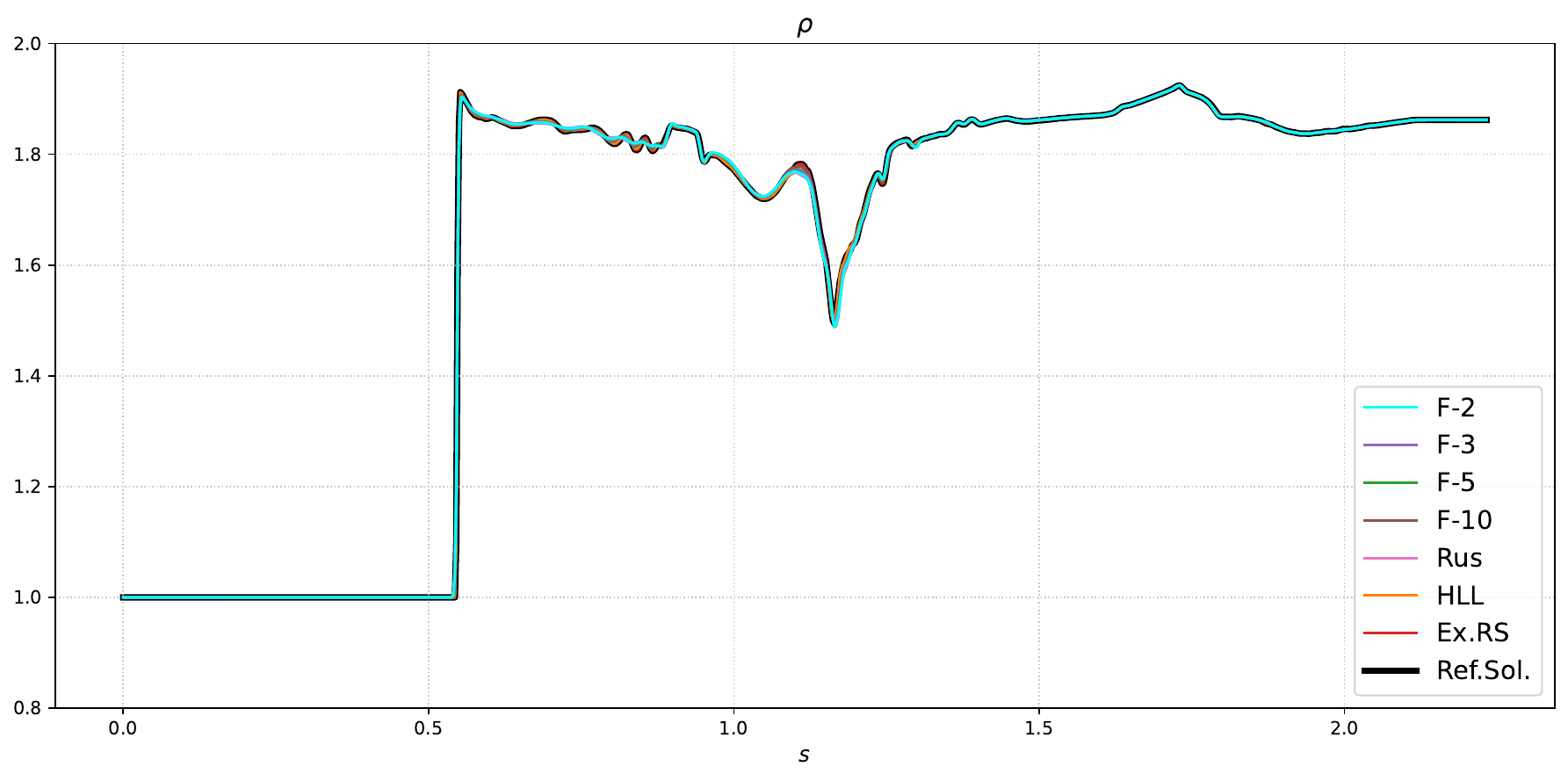}
		\caption{Order 7}
	\end{subfigure}
	\caption{Shock--vortex interaction: Density profile along the domain diagonal joining lower left and upper right corners obtained over a mesh with $800\times 401$ elements with $\sigma_{CFL}:=0.9$}
	\label{fig:density_slices_diag1}
\end{figure}

\begin{figure}[htbp]
	\centering
	\includegraphics[width=0.9\linewidth]{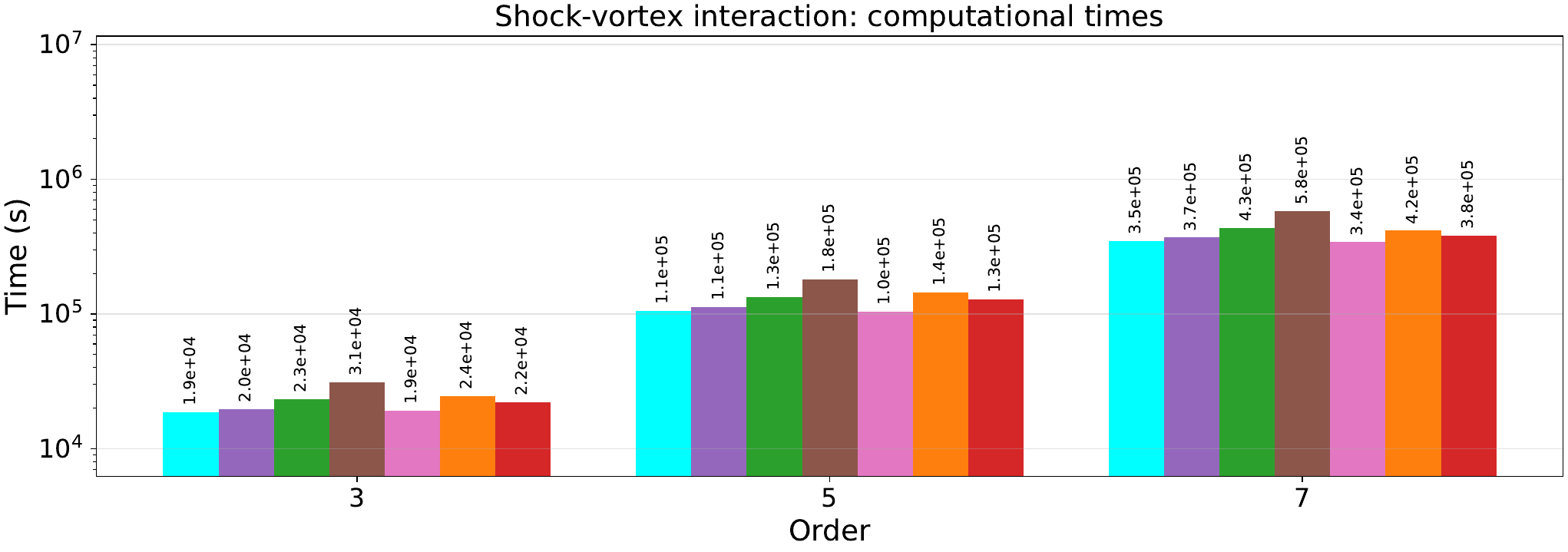}\\
	\includegraphics[width=0.9\linewidth]{figures_new/Euler_2D_flux_comparison_bar_plots_common_legend.pdf}	
	\caption{\RIIcolor{Shock--vortex interaction: Computational times in seconds for the simulation over meshes of $800\times 401$ elements with $\sigma_{CFL}:=0.9$.}}
	\label{fig:Euler_2d_barplots_shock_vortex}
\end{figure}

\RIIcolor{
Furthermore, in order to provide a better comparison, we report, in Figures~\ref{fig:density_slices_x1},~\ref{fig:density_slices_y0.5} and~\ref{fig:density_slices_diag1}, the density profiles obtained through all numerical fluxes and all orders along three slices of the computational domain, namely, the vertical line $x=1$, the horizontal line $y=0.5$ and the domain diagonal connecting the lower left and the upper right corners.
Due to the lack of an exact solution, we adopt the results obtained through exact RS with order 7 as reference.
The computational times are reported in Table~\ref{tab:Euler_2d_efficiency_shock_vortex} and graphically represented in Figure~\ref{fig:Euler_2d_barplots_shock_vortex}.

From the results, the influence of the numerical flux is clearly visible at order 3, where the exact RS provides the sharpest resolution of the structures generated by the shock--vortex interaction. 

The differences become much less pronounced when increasing the order of accuracy. For orders 5 and 7, all numerical fluxes provide a qualitatively similar description of the interaction, with only minor discrepancies in small-scale structures visible in the Schlieren plots and in the diagonal density profiles reported in Figure~\ref{fig:density_slices_diag1}. These differences are mostly related to the resolution of secondary waves, with FORCE-2 being slightly more diffusive. 

Concerning efficiency, FORCE-2 and FORCE-3 have computational times comparable to Rusanov, while FORCE-5 is comparable to HLL and exact RS. On the other hand, FORCE-10 appears excessively expensive for this test, without providing a sufficiently significant improvement in resolution to justify its cost. Overall, these results suggest that moderate values of $\alpha$, in particular FORCE-2 and FORCE-3, represent the most convenient choices. Although the computational cost of this two-dimensional test makes a systematic equal-cost comparison less practical than in the one-dimensional setting, the results clearly indicate that, in high order simulations, FORCE-$\alpha$ numerical fluxes can be employed without an excessive loss of accuracy with respect to classical upwind fluxes.
}

Let us finally remark that the results obtained with Rusanov and FORCE-$\alpha$ numerical fluxes feature a spurious acoustic wave propagating rightwards. This is a start-up error due to the fact that such numerical fluxes do not resolve exactly the stationary shock, and it can be fixed with standard techniques, e.g., running the test with the numerical stationary shock obtained after the spurious perturbation leaves the domain.

\begin{table}[htbp]
	\centering
	\begin{tabular}{|c||c|c|c|}
		\hline
		Numerical flux & Order 3 & Order 5 & Order 7\\
		\hline\hline
		FORCE--2  & 18662 & 105622 & 347501\\
		\hline
		FORCE--3  & 19751 & 113147 & 369561\\
		\hline
		FORCE--5  & 23252 & 133292 & 434373\\
		\hline
		FORCE--10 & 31005 & 181042 & 579572\\
		\hline
		Rusanov   & 19195 & 103892 & 341586\\
		\hline
		HLL       & 24458 & 143848 & 415805\\
		\hline
		exact RS  & 22168 & 128926 & 382664\\
		\hline
	\end{tabular}
	\caption{
		\RIIcolor{Shock--vortex interaction: Computational times in seconds for the simulation over meshes of $800\times 401$ elements with $\sigma_{CFL}:=0.9$.}
	}
	\label{tab:Euler_2d_efficiency_shock_vortex}
\end{table}

\subsubsection{Long--time evolution of smooth isentropic vortex}\label{sec:unsteady_vortex_longer_time}
The same vortex of the accuracy test, in Section~\ref{sec:unsteady_vortex}, is considered here on the computational domain $\Omega:=[-5,5]\times[-5,5]$, again with periodic boundary conditions.
This set-up will be used to test the performance of FORCE--$\alpha$ numerical fluxes in the context of problems involving long--time evolution.
Longer final times are associated with higher levels of diffusion, determining a smearing of the solution profile.
In such cases, increasing the order of accuracy and adopting more accurate numerical fluxes plays a crucial role~\cite{micalizzitoro2024}.

\begin{figure}[htbp]
	\centering
	\begin{subfigure}[b]{1.0\textwidth}
		\centering
		\includegraphics[width=0.32\textwidth]{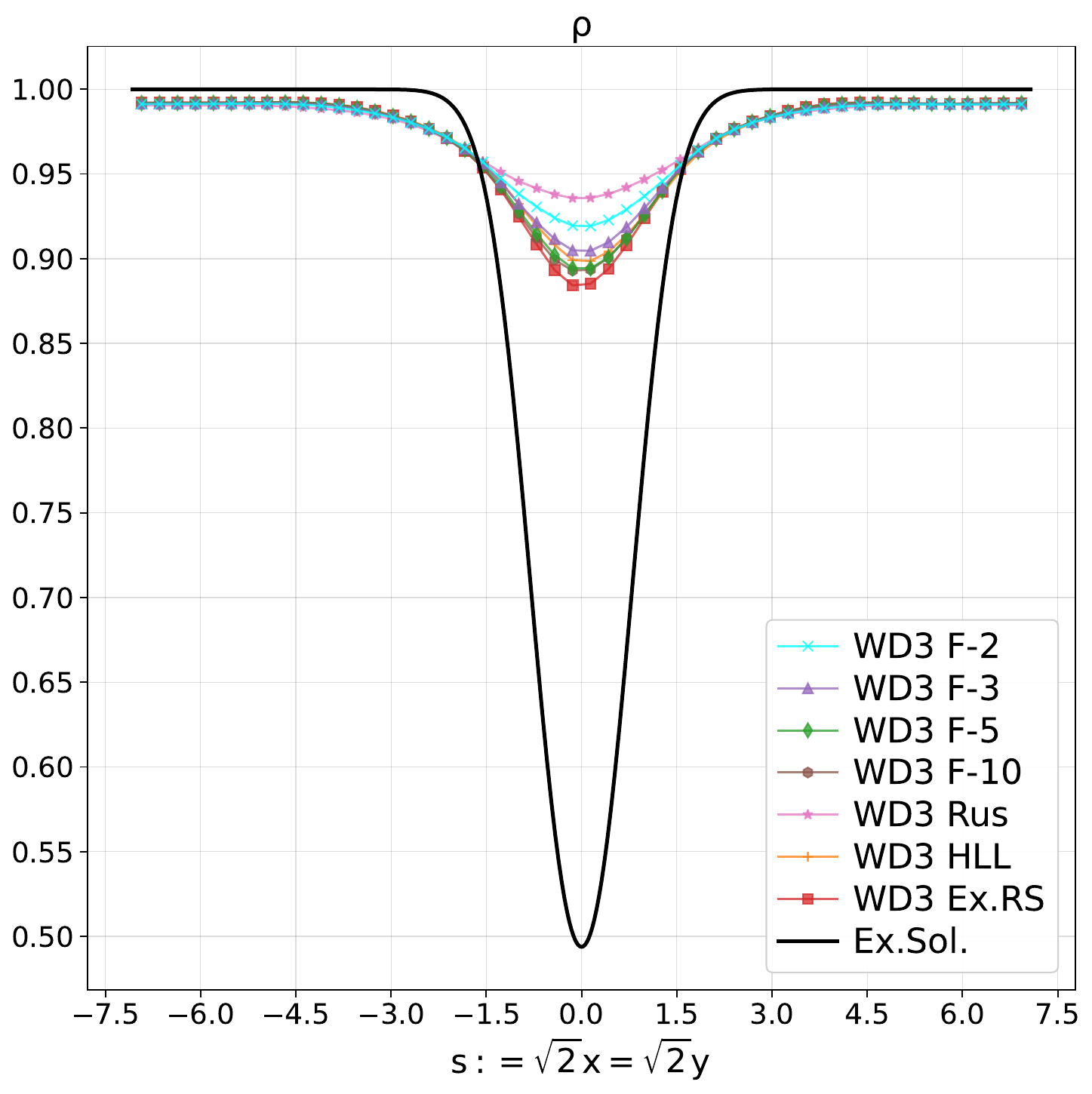}
		\includegraphics[width=0.32\textwidth]{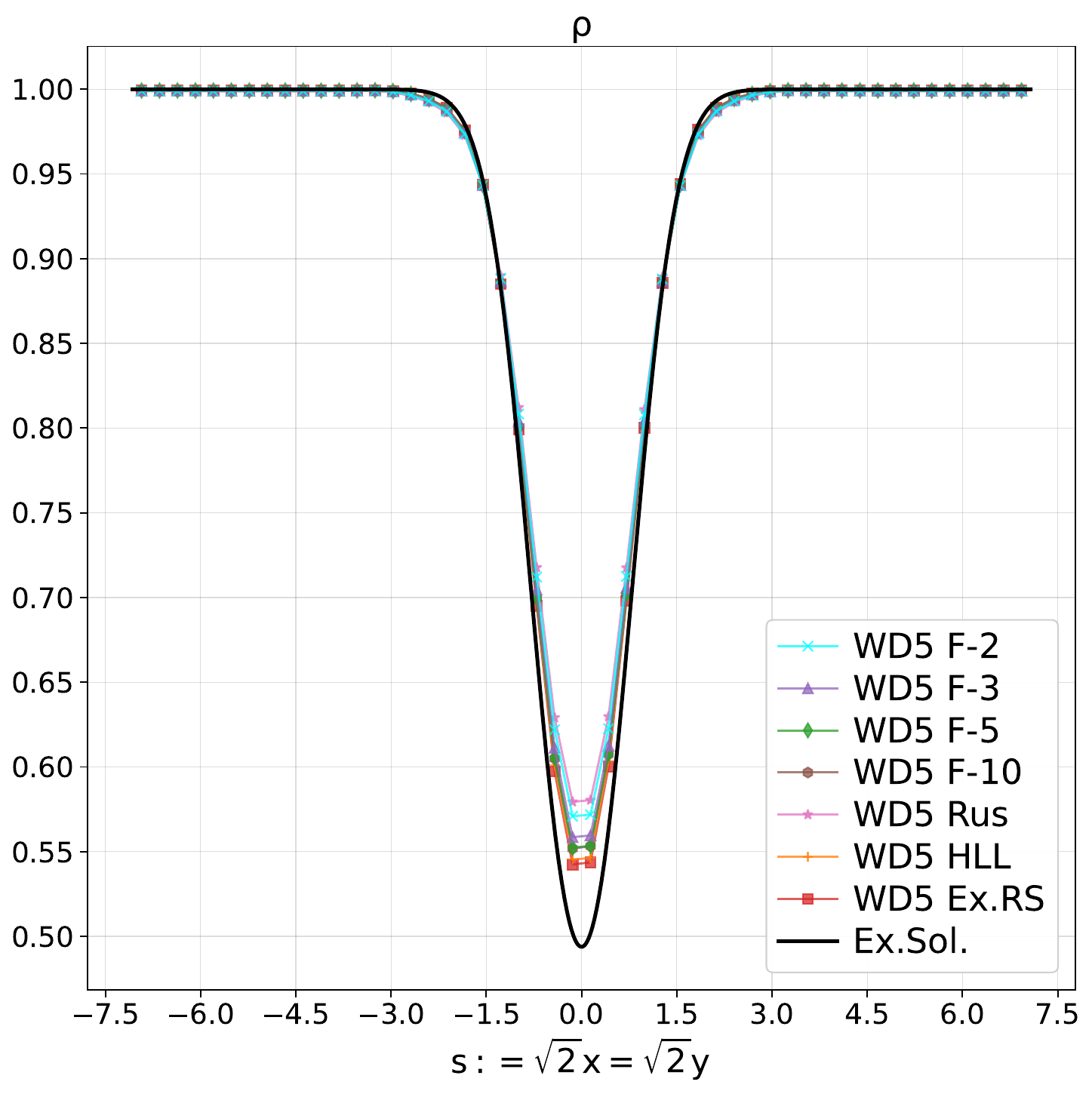}
		\includegraphics[width=0.32\textwidth]{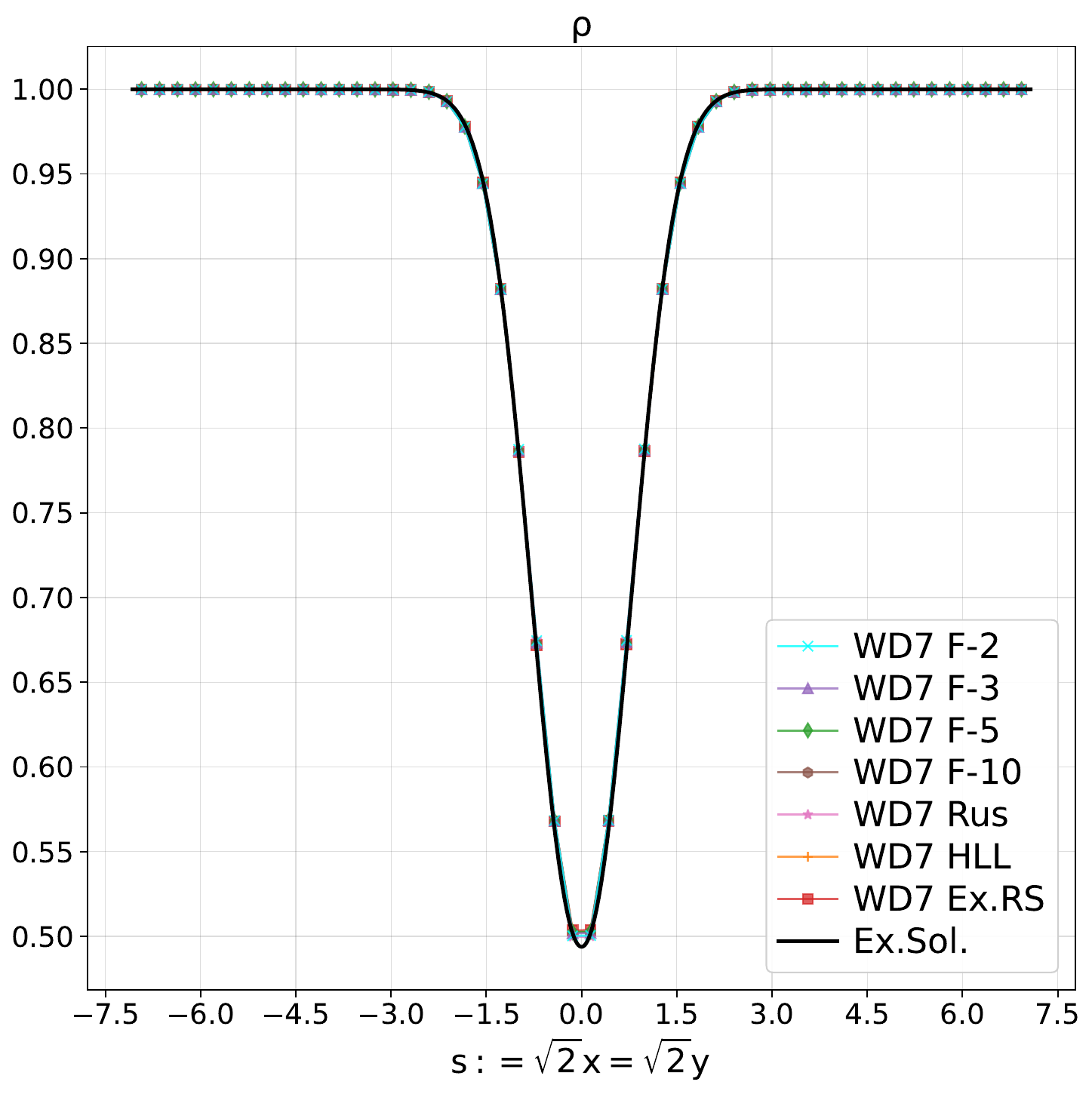}
		\caption{$T_f:=100$}
	\end{subfigure}\\
	\begin{subfigure}[b]{1.0\textwidth}
		\centering
		\includegraphics[width=0.32\textwidth]{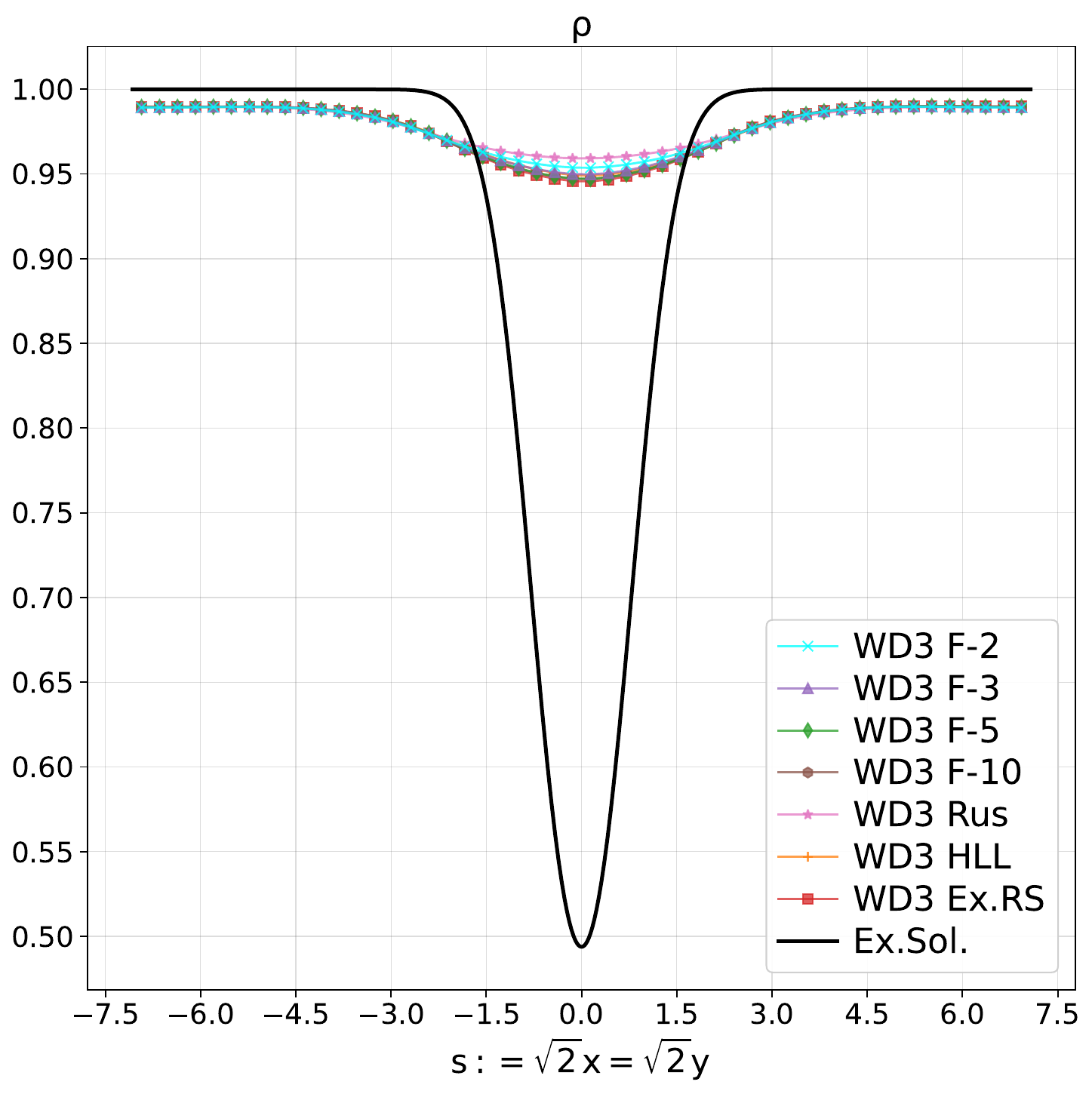}
		\includegraphics[width=0.32\textwidth]{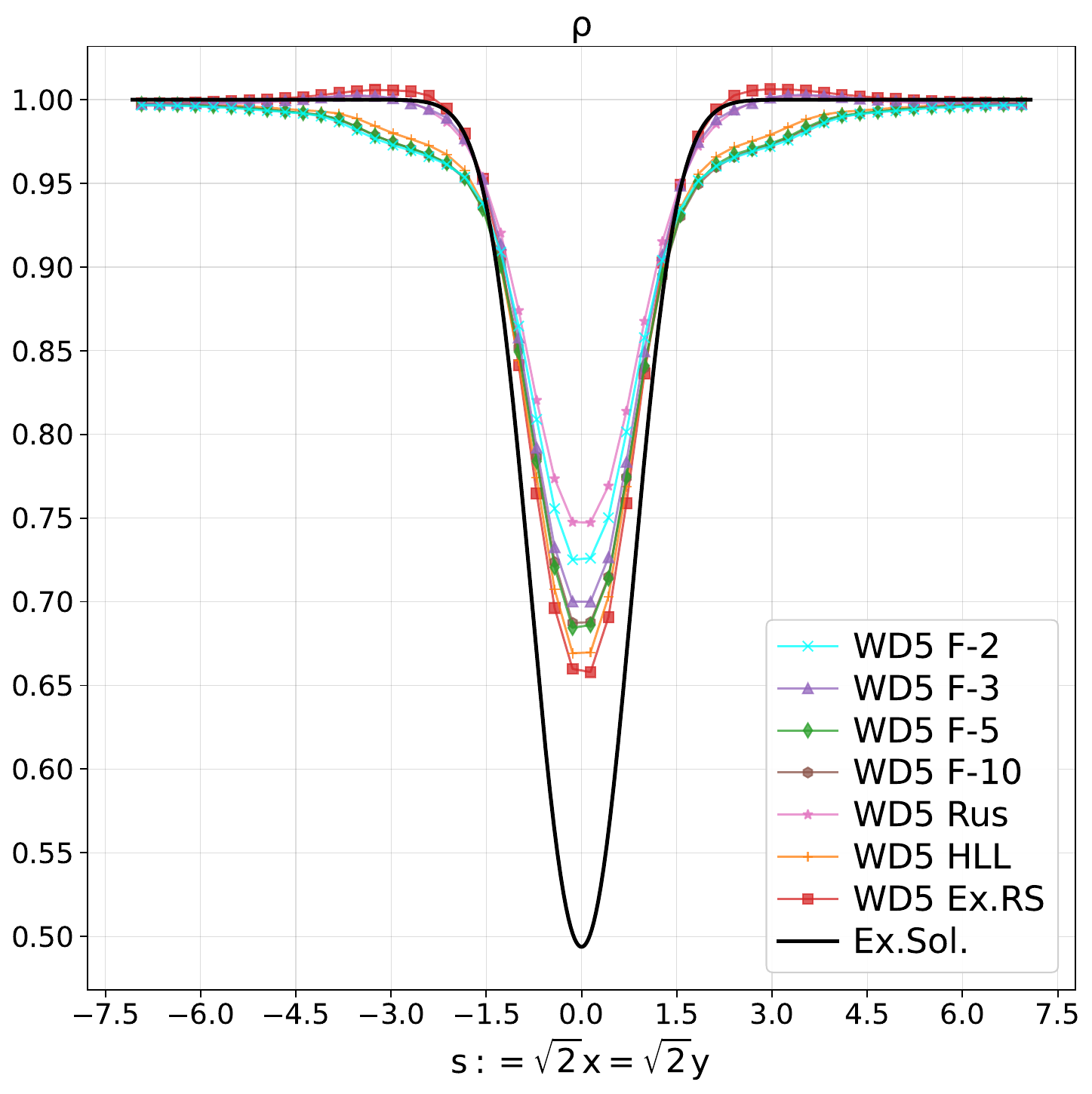}
		\includegraphics[width=0.32\textwidth]{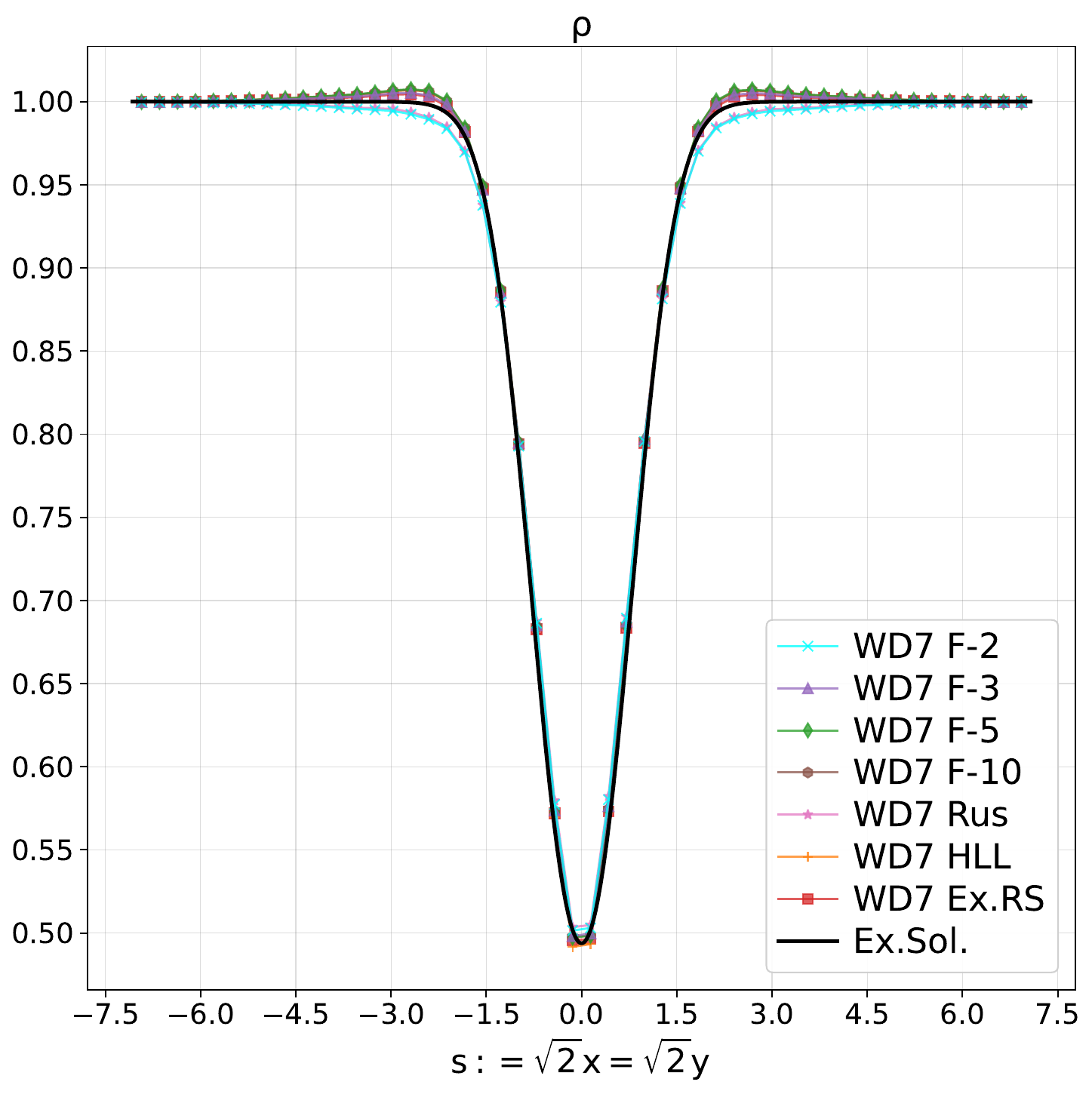}
		\caption{$T_f:=400$}
	\end{subfigure}
	\\
	\begin{subfigure}[b]{1.0\textwidth}
		\centering
		\includegraphics[width=0.32\textwidth]{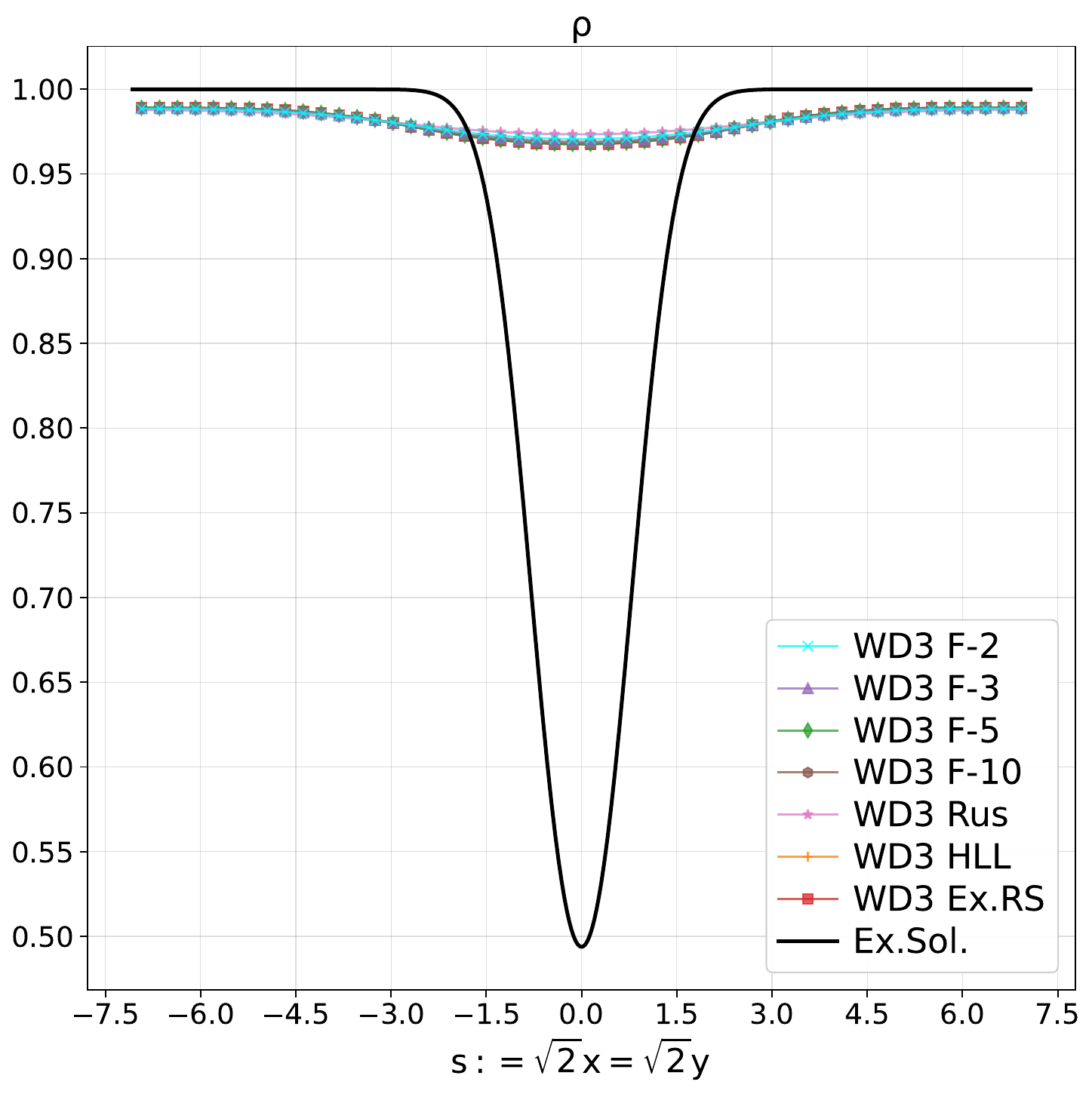}
		\includegraphics[width=0.32\textwidth]{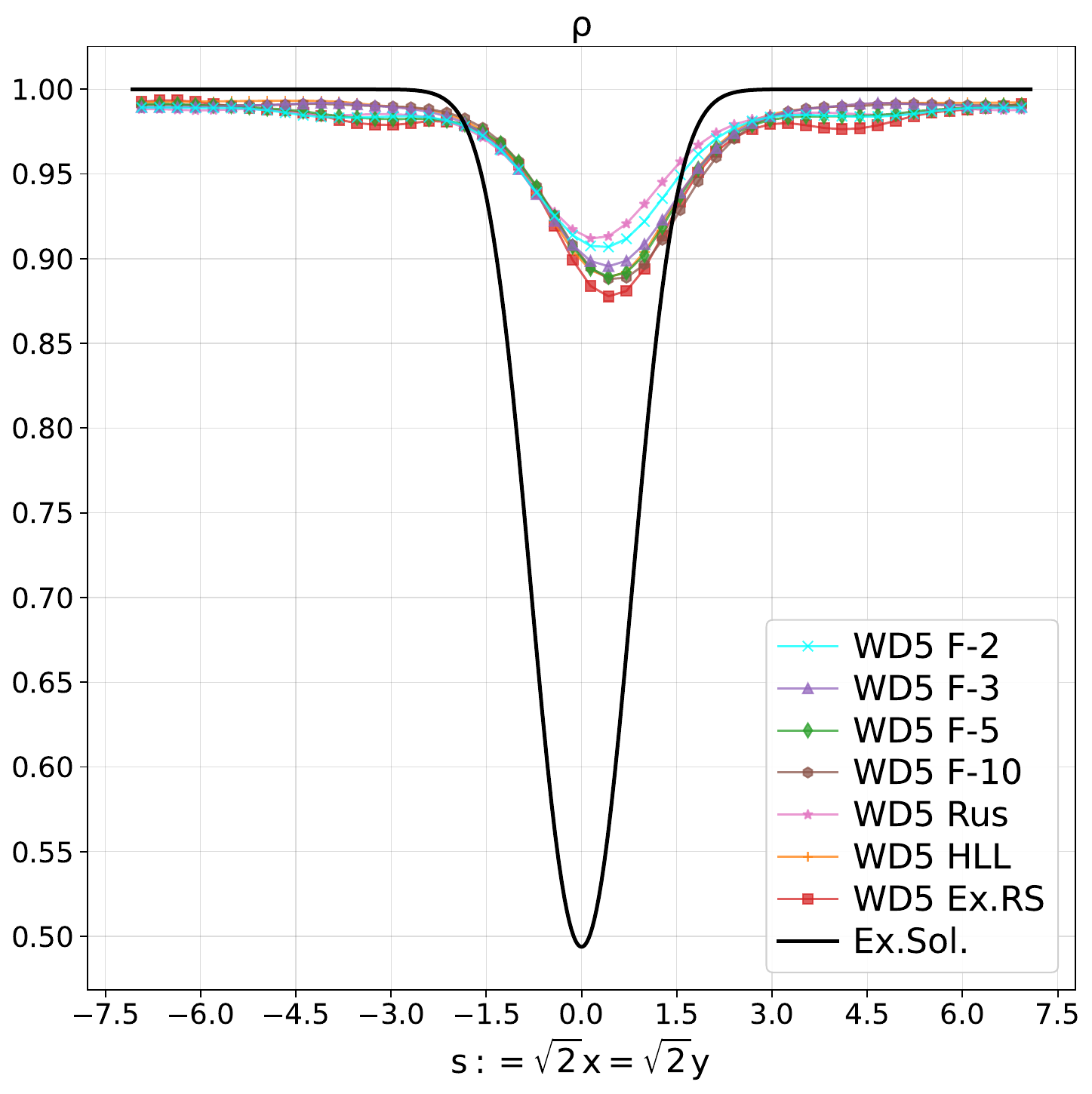}
		\includegraphics[width=0.32\textwidth]{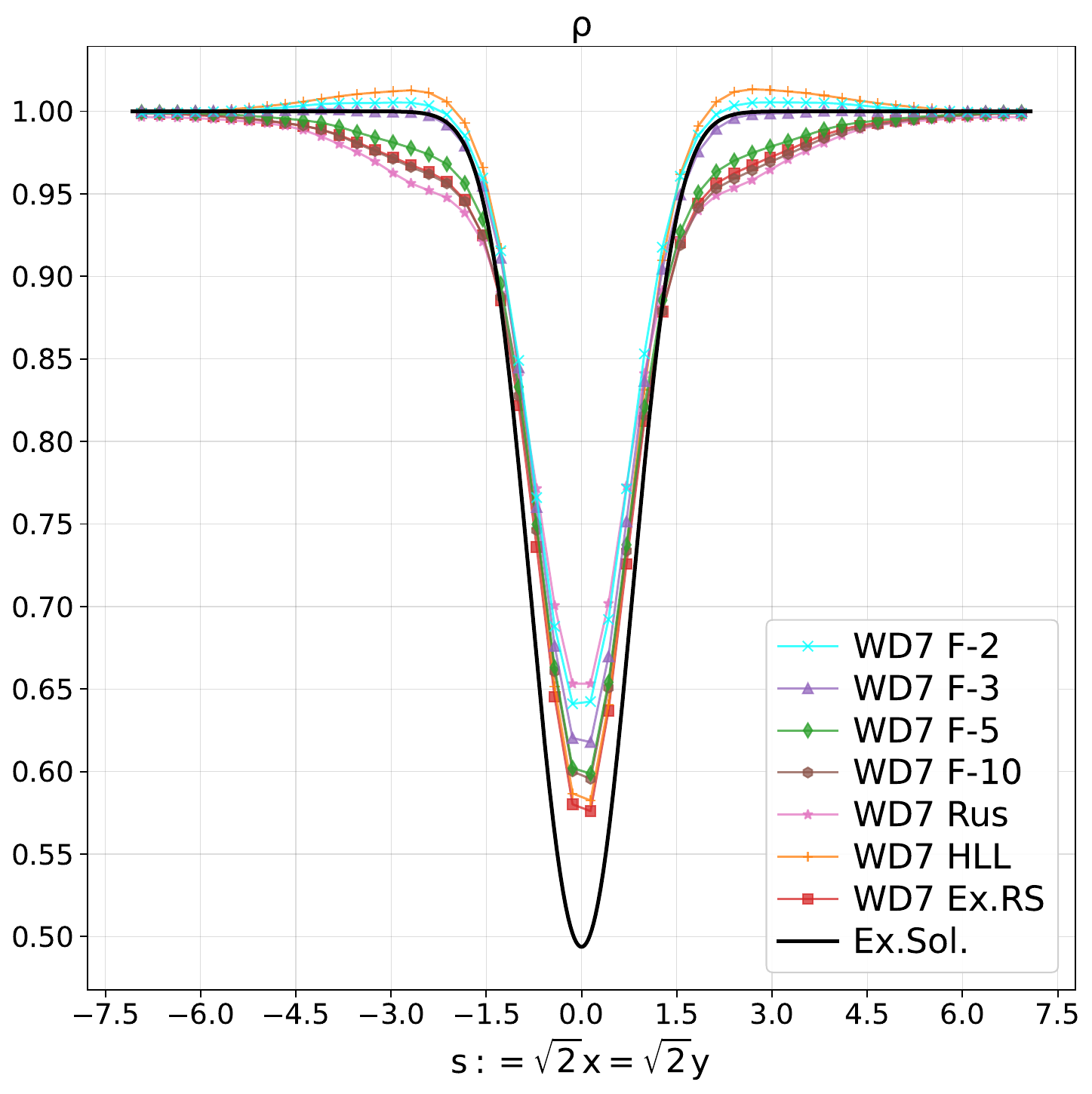}
		\caption{$T_f:=1600$}
	\end{subfigure}
	\caption{\RIcolor{Long--time evolution of smooth isentropic vortex: Density profile along the domain diagonal $y=x$ obtained over a mesh with $50\times 50$ elements with $\sigma_{CFL}:=0.9$ at different times. From left to right, orders 3, 5 and 7}}
	\label{fig:Euler_2d_unsteady_vortex_longer_time}
\end{figure}

\RIIcolor{Since the goal of this test is to assess the long-time dissipative behavior of the different numerical fluxes, rather than their computational efficiency, all comparisons are performed at fixed mesh resolution.}
We report the results obtained for the density along the domain diagonal $y=x$, over a mesh of $50\times 50$ elements with $\sigma_{CFL}:=0.9$, for all considered orders in Figure~\ref{fig:Euler_2d_unsteady_vortex_longer_time}. 
Three final times, $T_f:=100$, $400$ and $1600$, are considered.
For order 3, all results are basically unacceptable for all considered final times and all numerical fluxes, with smearing with respect to the exact amplitude of the vortex being at least around 80\% even for $T_f:=100$.
For order 5, results show little differences for $T_f:=100$, with 
exact RS being the best numerical flux followed by HLL, FORCE-5, FORCE-10, FORCE-3, FORCE-2 and Rusanov. The difference in the smearing with respect to the exact amplitude of the vortex between exact RS and Rusanov is only 8\% and the profile is well captured by all numerical fluxes.
There is indeed more variety for $T_f:=400$. In this case, the hierarchy of the results quality does not change but the difference between the best and the worst numerical fluxes, i.e., exact RS and Rusanov, is 18\%. Nonetheless, let us notice that the difference between FORCE-3 and exact RS is 8\%. Results are similar for FORCE-5 and FORCE-10, with a smearing which is higher by only 4\% with respect to HLL, and by 6\% with respect to the exact RS.
Let us also remark that the smearing obtained for such a final time with order 5 is considerable, being around 30\% with respect to the exact amplitude of the vortex for exact RS.
For order 5, results are unacceptable for all numerical fluxes with $T_f:=1600$.
For order 7, results are good for all final times and all numerical fluxes. In particular, they are optimal for $T_f:=100$ and $400$, with no visible smearing nor differences between the considered numerical fluxes.
Some differences arise for $T_f:=1600$.
Despite being an upwind flux, Rusanov gives the worst results, while, exact RS gives the best ones followed by HLL. Results of intermediate quality are obtained through all FORCE--$\alpha$ numerical fluxes with increasing quality as $\alpha$ increases.
Let us notice that for order 7 the smearing with respect to the exact solution obtained with FORCE-5 and FORCE-10 is only 4\% bigger than the one of HLL and exact RS.
The smearing obtained with FORCE-3 is 8\% bigger than the one of HLL and exact RS, while, with FORCE-2 it is 12\% bigger.
Clearly, the specific application determines whether such a loss in accuracy is acceptable or not.
Overall, also in this test, FORCE--$\alpha$ numerical fluxes prove to be a valid alternative to upwind fluxes within a very high order setting.

This example shows that long--time evolution problems are challenging and useful to assess the quality of numerical methods, especially in the context of very high order of accuracy, even though problems of such kind are rarely considered in existing literature.

\subsubsection{Explosion problems}\label{sec:explosion}
\RIIcolor{In this section, we test the robustness and the shock-capturing features of FORCE--$\alpha$ numerical fluxes on explosion problems.}

\subsubsection*{Radial explosion}\label{sec:sod_2d}
In this explosion problem, presented in~\cite{ToroBook}, we consider an initial condition, consisting of a circular region with high density and pressure in the middle of the domain $\Omega := [-1,1]\times[-1,1]$, which reads
\begin{equation}
	\begin{pmatrix}
		\rho\\
		u\\
		v\\
		p
	\end{pmatrix}(x,y,0) :=
	\begin{cases}
		(1,0,0,1)^T, & \text{if}~\sqrt{x^2+y^2} < 0.4, \\ 
		(0.125,0,0,0.1)^T, & \text{otherwise}.
	\end{cases}
	\label{eq:Euler_2d_sod_IC}
\end{equation}
We adopt transmissive boundary conditions, and final time $T_f := 0.25$.
The solution features a shock wave and a contact wave propagating towards the domain boundary, along with a circular rarefaction propagating towards the center of the domain.

\begin{figure}[htbp]
	\centering
	\begin{subfigure}[b]{1.0\textwidth}
		\centering
		\includegraphics[width=0.32\textwidth,trim={0 0 8.2cm 0},clip]
		{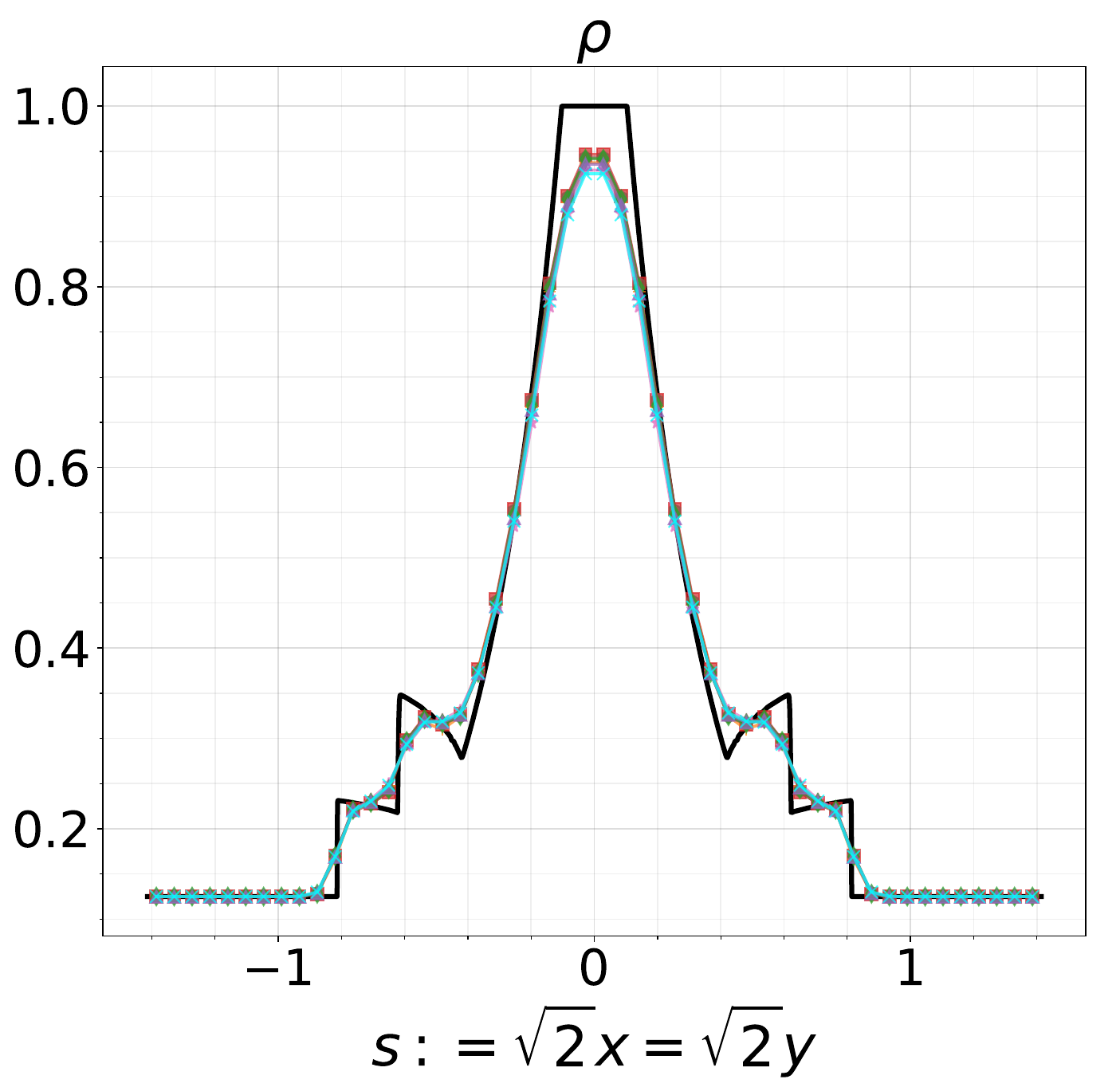}
		\includegraphics[width=0.32\textwidth,trim={0 0 8.2cm 0},clip]
		{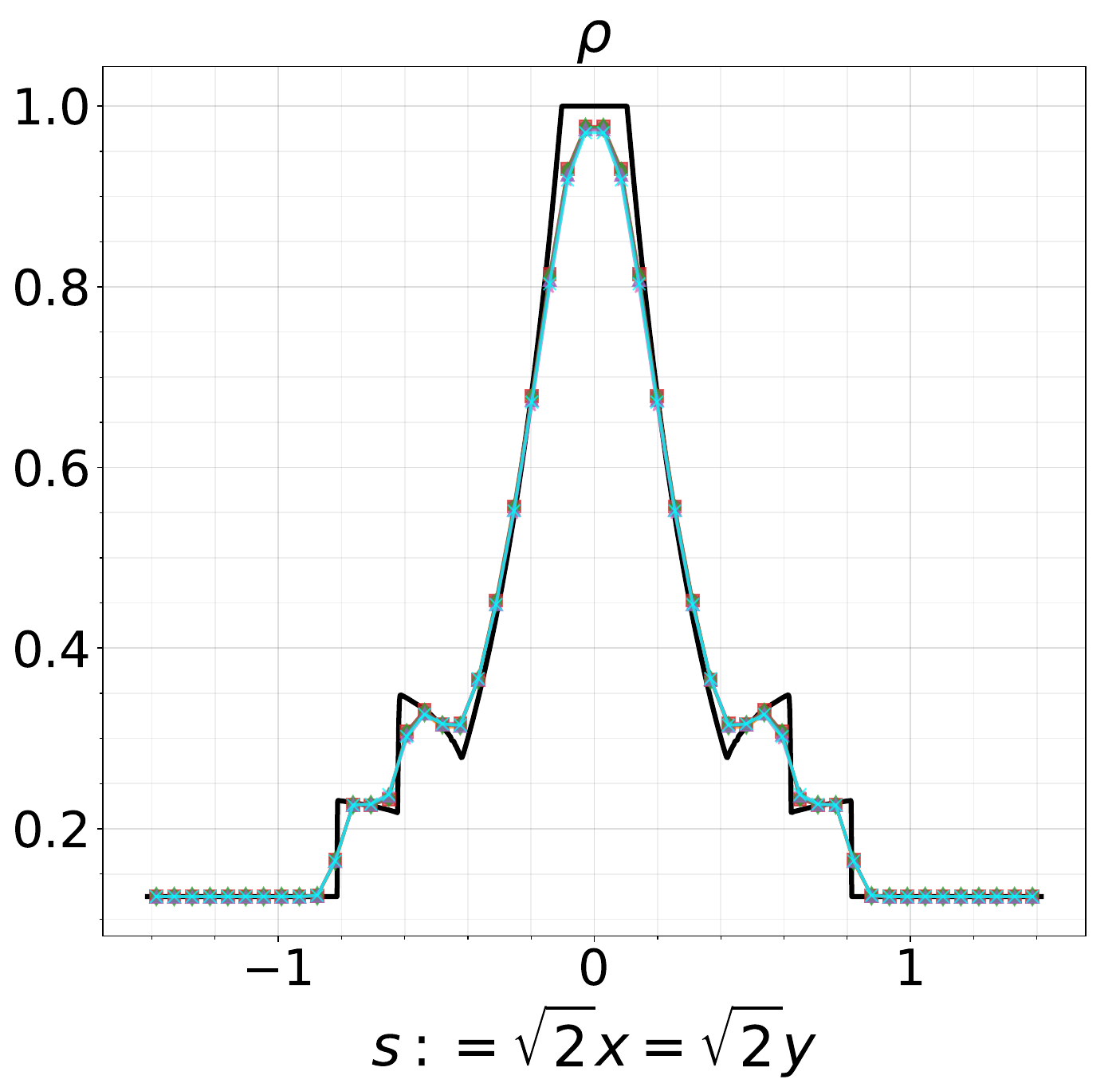}
		\includegraphics[width=0.32\textwidth,trim={0 0 8.2cm 0},clip]
		{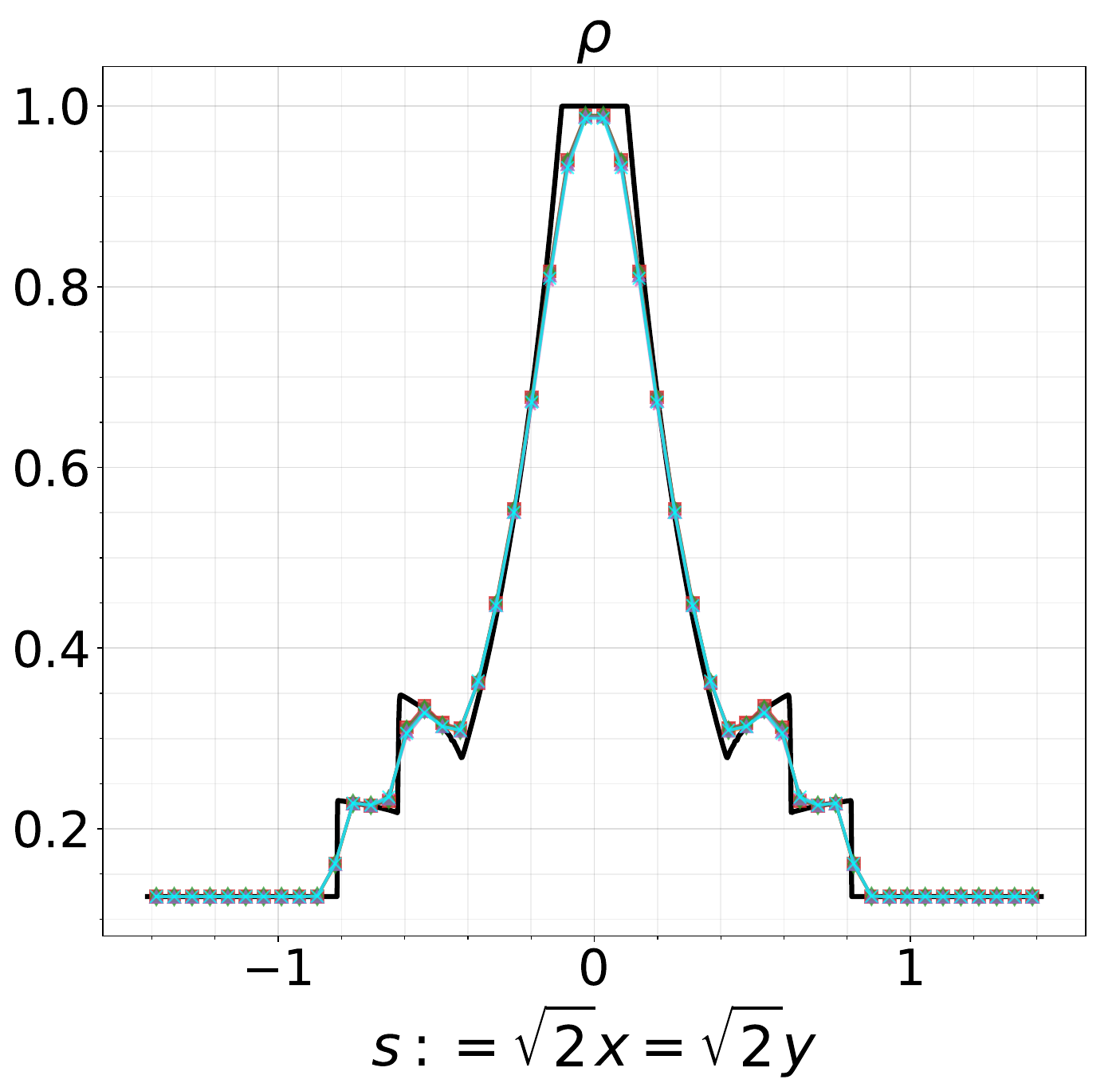}
		\caption{$50\times 50$ elements}
	\end{subfigure}\\
	\begin{subfigure}[b]{1.0\textwidth}
		\centering
		\includegraphics[width=0.32\textwidth,trim={0 0 8.2cm 0},clip]
		{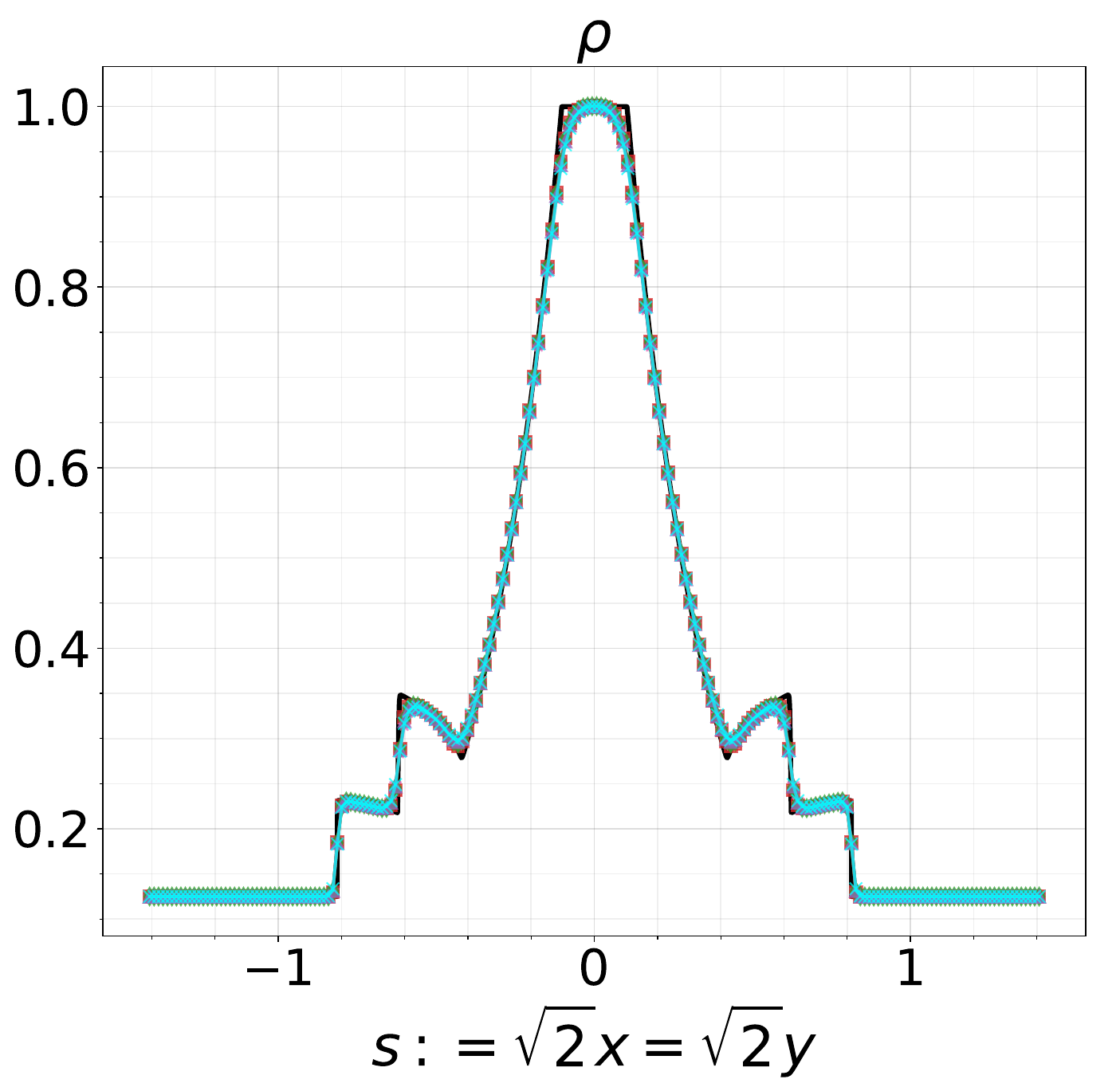}
		\includegraphics[width=0.32\textwidth,trim={0 0 8.2cm 0},clip]
		{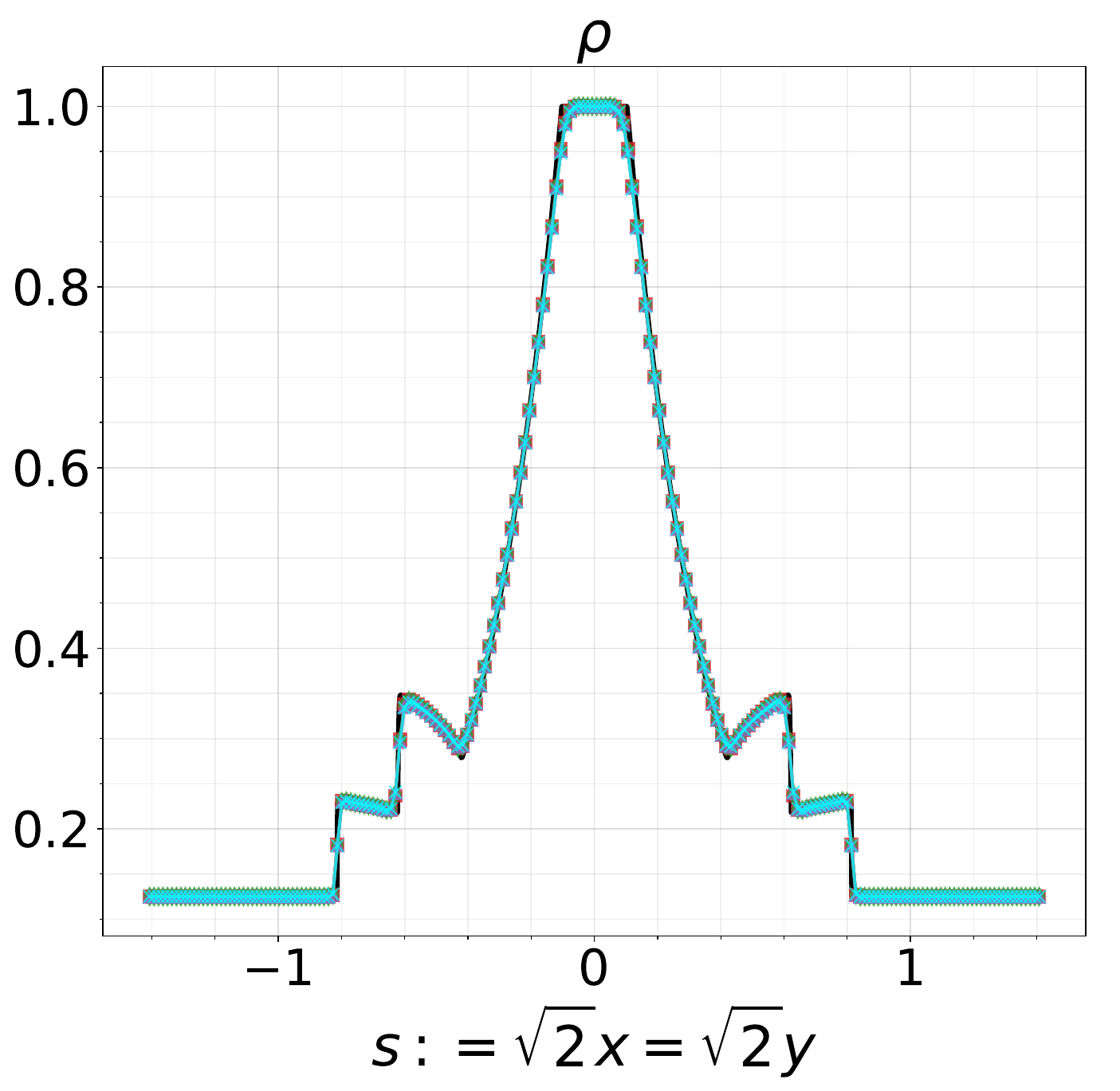}
		\includegraphics[width=0.32\textwidth,trim={0 0 8.2cm 0},clip]
		{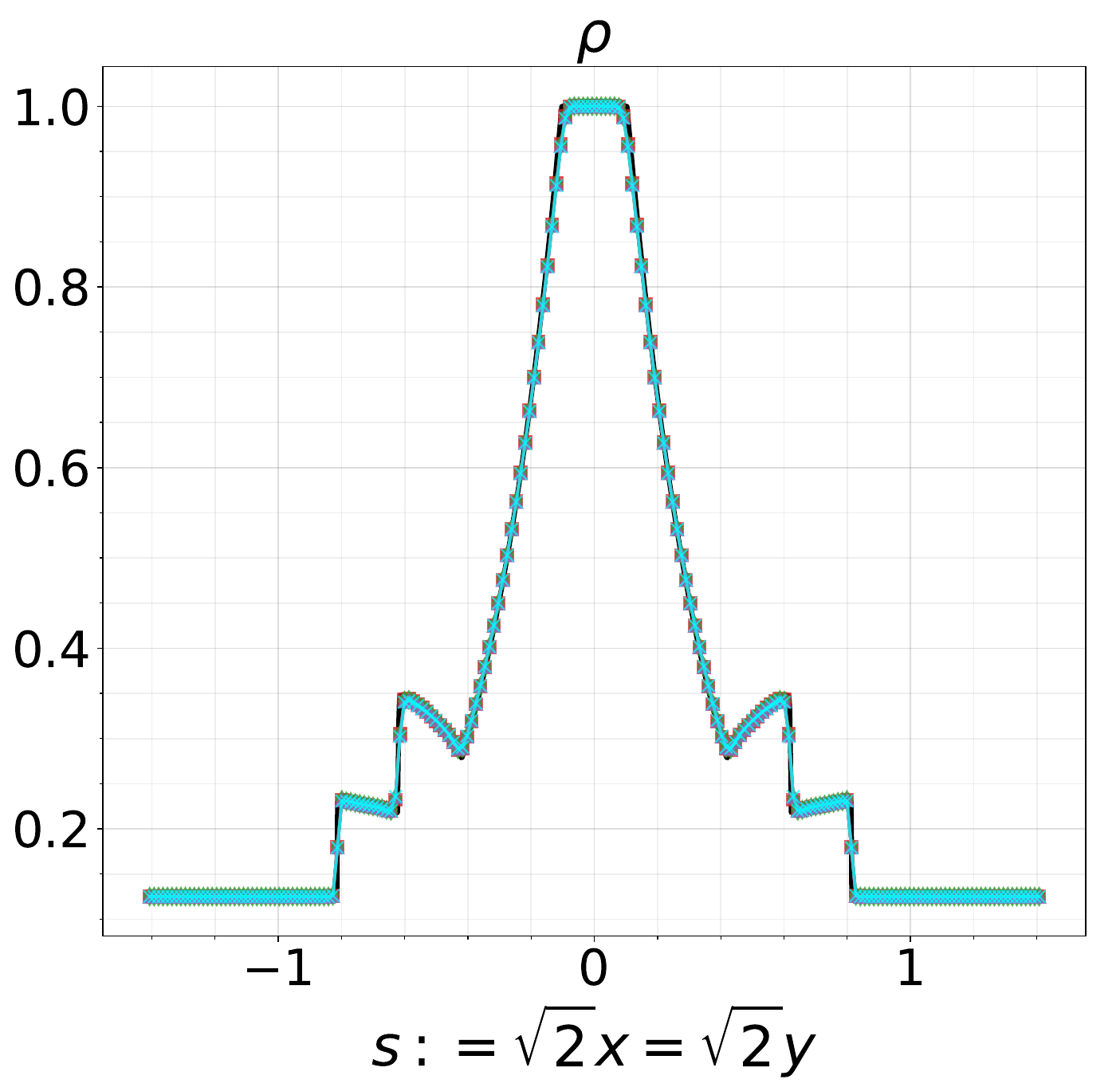}
		\caption{$200\times 200$ elements}
	\end{subfigure}
	
	\vspace{0.5em}
	
	\includegraphics[width=0.85\textwidth]{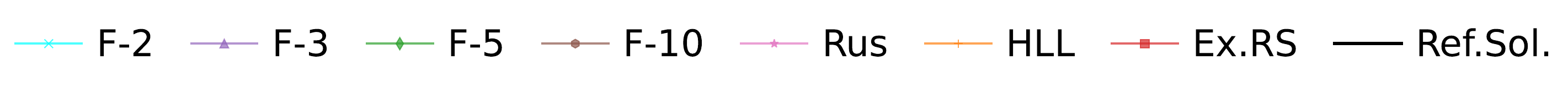}
	
	\caption{\RIIcolor{Radial explosion: Density profile along the domain diagonal $y=x$ obtained with $\sigma_{CFL}:=0.8$ over meshes with $50\times 50$ and $200\times 200$ elements. From left to right, orders 3, 5 and 7.}}
	\label{fig:Euler_2d_explosion_problem}
\end{figure}

\RIIcolor{The results obtained over meshes with $50\times 50$ and $200\times 200$ cells with $\sigma_{CFL}:=0.8$ are reported in Figure~\ref{fig:Euler_2d_explosion_problem}.
In particular, we report the density results along the domain diagonal $y=x$.}
The reference solution has been computed, over a mesh of $3000\times3000$ elements, through a second order FV scheme with reconstruction of characteristic variables, van Leer's minmod limiter~\cite{AbgrallMishranotes}, exact RS numerical flux, SSPRK2 time discretization and $C_{CFL}:=0.25$. Another possible and effective way to compute the reference solution would have been to reformulate the problem as a radial one--dimensional one with source term to be solved with high accuracy through a one--dimensional scheme, see~\cite{ToroBook}.

\RIIcolor{We remark that all order/flux configurations were able to reach the final time with $\sigma_{CFL}:=0.9$, except for FORCE-10 of order 7 on the $50\times50$ mesh and FORCE-10 of order 5 on the $200\times200$ mesh, for which $\sigma_{CFL}$ had to be reduced to $0.8$.
Concerning the results obtained on the $50\times50$ mesh, very little differences can be appreciated only for order 3, while, for orders 5 and 7 all numerical fluxes yield essentially indistinguishable results.
Instead, on the $200\times200$ mesh, no differences between the numerical fluxes can be appreciated at all for any order, with the reference solution being sharply captured by all schemes.
}

The results obtained on this test suggest the possibility to use FORCE--$\alpha$ numerical fluxes in place of upwind ones without appreciable loss of accuracy within the considered high order framework.

\subsubsection*{Long-time Liska--Wendroff explosion}
\begin{figure}[htbp]
	\centering
	
	\begin{subfigure}[b]{0.32\textwidth}
		\centering
		\includegraphics[width=\textwidth,trim={0 0 2.65cm 0},clip]
		{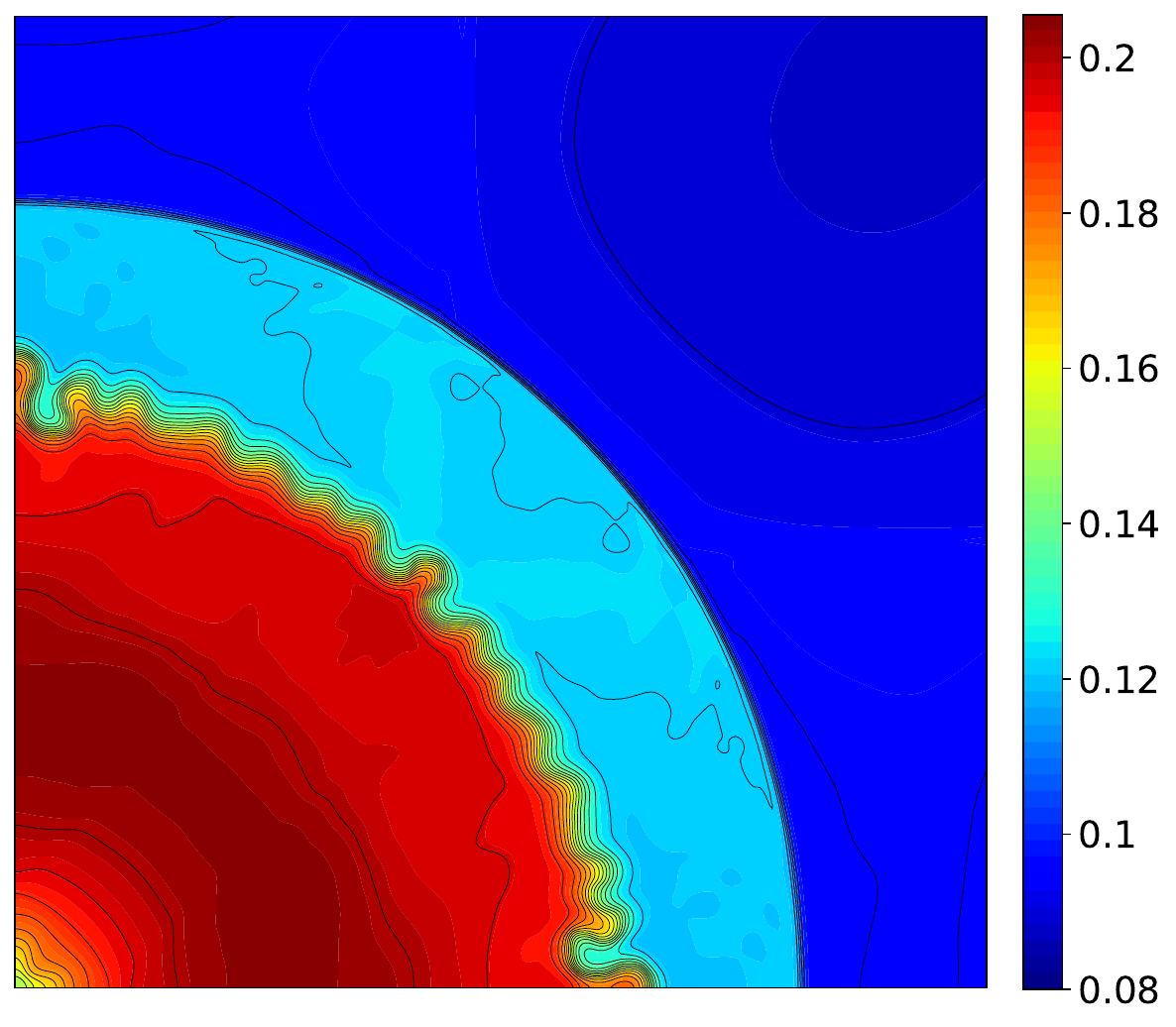}
		\caption{FORCE-2}
	\end{subfigure}
	\begin{subfigure}[b]{0.32\textwidth}
		\centering
		\includegraphics[width=\textwidth,trim={0 0 2.65cm 0},clip]
		{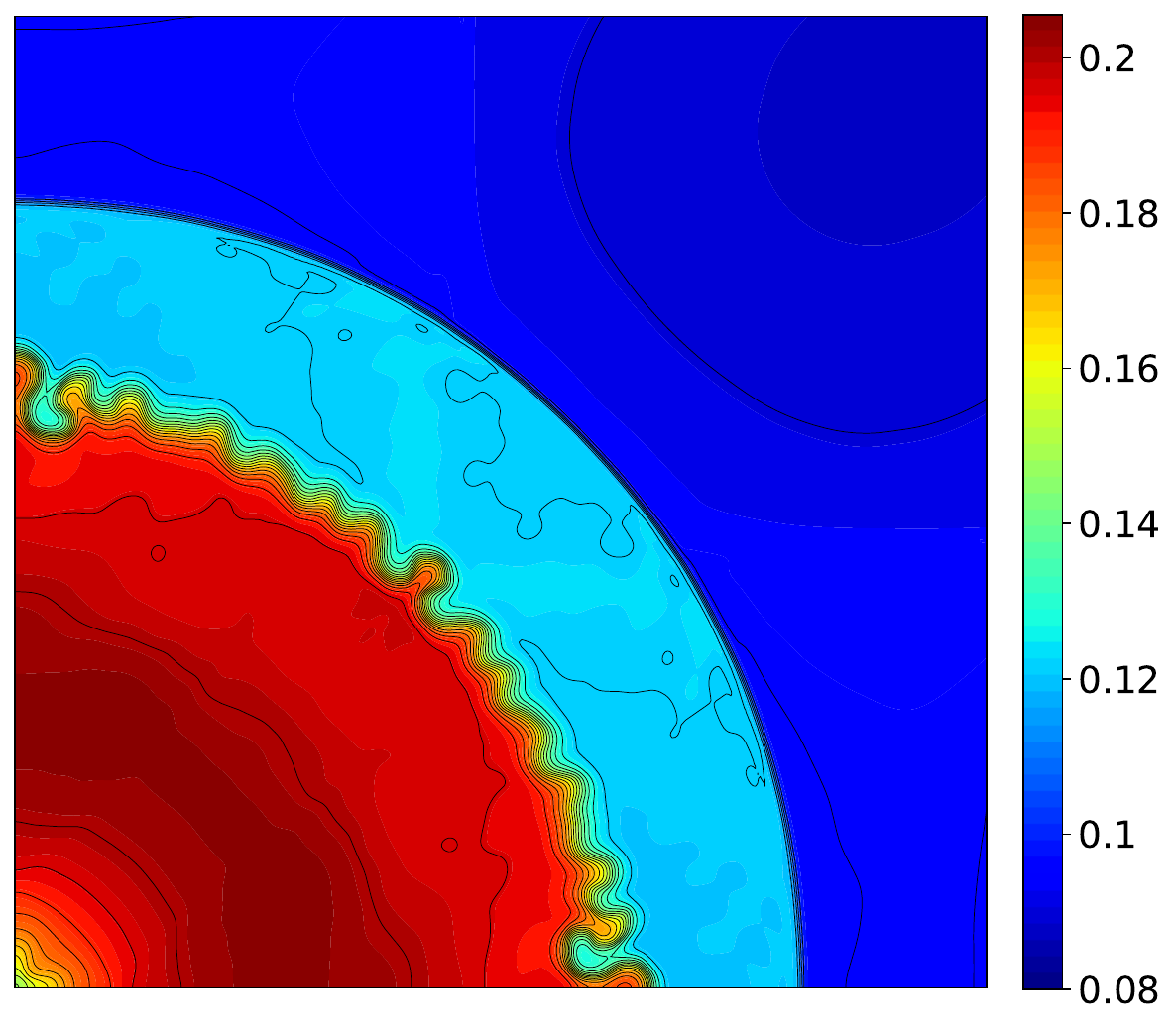}
		\caption{FORCE-3}
	\end{subfigure}\\
	\begin{subfigure}[b]{0.32\textwidth}
		\centering
		\includegraphics[width=\textwidth,trim={0 0 2.65cm 0},clip]
		{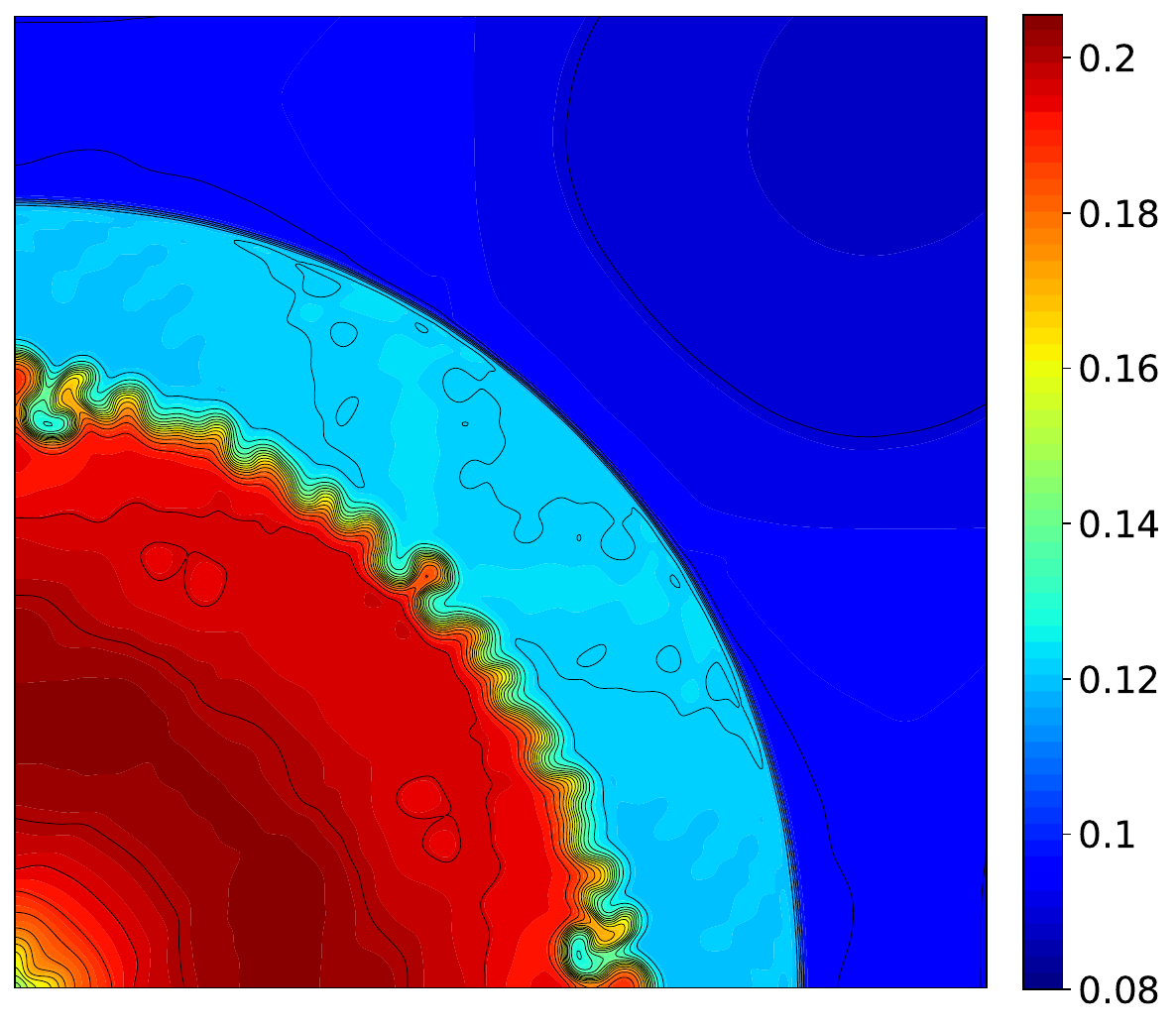}
		\caption{FORCE-5}
	\end{subfigure}
	\begin{subfigure}[b]{0.32\textwidth}
		\centering
		\includegraphics[width=\textwidth,trim={0 0 2.65cm 0},clip]
		{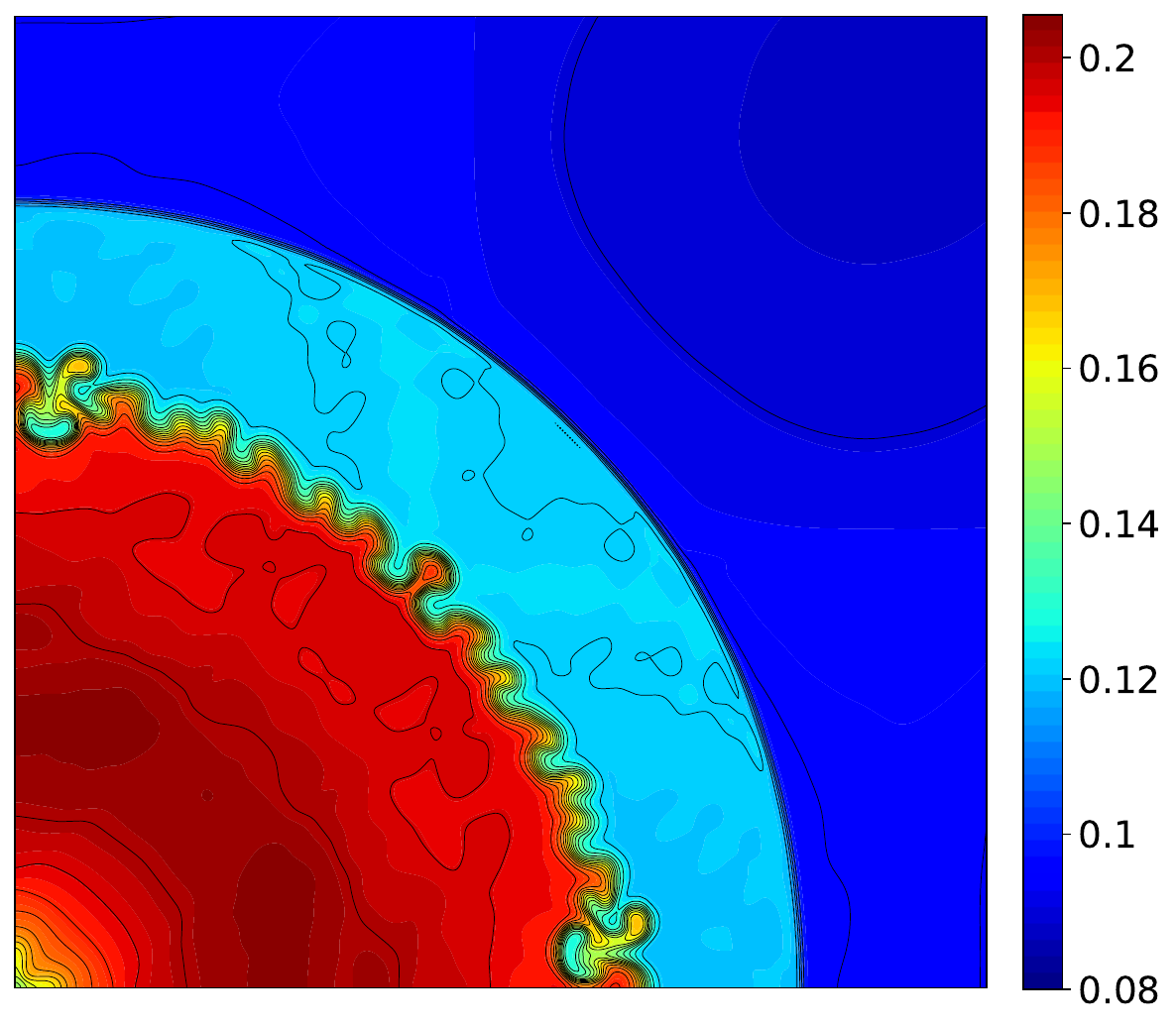}
		\caption{FORCE-10}
	\end{subfigure}
	
	\vspace{0.25em}
	
	\makebox[\textwidth][c]{%
		\begin{subfigure}[b]{0.32\textwidth}
			\centering
			\includegraphics[width=\textwidth,trim={0 0 2.65cm 0},clip]
			{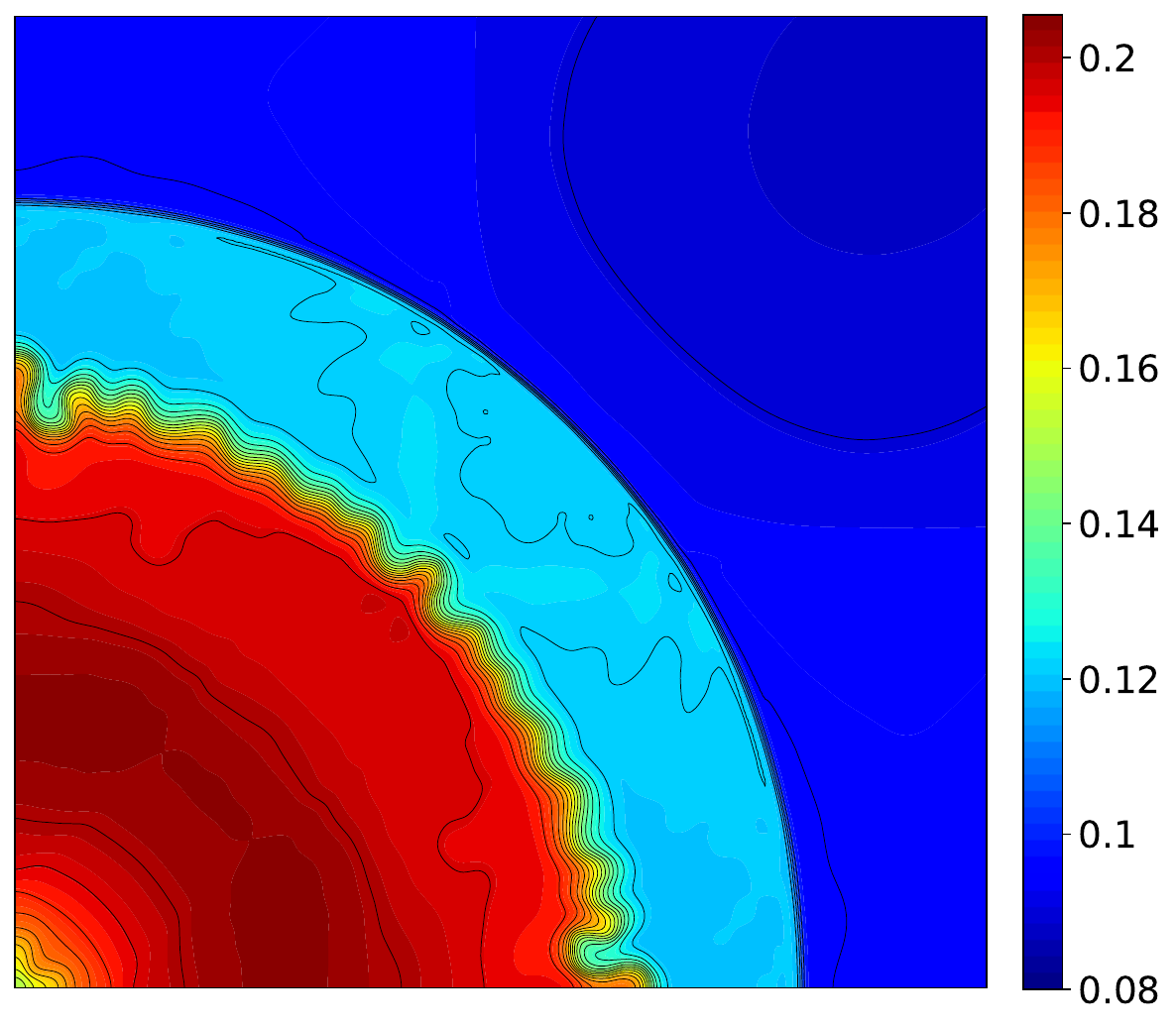}
			\caption{Rusanov}
		\end{subfigure}
		\begin{subfigure}[b]{0.32\textwidth}
			\centering
			\includegraphics[width=\textwidth,trim={0 0 2.65cm 0},clip]
			{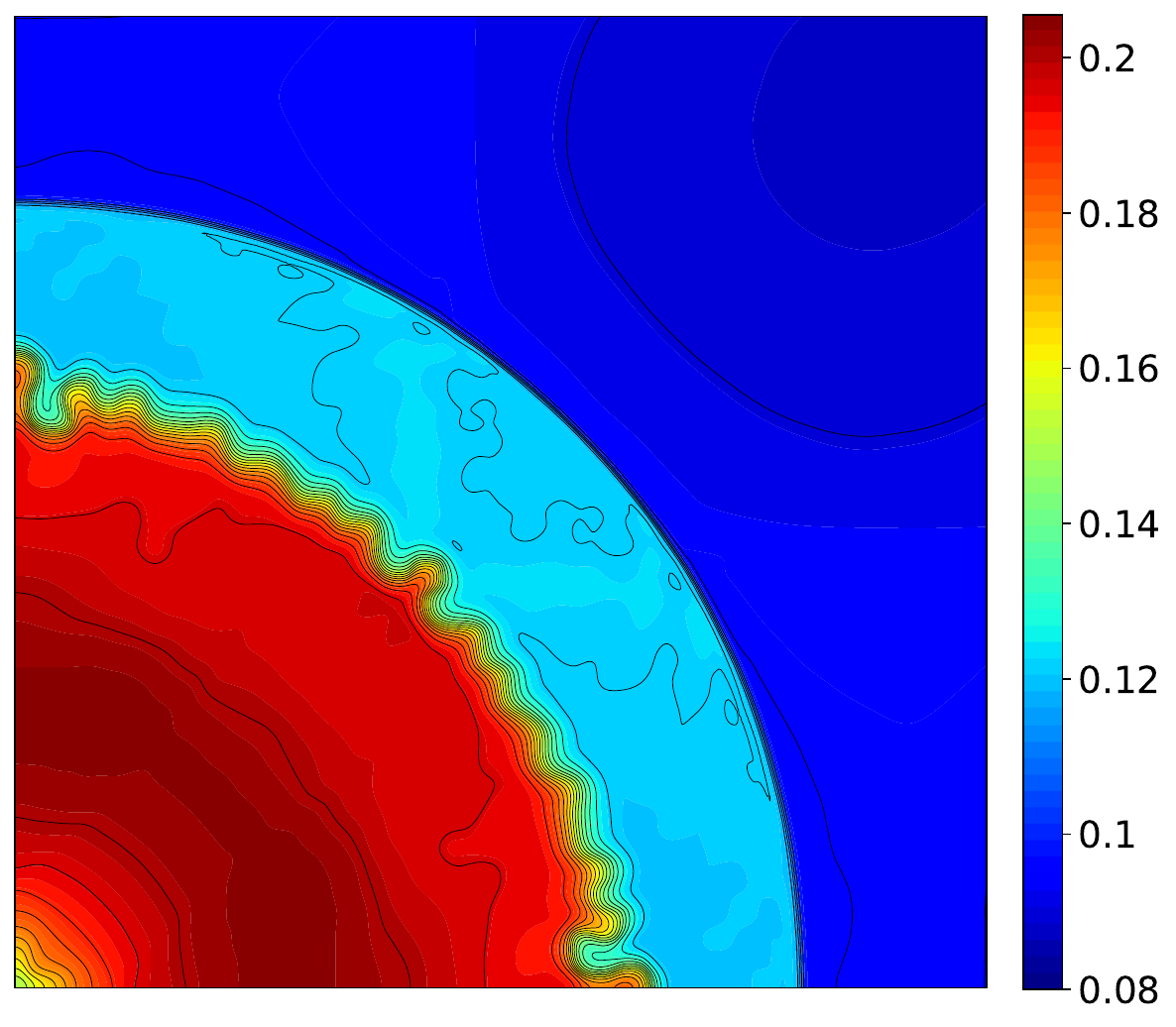}
			\caption{HLL}
		\end{subfigure}
		\begin{subfigure}[b]{0.32\textwidth}
			\centering
			\includegraphics[width=\textwidth,trim={0 0 2.65cm 0},clip]
			{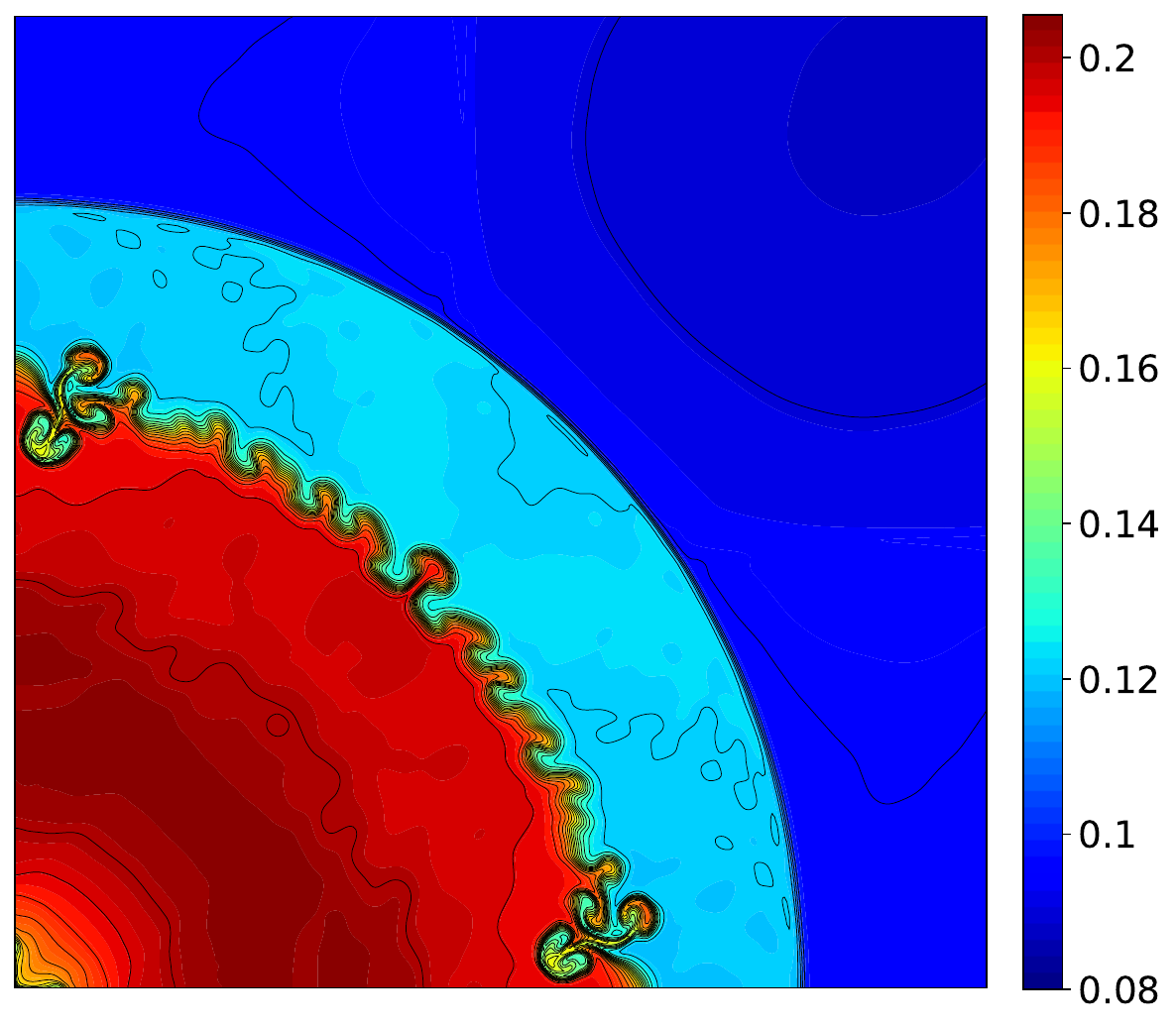}
			\caption{Exact RS}
		\end{subfigure}
	}
	
	\vspace{0.3em}
	\includegraphics[width=0.50\textwidth]
	{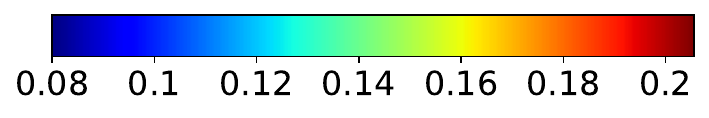}
	
	\caption{\RIIcolor{Long-time Liska--Wendroff explosion: Density obtained with order 3 over a mesh with $400\times 400$ elements and $\sigma_{CFL}:=0.9$.}}
	\label{fig:liska_wendroff_order3}
\end{figure}

\begin{figure}[htbp]
	\centering
	
	\begin{subfigure}[b]{0.32\textwidth}
		\centering
		\includegraphics[width=\textwidth,trim={0 0 2.65cm 0},clip]
		{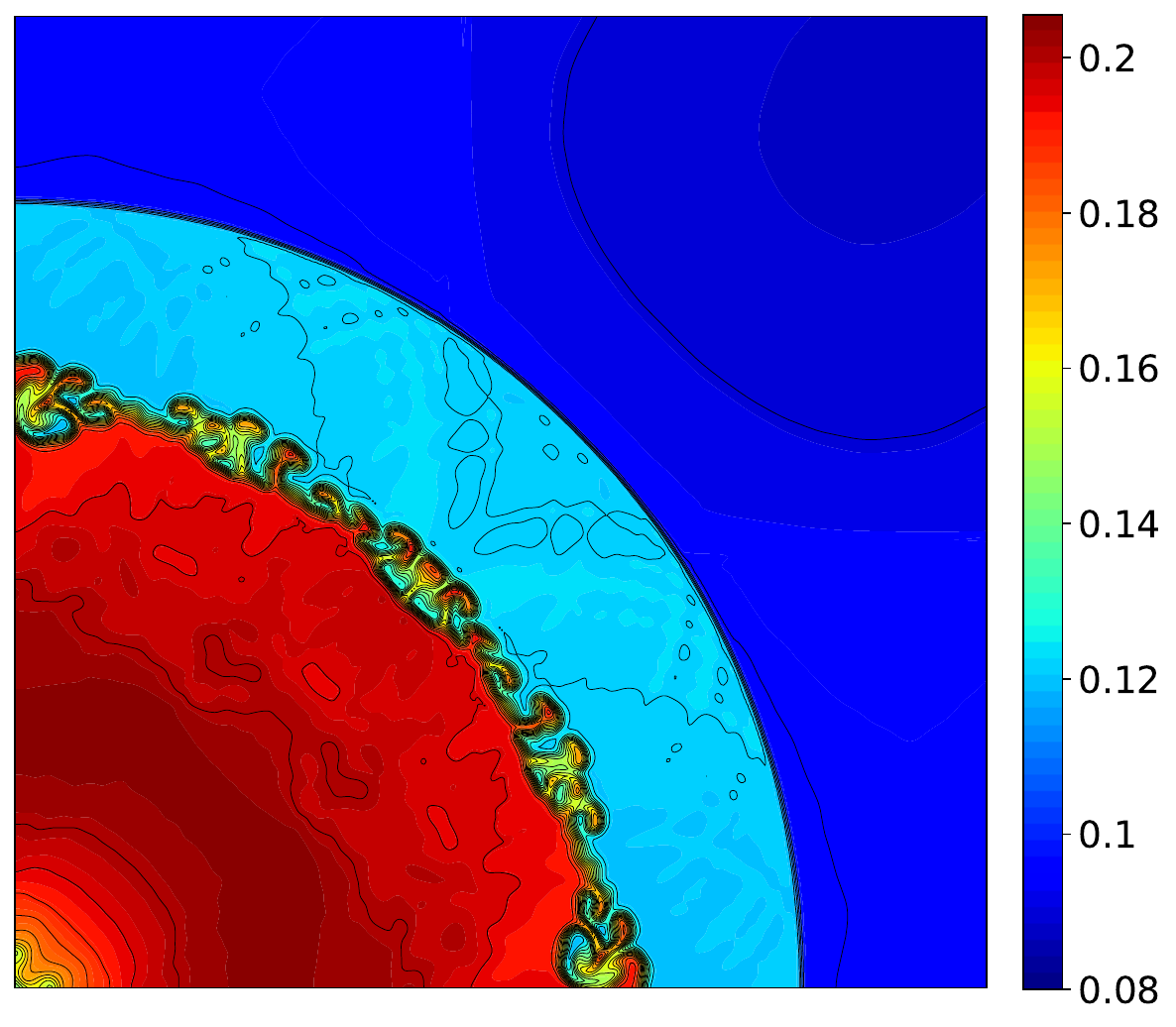}
		\caption{FORCE-2}
	\end{subfigure}
	\begin{subfigure}[b]{0.32\textwidth}
		\centering
		\includegraphics[width=\textwidth,trim={0 0 2.65cm 0},clip]
		{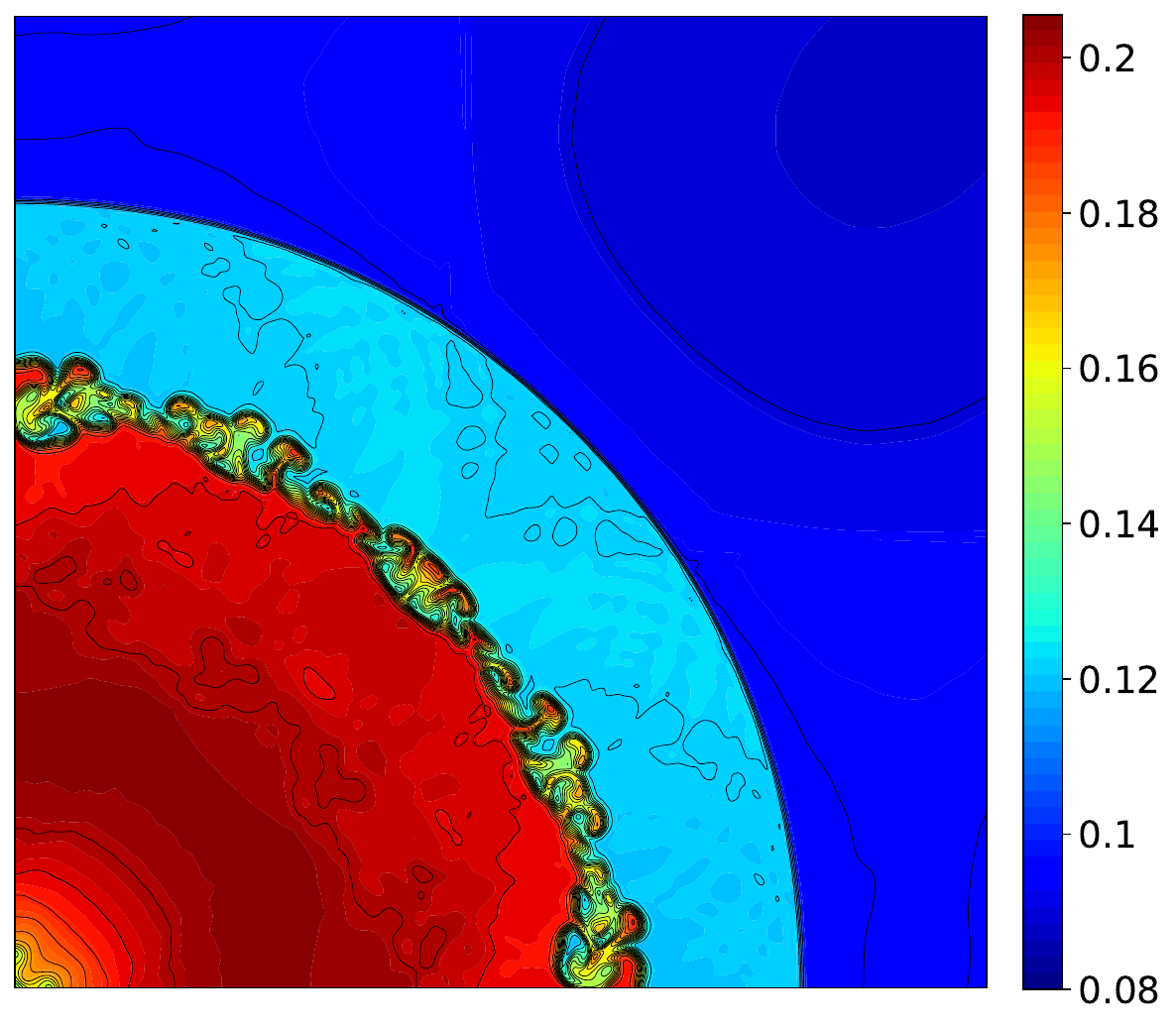}
		\caption{FORCE-3}
	\end{subfigure}\\
	\begin{subfigure}[b]{0.32\textwidth}
		\centering
		\includegraphics[width=\textwidth,trim={0 0 2.65cm 0},clip]
		{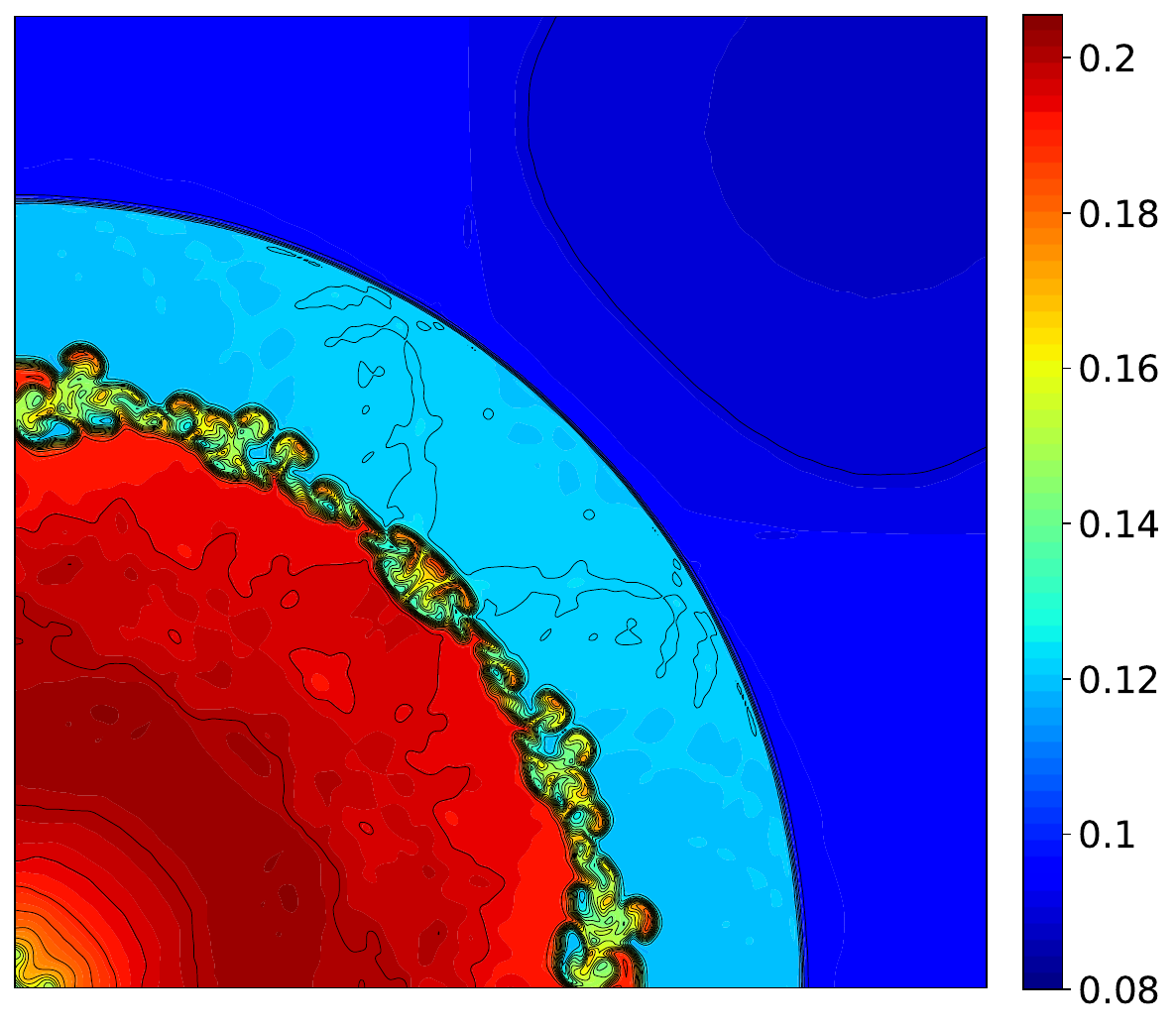}
		\caption{FORCE-5}
	\end{subfigure}
	\begin{subfigure}[b]{0.32\textwidth}
		\centering
		\includegraphics[width=\textwidth,trim={0 0 2.65cm 0},clip]
		{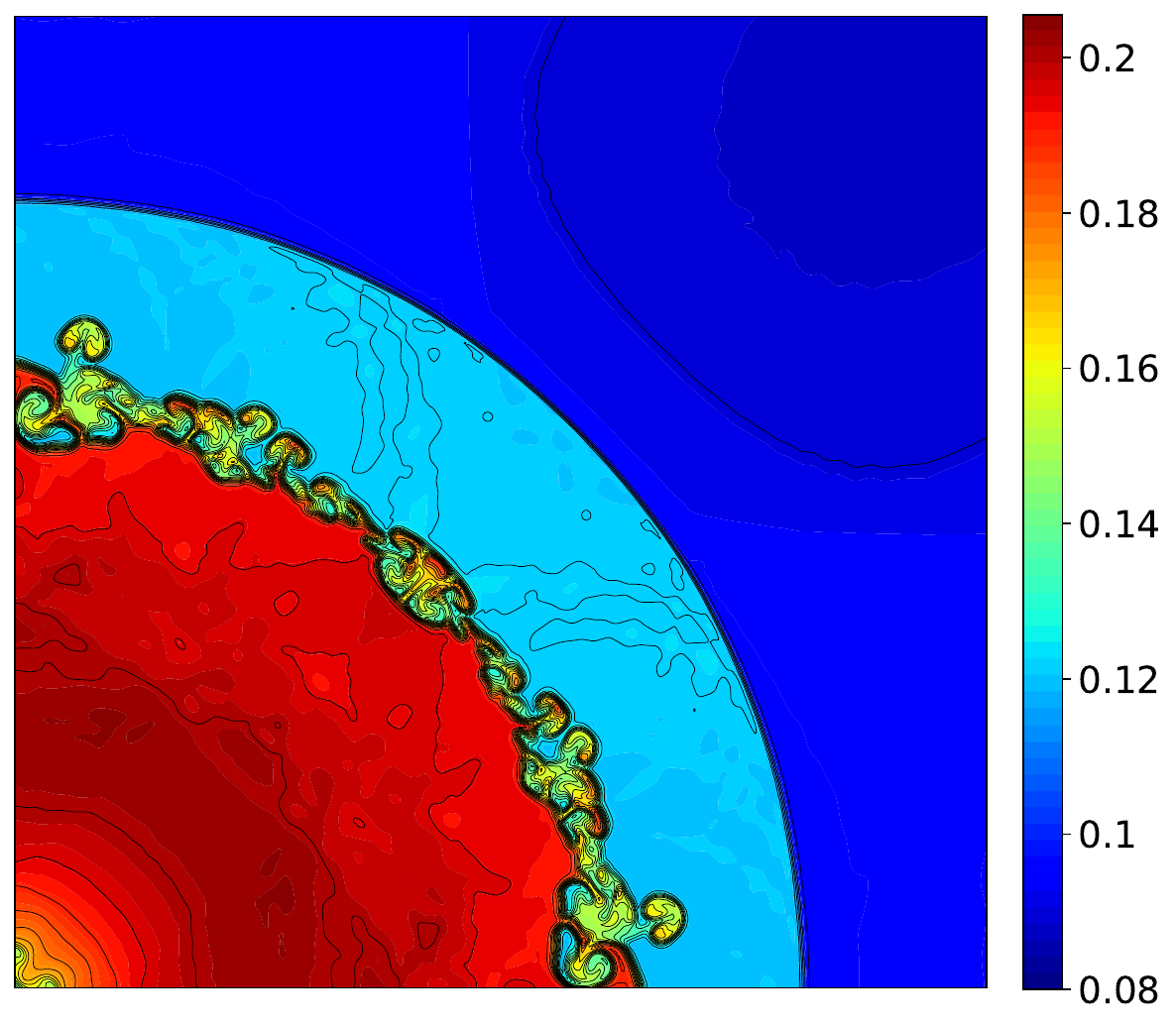}
		\caption{FORCE-10}
	\end{subfigure}
	
	\vspace{0.25em}
	
	\makebox[\textwidth][c]{%
		\begin{subfigure}[b]{0.32\textwidth}
			\centering
			\includegraphics[width=\textwidth,trim={0 0 2.65cm 0},clip]
			{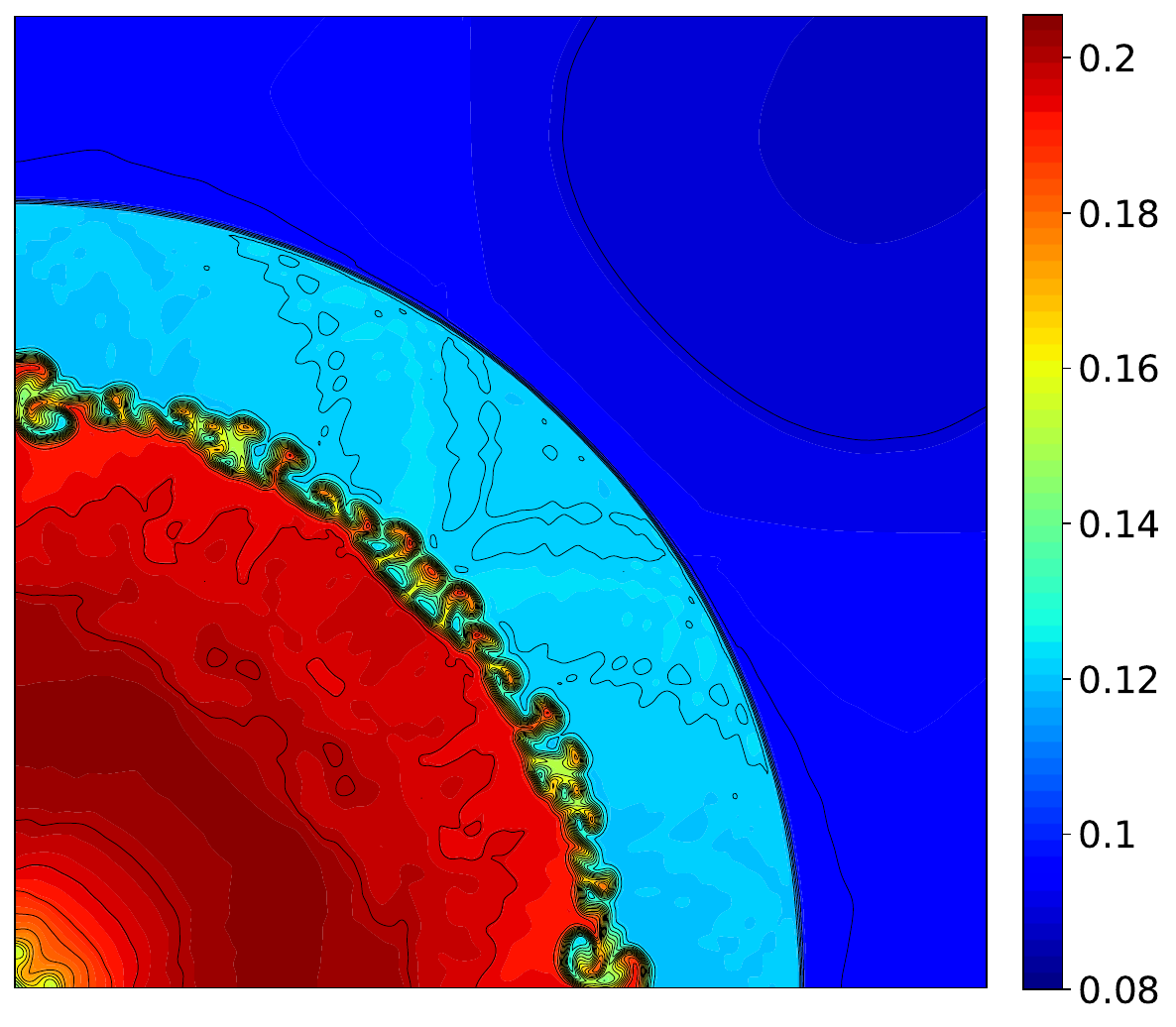}
			\caption{Rusanov}
		\end{subfigure}
		\begin{subfigure}[b]{0.32\textwidth}
			\centering
			\includegraphics[width=\textwidth,trim={0 0 2.65cm 0},clip]
			{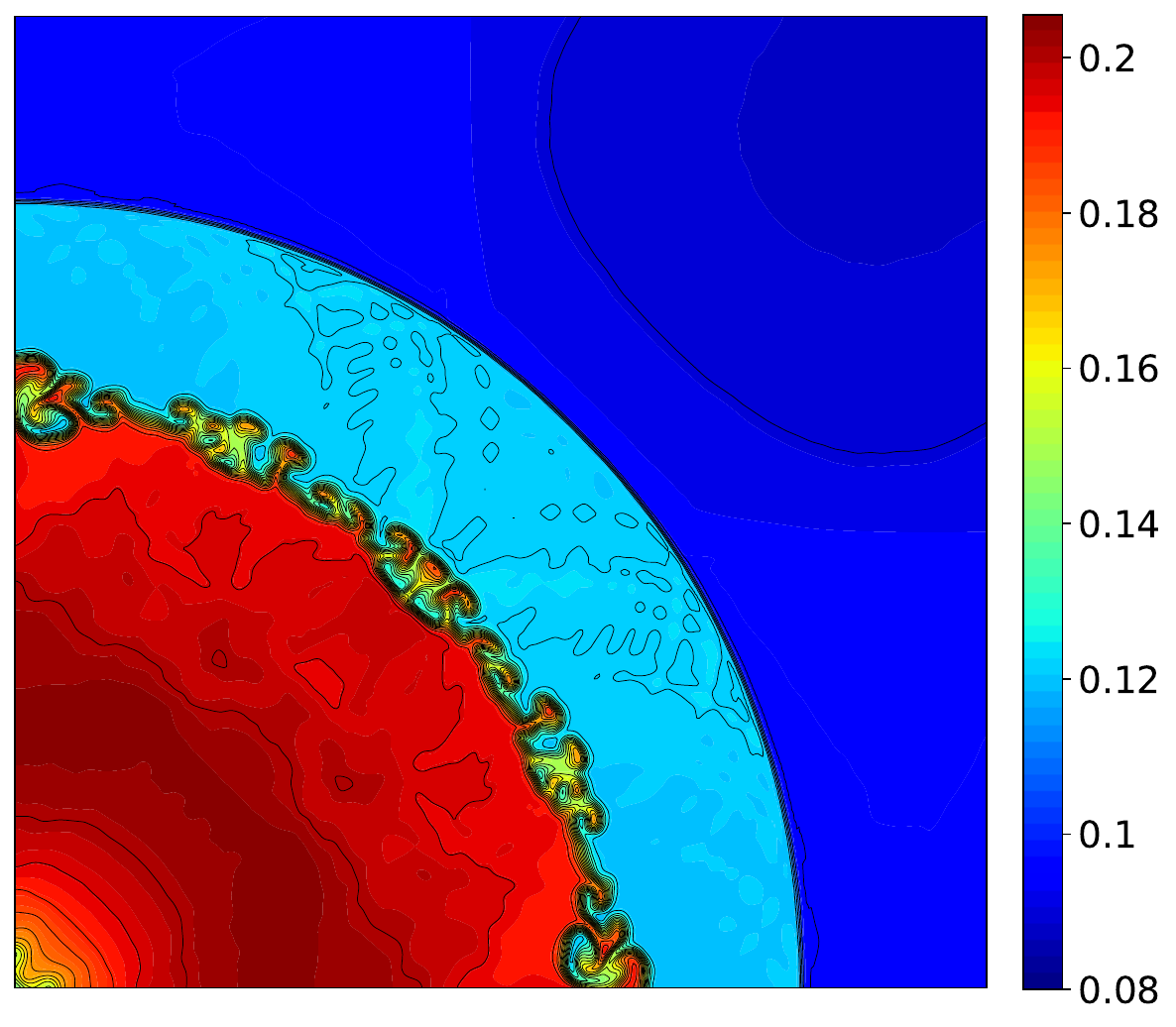}
			\caption{HLL}
		\end{subfigure}
		\begin{subfigure}[b]{0.32\textwidth}
			\centering
			\includegraphics[width=\textwidth,trim={0 0 2.65cm 0},clip]
			{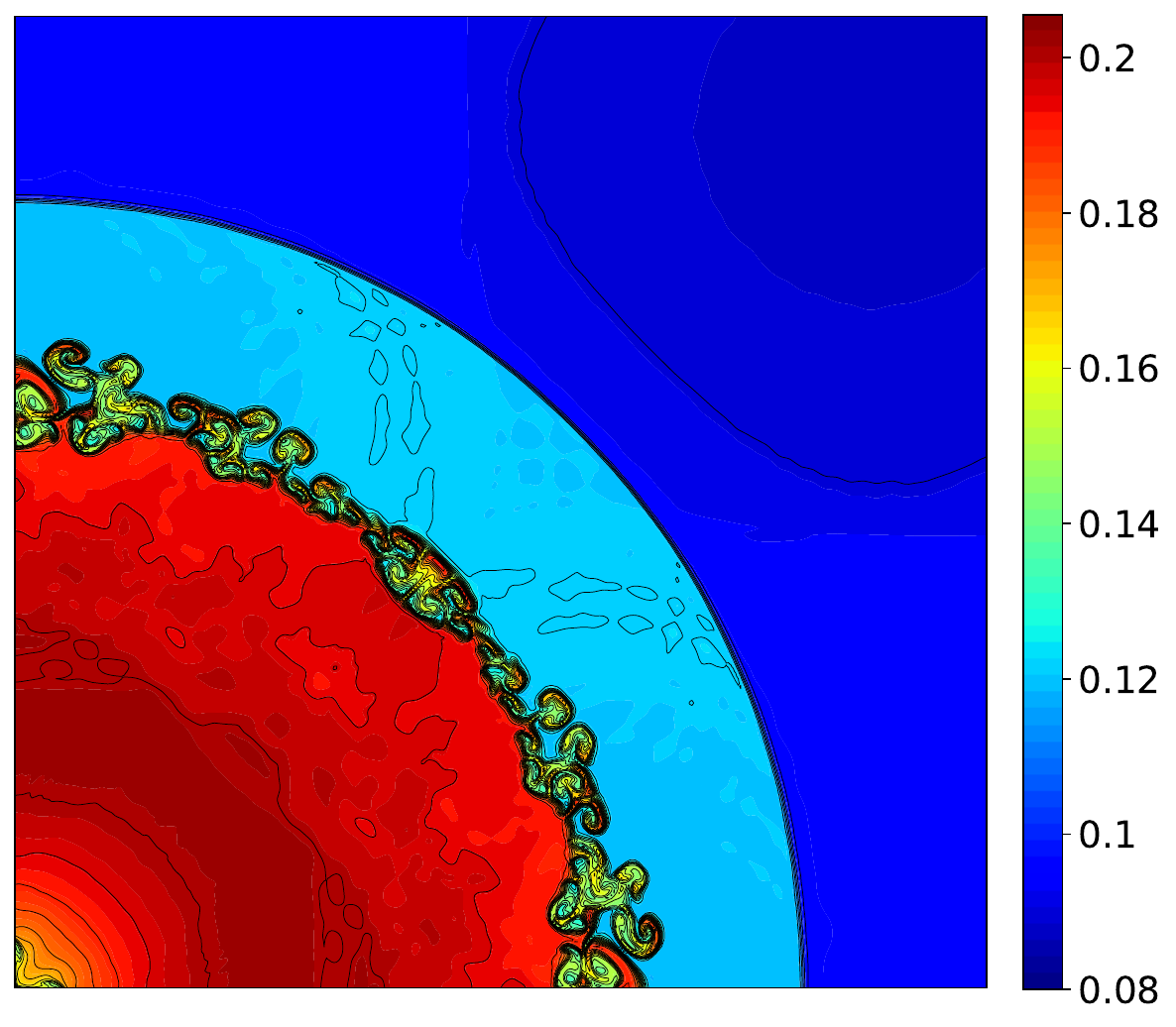}
			\caption{Exact RS}
		\end{subfigure}
	}
	
	\vspace{0.3em}
	\includegraphics[width=0.50\textwidth]
	{figures_new/liska_wendroff_common_density_colorbar_horizontal.pdf}
	
	\caption{\RIIcolor{Long-time Liska--Wendroff explosion: Density obtained with order 5 over a mesh with $400\times 400$ elements and $\sigma_{CFL}:=0.9$.}}
	\label{fig:liska_wendroff_order5}
\end{figure}

\begin{figure}[htbp]
	\centering
	
	\begin{subfigure}[b]{0.48\textwidth}
		\centering
		\includegraphics[width=\textwidth,trim={0 0 2.65cm 0},clip]
		{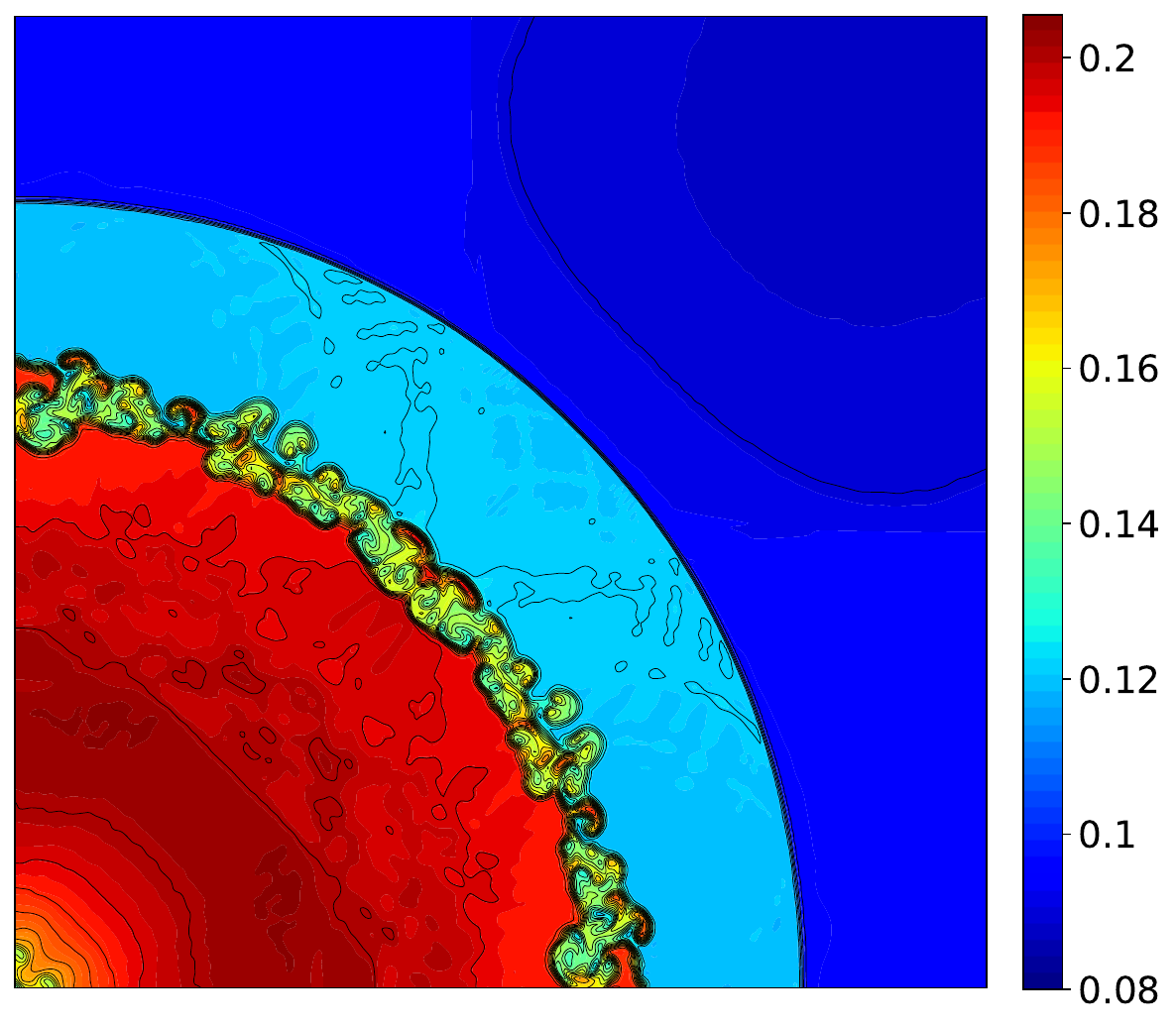}
		\caption{FORCE-2}
	\end{subfigure}
	\begin{subfigure}[b]{0.48\textwidth}
		\centering
		\includegraphics[width=\textwidth,trim={0 0 2.65cm 0},clip]
		{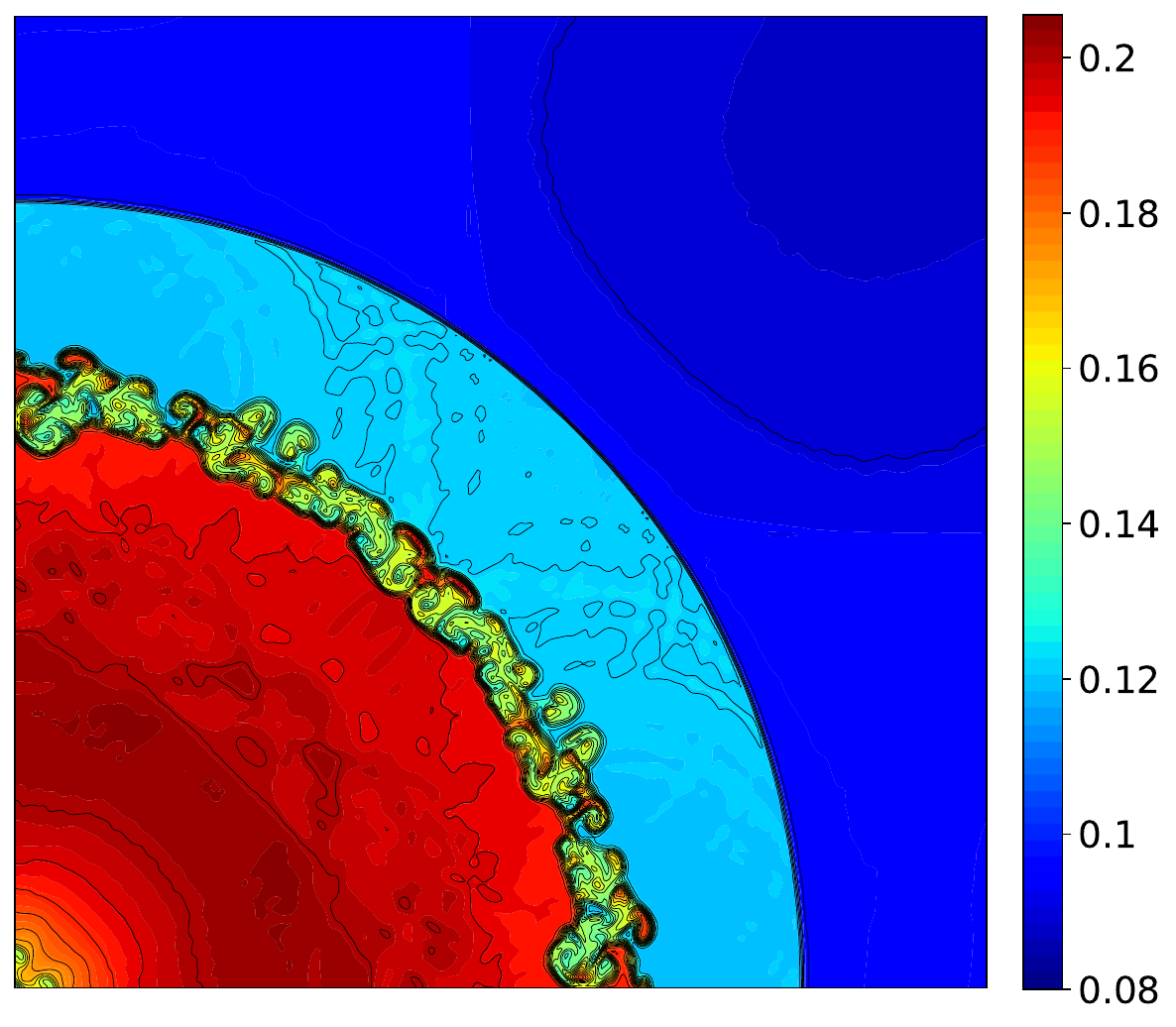}
		\caption{FORCE-3}
	\end{subfigure}
	
	\vspace{0.25em}
	
	\begin{subfigure}[b]{0.48\textwidth}
		\centering
		\includegraphics[width=\textwidth,trim={0 0 2.65cm 0},clip]
		{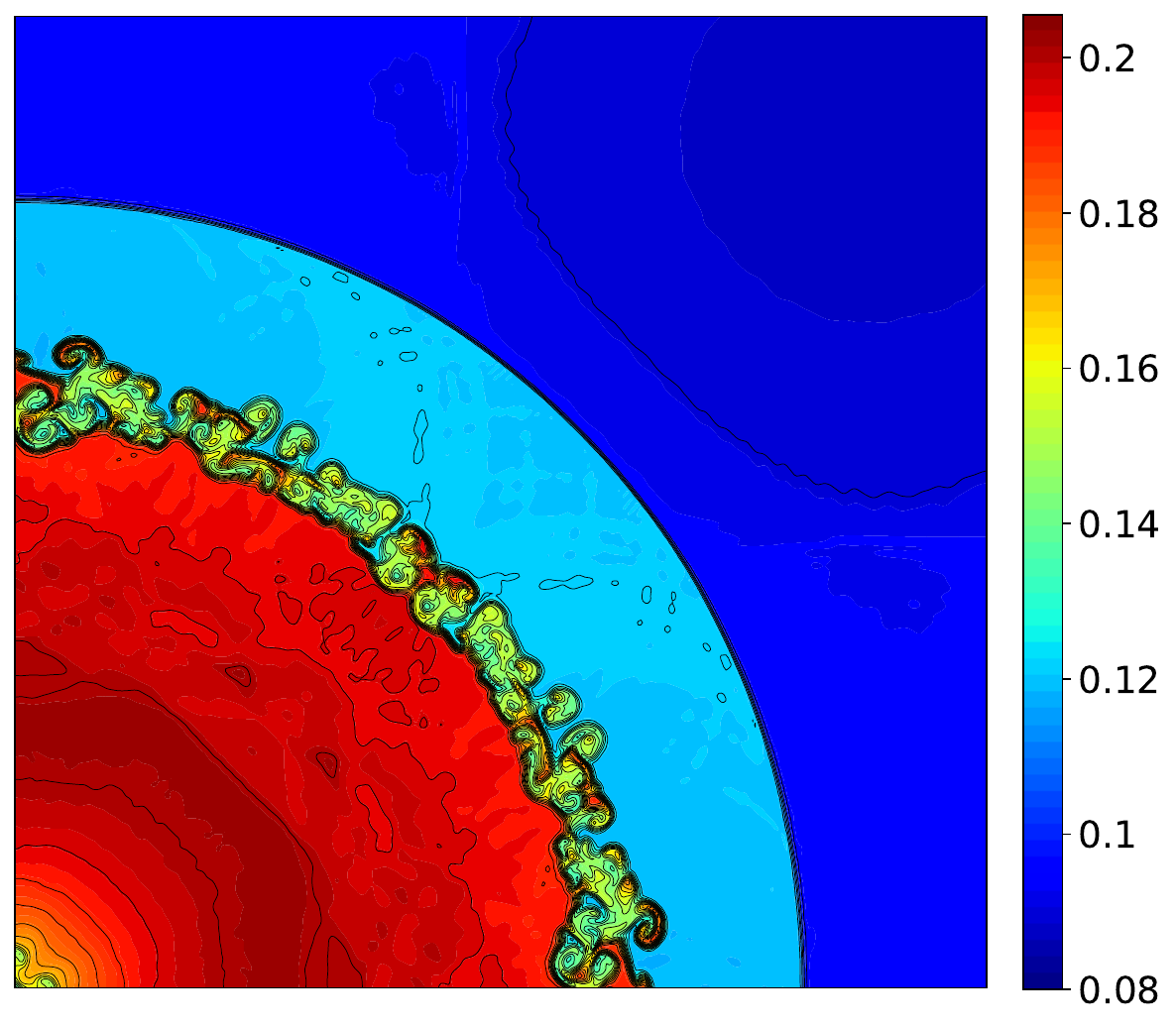}
		\caption{FORCE-5}
	\end{subfigure}
	\begin{subfigure}[b]{0.48\textwidth}
		\centering
		\includegraphics[width=\textwidth,trim={0 0 2.65cm 0},clip]
		{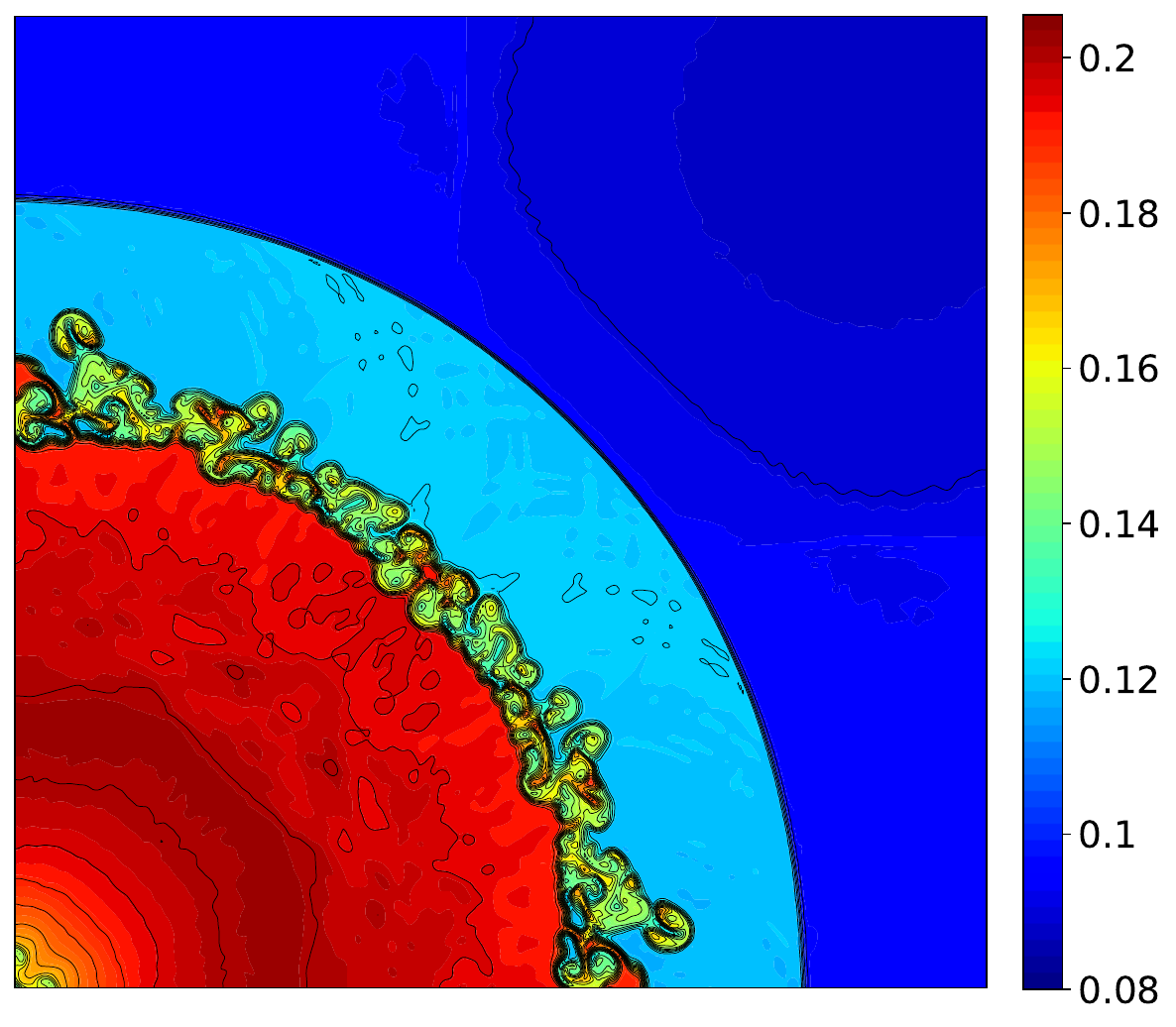}
		\caption{FORCE-10}
	\end{subfigure}
	
	\vspace{0.3em}
	\includegraphics[width=0.50\textwidth]
	{figures_new/liska_wendroff_common_density_colorbar_horizontal.pdf}
	
	\caption{\RIIcolor{Long-time Liska--Wendroff explosion: Density obtained with order 7 over a mesh with $400\times 400$ elements and $\sigma_{CFL}:=0.9$.}}
	\label{fig:liska_wendroff_order7}
\end{figure}

\RIIcolor{This test has been considered in~\cite{liska2003comparison,CKX_Ustar} and consists of a modification of the previous test from~\cite{ToroBook} involving the same initial condition but on a smaller computational domain, more specifically the square quadrant $[0,1.5]\times[0,1.5]$.
As in~\cite{CKX_Ustar}, we impose solid wall reflective boundary conditions at $x = 0$ and $y = 0$, while transmissive ones are imposed at $x = 1.5$ and $y = 1.5$.
The final time, set to $T_f:=3.2$, is much longer and determines the development of instabilities, whose handling can be challenging at the discrete level.
In fact, numerical schemes are required to capture such instabilities while guaranteeing a high level of robustness to prevent simulation crashes.

We ran our simulations on meshes with $400\times 400$ elements for $\sigma_{CFL}:=0.9$.
All configurations reached the final time, except for all upwind fluxes at order 7.
The contour plots of the density are reported in Figures~\ref{fig:liska_wendroff_order3}, \ref{fig:liska_wendroff_order5} and \ref{fig:liska_wendroff_order7}.

For order 3, exact RS is the sharpest numerical flux.
All other numerical fluxes follow with similar results.
Among them, FORCE-10 provides the best quality, followed by FORCE-3 and FORCE-5. 
Rusanov, HLL and FORCE-2 are slightly more diffusive and produce similar large-scale structures.

The situation changes when considering order 5.
The results obtained through different numerical fluxes become much more consistent among themselves. Still, exact RS, and to some extent FORCE-10, display slightly richer small-scale structures.

For order 7, the quality of the results further improves and all considered FORCE--$\alpha$ numerical fluxes lead to comparable resolution of the vortical structures.
We remark once again that only the FORCE--$\alpha$ numerical fluxes were able to reach the final time for order 7, proving to be robust alternatives to upwind fluxes in this long-time, instability-dominated regime.
}


\subsubsection{Two--dimensional Riemann problem}\label{sec:2DRP}
\begin{figure}[htbp]
	\centering
	
	\begin{subfigure}[b]{0.32\textwidth}
		\centering
		\includegraphics[width=\textwidth,trim={0 0 2.65cm 0},clip]
		{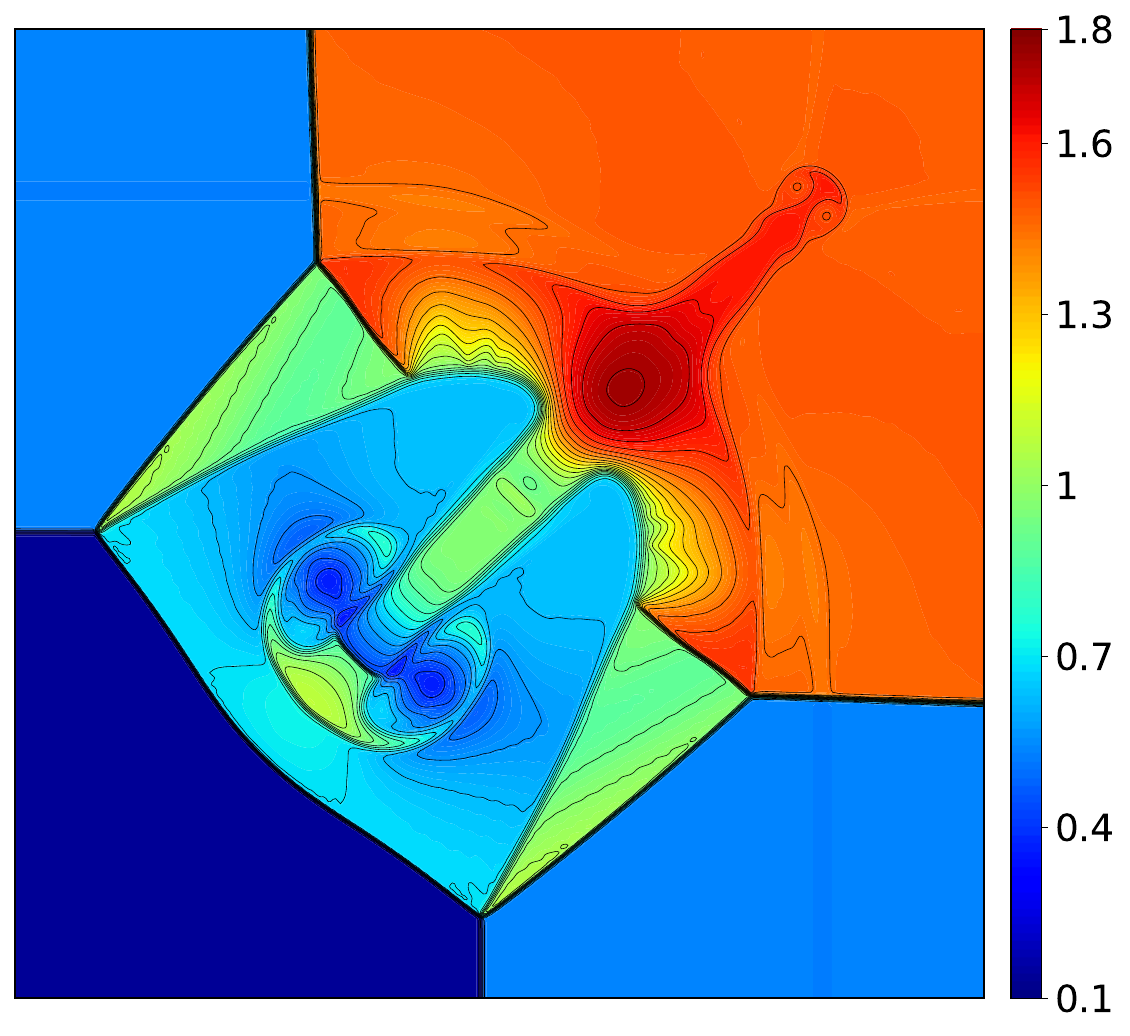}
		\caption{FORCE-2}
	\end{subfigure}
	\begin{subfigure}[b]{0.32\textwidth}
		\centering
		\includegraphics[width=\textwidth,trim={0 0 2.65cm 0},clip]
		{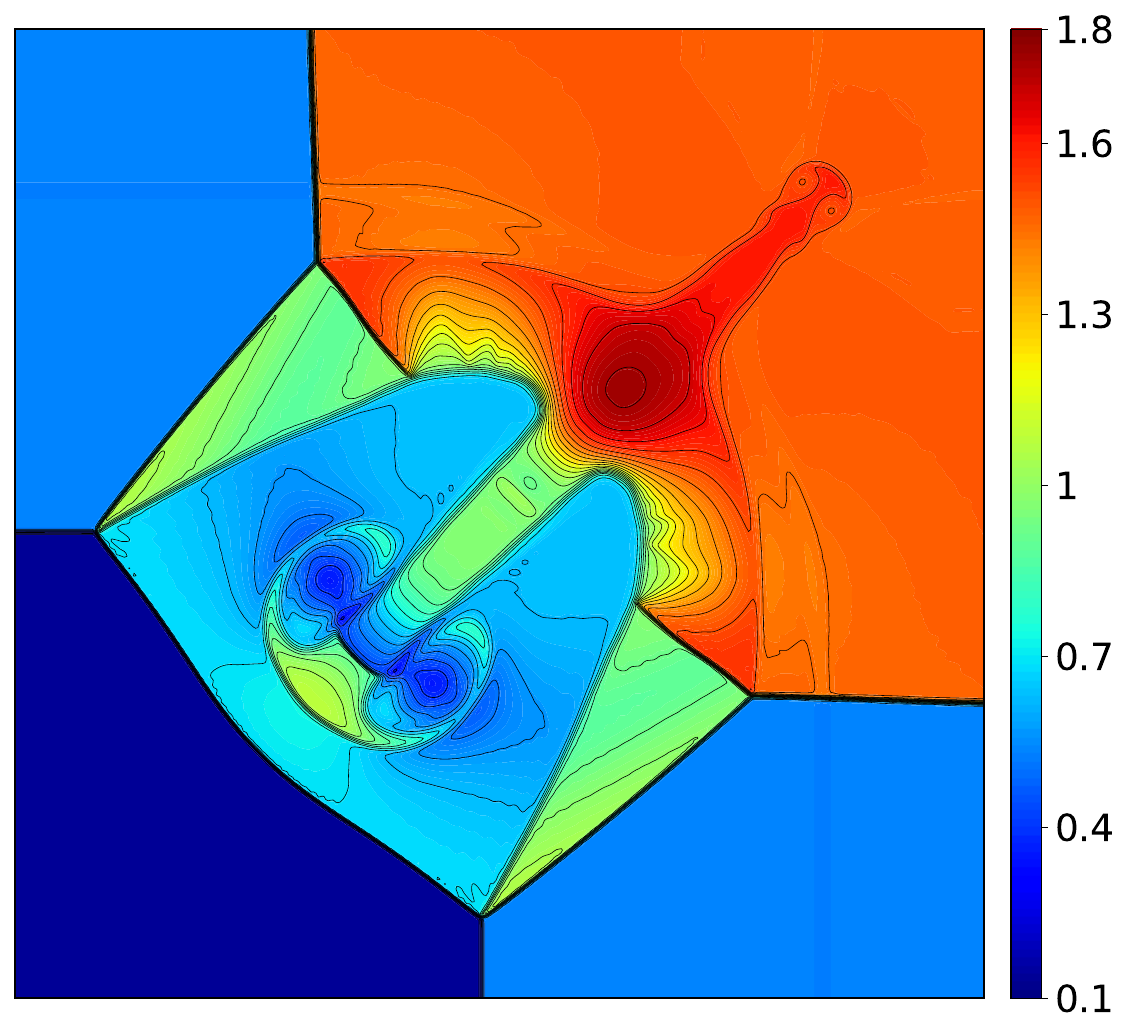}
		\caption{FORCE-3}
	\end{subfigure}
	
	\begin{subfigure}[b]{0.32\textwidth}
		\centering
		\includegraphics[width=\textwidth,trim={0 0 2.65cm 0},clip]
		{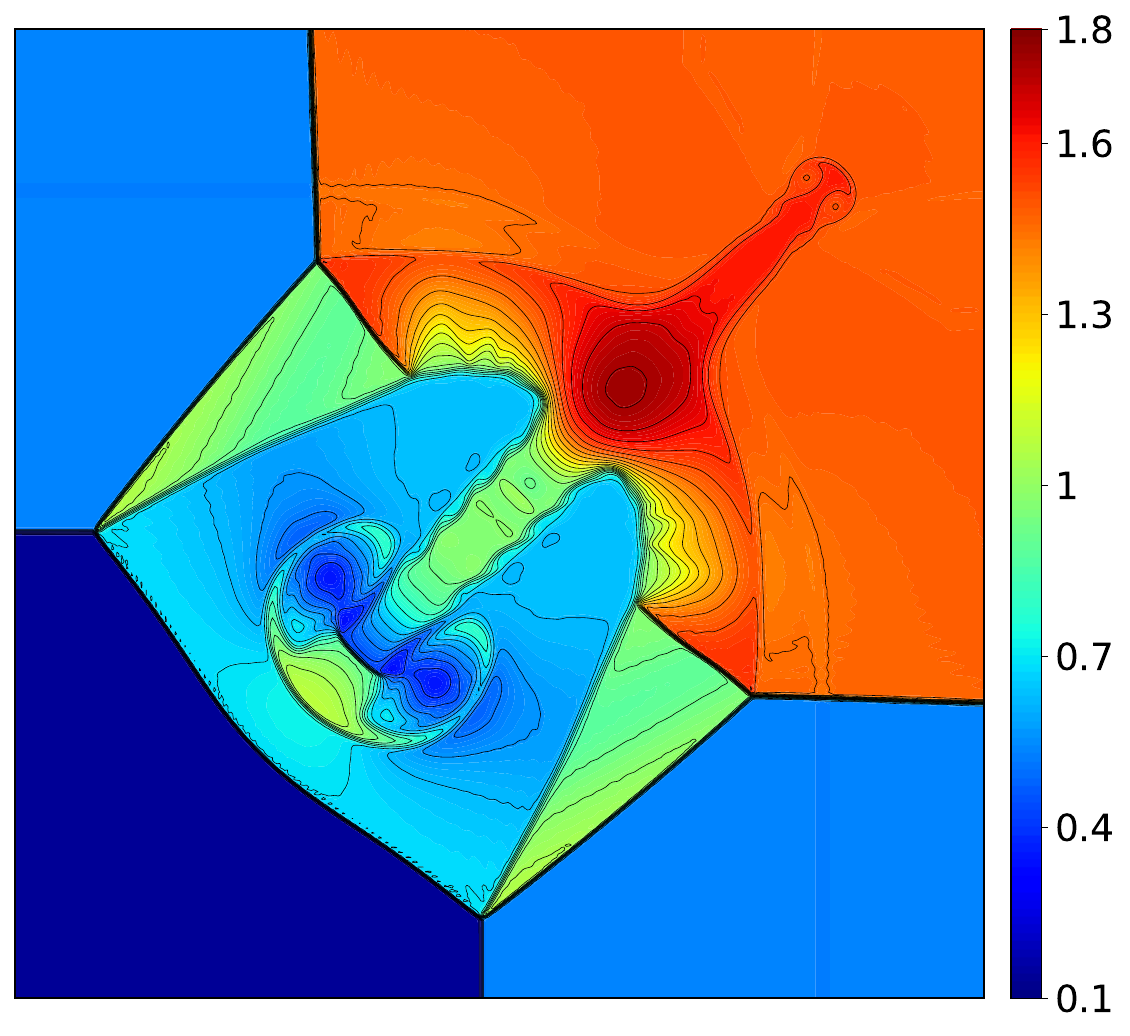}
		\caption{FORCE-5}
	\end{subfigure}
	\begin{subfigure}[b]{0.32\textwidth}
		\centering
		\includegraphics[width=\textwidth,trim={0 0 2.65cm 0},clip]
		{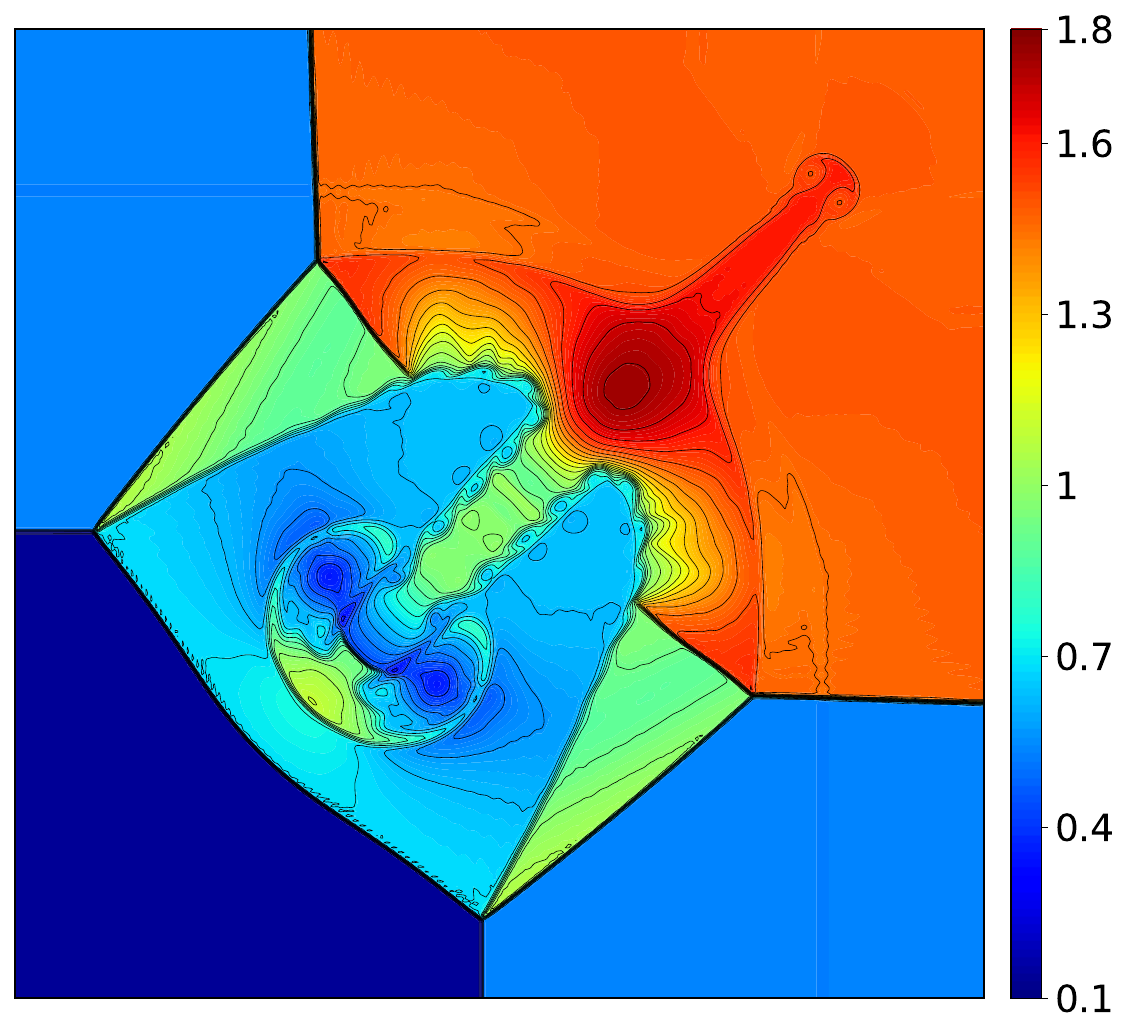}
		\caption{FORCE-10}
	\end{subfigure}
	
	\begin{subfigure}[b]{0.32\textwidth}
		\centering
		\includegraphics[width=\textwidth,trim={0 0 2.65cm 0},clip]
		{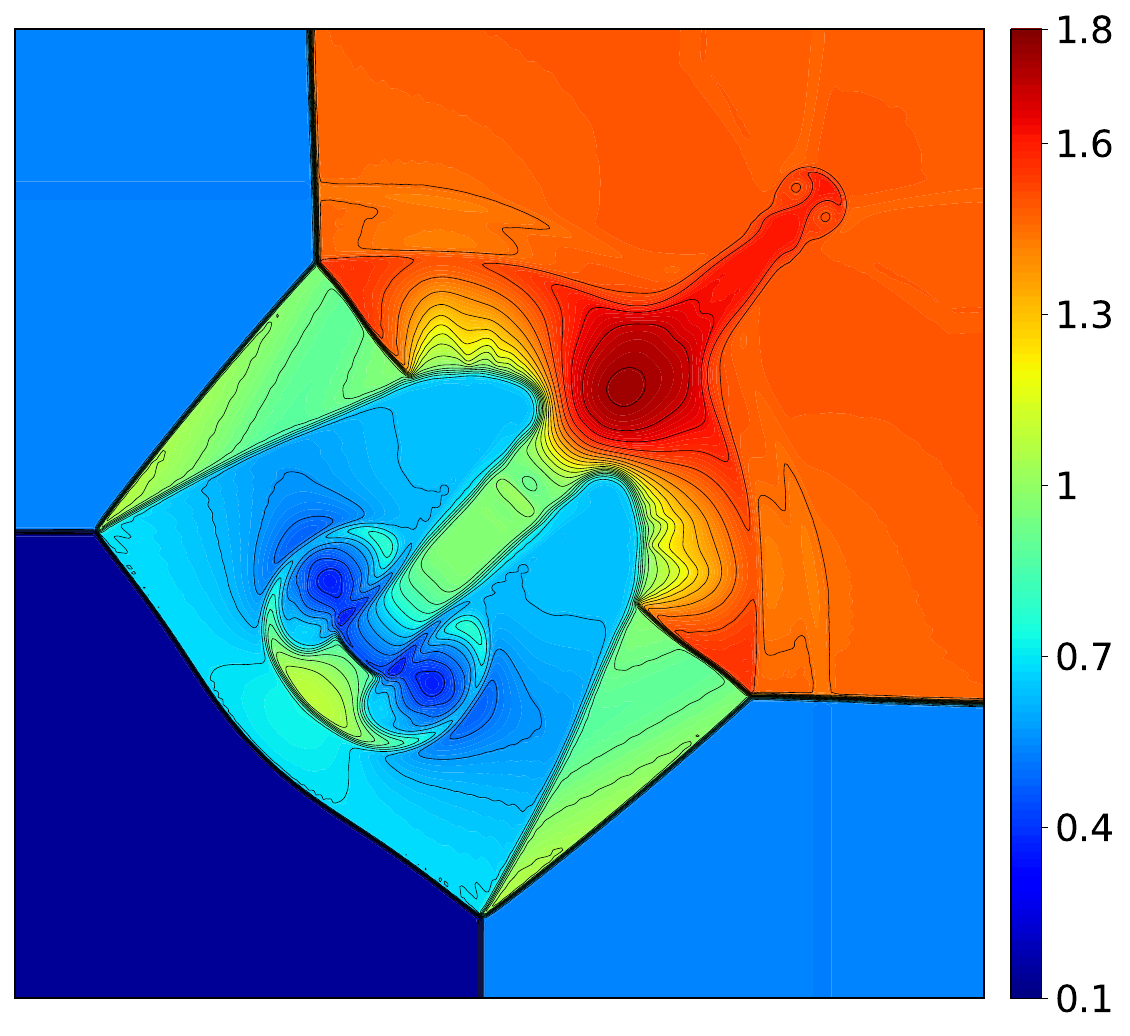}
		\caption{Rusanov}
	\end{subfigure}
	\hfill
	\begin{subfigure}[b]{0.32\textwidth}
		\centering
		\includegraphics[width=\textwidth,trim={0 0 2.65cm 0},clip]
		{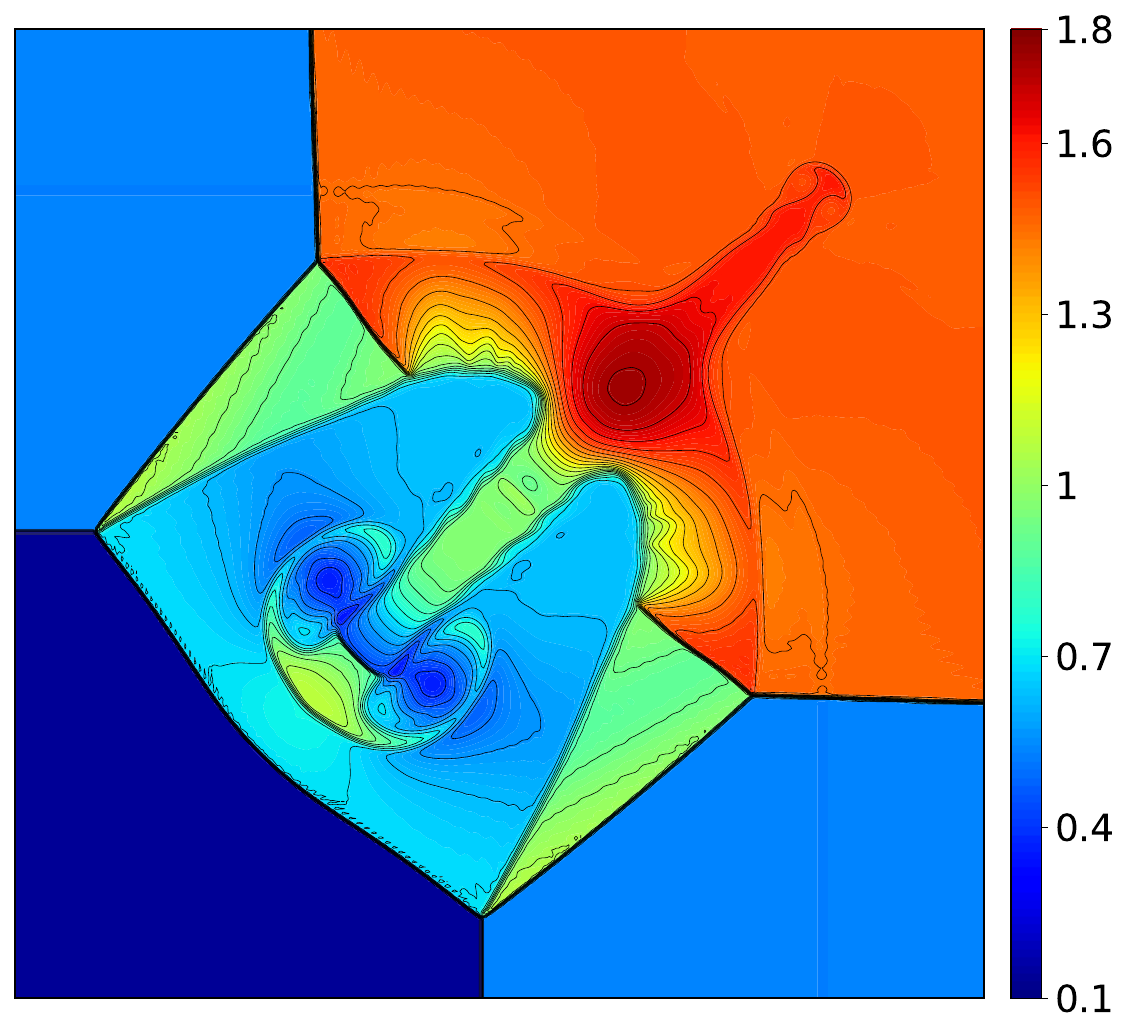}
		\caption{HLL}
	\end{subfigure}
	\hfill
	\begin{subfigure}[b]{0.32\textwidth}
		\centering
		\includegraphics[width=\textwidth,trim={0 0 2.65cm 0},clip]
		{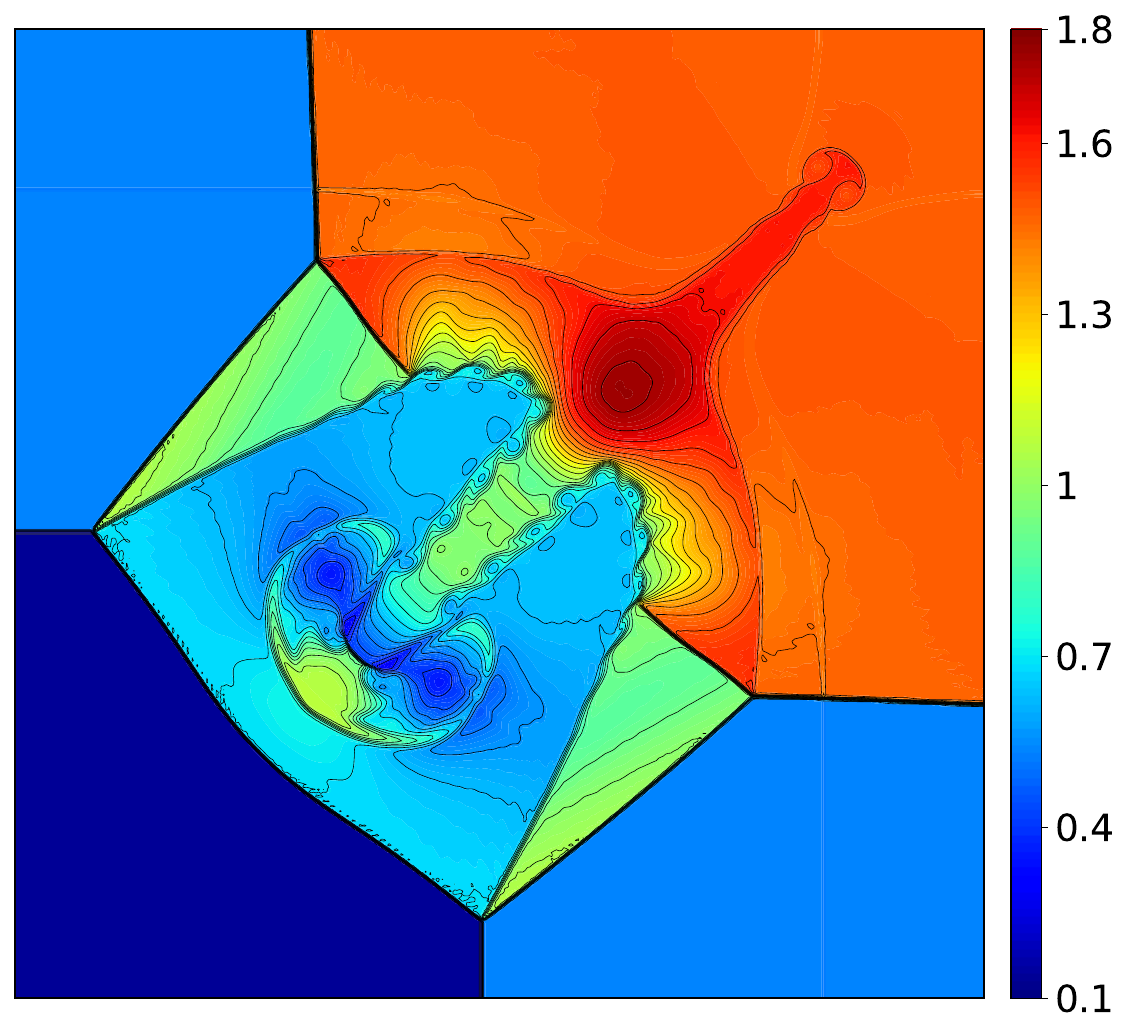}
		\caption{Exact RS}
	\end{subfigure}
	
	\vspace{0.4em}
	\includegraphics[width=0.42\textwidth]
	{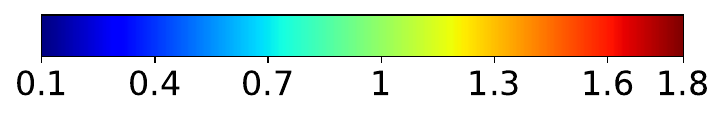}
	
	\caption{\RIIcolor{Two-dimensional Riemann problem: Density obtained with order 5 over a mesh with $400\times 400$ elements and $\sigma_{CFL}:=0.9$. The same density color scale is used in all panels.}}
	\label{fig:2DRP_order5}
\end{figure}

\begin{figure}[htbp]
	\centering
	
	\begin{subfigure}[b]{0.32\textwidth}
		\centering
		\includegraphics[width=\textwidth,trim={0 0 2.65cm 0},clip]
		{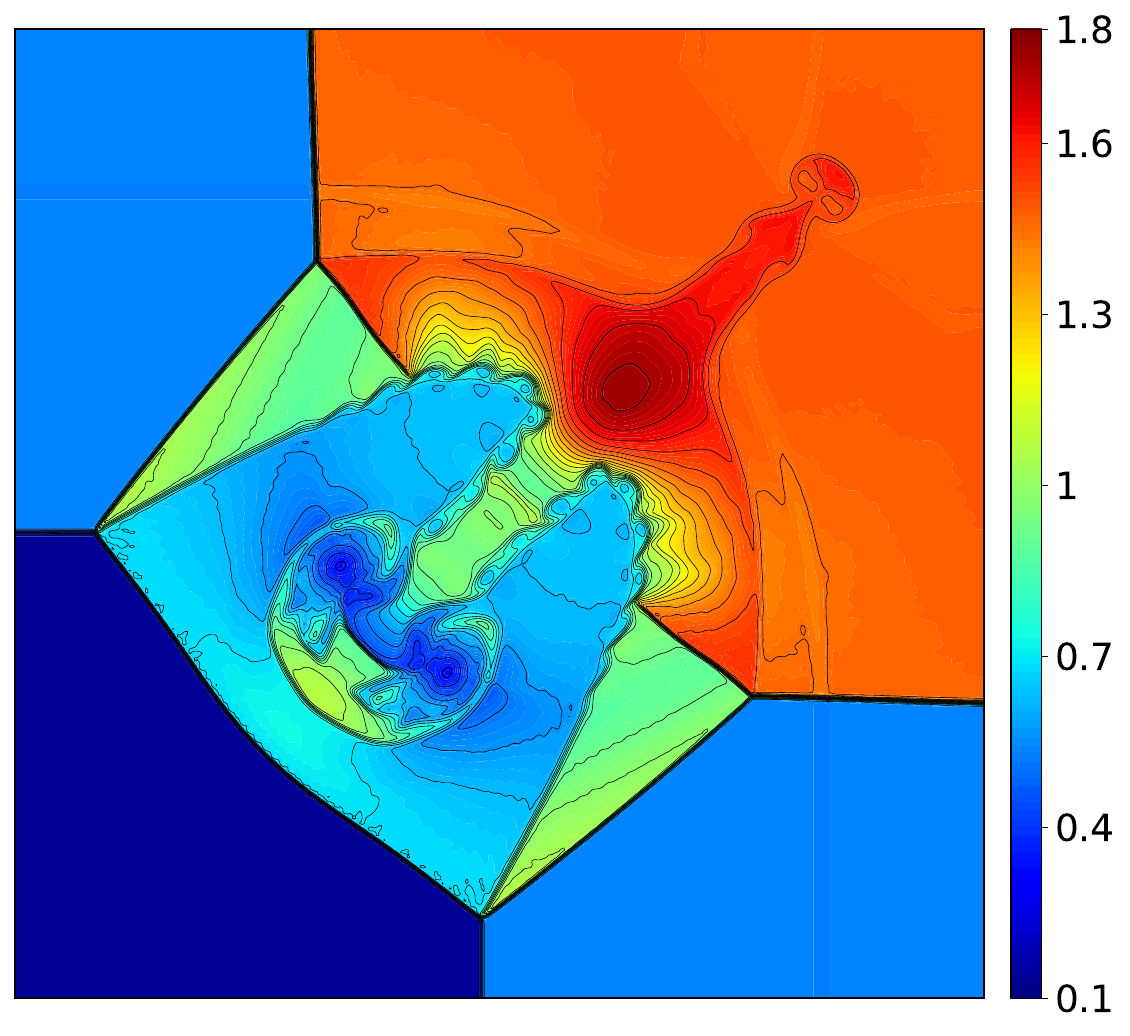}
		\caption{FORCE-2}
	\end{subfigure}
	\begin{subfigure}[b]{0.32\textwidth}
		\centering
		\includegraphics[width=\textwidth,trim={0 0 2.65cm 0},clip]
		{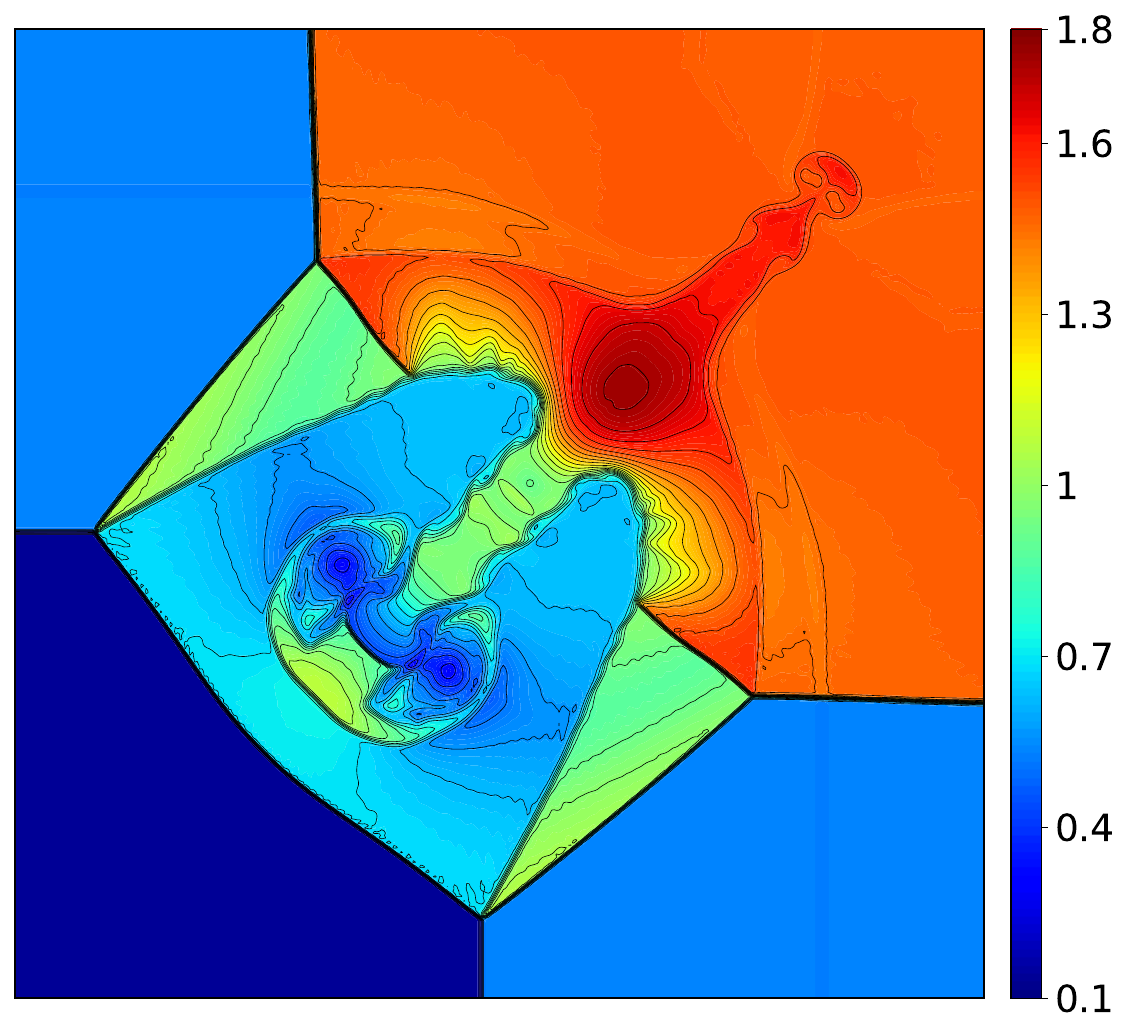}
		\caption{FORCE-3}
	\end{subfigure}
	
	\begin{subfigure}[b]{0.32\textwidth}
		\centering
		\includegraphics[width=\textwidth,trim={0 0 2.65cm 0},clip]
		{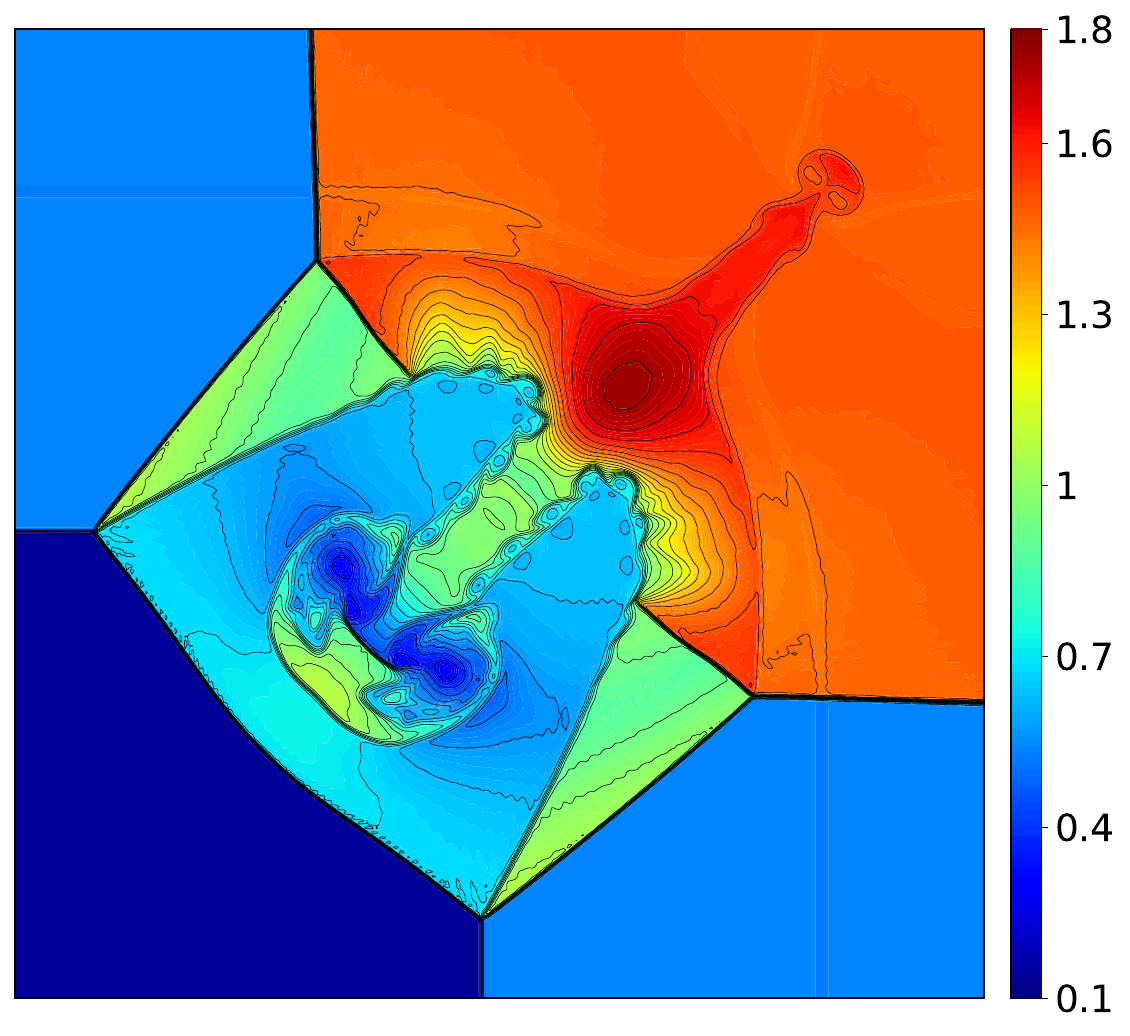}
		\caption{FORCE-5}
	\end{subfigure}
	\begin{subfigure}[b]{0.32\textwidth}
		\centering
		\includegraphics[width=\textwidth,trim={0 0 2.65cm 0},clip]
		{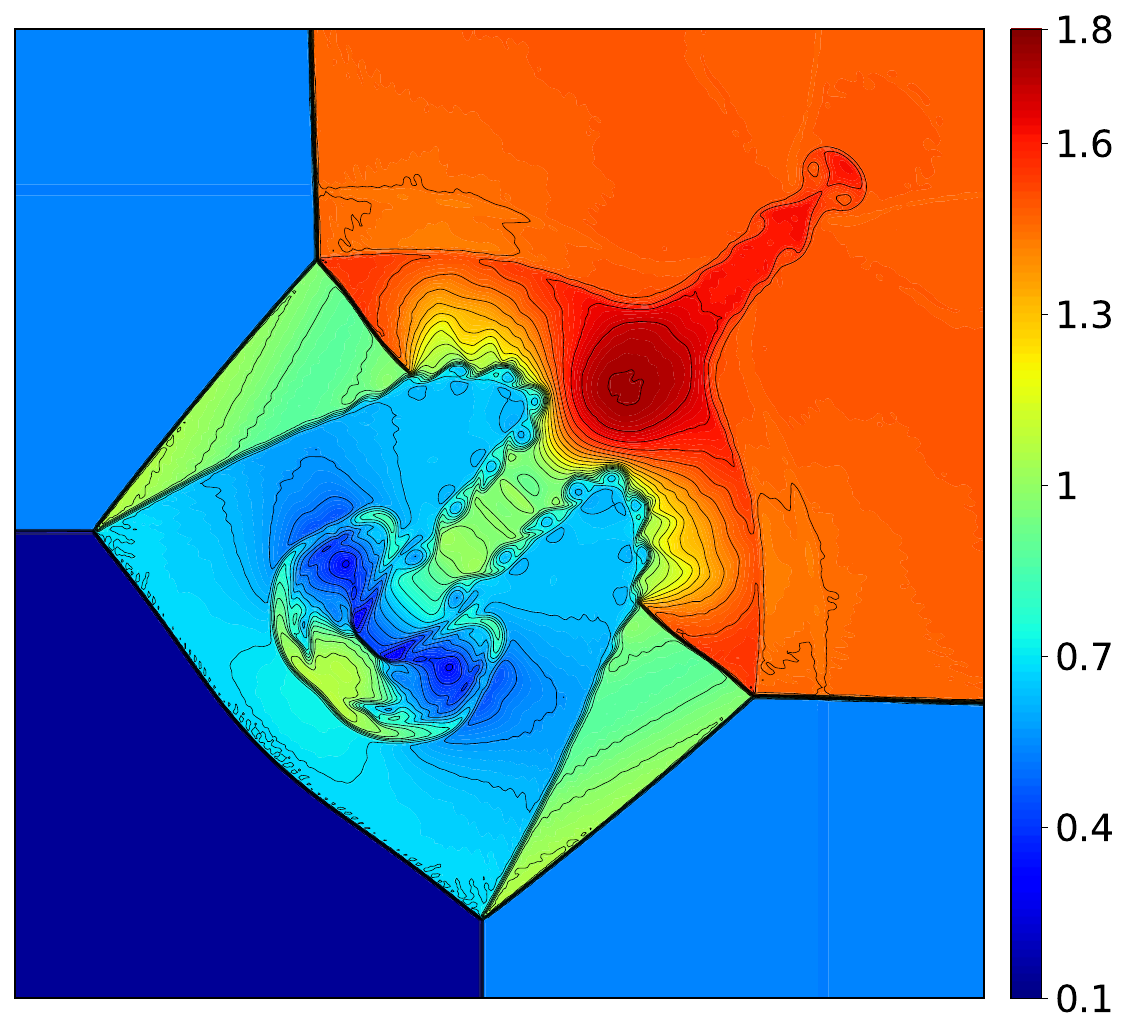}
		\caption{FORCE-10}
	\end{subfigure}
	
	\begin{subfigure}[b]{0.32\textwidth}
		\centering
		\includegraphics[width=\textwidth,trim={0 0 2.65cm 0},clip]
		{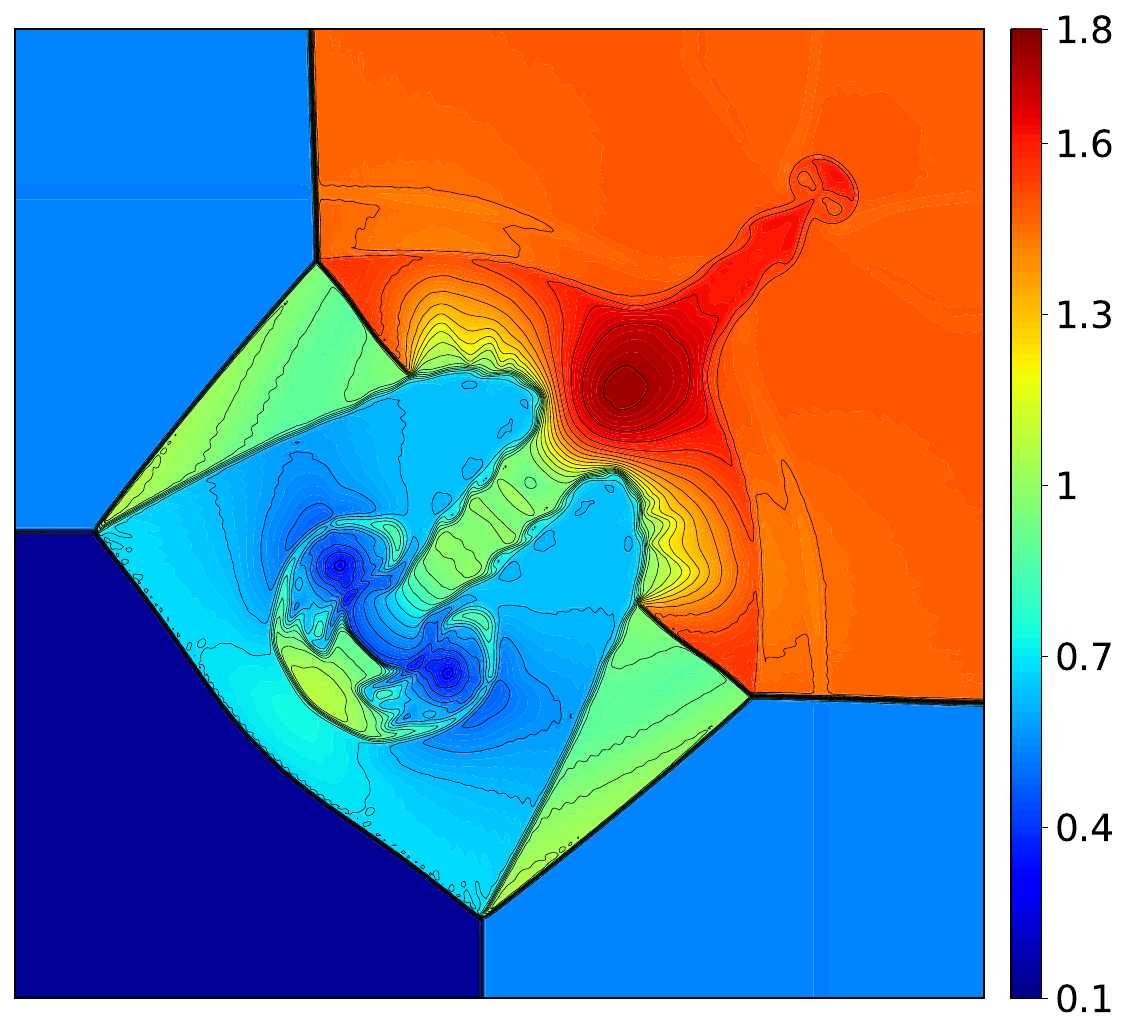}
		\caption{Rusanov}
	\end{subfigure}
	\hfill
	\begin{subfigure}[b]{0.32\textwidth}
		\centering
		\includegraphics[width=\textwidth,trim={0 0 2.65cm 0},clip]
		{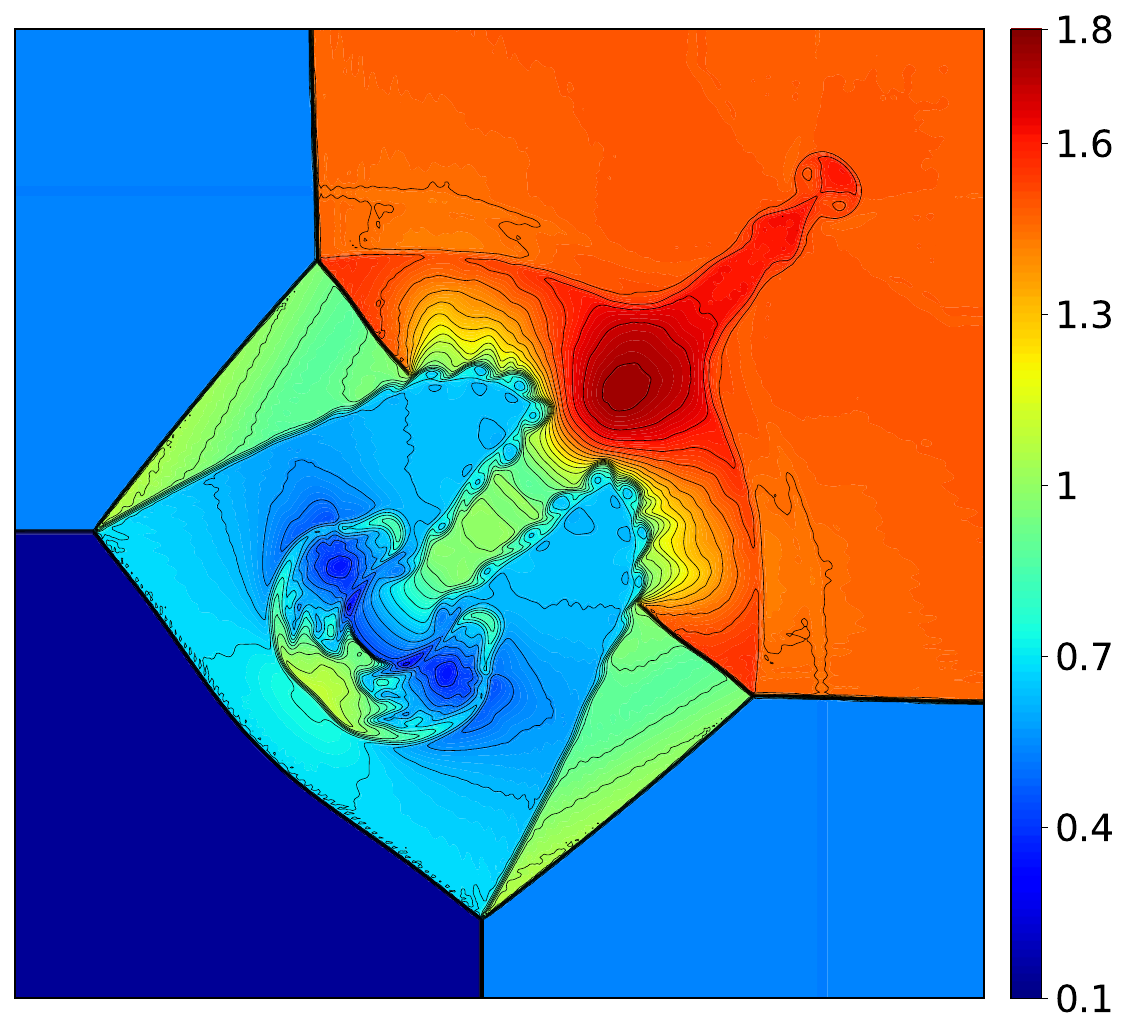}
		\caption{HLL}
	\end{subfigure}
	\hfill
	\begin{subfigure}[b]{0.32\textwidth}
		\centering
		\includegraphics[width=\textwidth,trim={0 0 2.65cm 0},clip]
		{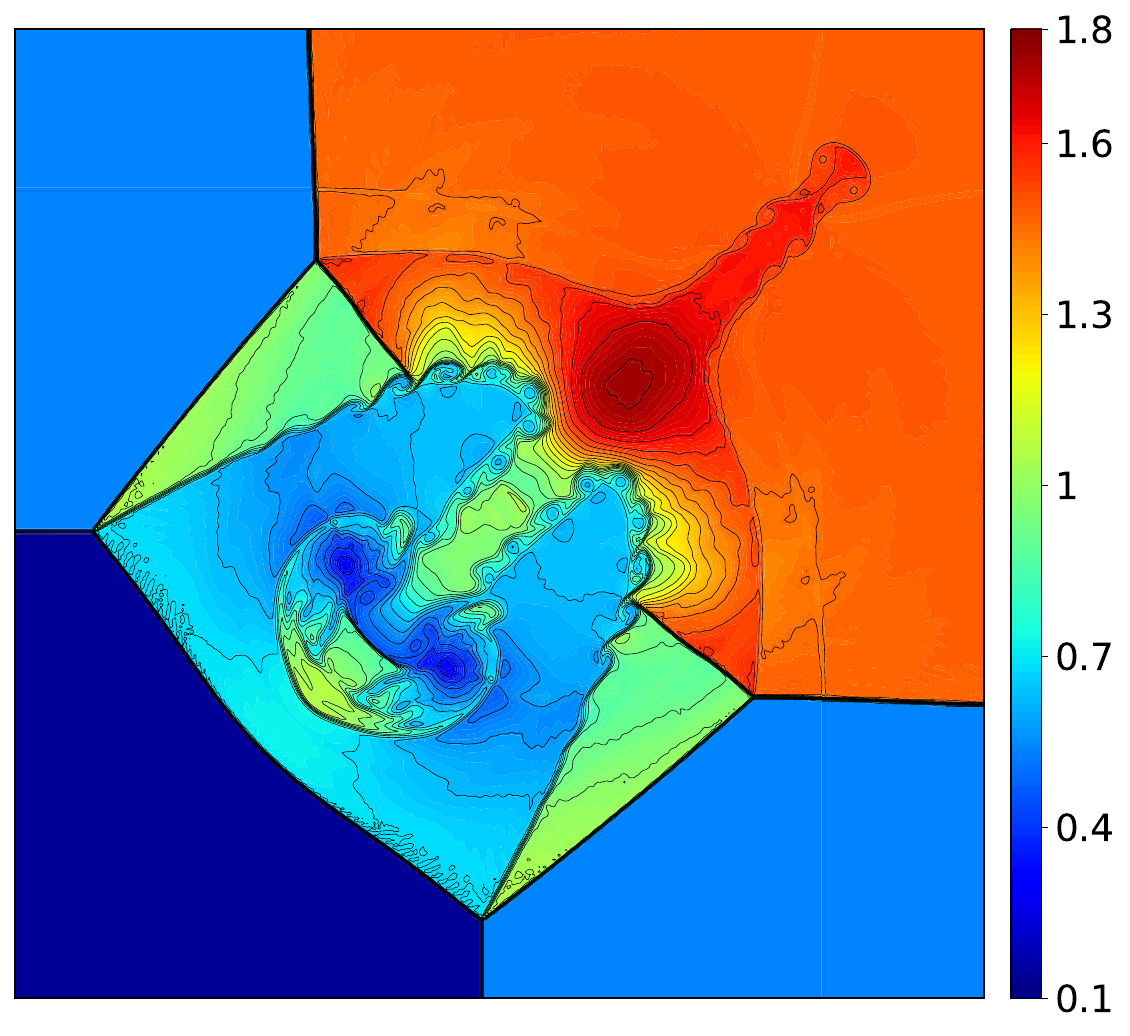}
		\caption{Exact RS}
	\end{subfigure}
	
	\vspace{0.4em}
	\includegraphics[width=0.42\textwidth]
	{figures_new/2DRP_common_density_colorbar_horizontal.pdf}
	
	\caption{\RIIcolor{Two-dimensional Riemann problem: Density obtained with order 7 over a mesh with $400\times 400$ elements and $\sigma_{CFL}:=0.9$. The same density color scale is used in all panels.}}
	\label{fig:2DRP_order7}
\end{figure}

\RIIcolor{In this last test, we consider the two--dimensional Riemann problem corresponding to Configuration 3 in \cite{kurganov2002solution}, see also \cite{CKX_Ustar}. 
On the computational domain $[0,1.2]\times[0,1.2]$ with transmissive boundary conditions, the initial condition reads
$$
\begin{pmatrix}
	\rho\\
	u\\
	v\\
	p
\end{pmatrix}(x,y,0) :=\begin{cases}(1.5,0,0,1.5)^\top,&x>1,~y>1,\\(0.5323,1.206,0,0.3)^\top,&x<1,~y>1,\\
	(0.138,1.206,1.206,0.029)^\top,&x<1,~y<1,\\(0.5323,0,1.206,0.3)^\top,&x>1,~y<1.
\end{cases}
$$

We ran our simulations on meshes with $400\times 400$ elements until the final time $T_f:=1$. 
All configurations could reach the final time with $\sigma_{CFL}:=0.9$.
The density results are reported in Figures~\ref{fig:2DRP_order5} and~\ref{fig:2DRP_order7} for orders 5 and 7 respectively.
For the considered mesh refinement, order 3 is rather diffusive for all numerical fluxes without fundamental differences to be appreciated, hence, it is omitted.

For order 5, all numerical fluxes still display a significant amount of numerical diffusion. 
Nevertheless, exact RS and FORCE-10 are able to capture the first visible vortical structures, while the remaining fluxes produce smoother results.

At order 7, the resolution significantly improves for all numerical fluxes. 
The global wave pattern is consistently reproduced, and the instability-driven vortical structures become clearly visible. 
As expected, exact RS provides the sharpest result.
FORCE-$\alpha$ fluxes with $\alpha=2$, 5 and 10 achieve a resolution comparable with HLL. 
Rusanov and FORCE-3 yield more diffusive similar results.

Also in this test, we can see how, within a very high order setting, FORCE-$\alpha$ fluxes, in particular FORCE-2, become competitive alternative to upwind fluxes.
}


\section{Conclusions and future perspectives}\label{sec:conclusions}


\RIIcolor{The main goal of this work was to systematically assess
the computational performance of the FORCE--$\alpha$ numerical fluxes introduced in~\cite{toro2020low} within the WENO--DeC framework~\cite{ciallella2022arbitrary,ciallella2023arbitrary,ciallella2025high,micalizzitoro2024,micalizzi2025algorithms}.} Such numerical fluxes are centred, meaning that they make no use of information coming from the structure of the solution of the Riemann problem at cell interfaces.
This aspect makes them particularly convenient for applications where such information is not directly and explicitly accessible, such as two (or more)--layer shallow-water systems~\cite{abgrall2009two,kim2009two}.
Moving from one of the main conclusions of~\cite{micalizzitoro2024,micalizzi2025algorithms}, which show that the impact of the numerical flux choice is less crucial as the order of accuracy increases towards very high values, we investigated FORCE--$\alpha$ numerical fluxes within a semidiscrete FV/DeC framework that allows the construction of schemes of arbitrarily high order accuracy, in both space and time. In the present paper, we presented results up to order 7.
%
%
\RIIcolor{In particular, the results indicate that moderate values of $\alpha$ ($1$--$2$ and $2$--$3$ in one and two space dimensions respectively) yield the best compromise between resolution, robustness, and computational cost. 
Moreover, FORCE--$\alpha$ schemes are shown to be a competitive alternative to classical upwind fluxes, e.g., exact RS, HLL and Rusanov, in the framework of very high order schemes.}
Future investigations in mind include the application of the present computational framework to systems for which exact Riemann solvers are very complex or simply  not available,  and complete approximate Riemann solvers might still be too expensive.
Other possible directions of investigation concern other types of FORCE--$\alpha$ numerical fluxes. In fact, alternative versions are suggested in the original reference~\cite{toro2020low}, while, even more sophisticated versions are currently under study.

\section*{Acknowledgments}
Lorenzo Micalizzi has been funded by the LeRoy B. Martin, Jr. Distinguished Professorship Foundation. 
The authors would like to acknowledge Prof. Rémi Abgrall from University of Zurich, Prof. Alina Chertock from North Carolina State University, Prof. Alexander Kurganov, Qingcheng Fu and Haoan Yi from Southern University of Science and Technology, members of the organizing committee of the Workshop ``Active Flux Methods: Developments and Applications'', held in Shenzhen, December 6-8 2025.´
A consistent part of this work was developed during such a conference.


\bibliography{sn-bibliography}
\bibliographystyle{abbrv}

\end{document}